\numberwithin{equation}{section}

\newcommand{\set}[1]{\left\{#1\right\}}

\newcommand{\To}{\longrightarrow}

\newcommand{\C}[1]{\mathbf{#1}} 

\def\r{\rightarrow} 
\def\l{\leftarrow} 
\def\rr{\Rightarrow} 

\def\hom{\operatorname{Hom}}

\def\aut{\operatorname{Aut}}

\def\im{\operatorname{Im}}

\def\ho{\operatorname{Ho}}

\def\F{\mathcal{F}} 

\def\P{\mathcal{P}} 

\def\st{\stackrel} 

\newcommand{\Ab}{\mathrm{\mathbf{Ab}}}

\def\Z{\mathbb{Z}}

\newcommand{\grupo}[1]{\langle #1\rangle}

\newcommand{\vc}{\square}    
\newcommand{\vi}{\boxminus}  

\newcommand{\pf}{\rightsquigarrow}  

\newcommand{\homm}[1]{[\![#1]\!]} 

\newcommand{\I}{\mathbb{I}} 

\newcommand{\J}{\mathbb{J}} 

\newcommand{\D}{\mathbb{D}} 

\numberwithin{equation}{section}

\makeatletter
\def\blfootnote{\xdef\@thefnmark{}\@footnotetext}
\makeatletter

\makeindex

\begin{document}

\title{{[III]}\\\hspace{5pt}\\The homotopy category of pseudofunctors and translation cohomology}%
\author{H.-J. Baues and F. Muro}%

\blfootnote{\textit{$2000$ Mathematics Subject Classification:} 18D05, 18G60, 55S35}%
\blfootnote{The second author was partially supported by the MEC
grant MTM2004-01865 and postdoctoral fellowship EX2004-0616}

\maketitle

\begin{abstract}
We develop the obstruction theory of the $2$-category of abelian
track categories, pseudofunctors and pseudonatural transformations
by using the cohomology of categories. We introduce
translation cohomology to classify endomorphisms in this category.
\end{abstract}


\section*{Introduction}

A homotopy category $\C{A}$ is derived from a groupoid enriched
category $\C{B}$ in which the morphism groupoids $\homm{X,Y}$ are
given by maps $f,g\colon X\r Y$ and tracks $\alpha\colon f\rr g$
which are homotopy classes of homotopies. The set of connected
components $$[X,Y]=\pi_0\homm{X,Y}$$ is the set of morphisms in
$\C{A}=\pi_0\C{B}$. Now let $\C{A}$ be the homotopy category of a
pointed stable Quillen model category. The suspension functor
$\Sigma$ is an endofunctor of $\C{A}$ which on $\C{B}$, however,
is only a pseudofunctor $\tilde{\Sigma}\colon\C{B}\pf\C{B}$. If
$\C{B}$ is enriched in abelian groupoids with $\pi_0\C{B}=\C{A}$
and $\pi_1\C{A}=D$ then the equivalence class of $\C{B}$ is
determined by an element in the cohomology group
\begin{equation}\label{intro1}
\grupo{\C{B}}\in H^3(\C{A},D),
\end{equation}
as proved in \cite{ccglg}, \cite{mhtcc}. In this paper we consider
the equivalence class of the pair $(\C{B},\tilde{\Sigma})$. For
this we introduce the \emph{translation cohomology}
$H^*(\C{A},\Sigma)$ and we show as a main result that the
equivalence class of $(\C{B},\tilde{\Sigma})$ is determined by the
cohomology class
\begin{equation}\label{intro2}
\grupo{\C{B},\tilde{\Sigma}}\in H^3(\C{A},\Sigma),
\end{equation}
generalizing the result (\ref{intro1}).

In order to achieve the result (\ref{intro2}) we study the homotopy
category of pseudofunctors $\C{Pseudo}^{ab}_\simeq$ which is a
$2$-category via pseudonatural transformations. We achieve a useful
obstruction theory for this category with obstructions in the
cohomology of categories. As the motivating application we show in
\cite{ccctc1} that for a stable homotopy category $(\C{A},\Sigma)$
the cohomology class $\grupo{\C{B},\tilde{\Sigma}}$ determines the
triangulation of $(\C{A},\Sigma)$, so that triangulated categories
have a cohomological description. On the other hand there is an
 application concerning the category of topological
$2$-types, which can be considered as a full subcategory of
$\C{Pseudo}^{ab}_\simeq$. This leads to the \emph{translation
cohomology of groups} which classifies pairs $(X,h)$ where $X$ is
a $2$-type and $h\colon X\r X$ is an endomorphism in the homotopy
category of spaces, see \cite{tecc} and \cite{e2tc}.

Some proofs in this paper are somewhat technical, for this reason
we  present the results in the first sections, \ref{homocat},
\dots, \ref{c2c}, and  concentrate most of the proofs in the final
ones, \ref{81}, \dots, \ref{gs}. Section \ref{ctc} contains some
background material about track categories which is completed by
the Appendix, where we recall the classical concepts of
pseudofunctor and pseudonatural transformation. In Section
\ref{homocat} we define the homotopy category of pseudofunctors
$\C{Pseudo}_\simeq$ as a quotient category, which is the main
object of study in this paper. The full subcategory
$\C{Pseudo}^{ab}_\simeq$ given by abelian track categories is the
localization of the category of pseudofunctors and of track
functors with respect to weak equivalences, as we prove in Section
\ref{atc} (Theorem \ref{local}).
In Section \ref{cc}
we recall the definition of the cohomology of small categories in
the sense of \cite{csc} as well as some of its $2$-functorial
structure developed in \cite{fcc}. Translation cohomology of
categories is defined in Section \ref{ct}, and the classification
result (\ref{intro2}) is proved in Section \ref{catctc} (Theorem
\ref{clastrans}).

Section \ref{os} is the core of the paper, it is devoted to the
obstruction theory in the homotopy category of pseudofunctors
$\C{Pseudo}^{ab}_\simeq$. In Section \ref{hcple} we give an
interpretation of this obstruction theory in terms of exact
sequences for functors and linear extensions of categories in the
sense of \cite{ah}. As a consequence we show in Section \ref{c2c}
that the homotopy category of pseudofunctors
$\C{Pseudo}^{ab}_\simeq$ is equivalent to a category defined
only in terms of cocycles in the cohomology of categories.

\section{Categories and track categories}\label{ctc}

\subsection*{Categories, functors and natural transformations}

The category $\C{Cat}$ of all categories is a well-known
$2$-category. The objects or $0$-cells of $\C{Cat}$ are the
categories which are small with respect to a fixed universe. The
$1$-cells are functors, and the $2$-cells are natural
transformations. The category $\C{Cat}$ has products denoted by
$\C{A}\times\C{B}$. A natural transformation can be regarded as a
functor from the \emph{cylinder category}\index{cylinder category}
$\C{A}\times\I\r\C{B}$. Here $\I$ is the \emph{interval
category}\index{interval category}, it has two objects $0$, $1$
and a unique non-trivial morphism $0\r 1$. There are two obvious
functors from the trivial category $*$, which has only one
morphism, $i_0,i_1\colon*\r\I,$ and a functor $p\colon\I\r*$.
These functors induce the structural functors of a cylinder
\begin{equation}\label{cyl}
\xymatrix{\C{A}\ar@<.5ex>[r]^{i_0}\ar@<-.5ex>[r]_{i_1}&\C{A}\times\I\ar[r]^p&\C{A}}.
\end{equation}

The \emph{double interval category}\index{double interval
category} $\J$ has three objects $0$, $1$, $2$ and two essential
morphisms $0\r1\r2$. It fits into a push-out diagram in $\C{Cat}$
\begin{equation}\label{dc}
\xymatrix{{*}\ar[r]^{i_1}\ar[d]_{i_0}\ar@{}[rd]|{\mathrm{push}}&\I\ar[d]\\\I\ar[r]&\J}
\end{equation}
There is a unique functor $j\colon\I\r\J$ with
$j(0)=0$ and $j(1)=2$. It can be used to define vertical composition of natural transformations regarded as functors from cylinders
since the product functor $\C{A}\times-$ preserves colimits in $\C{Cat}$.

In this paper a single arrow $\rightarrow$ will denote a morphism, a functor, or
a $1$-cell in a $2$-category, while double arrows $\rr$
will be kept for natural transformations and $2$-cells in $2$-categories. The composite of
morphisms $f, g$ will be indicated by juxtaposition $fg$, as well
as horizontal compositions in $2$-categories. The
symbol $\square$ will be used for vertical composition of
$2$-cells $\alpha\vc\beta$, and the inverse of an invertible
$2$-cell $\alpha$ will be denoted by $\alpha^\vi$.

\subsection*{Track categories, track functors and track natural transformations}

A \emph{track category}\index{track category} $\C{A}$ is a
category enriched in groupoids. Hence for objects $X, Y$ in
$\C{A}$ one has the hom-groupoid $\homm{X,Y}_{\C{A}}$ in $\C{A}$.
The objects $f\colon X\r Y$ of this groupoid are called
\emph{maps}, and the morphisms $\alpha\colon f\rr g$ are called
\emph{tracks}. Let $\homm{X,Y}_\C{A}(f,g)$ be the set of all
tracks $f\rr g$. The symbol for the identity track in $f$ is
$0^\vc_f$. Equivalently, a track category is a $2$-category all of
whose $2$-cells are invertible. Ordinary categories can be
regarded as track categories with only the identity tracks.

A track category $\C{A}$ will
also be depicted as $\C{A}_1\rightrightarrows\C{A}_0$. Here
$\C{A}_0$ is the underlying ordinary category; $\C{A}_1$ has the
same objects as $\C{A}$, morphisms form $X$ to $Y$ in $\C{A}_1$
are tracks $\alpha\colon f\rr g$ between maps $f,g\colon X\rightarrow Y$
in $\C{A}$, and the composition of $\alpha$ and $\beta$ as in the following  diagram
\begin{equation}\label{horizontal}
\xymatrix@C=35pt{Z&Y\ar@/_10pt/[l]_f^{\;}="a"\ar@/^10pt/[l]^{f'}_{\;}="b"&
X\ar@/_10pt/[l]_g^{\;}="c"\ar@/^10pt/[l]^{g'}_{\;}="d"\ar@{=>}^{\;\alpha}"a";"b"\ar@{=>}^{\;\beta}"c";"d"}
\end{equation}
is horizontal composition
$$\alpha\beta=(\alpha g')\vc(f\beta)=(f'\beta)\vc(\alpha g).$$
The arrows $\rightrightarrows$ from $\C{A}_1$ to $\C{A}_0$ denote
the obvious ``source'' and ``target'' functors.

A \emph{track functor}\index{track functor} is an enriched functor
between track categories, that is a $2$-functor, or equivalently,
an assignment of objects, maps and tracks which preserves all
composition laws and identities. The category $\C{Track}$ of track
categories and track functors also has colimit-preserving
products.

A \emph{track natural transformation}\index{track natural
transformation} $\alpha\colon\varphi\rr\psi$ between track
functors $\varphi,\psi\colon\C{A}\r\C{B}$. is a collection of maps
$\alpha_X\colon\varphi(X)\r\psi(X)$ in $\C{B}$ indexed by the
objects of $\C{A}$ satisfying certain axioms described in the
Appendix. One can alternatively define such a track natural
transformation $\alpha$ as a track functor
$\alpha\colon\C{A}\times\I\r\C{B}$. Track natural transformations
endow $\C{Track}$ with a $2$-category structure. Moreover,
(\ref{dc}) can be used to define vertical composition of track
natural transformations.


\subsection*{Pseudofunctors and pseudonatural transformations}\label{pseudo}

A \emph{pseudofunctor}\index{pseudofunctor} between track
categories $\varphi\colon\C{A}\pf\C{B}$ is an assignment of
objects, maps and tracks which preserves composition and
identities only up to certain given tracks
$$\varphi_{f,g}\colon\varphi(f)\varphi(g)\rr\varphi(fg)\text{ and
} \varphi_X\colon\varphi(1_X)\rr1_{\varphi(X)}.$$ These tracks
must satisfy well-known coherence and naturality properties which
are recalled in the Appendix. Moreover, the category $\C{Pseudo}$
of small track categories and pseudofunctors is defined.

Track functors can be regarded as a special case of
pseudofunctors. The inclusion functor $\C{Track}\r\C{Pseudo}$ has
a left adjoint
\begin{equation}\label{la}
\P\colon\C{Pseudo}\To\C{Track}.
\end{equation}
This adjoint is defined in \cite{fct} I,4.23 where
$\P(\C{A})=\tilde{\C{A}}$. A more handy description that we now
recall can be found in \cite{stt} 2.4. The category $\P(\C{A})_0$
is the free category on $\C{A}_0$ regarded as a graph, that is
objects are the same as in $\C{A}_0$ and a morphism
$(\sigma_1,\cdots,\sigma_n)\colon X\r Y$ in $\P(\C{A})_0$ is a
sequence of composable morphisms
$Y\st{\sigma_1}\l\cdots\st{\sigma_n}\l X$ in $\C{A}$, the identity
morphism is the empty sequence $()\colon X\r X$ and composition is
defined by concatenation. Track sets in $\P(\C{A})$ are defined by
the formula
$$\homm{X,Y}_{\P(\C{A})}((\sigma_1,\dots,\sigma_n),(\tau_1,\dots,\tau_m))=\homm{X,Y}_{\C{A}}(\sigma_1\cdots\sigma_n,\tau_1\cdots\tau_m),$$
and their composition laws are induced by the laws in $\C{A}$.
The adjointness yields natural bijections
\begin{equation}\label{adj}
ad\colon\C{Pseudo}(\C{A},\C{B})\cong\C{Track}(\P(\C{A}),\C{B}).
\end{equation}
Here $\C{Track}$ is a $2$-category so that the right hand side of
(\ref{adj}) is a category consisting of track functors as objects
and track natural transformations as morphisms. If we endow
$\C{Pseudo}$ with the $2$-category structure given by
pseudonatural transformations then (\ref{adj}) is an isomorphism
of categories.

Here a \emph{pseudonatural transformation}\index{pseudonatural
transformation} $\alpha\colon\varphi\rr\psi$ between two
pseudofunctors $\varphi,\psi\colon\C{A}\pf\C{B}$ is a collection
of maps $\alpha_X\colon\varphi(X)\r\psi(X)$ in $\C{B}$ indexed by
the objects of $\C{A}$ such that the usual squares only commute up
to given tracks
$\alpha_f\colon\alpha_Y\varphi(f)\rr\psi(f)\alpha_X$ indexed by
the maps $f\colon X\r Y$ in $\C{A}$. These tracks must satisfy
certain axioms described in the Appendix. By using the
isomorphisms of categories (\ref{adj}) one can also regard the
pseudonatural transformation $\alpha\colon\varphi\rr\psi$ as a
track functor $\alpha\colon\P(\C{A})\times\I\r\C{B}$.

There is a $2$-functor
\begin{equation}\label{pi0}
\pi_0\colon\C{Pseudo}\To\C{Cat}
\end{equation}
defined on objects by sending a track category $\C{A}$ to its
\emph{homotopy category} $\pi_0\C{A}=\C{A}_\simeq$ obtained from
$\C{A}_0$ by identifying those maps which are connected by tracks.
The equivalence class in $\pi_0\C{A}$ of a map $f$ in $\C{A}$ is
denoted by $\set{f}$. The quotient category $\pi_0\C{A}$ of
$\C{A}_0$ is a coequalizer in $\C{Cat}$ of the form
$$\C{A}_1\rightrightarrows\C{A}_0\r\pi_0\C{A}.$$

\section{The homotopy category of pseudofunctors}\label{homocat}

A \emph{homotopy}\index{homotopy between pseudofunctors}
$\xi\colon\varphi\rr\psi$ between two pseudofunctors
$\varphi,\psi\colon\C{A}\pf\C{B}$ in $\C{Pseudo}$ is a
pseudonatural transformation $\xi$ such that the maps
$\xi_X\colon\varphi(X)\r\psi(X)$
are identities, that is $\varphi(X)=\psi(X)$ and 
$\xi_X=1_{\varphi(X)}=1_{\psi(X)}$,
in particular homotopic pseudofunctors induce the same functor on $\pi_0$.

The homotopy relation $\simeq$ 
is a natural equivalence relation 
since homotopies
are invertible $2$-cells in this $2$-category. Therefore the
homotopy category of pseudofunctors is well defined
$$\C{Pseudo}_\simeq.$$
This category, which is of main interest in this paper, can be endowed
with a $2$-category structure in the following way: let us write $[\varphi]$ for the homotopy class of a pseudofunctor $\varphi$,
a $2$-cell
$[\alpha]\colon[\varphi]\rr[\psi]$ is represented by a pseudonatural
transformation $\alpha\colon\varphi\rr\psi$ between two
pseudofunctors $\varphi,\psi\colon\C{A}\pf\C{B}$, which is
equivalent to $\alpha'\colon\varphi'\rr\psi'$, $[\alpha]=[\alpha']$, if $\varphi\simeq\varphi'$,
$\psi\simeq\psi'$, and there exist tracks $\alpha_X\rr\alpha'_X$ for all objects $X$ in $\C{A}$. The horizontal composition of $2$-cells
is the composition of representatives
$$\xymatrix@C=40pt{\C{A}\ar@/^10pt/[r]^{[\varphi]}_{\;}="a"\ar@/_10pt/[r]_{[\varphi']}^{\;}="b"&
\C{B}\ar@/^10pt/[r]^{[\psi]}_{\;}="c"\ar@/_10pt/[r]_{[\psi']}^{\;}="d"&\C{C}\ar@{=>}^{[\alpha]}"a";"b"\ar@{=>}^{[\beta]}"c";"d"}
=\xymatrix@C=45pt{\C{A}\ar@/^10pt/[r]^{[\psi\varphi]}_{\;}="a"\ar@/_10pt/[r]_{[\psi'\varphi']}^{\;}="b"&
\C{C}\ar@{=>}^{[\beta\alpha]}"a";"b"}.$$
For the vertical composition of $2$-cells 
represented by pseudonatural transformations $\alpha\colon\varphi\rr\psi$ and $\beta\colon\psi'\rr\phi$ we have to choose
a homotopy $\xi\colon\psi\rr\psi'$
$$\xymatrix@C=60pt{\C{A}\ar@/^25pt/[r]^{[\varphi]}_{\;}="a"\ar[r]|{[\psi]=[\psi']}^{\;}="b"_{\;}="c"\ar@/_25pt/[r]_{[\phi]}^{\;}="d"&
\C{C}\ar@{=>}^{[\alpha]}"a";"b"\ar@{=>}^{[\beta]}"c";"d"}
=\xymatrix@C=50pt{\C{A}\ar@/^15pt/[r]^{[\varphi]}_{\;}="a"\ar@/_15pt/[r]_{[\phi]}^{\;}="b"&
\C{C}\ar@{=>}^{[\beta\xi\alpha]}"a";"b"}.$$

The $2$-functor $\pi_0$ in (\ref{pi0}) clearly factors through the
homotopy $2$-category of pseudo functors
\begin{equation}\label{pi1}
\pi_0\colon\C{Pseudo}_{\simeq}\To\C{Cat}.
\end{equation}
Moreover, the following result is an immediate consequence of the equivalence relation defining $2$-cells in $\C{Pseudo}_\simeq$.

\begin{prop}\label{faith}
The $2$-functor $\pi_0$ in (\ref{pi1}) is faithful on $2$-cells.
\end{prop}

\section{Abelian track categories}\label{atc}

In this paper we will concentrate on \emph{abelian track
categories}\index{track category!abelian}, which are the track
categories $\C{A}$ for which the automorphism group
$\aut_{\homm{X,Y}_\C{A}}(f)$ of any map $f\colon X\r Y$ in the
groupoid $\homm{X,Y}_\C{A}$ is an abelian group. Equivalently an
abelian track category is a category enriched in abelian
groupoids. We write $\C{Pseudo}^{ab}$ for the $2$-category of
small abelian track categories, pseudofunctors and pseudonatural
transformations, and $\C{Track}^{ab}$ for the $2$-category of
small abelian track categories, track functors and track natural
transformations.

The $\pi_0$ of a track category was defined in
Section \ref{pseudo}. In order to define the $\pi_1$ of an abelian track
category we need to recall some concepts from \cite{csc}.

For any category $\C{C}$ the \emph{factorization
category}\index{factorization category} $\F\C{C}$ is given as
follows. Objects are the morphisms in $\C{C}$, and a morphism
$(h,k)\colon f\r g$ is a pair of morphisms in $\C{C}$ with
$kfh=g$. This gives rise to a functor
\begin{equation}\label{fac1}
\F\colon\C{Cat}\To\C{Cat}.
\end{equation}
%
%
%
%
A \emph{natural system}\index{natural system} $D$ (of abelian
groups) on $\C{C}$ is just a functor $D\colon \F\C{C}\r\Ab$ from
the factorization category to the category of abelian groups. We
usually denote $D_f=D(f)$, $h^*=D(h,1)$ and $k_*=D(1,k)$.

For any abelian track category $\C{A}$ there is a well-defined
natural system $\pi_1\C{A}$ on $\pi_0\C{A}$, see \cite{catc} 2.4. 
This natural system is up to isomorphism determined by the existence of
isomorphisms in $\C{Ab}$
\begin{equation}
\sigma_f\colon(\pi_1\C{A})_{\set{f}}\cong\aut_{\homm{X,Y}_\C{A}}(f)
\end{equation}
such that given a track $\alpha\colon f\rr g$ in
$\C{A}$,
\begin{equation}\label{lec1}
\alpha\vc\sigma_f(x)=\sigma_g(x)\vc\alpha;
\end{equation}
and given composable maps
$\bullet\st{h}\r\bullet\st{g}\r\bullet\st{f}\r\bullet$ in $\C{A}$
\begin{equation}\label{lec2}
f\sigma_g(x)=\sigma_{fg}(\set{f}_*x) \text{ and }
\sigma_g(x)h=\sigma_{gh}(\set{h}^*x).
\end{equation}
This means that any abelian track category $\C{A}$ is endowed with
the structure of a \emph{linear track extension}\index{linear
track extension} denoted by
$$\pi_1\C{A}\r\C{A}_1\rightrightarrows\C{A}_0\r\pi_0\C{A},$$ compare \cite{ccglg}.

The functorial properties of $\pi_1$ are described in terms of the following $2$-category $\C{nat}$. The $0$-cells
of $\C{nat}$ are pairs $(\C{C},D)$
where $D$ is a natural system on a small category $\C{C}$, $1$-cells
$(\varphi,\lambda)\colon(\C{C},D)\r(\C{D},E)$ are given by a functor
$\varphi\colon\C{C}\r\C{D}$ together with a natural transformation $\lambda\colon D\rr E\F(\varphi)$, and $2$-cells $\alpha\colon
(\varphi,\lambda)\rr(\psi,\zeta)$ are natural transformations $\alpha\colon\varphi\rr\psi$ such that for any morphism
$f\colon X\r Y$ in the source of $\varphi$ and $\psi$ the following diagram of abelian groups commutes
\begin{equation}\label{hdc}
\xymatrix{D_f\ar[r]^{\lambda_f}\ar[d]_{\zeta_f}&E_{\varphi(f)}\ar[d]^{{\alpha_Y}_*}\\E_{\psi(f)}\ar[r]_{\alpha_X^*}&
**[r] E_{\alpha_Y\varphi(f)}=E_{\psi(f)\alpha_X}}
\end{equation}
Composition of $1$-cells is defined as
$$(\varphi,\lambda)(\psi,\zeta)=(\varphi\psi,(\lambda\F(\psi))\vc\zeta),$$
and the composition laws involving $2$-cells in $\C{nat}$ are induced by those in $\C{Cat}$.

There is a $2$-functor
\begin{equation}\label{pi}
\pi\colon\C{Pseudo}^{ab}\To\C{nat}
\end{equation}
defined on cells of dimensions $0$, $1$ and $2$  as
$$\pi\C{A}=(\pi_0\C{A},\pi_1\C{A}),\;\;\;\pi\varphi=(\pi_0\varphi,\pi_1\varphi),\text{ and }\pi\tau=\pi_0\tau,$$
respectively. Here given a pseudofunctor
$\varphi\colon\C{A}\pf\C{B}$ between abelian track categories the
natural transformation
$$\pi_1\varphi\colon\pi_1\C{A}\rr(\pi_1\C{B})\F(\pi_0\varphi)$$
is determined by the formula
\begin{equation}\label{lec3}
(\pi_1\varphi)_{\set{f}}\sigma_f^{-1}(\alpha)=\sigma_{\varphi(f)}^{-1}(\varphi(\alpha)),
\end{equation}
where $\alpha\colon f\rr f$
is any self-track in $\C{A}$.
In fact this $2$-functor factors through the homotopy $2$-category of small abelian track categories
\begin{equation}\label{pi2}
\pi\colon\C{Pseudo}^{ab}_\simeq\;\To\C{nat}.
\end{equation}


A pseudofunctor $\varphi\colon\C{A}\pf\C{B}$ between abelian track
categories is said to be a \emph{homotopy equivalence} if
$[\varphi]$ is an isomorphism in the homotopy category of
pseudofunctors, and it is a \emph{weak equivalence}\index{weak
equivalence of track categories} if $\pi\varphi$ is an isomorphism
in $\C{nat}$. Obviously homotopy equivalences are weak
equivalences. We write $\ho(\C{Pseudo}^{ab})$ and
$\ho(\C{Track}^{ab})$ for the corresponding localized categories
with respect to weak equivalences.

\begin{thm}\label{local}
Weak equivalences are homotopy equivalences. Moreover, the
quotient functor induces a functor
$$\ho(\C{Pseudo}^{ab})\cong\C{Pseudo}^{ab}_\simeq$$ which is an isomorphism of categories. Furthermore, the inclusion
functor induces an isomorphism of categories
$$\ho(\C{Track}^{ab})\cong\ho(\C{Pseudo}^{ab}).$$
\end{thm}

This result underlines the natural significance of he homotopy
category $\C{Pseudo}^{ab}_\simeq$ as the localization of
$\C{Track}^{ab}$.

\begin{proof}[Proof of Theorem \ref{local}]
Weak equivalences in $\C{Pseudo}^{ab}$ are homotopy equivalences
because $\pi$ in (\ref{pi2}) fits into an exact sequence for
functors, as we show in (\ref{esf}) below, see
\cite{ah} IV.4.11 and IV.1.3 (a).

Let us check that the quotient functor
$\varrho\colon\C{Pseudo}^{ab}\r\C{Pseudo}^{ab}_\simeq$ satisfies
the universal property of the localization $\ho(\C{Pseudo}^{ab})$.
We just have to show that a functor
$\varsigma\colon\C{Pseudo}^{ab}\r\C{C}$ carrying weak equivalences
to isomorphisms does not distinguish between homotopic
pseudofunctors.

In the Appendix we consider the reduced cylinder track category
$\D$. This is an abelian track category together with a diagram of
weak equivalences
\begin{equation*}
\xymatrix{{*}\ar@<.5ex>[r]^{j^0}\ar@<-.5ex>[r]_{j^1}&\D\ar[r]^q&{*}},
\end{equation*}
where $*$ is the final track category. We can multiply
this diagram by any other abelian track category $\C{A}$ to obtain
a new diagram of weak equivalences
\begin{equation*}
\xymatrix{\C{A}\ar@<.5ex>[r]^<(.3){j^0}\ar@<-.5ex>[r]_<(.3){j^1}&\C{A}\times\D\ar[r]^q&\C{A}}.
\end{equation*}
Suppose that we have a homotopy $\xi\colon\varphi\rr\psi$ between
pseudofunctors $\varphi,\psi\colon\C{A}\pf\C{B}$ in
$\C{Pseudo}^{ab}$. By Proposition \ref{carhomo} there is a
pseudofunctor $\bar{\xi}\colon\C{A}\times\D\pf\C{B}$ with
$\varphi=\bar{\xi} j^0$ and $\psi=\bar{\xi} j^1$. The equalities
$q j^0=1=q j^1$ hold in $\C{Pseudo}^{ab}$, therefore
$\varsigma(j^0)=\varsigma(q)^{-1}=\varsigma(j^1)$ in $\C{C}$, and
so
$\varsigma(\varphi)=\varsigma(\bar{\xi})\varsigma(q)^{-1}=\varsigma(\psi)$.

Let us now check that the functor
$\bar{\varrho}\colon\C{Track}^{ab}\r\ho(\C{Pseudo}^{ab})$ given by
the inclusion $\C{Track}^{ab}\subset\C{Pseudo}^{ab}$ and the
quotient functor $\varrho$ above satisfies the universal property
of the localization $\ho(\C{Track}^{ab})$. Obviously
$\bar{\varrho}$ sends weak equivalences to isomorphisms.

Recall from (\ref{la}) and (\ref{adj}) that there is a natural
track functor $ad(1_\C{A})\colon\P(\C{A})\r\C{A}$, namely the
adjoint of the identity. The track category $\P(\C{A})$ is abelian
provided $\C{A}$ is, and in this case $ad(1_\C{A})$ is a week
equivalence, see \cite{stt} 2.4. The track functor $ad(1_\C{A})$
is a retraction in $\C{Pseudo}^{ab}$ with splitting section
$\epsilon\colon\C{A}\pf\P(\C{A})$. Here $\epsilon$ is the unit of
the adjunction between the inclusion $\C{Track}\subset\C{Pseudo}$
and $\P$. In particular the section $\epsilon$ is necessarily a
weak equivalence. Actually any pseudofunctor $\varphi$ in $\C{Pseudo}^{ab}$
factors as $\varphi=ad(\varphi)\epsilon$

Suppose now that we have a functor
$\varsigma\colon\C{Track}^{ab}\r\C{C}$ sending weak equivalences
to isomorphisms, then we define a factorization
$\bar{\varsigma}\colon\C{Pseudo}^{ab}\r\C{C}$ through the
inclusion as follows, given a pseudofunctor
$\varphi\colon\C{A}\pf\C{B}$ in $\C{Pseudo}^{ab}$,
$$\bar{\varsigma}\varrho(\varphi)=\varsigma(ad(\varphi))\varsigma(ad(1_\C{A}))^{-1}.$$
This actually defines an factorization because we know that
$\varrho$ is full, since we have already established the first
isomorphism of the statement. This factorization is the unique
possible one because of the equalities
$ad(1_\C{A})\epsilon=1_\C{A}$ and $\varphi=ad(\varphi)\epsilon$,
therefore we will be done if we manage to prove that
$\bar{\varsigma}$ also sends weak equivalences to isomorphisms.
But if $\varphi$ is a weak equivalence then $ad(\varphi)$ is also
a weak equivalence because the equality
$\varphi=ad(\varphi)\epsilon$ holds and in this case
$\bar{\varsigma}\varrho(\varphi)$ is an isomorphism. The proof is
now finished.
\end{proof}

\section{Cohomology of categories}\label{cc}

The cohomology $H^*(\C{C},D)$ of a small category $\C{C}$ with
coefficients in a natural system $D$ is the cohomology of the
Baues-Wirsching complex defined in \cite{csc}. This complex
$F^*(\C{C},D)$ is a cochain complex of abelian groups concentrated
in dimensions $\geq 0$. Its group in dimension $n$ is the
following product indexed by all sequences of morphisms of length
$n$ in $\C{C}$
\begin{equation*}
F^n(\C{C},D)=\prod_{\bullet\st{\sigma_1}\l\cdots\st{\sigma_n}\l\bullet}D_{\sigma_1\cdots\sigma_n}.
\end{equation*}
Here we assume that a sequence of length $0$ is an object $X$ in
$\C{C}$ which we also identify with the identity morphism $1_X$.
The coordinate of $c\in F^n(\C{C},D)$ in
$\bullet\st{\sigma_1}\l\cdots\st{\sigma_n}\l\bullet$ will be
denoted by $c(\sigma_1,\dots,\sigma_n)\in D_{\sigma_1\cdots\sigma_n}$. The differential $\delta$
is defined as
\begin{eqnarray*}
  \delta(c)(\sigma_1,\dots,\sigma_{n+1}) &=& {\sigma_1}_*c(\sigma_2,\dots,\sigma_{n+1}) \\
   && +\sum_{i=1}^{n}(-1)^ic(\sigma_1,\dots,\sigma_i\sigma_{i+1},\dots,\sigma_{n+1}) \\
   && +(-1)^{n+1}\sigma_{n+1}^*c(\sigma_1,\dots,\sigma_n)
\end{eqnarray*}
over an $n$-cochain $c$ for $n\geq 1$, and
$\delta(c)(\sigma)=\sigma_*c(X)-\sigma^*c(Y)$ for $n=0$ and
$\sigma\colon X\r Y$.

A functor $\varphi\colon\C{D}\r\C{C}$ and a natural transformation
$\lambda\colon D\rr E$ induce cochain homomorphisms
$$\varphi^*\colon F^*(\C{C},D)\To F^*(\C{D},D\F(\varphi)),\;\; \varphi^*(c)(\sigma_1,\dots,\sigma_n)=
c(\varphi(\sigma_1),\dots,\varphi(\sigma_n));$$
$$\lambda_*\colon F^*(\C{C},D)\To F^*(\C{C},E),\;\; \lambda_*(c)(\sigma_1,\dots,\sigma_n)=
\lambda_{\sigma_1\cdots\sigma_n}(c(\sigma_1,\dots,\sigma_n));$$
and the corresponding homomorphisms in cohomology
$$\varphi^*\colon H^*(\C{C},D)\To H^*(\C{D},D\F(\varphi)),$$
$$\lambda_*\colon H^*(\C{C},D)\To H^*(\C{C},E).$$
See \cite{csc} for further details on the functoriality of cohomology of categories.

A cochain $c\in F^n(\C{C},D)$ is said to be
\emph{normalized}\index{normalized cochain} (at identities) if
$c(\sigma_1,\dots,\sigma_n)=0$ provided $\sigma_i$ is an identity
morphism for some $1\leq i\leq n$. The inclusion of the subcomplex
$\bar{F}^*(\C{C},D)\subset F^*(\C{C},D)$ of normalized cochains
induces an isomorphism in cohomology, see \cite{ccglg}. Moreover,
the induced cochain homomorphisms above restrict to normalized
cochains.


The factorization functor $\F$ in (\ref{fac1}) can be extended to
a $2$-functor if we consider the $2$-category $\C{Cat}_\F$, see
\cite{fcc} 3.2. This $2$-category is obtained by applying the
product-preserving functor $\F$ to morphism categories in the
$2$-category $\C{Cat}$. This $2$-functor will be used in the next section to describe the functorial properties of translation cohomology
as well as
from Section \ref{os} on in the context of the obstruction theory
for pseudonatural transformations. More precisely, $0$-cells in
$\C{Cat}_\F$ are small categories, a $1$-cell
$\alpha\colon\C{C}\r\C{D}$ is a natural transformation
$\alpha\colon\varphi\rr\psi$ between functors
$\varphi,\psi\colon\C{C}\r\C{D}$, composition is horizontal
composition, so that composition is defined in the same way as in
(\ref{horizontal}). A $2$-cell
$(\varepsilon,\gamma)\colon\alpha\rr\beta$ between morphisms
$\alpha,\beta\colon\C{C}\r\C{D}$ is a pair of natural
transformations such that $\gamma\vc\alpha\vc\varepsilon=\beta$.
We can place $\C{Cat}$ as an ordinary category inside $\C{Cat}_\F$
by using the ordinary functor $\C{Cat}\r\C{Cat}_\F$ which is the
identity on objects and sends a functor
$\varphi\colon\C{C}\r\C{D}$ to the identity natural transformation
$0_{\varphi}^\square\colon\varphi\rr\varphi$ regarded as a
morphism $0^\square_\varphi\colon\C{C}\r\C{D}$ in $\C{Cat}_\F$.

The $2$-functor
\begin{equation}\label{fac2}
\F\colon\C{Cat}_\F\To\C{Cat}
\end{equation}
extending (\ref{fac1}) is defined in the same way on $0$-cells, on a $1$-cell as above $\F(\alpha)\colon\C{C}\r\C{D}$ is the functor
with $\F(\alpha)(f)=\alpha_Y\varphi(f)=\psi(f)\alpha_X$ for any morphism $f\colon X\r Y$ in $\C{C}$ regarded as an object in $\F(\C{C})$ and
$\F(\alpha)(h,k)=(\varphi(h),\psi(k))$ for any morphism $(h,k)$ in $\F(\C{C})$, and on a $2$-cell as above $\F(\varepsilon,\gamma)\colon
\F(\alpha)\rr\F(\beta)$ is the natural transformation given by $\F(\varepsilon,\gamma)_f=(\varepsilon_X,\gamma_Y)$.

This formalism can be used to restate the condition defining
$2$-cells in $\C{nat}$, more precisely, (\ref{hdc}) commutes for
all $f$ if and only if the following equality holds
\begin{equation}\label{hdc2}
(E\F(0^\vc_\varphi,\alpha))\vc\lambda=(E\F(\alpha,0^\vc_\psi))\vc\zeta.
\end{equation}

A natural transformation $\alpha\colon\varphi\rr\psi$
between functors $\varphi,\psi\colon\C{D}\r\C{C}$ induces a cochain homomorphism
$$\alpha^*\colon F^*(\C{C},D)\To F^*(\C{D},D\F(\alpha))$$ defined as $$\alpha^*(c)(\sigma_1,\cdots,\sigma_n)=
(D\F(0^\square_\varphi,\alpha))_{\sigma_1\cdots\sigma_n}
(c(\varphi(\sigma_1),\dots,\varphi(\sigma_n))).$$ This cochain
homomorphism is called $F^*(\alpha,0^\vc_{D\F(\alpha)})$ in
\cite{fcc} 4.1. The induced homomorphism in cohomology is denoted
in the same way
$$\alpha^*\colon H^*(\C{C},D)\To H^*(\C{D},D\F(\alpha)).$$
Notice that the cochain homomorphism defined by an identity
natural transformation coincides with the cochain homomorphism
defined by the corresponding functor
$$(0^\square_\varphi)^*=\varphi^*.$$ The natural transformation
$\alpha\colon\varphi\rr\psi$
 also induces cochain homotopies given by degree $-1$ homomorphisms
$$\alpha_\#,\alpha^\#\colon F^*(\C{C},D)\To F^*(\C{D},D\F(\alpha)),$$
with
\begin{equation}\label{homo11}
\delta\alpha_\#+\alpha_\#\delta=-(D\F(0^\square_\varphi,\alpha))_*\varphi^*+\alpha^*,
\end{equation}
\begin{equation}\label{homo12}
\delta\alpha^\#+\alpha^\#\delta=-(D\F(\alpha,0^\square_\psi))_*\psi^*+\alpha^*.
\end{equation}
The cochain homotopies $\alpha_\#$ and $\alpha^\#$ are
$h_{0^\vc_\varphi,\alpha}$ and $h_{\alpha,0^\vc_\psi}$ in
\cite{fcc} 4.1, respectively. They are defined by the following
formulas
\begin{flushleft}
$\alpha_\#(c)(\sigma_1,\dots,\sigma_n)=$
\end{flushleft}
\begin{flushright}
$\sum\limits_{i=0}^n(-1)^i(D\F(0^\square_\varphi,\alpha))_{\sigma_1\cdots\sigma_n}(
c(\varphi(\sigma_1),\dots,\varphi(\sigma_i),1_{\varphi(X_i)},\varphi(\sigma_{i+1}),\dots,\varphi(\sigma_n))),$
\end{flushright}

\begin{equation}\label{almo}
\alpha^\#(c)(\sigma_1,\dots,\sigma_n)=\sum_{i=0}^n(-1)^i
c(\psi(\sigma_1),\dots,\psi(\sigma_i),\alpha_{X_i},\varphi(\sigma_{i+1}),\dots,\varphi(\sigma_n)).
\end{equation}
Here $X_i$ is the source of $\sigma_i$ and/or the target of
$\sigma_{i+1}$. These homotopies are part of the $2$-functorial
properties of the Baues-Wirsching complex studied in \cite{fcc}.
Notice that these induced cochain homomorphisms and homotopies
restrict to the subcomplex of normalized cochains. Moreover, as
one can readily notice on the normalized cochain complex the
following formulas hold
\begin{equation}\label{alpha0}
\alpha_\#=0,
\end{equation}
\begin{equation}\label{00}
(0^\vc_\varphi)^\#=0,
\end{equation}
\begin{equation}\label{homo12'}
\alpha^*=(D\F(0^\square_\varphi,\alpha))_*\varphi^*,
\end{equation}
in particular by (\ref{homo11}) and (\ref{homo12})
\begin{equation}\label{homo13}
\delta\alpha^\#+\alpha^\#\delta=-(D\F(\alpha,0^\square_\psi))_*\psi^*+(D\F(0^\square_\varphi,\alpha))_*\varphi^*.
\end{equation}
Observe that for any natural transformation $\lambda\colon D\r E$ the following equality holds,
\begin{equation}\label{homo14}
\alpha^\#\lambda_*=(\lambda\F(\alpha))_*\alpha^\#.
\end{equation}

We finally recall the behaviour of these cochain homotopies with
respect to horizontal and vertical composition, see \cite{fcc}
5.1. Given natural transformations
$\phi\st{\beta}\rr\varphi\st{\alpha}\rr\psi$ there is a degree
$-2$ homomorphism
$$(\alpha,\beta)^\%\colon F^*(\C{C},D)\To F^*(\C{D},D\F(\alpha\vc\beta))$$
such that
\begin{eqnarray}
\label{homo21} \delta(\alpha,\beta)^\%-(\alpha,\beta)^\%\delta&=&-(D\F(\beta,0^\vc_{\psi}))_*\alpha^\#\\
\nonumber&&-(D\F(0^\vc_{\phi},\alpha))_*\beta^\#+(\alpha\vc\beta)^\#.
\end{eqnarray}
This order $2$ cochain homotopy coincides with
$r_{(\beta,0^\vc_\varphi);(\alpha,0^\vc_\psi)}$ in \cite{fcc} 4.1.
If we have functors $\varphi,\varphi'\colon\C{D}\r\C{C}$,
$\psi,\psi'\colon\C{E}\r\C{D}$ and natural transformations
$\alpha\colon\varphi\rr\varphi'$, $\beta\colon\psi\rr\psi'$, there
is a degree $-2$ homomorphism
$$(\alpha,\beta)^\$\colon
F^*(\C{C},D)\To F^*(\C{E},D\F(\alpha\beta)),$$ such
that
\begin{eqnarray}
\label{homo22}\delta(\alpha,\beta)^\$-(\alpha,\beta)^\$\delta&=&-\beta^\#\alpha^*-\beta^*\alpha^\#+(\alpha\beta)^\#.
\end{eqnarray}
The order $2$ cochain homotopy $(\alpha,\beta)^\$$ is
$r'_{(\beta,0^\vc_{\psi'});(\alpha,0^\vc_{\varphi'})}$ in
\cite{fcc} 4.1. We shall not give here the explicit formulas of
$(\alpha,\beta)^\%$ and $(\alpha,\beta)^\$$ since they are not
needed in this paper, however we point out that they restrict to
the subcomplex of normalized cochains. Moreover, on normalized
cochains the following equalities holds
\begin{equation}\label{homo23}
(\alpha,0^\vc_\psi)^\$=0,
\end{equation}
\begin{equation}\label{homo24}
(0^\vc_\varphi,\beta)^\$=0,
\end{equation}
and therefore by (\ref{homo22}) the equalities
\begin{equation}\label{homo25}
\psi^*\alpha^\#=(\alpha\psi)^\#,
\end{equation}
\begin{equation}\label{homo26}
\beta^\#\varphi^*=(\varphi\beta)^\#,
\end{equation}
are also satisfied on the normalized cochain complex.


\section{Translation cohomology}\label{ct}

Let $\C{A}$ be a category, $D$ a natural system on $\C{A}$, and
$(t,\bar{t})$ an endomorphism of $(\C{A},D)$ in $\C{nat}$. Such an
endomorphism is termed a \emph{translation}\index{translation} of
$(\C{A},D)$. The \emph{translation cohomology}\index{translation
category}
$$H^*(t,\bar{t})$$
is the cohomology of the homotopy fiber of the cochain homomorphism
$$\bar{t}_*-t^*\colon\bar{F}^*(\C{A},D)\To\bar{F}^*(\C{A},D\F(t)),$$
which will be denoted by $\bar{F}^*(t,\bar{t})$. In dimension $n$
this cochain complex is
$$\bar{F}^n(t,\bar{t})=\bar{F}^n(\C{A},D)\oplus\bar{F}^{n-1}(\C{A},D\F(t))$$
and the differential is the matrix
$$\left(%
\begin{array}{cc}
  \delta & 0 \\
  \bar{t}_*-t^* & -\delta \\
\end{array}%
\right).$$ In particular there is a long exact sequence $(n\in\Z)$
\begin{equation}\label{exacta}
\cdots\r H^n(t,\bar{t})\st{j}\r H^n(\C{A},D)\st{\bar{t}_*-t^*}\To
H^n(\C{A},D\F(t))\st{\partial}\r H^{n+1}(t,\bar{t})\r\cdots.
\end{equation}

\begin{rem}
This is a general concept of \emph{cohomology of diagrams} as for
example described in \cite{slds} or \cite{cfhh} A.2. Then
translation cohomology corresponds to the cohomology of the
diagram
$$
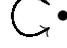
\xymatrix{\bullet\ar@(ul,dl)[]}$$
with one vertex and one arrow.
\end{rem}

Given functors $s\colon\C{B}\r\C{B}$ and $\varphi\colon\C{B}\r\C{A}$ such that the square
$$\xymatrix{\C{B}\ar[r]^s\ar[d]_\varphi&\C{B}\ar[d]^\varphi\\\C{A}\ar[r]_t&\C{A}}$$
commutes up to a given natural isomorphism $\alpha\colon t\varphi\rr\varphi s$
there is an induced translation
$$(s,(D\F(\alpha^\vi,\alpha))\vc(\bar{t}\F(\varphi)))\colon(\C{B},D\F(\varphi))\To(\C{B},D\F(\varphi))$$
and an induced cochain homomorphism
$$(\varphi,\alpha)^*\colon\bar{F}^*(t,\bar{t})\To\bar{F}^*(s,(D\F(\alpha^\vi,\alpha))\vc(\bar{t}\F(\varphi)))$$
given by the matrix
$$(\varphi,\alpha)^*=\left(\begin{array}{cc}
\varphi^*&0\\
  (D\F(\alpha^\vi,0^\vc_{\varphi s}))_*\alpha^\# & (D\F(\alpha^\vi,\alpha))_*\varphi^*
\end{array}%
\right).$$
Here we use (\ref{homo13}).
The induced homomorphisms in translation cohomology will be denoted in the same way
$$(\varphi,\alpha)^*\colon H^*(t,\bar{t})\To H^*(s,(D\F(\alpha^\vi,\alpha))\vc(\bar{t}\F(\varphi))).$$
Similarly two natural transformations $\lambda\colon D\rr E$ and $\bar{s}\colon E\rr E\F(t)$ between natural systems on $\C{A}$
such that $(\lambda\F(t))\vc\bar{t}=\bar{s}\vc\lambda$ give rise to a translation
$$(t,\bar{s})\colon(\C{A},E)\To(\C{A},E)$$
and a cochain homomorphism
$$\lambda_*\colon\bar{F}^*(t,\bar{t})\To\bar{F}^*(t,\bar{s})$$
given by the matrix
$$\lambda_*=\left(\begin{array}{cc}
\lambda_*&0\\
0 & (\lambda\F(t))_*
\end{array}%
\right).$$
The corresponding homomorphism in translation cohomology will also be denoted by
$$\lambda_*\colon H^*(t,\bar{t})\To H^*(t,\bar{s}).$$

These homomorphisms in translation cohomology can be used to describe the functorial properties of this cohomology.
More precisely, let $\C{Natt}$ be the category whose objects are all translations. A morphism $(\varphi,\alpha,\lambda)\colon
(t,\bar{t})\r(s,\bar{s})$ between two translations $(t,\bar{t})\colon(\C{A},D)\r(\C{A},D)$ and $(s,\bar{s})\colon(\C{B},E)\r(\C{B},E)$
is given by a functor $\varphi\colon\C{B}\r\C{A}$, a natural isomorphism $\alpha\colon t\varphi\rr\varphi s$, and a
natural transformation $\lambda\colon D\F(\varphi)\rr E$ such that $$(\lambda\F(s))\vc(D\F(\alpha^\vi,\alpha))\vc(\bar{t}\F(\varphi))=
\bar{s}\vc\lambda.$$
Composition of morphisms
$$(t,\bar{t})\st{(\varphi,\alpha,\lambda)}\To(s,\bar{s})\st{(\psi,\beta,\zeta)}\To(r,\bar{r})$$
in $\C{Natt}$ is defined by the formula
$$(\psi,\beta,\zeta)(\varphi,\alpha,\lambda)=(\varphi\psi,(\varphi\beta)\vc(\alpha\psi) ,\zeta\vc(\lambda\F(\psi))).$$

\begin{prop}\label{funw}
Translation cohomology groups are functors
$$H^*\colon\C{Natt}\To\C{Ab}$$ given on morphisms by $H^*(\varphi,\alpha,\lambda)=\lambda_*(\varphi,\alpha)^*$.
\end{prop}

This is just a simplified version of a more extended functoriality
of translation cohomology that we now describe. Suppose that we
have a translation $(t,\bar{t})\colon(\C{A},D)\r(\C{A},D)$, a
functor $s\colon\C{B}\r\C{B}$, a natural transformation
$\gamma\colon\varphi\rr\psi$ between functors
$\varphi,\psi\colon\C{B}\r\C{A}$, and two natural isomorphisms
$\alpha\colon t\varphi\rr\varphi s$ and $\beta\colon t\psi\rr\psi
s$ such that
\begin{equation}\label{laeq}
\beta\vc(t\gamma)=(\gamma s)\vc\alpha.
\end{equation}
In this situation there is an induced translation
$$(s,(D\F(\alpha^\vi,\beta))\vc(\bar{t}\F(\gamma)))\colon(\C{B},D\F(\gamma))\To(\C{B},D\F(\gamma))$$
and an induced cochain homomorphism
$$(\gamma,\alpha,\beta)^*\colon\bar{F}^*(t,\bar{t})\To\bar{F}^*(s,(D\F(\alpha^\vi,\beta))\vc(\bar{t}\F(\gamma)))$$
given by the matrix
$$(\gamma,\alpha,\beta)^*=\left(\begin{array}{cc}
\gamma^*&0\\
  (D\F(\alpha^\vi,\gamma s))_*\alpha^\# & (D\F(\alpha^\vi,\beta))_*\gamma^*
\end{array}%
\right).$$
The reader can check that this is indeed a cochain homomorphism by using (\ref{homo12'}), (\ref{homo13}) and (\ref{laeq}).

Let $\C{Natt}_1$ be the category whose objects are all
translations. A morphism $(\gamma,\alpha,\beta,\lambda)\colon
(t,\bar{t})\r(s,\bar{s})$ between two translations
$(t,\bar{t})\colon(\C{A},D)\r(\C{A},D)$ and
$(s,\bar{s})\colon(\C{B},E)\r(\C{B},E)$ is given by a natural
transformation $\gamma\colon\varphi\rr\psi$ between functors
$\varphi,\psi\colon\C{B}\r\C{A}$, two natural isomorphisms
$\alpha\colon t\varphi\rr\varphi s$ and $\beta\colon t\psi\rr\psi
s$ satisfying (\ref{laeq}), and a natural transformation
$\lambda\colon D\F(\gamma)\rr E$ such that
$$(\lambda\F(s))\vc(D\F(\alpha^\vi,\beta))\vc(\bar{t}\F(\gamma))=
\bar{s}\vc\lambda.$$ Composition of morphisms
$$(t,\bar{t})\st{(\gamma,\alpha,\beta,\lambda)}\To(s,\bar{s})\st{(\varepsilon,\mu,\nu,\zeta)}\To(r,\bar{r})$$
in $\C{Natt}_1$ with $\gamma\colon\varphi\rr\psi$ and $\varepsilon\colon\bar{\varphi}\rr\bar{\psi}$ is defined by the formula
$$(\varepsilon,\mu,\nu,\zeta)(\gamma,\alpha,\beta,\lambda)=
(\gamma\varepsilon,(\varphi\mu)\vc(\alpha\bar{\varphi}), (\psi\nu)\vc(\beta\bar{\psi}),
\zeta\vc(\lambda\F(\varepsilon))).$$
Compare with the definition of $\C{Cat}_\F$ in Section \ref{cc}.

\begin{prop}\label{fun}
Translation cohomology groups are functors
$$H^*\colon\C{Natt}_1\To\C{Ab}$$ given on morphisms by $H^*(\gamma,\alpha,\beta,\lambda)=\lambda_*(\gamma,\alpha,\beta)^*$.
\end{prop}

\begin{proof}
All we have to check is that composition is preserved. For this we are going to prove that the cochain homomorphisms
\begin{small}
$$\zeta_*(\varepsilon,\mu,\nu)^*\lambda_*(\gamma,\alpha,\beta)^*,
(\zeta\vc(\lambda\F(\varepsilon)))_*(\gamma\varepsilon,(\varphi\mu)\vc(\alpha\bar{\varphi}),(\psi\nu)\vc(\beta\bar{\psi}))^*
\colon\bar{F}^*(t,\bar{t})\r\bar{F}^*(s,\bar{s})$$
\end{small}
are homotopic. The first of these cochain homomorphisms is given by the matrix
$$\left(\begin{array}{cc}
a_{11}&a_{12}\\
a_{21} & a_{22}
\end{array}%
\right)$$
with entries
\begin{eqnarray*}a_{11}&=&
\zeta_*\varepsilon^*\lambda_*\gamma^*,\\a_{12}&=&0,\\
a_{21}&=&(\zeta\F(r))_*(E\F(\mu^\vi,\varepsilon r))_*\mu^\#\lambda_*\gamma^*\\&&+
(\zeta\F(r))_* (E\F(\mu^\vi,\nu))_*\varepsilon^*(\lambda\F(s))_*
(D\F(\alpha^\vi,\gamma s))_*\alpha^\#,\\
a_{22}&=&(\zeta\F(r))_*(E\F(\mu^\vi,\nu))_*\varepsilon^*(\lambda\F(s))_*(D\F(\alpha^\vi,\beta))_*\gamma^*,
\end{eqnarray*}%
and the second one by the matrix
$$\left(\begin{array}{cc}
b_{11}&b_{12}\\
b_{21} & b_{22}
\end{array}%
\right)$$
with
\begin{eqnarray*}
b_{11}&=&(\zeta\vc(\lambda\F(\varepsilon)))_*(\gamma\varepsilon)^*,\\
b_{12}&=&0,\\
b_{21}&=&((\zeta\vc(\lambda\F(\varepsilon)))\F(r))_*(D\F((\varphi\mu^\vi)\vc(\alpha^\vi\bar{\varphi}),0^\vc_{\gamma\varepsilon r}))_*
((\varphi\mu)\vc(\alpha\bar{\varphi}))^\#,\\
b_{22}&=&((\zeta\vc(\lambda\F(\varepsilon)))\F(r))_*(D\F((\varphi\mu^\vi)\vc(\alpha^\vi\bar{\varphi}),(\psi\nu)\vc(\beta\bar{\psi})))_*
(\gamma\varepsilon)^*.
\end{eqnarray*}

One can easily check that the equalities $a_{ij}=b_{ij}$ hold for $(i,j)\neq (2,1)$.
Moreover, by using (\ref{homo12'}), (\ref{homo14}) and (\ref{laeq}) the reader can check that
\begin{eqnarray*}
a_{21}&=&((\zeta\vc(\lambda\F(\varepsilon)))\F(r))_*(D\F(\varphi\mu^\vi,\gamma\varepsilon r))\mu^\#\varphi^*\\
&&+((\zeta\vc(\lambda\F(\varepsilon)))\F(r))_*(D((\alpha^\vi\bar{\varphi})\vc(\varphi\mu^\vi),(\gamma\varepsilon r)\vc(\varphi\mu))_*\bar{\varphi}^*\alpha^\#
\end{eqnarray*}
 Furthermore, by using (\ref{homo25}), (\ref{homo26}),
and
(\ref{homo13})  one can now check that the matrix
$$\left(\begin{array}{cc}
0&0\\
(\varphi\mu,\alpha\bar{\varphi})^\% & 0
\end{array}%
\right)$$
yields the desired cochain homotopy.
\end{proof}

Now Proposition \ref{funw} follows from Proposition \ref{fun} since there is a commutative diagram of functors
$$\xymatrix@R=10pt{\C{Natt}\ar[rd]^{H^*}\ar[dd]&\\&\C{Ab}\\\C{Natt}_1\ar[ru]_{H^*}}$$
where the vertical arrow is the identity on objects and sends a morphism $(\varphi,\alpha)$ to $(0^\vc_\varphi,\alpha,\alpha)$.
Compare with the definition of the functor $\C{Cat}\r\C{Cat}_\F$ in Section \ref{cc}.

The following identities between induced homomorphisms in translation cohomology will be useful for future applications.

\begin{prop}\label{uni}
Let be a morphism $(\gamma,\alpha,\beta,\lambda)\colon
(t,\bar{t})\r(s,\bar{s})$ between two translations $(t,\bar{t})\colon(\C{A},D)\r(\C{A},D)$ and $(s,\bar{s})\colon(\C{B},E)\r(\C{B},E)$
in $\C{Natt}_1$. The following homomorphisms $H^*(t,\bar{t})\r H^*(s,\bar{s})$ coincide
$$\lambda_*(\gamma,\alpha,\beta)^*=(\lambda\vc D(0^\vc_\varphi,\gamma))_*(\varphi,\alpha)^*=
(\lambda\vc D(\gamma,0^\vc_\psi))_*(\psi,\beta)^*$$
\end{prop}

\begin{proof}
The first equality can be easily checked at the level of cochains by using (\ref{homo12'}) and (\ref{laeq}). For the second one
we only need to use the cochain homotopy in the translation cochain complex given by the matrix
$$\left(\begin{array}{cc}
\lambda_*\gamma^\#&0\\
0&((\lambda\F(s))\vc D\F(\alpha^{-1},\beta))_*\gamma^\#
\end{array}\right),$$
see (\ref{homo13}).
\end{proof}

\begin{rem}
As the reader may suspect Proposition \ref{uni} is a consequence
of a higher functoriality of the translation cochain complex. More
precisely, the category  $\C{Natt}_1$ can be endowed with a
$2$-category structure by using factorizations, compare with the
definition of the $2$-category $\C{Cat}_\F$ in Section \ref{cc}.
The translation cochain complex turns out to be a  pseudofunctor
on the $2$-category $\C{Natt}_1$. This pseudofunctor can be
described by using the $2$-functorial properties of the
Baues-Wirsching cochain complex described in \cite{fcc}.
\end{rem}

If $\C{A}$ is additive and $t\colon\C{A}\r\C{A}$ is an additive
endofunctor we can consider the natural system $D$ given by the
$\C{A}$-bimodule $\hom_\C{A}(t,-)$ and the natural transformation
$\bar{t}\colon D=\hom_\C{A}(t,-)\rr \hom_\C{A}(t^2,t)=D\F(t)$
which is given by $(-1)t$. In this particular case we denote the
translation cohomology as follows
\begin{equation}
H^*(\C{A},t)=H^*(t,\bar{t}).
\end{equation}
This cohomology is of particular interest in case $\C{A}$ is a
triangulated category and $t$ its translation functor. This example,
in fact, is the motivation of the theory in this paper. We shall use
translation cohomology in order to study the ``characteristic
cohomology class'' in $H^3(\C{A},t)$ of a triangulated category
$(\C{A},t)$, see \cite{ccctc1} and \cite{ccctc2}.

\section{Classification of abelian track categories and translation track categories}\label{catctc}


Let $D$ be a natural system over $\C{C}$. We introduce the
categories $\C{Track}(\C{C},D)$ and $\C{Pseudo}_\simeq(\C{C},D)$.
Both categories have the same objects $(\C{A},\chi)$ which are
given by an abelian track category $\C{A}$ together with an
isomorphism $\chi=(\chi_0,\chi_1)\colon\pi\C{A}\cong(\C{C},D)$ in
$\C{nat}$. We often use $\chi$ as an identification and drop it
from notation. Morphisms
$\varphi\colon(\C{A},\chi)\r(\C{B},\kappa)$ in
$\C{Track}(\C{C},D)$ are track functors
$\varphi\colon\C{A}\r\C{B}$ such that $\chi=\kappa(\pi\varphi)$,
and morphisms $[\psi]\colon(\C{A},\chi)\r(\C{B},\kappa)$ in
$\C{Pseudo}_\simeq(\C{C},D)$ are homotopy classes of
pseudofunctors $[\psi]\colon\C{A}\pf\C{B}$ such that
$\chi=\kappa(\pi[\psi])$.

Given an ordinary category $\C{G}$ we write $\pi_0\C{G}$ for the
set of connected components in $\C{G}$. Two objects $X$, $Y$
belong to the same connected component iff there is a chain of
morphisms in $\C{G}$ as follows
$$X\l\bullet\r\cdots\l\bullet\r Y.$$
This equivalence relation simplifies when $\C{G}$ is a groupoid,
in this case two objects belong to the same connected component if
and only if they are isomorphic.

The next result classifies connected components in
$\C{Track}(\C{C},D)$ by the cohomology in Section \ref{cc}.

\begin{prop}[\cite{ccglg} 4.6, \cite{mhtcc}]
There is a bijection $$\pi_0\C{Track}(\C{C},D)\cong H^3(\C{C},D)$$
which carries $(\C{A},\chi)$ to the characteristic class
$\grupo{\C{A}}\in H^3(\C{C},D)$ given by the universal Toda
bracket of the linear track extension $(\C{A},\chi)$.
\end{prop}

A $3$-cocycle $c_\C{A}$ representative of the cohomology class
$\grupo{\C{A}}$ is defined in Section \ref{gs}.


The category $\C{Track}(\C{C},D)$ is far from being a groupoid,
however $\C{Pseudo}^{ab}_\simeq(\C{C},D)$ is indeed a groupoid as
a consequence of Theorem \ref{local}, therefore its set of
connected components can be described in a much easier way.
Moreover, we have the following result.

\begin{prop}\label{clasi1}
There is a bijection $$\pi_0\C{Pseudo}^{ab}_\simeq(\C{C},D)\cong
H^3(\C{C},D)$$ which carries $(\C{A},\chi)$ to $\grupo{\C{A}}$.
\end{prop}

\begin{proof}
This follows from the existence of the exact sequence for functors
(\ref{esf}), the definition of its obstruction operator $\theta$
in (\ref{teobs}), and \cite{ah} IV.4.12.
\end{proof}

A \emph{translation track category}\index{translation track
category} $(\C{A},[\varphi])$ is an abelian track category
equipped with a homotopy class of pseudofunctors
$[\varphi]\colon\C{A}\pf\C{A}$. Notice that $\pi[\varphi]$ is a
translation in the sense of Section \ref{ct}.



Let $(t,\bar{t})\colon(\C{C},D)\r(\C{C},D)$ be a translation. We
now introduce the category $\C{Trans}(t,\bar{t})$ whose objects
are triples $(\C{A},[\varphi],\chi)$ where $(\C{A},[\varphi])$ is
a translation track category and
$\chi\colon\pi\C{A}\cong(\C{C},D)$ is an isomorphism with
$(t,\bar{t})=\chi(\pi[\varphi])\chi^{-1}$. A morphism
$[\phi]\colon(\C{A},[\varphi],\chi)\r(\C{B},[\psi],\kappa)$ is a
homotopy class of pseudofunctors $[\phi]\colon\C{A}\pf\C{B}$ such
that $[\phi][\varphi]=[\psi][\phi]$ and $\chi=\kappa(\pi[\phi])$.
Again by Theorem \ref{local} we see that $\C{Trans}(t,\bar{t})$ is
in fact a groupoid. The abelian group $H^2(\C{C},D\F(t))$ acts on
this groupoid by using the action defined in Theorem \ref{2act}.
More precisely, this action is given by
$(\C{A},[\varphi],\xi)+\omega=(\C{A},[\varphi]+\omega,\xi)$ on
objects and the identity on morphisms. Moreover, there is a
functor
$$\bar{\jmath}\colon\C{Trans}(t,\bar{t})\To\C{Pseudo}_\simeq(\C{C},D),\;\;\bar{\jmath}(\C{A},[\varphi],\chi)=(\C{A},\chi).$$


\begin{thm}\label{clastrans}
There is a bijection $$\pi_0\C{Trans}(t,\bar{t})\cong
H^3(t,\bar{t})$$ where we use the translation cohomology in
Section \ref{ct}. This bijection is equivariant with respect to
the action of $H^2(\C{C},D\F(t))$ on $H^3(t,\bar{t})$ given by
$\partial$ in (\ref{exacta}). Moreover, if we use this bijection
and the bijection in Proposition \ref{clasi1} as identifications
then the following diagram commutes
$$\xymatrix{\pi_0\C{Trans}(t,\bar{t})\ar[r]^<(.15){\pi_0\bar{\jmath}}\ar@{=}[d]&\pi_0\C{Pseudo}^{ab}_\simeq(\C{C},D)\ar@{=}[d]\\
H^3(t,\bar{t})\ar[r]_j&H^3(\C{C},D)}$$
Here $j$ is the homomorphism in the exact sequence (\ref{exacta}).
\end{thm}

\begin{proof}
The characteristic class of a translation track category
$(\C{A},[\varphi])$ is the translation cohomology class
$$\grupo{\C{A},[\varphi]}\in H^3(\pi\varphi)$$
represented by the cocycle
$$(c_\C{A},b_\varphi).$$
Here we use the cochains defined in Section \ref{gs} obtained from
a fixed global section. This is indeed a cocycle by Proposition
\ref{fmla1} (1). Moreover, its cohomology class does not depend on
the choice of the global section by Proposition \ref{fmla3}.
Notice also that the equality
$j\grupo{\C{A},[\varphi]}=\grupo{\C{A}}$ holds from the very
definition of $j$ and these characteristic classes.

Let us check that the map
$$\C{Trans}(t,\bar{t})\To
H^3(t,\bar{t}),\;\;(\C{A},[\varphi],\chi)\mapsto
\grupo{\C{A},[\varphi]}$$ induces the desired bijection. This map
is equivariant with respect to the action of $H^2(\C{C},D\F(t))$,
see Proposition \ref{suma} and the proof of Theorem \ref{2act} in
Section \ref{p2a}. Let us check that it factors through the set of
connected components $\pi_0\C{Trans}(t,\bar{t})$. Let
$[\phi]\colon(\C{A},[\varphi],\chi)\r(\C{B},[\psi],\kappa)$ be an
isomorphism in $\C{Trans}(t,\bar{t})$ and
$\xi\colon\phi\varphi\rr\psi\phi$ a homotopy. By Proposition
\ref{fmla1} the coboundary of $(b_\phi,o_{\phi,\varphi}
-o_{\psi,\phi}+e_\xi)$ in the translation cochain complex is
$(c_\C{A}-c_\C{B},b_\varphi-b_\psi)$, therefore
$\grupo{\C{A},[\varphi]}=\grupo{\C{B},[\psi]}$.

We claim that the isotropy groups of the action of
$H^2(\C{C},D\F(t))$ on $\pi_0\C{Trans}(t,\bar{t})$ are the image
of $\bar{t}_*-t^*\colon H^2(\C{C},D)\r H^2(\C{C},D\F(t))$ in the
exact sequence (\ref{exacta}). By Theorem \ref{2act} two
translation track categories $(\C{A},[\varphi],\chi)$ and
$(\C{A},[\varphi]+\omega,\chi)$ are isomorphic in
$\C{Trans}(t,\bar{t})$ if and only if there exists $\varpi\in
H^2(\C{C},D)$ such that
$$(1_{(\C{A},\chi)}+\varpi)[\varphi]=([\varphi]+\omega)(1_{(\C{A},\chi)}+\varpi).$$
Now by the effectivity of the action in Theorem \ref{2act} and the
linear distributivity law in Proposition \ref{2ldl} we obtain that
$$\omega=(\bar{t}_*-t^*)(-\varpi).$$

In order to complete the proof it is enough to check that the
image of the map $\pi_0\bar{\jmath}$ coincides with the image of
$j$ under the identification in Proposition \ref{clasi1}. By
Proposition \ref{clasi1} any cohomology class $j(c)\in
H^3(\C{C},D)$ is the characteristic class of some $(\C{A},\chi)$
in $\C{Pseudo}^{ab}_\simeq(\C{C},D)$. We claim that there exists a
homotopy class $[\varphi]\colon\C{A}\pf\C{A}$ such that
$(\C{A},[\varphi],\chi)$ is a translation. By Theorem \ref{obs1}
the existence of such a $[\varphi]$ is equivalent to the vanishing
of $\theta_{\C{A},\C{A}}(t,\bar{t})$, and by applying
(\ref{teobs}) and the exactness of (\ref{exacta}) we obtain
$$\theta_{\C{A},\C{A}}(t,\bar{t})=\bar{t}_*j(c)-t^*j(c)=0.$$

\end{proof}

\section{Obstruction theory for the $2$-functor $\pi$}\label{os}


This section is devoted to the realizability of $0$-, $1$- and $2$-cells through the $2$-functor
$$\pi\colon\C{Pseudo}^{ab}_\simeq\;\To\C{nat}$$
in (\ref{pi2}). All $0$-cells are known to be realizable since
there is always the trivial linear track extension for
$(\C{C},D)$, see \cite{ccglg}. We here define classes in
cohomology of categories which represent obstructions to the
realizability of $1$- and $2$-cells. We also consider their
behaviour with respect to the various composition laws.

Given two abelian track categories $\C{A}$, $\C{B}$ and a morphism $(\varphi,\lambda)\colon\pi\C{A}\r\pi\C{B}$
in $\C{nat}$ we define the cohomology class
\begin{equation}\label{teobs}
\theta_{\C{A},\C{B}}(\varphi,\lambda)=\lambda_*\grupo{\C{A}}-\varphi^*\grupo{\C{B}}\in H^3(\C{A},(\pi_1\C{B})\F(\varphi)).
\end{equation}

\begin{thm}\label{obs1}
The element $\theta_{\C{A},\C{B}}(\varphi,\lambda)$ vanishes if and only if there exists a pseudofunctor
$\psi\colon\C{A}\pf\C{B}$ such that $\pi\psi=(\varphi,\lambda)$.
\end{thm}

The proof of this result will be given in Section \ref{81}.

\begin{prop}\label{1ldl}
Given abelian track categories $\C{A}$, $\C{B}$, $\C{C}$ and
morphisms $\pi\C{A}\st{(\varphi,\lambda)}\r
\pi\C{B}\st{(\psi,\zeta)}\r\pi\C{C}$ in $\C{nat}$ the following
formula holds
$$\theta_{\C{A},\C{C}}((\psi,\zeta)(\varphi,\lambda))=(\zeta\F(\varphi))_*\theta_{\C{A},\C{B}}(\varphi,\lambda)
+\varphi^*\theta_{\C{B},\C{C}}(\psi,\zeta).$$
\end{prop}

This proposition follows from (\ref{teobs}) and the composition
law for $1$-cells in $\C{nat}$.

\begin{prop}\label{con2}
Given abelian track categories $\C{A}$, $\C{B}$ and a $2$-cell
$\alpha\colon(\varphi,\lambda)\rr(\psi,\zeta)$ in $\C{nat}$
between $1$-cells $(\varphi,\lambda),(\psi,\zeta)\colon\pi\C{A}\r
\pi\C{B}$ the following equality holds
$$((\pi_1\C{B})\F(0^\vc_{\pi_0\varphi},\alpha))_*\theta_{\C{A},\C{B}}(\varphi,\lambda)=
((\pi_1\C{B})\F(\alpha,0^\vc_{\pi_0\psi}))_*\theta_{\C{A},\C{B}}(\psi,\zeta).$$
\end{prop}

This proposition is a consequence of (\ref{teobs}), (\ref{hdc2})
and (\ref{homo13}).

\begin{prop}\label{actob}
Given an abelian track category $\C{A}$ and $\omega\in
H^3(\pi_0\C{A},\pi_1\C{A})$ there exists another abelian track
category $\C{B}$ with $\pi\C{A}=\pi\C{B}$ such that
$\theta_{\C{A},\C{B}}(1_{\pi\C{A}})=\omega$. This abelian track
category is unique up to isomorphism in $\C{Pseudo}^{ab}_\simeq$
and we write $\C{A}=\C{B}+\omega$.
\end{prop}

This result is a consequence of the definition of $\theta$ and \cite{catc} 3.3.

Let $\varphi,\psi\colon\C{A}\pf\C{B}$ be pseudofunctors between
small abelian track categories and
$\alpha\colon\pi\varphi\rr\pi\psi$ a $2$-cell in $\C{nat}$. We
define the element
$$\vartheta_{\varphi,\psi}(\alpha)\in H^2(\pi_0\C{A},(\pi_1\C{B})\F(\alpha))$$
as the cohomology class represented by
\begin{equation}\label{teobs2}
-((\pi_1\C{B})\F(0^\square_{\pi_0\varphi},\alpha))_*b_\varphi-\alpha^\#c_\C{B}
+((\pi_1\C{B})\F(\alpha,0^\square_{\pi_0\psi}))_*b_\psi.
\end{equation}
Here we use the cochains defined in Section \ref{gs} for a fixed
global section. This is indeed a cocycle by Proposition
\ref{fmla2}. Moreover, by using Proposition \ref{fmla4} we see
that its cohomology class does not depend on the global section
chosen for its definition.

\begin{thm}\label{obs2}
The element $\vartheta_{\varphi,\psi}(\alpha)$ vanishes if and only if there exists a pseudonatural transformation $\tau\colon\varphi\rr\psi$
such that $\pi_0\tau=\alpha$.
\end{thm}

This proposition will be proved in Section \ref{84}.

\begin{prop}\label{homotono}
Two pseudofunctors between abelian track categories
$\varphi,\psi\colon\C{A}\pf\C{B}$ are homotopic if and only if
$\pi\varphi=\pi\psi$ and
$\vartheta_{\varphi,\psi}(0^\vc_{\pi_0\varphi})=0$.
\end{prop}

\begin{proof}
Obviously if $\xi\colon\varphi\rr\psi$ is a homotopy then
$\pi\varphi=\pi\psi$ and $\pi_0\xi=0^\vc_{\pi_0\varphi}$, hence by
Theorem \ref{obs2}
$\vartheta_{\varphi,\psi}(0^\vc_{\pi_0\varphi})=0$. On the other
hand if we suppose that $\pi\varphi=\pi\psi$ and
$\vartheta_{\varphi,\psi}(0^\vc_{\pi_0\varphi})=0$ the
pseudonatural transformation $\tau\colon\varphi\rr\psi$ with
$\pi_0\tau=0^\vc_{\pi_0\varphi}$ constructed in the proof of
Theorem \ref{obs2}, given in Section \ref{84}, satisfies
$\tau_X=t(0^\vc_{\pi_0\varphi})_X=t\set{1_{\varphi(X)}}=1_{\varphi(X)}$,
here we use (\ref{gs1}), therefore $\tau$ is a homotopy.
\end{proof}

This result proves that, if we restrict to the category of abelian
track categories, we could have defined homotopies as
pseudonatural transformations $\xi\colon\varphi\rr\psi$ such that
$\pi_0\xi$ is and identity natural transformation, obtaining in
this way the same homotopy category $\C{Pseudo}^{ab}_\simeq$.
Therefore homotopies in the sense of Section \ref{homocat} would
be a reduced version of these more general ones, however the
existence of $\C{Pseudo}^{ab}_\simeq$ is not so easy to establish
if we try to use non-reduced homotopies from the beginning.

\begin{prop}\label{deriv2}
Given pseudofunctors $\phi,\varphi,\psi\colon\C{A}\pf\C{B}$ and
$2$-cells $\pi\phi\st{\beta}\rr\pi\varphi\st{\alpha}\rr\pi\psi$ in
$\C{nat}$ the following equality holds
$$\vartheta_{\phi,\psi}(\alpha\vc\beta)=((\pi_1\C{B})\F(0^\vc_{\pi_0\phi},\alpha))_*\vartheta_{\phi,\varphi}(\beta)+
((\pi_1\C{B})\F(\beta,0^\vc_{\pi_0\psi}))_*\vartheta_{\varphi,\psi}(\alpha).$$
\end{prop}

\begin{proof}
The result follows from the equations:
\begin{eqnarray*}
((\pi_1\C{B})\F(0^\square_{\pi_0\phi},\alpha))_*(-((\pi_1\C{B})\F(0^\square_{\pi_0\phi},\beta))_*b_\phi&&\\-\beta^\#c_\C{B}
+((\pi_1\C{B})\F(\beta,0^\square_{\pi_0\varphi}))_*b_\varphi)&&\\
+((\pi_1\C{B})\F(\beta,0^\square_{\pi_0\psi}))_*(-((\pi_1\C{B})\F(0^\square_{\pi_0\varphi},\alpha))_*b_\varphi&&\\-\alpha^\#c_\C{B}
+((\pi_1\C{B})\F(\alpha,0^\square_{\pi_0\psi}))_*b_\psi)&=&
-((\pi_1\C{B})\F(0^\square_{\pi_0\phi},\alpha\vc\beta))_*b_\phi\\&&-((\pi_1\C{B})\F(0^\square_{\pi_0\phi},\alpha))_*\beta^\#c_\C{B}\\&&
+((\pi_1\C{B})\F(\beta,\alpha))_*b_\varphi\\&&-((\pi_1\C{B})\F(\beta,\alpha))_*b_\varphi\\&&
-((\pi_1\C{B})\F(\beta,0^\square_{\pi_0\psi}))_*\alpha^\#c_\C{B}\\&&
+((\pi_1\C{B})\F(\alpha\vc\beta,0^\square_{\pi_0\psi}))_*b_\psi\\&=&
-((\pi_1\C{B})\F(0^\square_{\pi_0\phi},\alpha\vc\beta))_*b_\phi\\&&
-(\alpha\vc\beta)^\#c_\C{B}\\&&
+((\pi_1\C{B})\F(\alpha\vc\beta,0^\square_{\pi_0\psi}))_*b_\psi\\&&
+\delta(\alpha,\beta)^\% c_\C{B}.
\end{eqnarray*}
Here we use (\ref{homo21}) and the fact that $c_\C{B}$ is a cocycle.
\end{proof}

\begin{prop}\label{deriv3}
Given pseudofunctors $\varphi,\varphi'\colon\C{B}\pf\C{A}$,
$\psi,\psi'\colon\C{C}\pf\C{B}$ and $2$-cells
$\alpha\colon\pi\varphi\rr\pi\varphi'$,
$\beta\colon\pi\psi\rr\pi\psi'$ in $\C{nat}$ the following
equalities hold
\begin{eqnarray*}
\vartheta_{\varphi\psi,\varphi'\psi'}(\alpha\beta)&=&
((((\pi_1\C{A})\F(0^\vc_{\pi_0\varphi},\alpha))\vc\pi_1\varphi)\F(\beta))_*\vartheta_{\psi,\psi'}(\beta)+
\beta^*\vartheta_{\varphi,\varphi'}(\alpha)
\\&=&((((\pi_1\C{A})\F(\alpha,0^\vc_{\pi_0\varphi'}))\vc\pi_1\varphi')\F(\beta))_*\vartheta_{\psi,\psi'}(\beta)+
\beta^*\vartheta_{\varphi,\varphi'}(\alpha).
\end{eqnarray*}
\end{prop}

\begin{proof}
The result follows from the next equalities
\begin{small}
\begin{eqnarray*}
((((\pi_1\C{A})\F(0^\vc_{\pi_0\varphi},\alpha))\vc\pi_1\varphi)\F(\beta))_*(&&\\
-((\pi_1\C{B})\F(0^\square_{\pi_0\psi},\beta))_*b_\psi&&\\-\beta^\#c_\C{B}
+((\pi_1\C{B})\F(\beta,0^\square_{\pi_0\psi'}))_*b_{\psi'})&&\\
+\beta^*(-((\pi_1\C{A})\F(0^\square_{\pi_0\varphi},\alpha))_*b_\varphi&&\\
-\alpha^\#c_\C{A}
+((\pi_1\C{A})\F(\alpha,0^\square_{\pi_0\varphi'}))_*b_{\varphi'})
&\st{\mathrm{(a)}}=&
-((\pi_1\C{A})\F(0^\vc_{\pi_0\varphi\psi},\alpha\beta))_*((\pi_1\varphi)\F(\pi_0\psi))_*b_\psi\\&&
-\beta^\#((\pi_1\C{A})\F(0^\vc_{\pi_0\varphi},\alpha))_*(\pi_1\varphi)_*c_\C{B}\\&&
+((\pi_1\C{A})\F(\alpha\beta,0^\vc_{\pi_0\varphi'\psi'}))_*((\pi_1\varphi')\F(\pi_0\psi'))_*b_{\psi'}\\&&
-((\pi_1\C{A})\F(0^\vc_{\pi_0\varphi\psi},\alpha\beta))_*(\pi_0\varphi)^*b_\varphi\\&&
-\beta^*\alpha^\#c_\C{A}\\&&
+((\pi_1\C{A})\F(\alpha\beta,0^\vc_{\pi_0\varphi'\psi'}))_*(\pi_0\psi')^*b_{\varphi'}\\&&
+\delta\beta^\#((\pi_1\C{A})\F(\alpha,0^\square_{\pi_0\varphi'}))_*b_{\varphi'}
\\&&+\beta^\#((\pi_1\C{A})\F(\alpha,0^\square_{\pi_0\varphi'}))_*\delta b_{\varphi'}\\
&\st{\mathrm{(b)}}=&
-(\alpha\beta)^*(b_{\varphi\psi}-\delta(o_{\varphi,\psi}))\\&&
-\beta^\#((\pi_1\C{A})\F(0^\vc_{\pi_0\varphi},\alpha))_*(\pi_1\varphi)_*c_\C{B}\\&&
-\beta^*\alpha^\#c_\C{A}\\&&
+\beta^\#((\pi_1\C{A})\F(\alpha,0^\square_{\pi_0\varphi'}))_*((\pi_1\varphi')_*c_\C{B}\\&&-(\pi_0\varphi')^*c_\C{A})\\&&
+((\pi_1\C{A})\F(\alpha\beta,0^\vc_{\pi_0\varphi'\psi'}))_*(b_{\varphi'\psi'}-\delta(o_{\varphi',\psi'}))\\&&
+\delta\beta^\#((\pi_1\C{A})\F(\alpha,0^\square_{\pi_0\varphi'}))_*b_{\varphi'}\\
&\st{\mathrm{(c)}}=&
-(\alpha\beta)^*b_{\varphi\psi}-\beta^*\alpha^\#c_\C{A}-\beta^\#\alpha^*c_\C{A}\\&&
+((\pi_1\C{A})\F(\alpha\beta,0^\vc_{\pi_0\varphi'\psi'}))_*b_{\varphi'\psi'}\\&&
+\delta((\alpha\beta)^*o_{\varphi,\psi}-((\pi_1\C{A})\F(\alpha\beta,0^\vc_{\pi_0\varphi'\psi'}))_*o_{\varphi',\psi'}
\\&&+\beta^\#((\pi_1\C{A})\F(\alpha,0^\square_{\pi_0\varphi'}))_*b_{\varphi'})\\&\st{\mathrm{(d)}}=&
-(\alpha\beta)^*b_{\varphi\psi}-(\alpha\beta)^\#c_\C{A}\\&&
+((\pi_1\C{A})\F(\alpha\beta,0^\vc_{\pi_0\varphi'\psi'}))_*b_{\varphi'\psi'}\\&&
+\delta((\alpha\beta)^*o_{\varphi,\psi}-((\pi_1\C{A})\F(\alpha\beta,0^\vc_{\pi_0\varphi'\psi'}))_*o_{\varphi',\psi'}
\\&&+\beta^\#((\pi_1\C{A})\F(\alpha,0^\square_{\pi_0\varphi'}))_*b_{\varphi'}+(\alpha,\beta)^\$ c_\C{A})
\end{eqnarray*}
\end{small}
Here for (a) we use (\ref{hdc2}), (\ref{homo12'}), (\ref{homo13}) and (\ref{homo14}); in (b) we apply
(\ref{homo12'}) as well as (\ref{fmla1}) (1) and (3); for (c) we use (\ref{hdc2}) and (\ref{homo12'}); and finally
for (d) we apply (\ref{homo22}) and the fact that $c_\C{A}$ is a cocycle.
\end{proof}

\begin{prop}\label{factorhomo}
Given pseudofunctors
$\varphi,\varphi',\psi,\psi'\colon\C{A}\pf\C{B}$ between abelian
track categories such that $\varphi\simeq\varphi'$,
$\psi\simeq\psi'$ and a $2$-cell
$\alpha\colon\pi\varphi=\pi\varphi'\rr\pi\psi=\pi\psi'$ in
$\C{nat}$ then
$\vartheta_{\varphi,\psi}(\alpha)=\vartheta_{\varphi',\psi'}(\alpha)$.
\end{prop}

\begin{proof}
This is a consequence of the following equalities
\begin{eqnarray*}
\vartheta_{\varphi,\psi}(\alpha)&=&((\pi_1\C{B})\F(\alpha,0^\vc_{\pi_0\psi}))_*\vartheta_{\psi',\psi}(0^\vc_{\pi_0\psi})
\\&&+\vartheta_{\varphi',\psi'}(\alpha)
+((\pi_1\C{B})\F(0^\vc_{\pi_0\varphi},\alpha))_*\vartheta_{\varphi,\varphi'}(0^\vc_{\pi_0\varphi})\\
&=&\vartheta_{\varphi',\psi'}(\alpha).
\end{eqnarray*}
In the first equality we use Proposition \ref{deriv2} applied to
the vertical composite $\alpha=0^\vc_{\pi_0\psi}\vc\alpha\vc
0^\vc_{\pi_0\varphi}$, for the second one we use Proposition
\ref{homotono}.
\end{proof}

This proposition proves that the cohomology class $\vartheta_{\varphi,\psi}(\alpha)$ only depends on the homotopy classes of
$\varphi$ and $\psi$, therefore we can denote it by
$$\vartheta_{[\varphi],[\psi]}(\alpha).$$

\begin{thm}\label{2act}
Given a homotopy class of pseudofunctors between abelian track categories $[\varphi]\colon\C{A}\pf\C{B}$
and $\omega\in H^2(\pi_0\C{A},(\pi_1\C{B})\F(\pi_0[\varphi]))$ there exists a unique homotopy class $[\varphi]+\omega\colon\C{A}\pf\C{B}$
such that $\vartheta_{[\varphi]+\omega,[\varphi]}(0^\vc_{\pi_0[\varphi]})=\omega$. This induces a transitive and effective action
of the abelian group $H^2(\pi_0\C{A},(\pi_1\C{B})\F(\pi_0[\varphi]))$ on the morphism set $\pi^{-1}(\psi,\lambda)\subset\C{Pseudo}^{ab}_\simeq
(\C{A},\C{B})$ for any $1$-cell in
$(\psi,\lambda)\colon\pi\C{A}\r\pi\C{B}$ in $\C{nat}$ such that $\theta_{\C{A},\C{B}}(\psi,\lambda)=0$.
\end{thm}

This proposition will be proved in Section \ref{p2a}.

\begin{prop}\label{2ldl}
The action in Theorem \ref{2act} satisfies the following linear
distributivity law: given composable morphisms
$\C{C}\st{[\psi]}\r\C{B}\st{[\varphi]}\r\C{A}$ in
$\C{Pseudo}^{ab}_\simeq$ and cohomology classes $\varpi\in
H^2(\pi_0\C{C},(\pi_1\C{B})\F(\pi_0[\psi]))$ and $\omega\in
H^2(\pi_0\C{B},(\pi_1\C{A})\F(\pi_0[\varphi]))$ the equality
$$([\varphi]+\omega)([\psi]+\varpi)=[\varphi][\psi]+((\pi_1[\varphi])\F(\pi_0[\psi]))_*\varpi+(\pi_0[\psi])^*\omega$$
is satisfied.
\end{prop}

\begin{proof}
The formula is a consequence of the following equalities,
\begin{eqnarray*}
\vartheta_{([\varphi]+\omega)([\psi]+\varpi),[\varphi][\psi]}(0^\vc_{\pi_0[\varphi][\psi]})&=&
((\pi_1[\varphi])\F(\pi_0[\psi]))_*\vartheta_{[\psi]+\varpi,[\psi]}(0^\vc_{\pi_0[\psi]})
\\&&+(\pi_0[\psi])^*\vartheta_{[\varphi]+\omega,[\varphi]}(0^\vc_{\pi_0[\varphi]})\\
&=&((\pi_1[\varphi])\F(\pi_0[\psi]))_*\varpi+(\pi_0[\psi])^*\omega.
\end{eqnarray*}
Here we apply Proposition \ref{deriv3} to the composite $0^\vc_{\pi_0[\varphi][\psi]}=0^\vc_{\pi_0[\varphi]}0^\vc_{\pi_0[\psi]}$.
\end{proof}

\section{The homotopy category of pseudofunctors as a linear extension}\label{hcple}

An \emph{exact sequence for a functor}\index{exact sequence for a
functor} $\lambda$
$$D\st{+}\r\C{A}\st{\lambda}\To\C{B}\st{\theta}\r H$$
is given by two natural systems $D$, $H$ on $\C{B}$ together with
an obstruction operator $\theta$ which determines whether a
morphism in $\C{B}$ is realizable through $\lambda$ or not, and an
action $+$ counting all possible realizations. We refer the reader
to \cite{ah} IV.4.10 for an explicit definition. \emph{Linear
extensions of categories}\index{linear extension of categories},
in the sense of \cite{ah} IV.3.2,
$$D\st{+}\rightarrowtail\C{A}\st{\lambda}\To\C{B}$$
can be regarded as special cases of exact sequences for functors
with $H=0$ and an effective action $+$ of $D$. A weak linear
extension of categories occurs when we change $\C{B}$ by an
equivalent category in a linear extension as above.

As a consequence of the results in Section \ref{os} the $2$-functor $\pi$ in (\ref{pi2})
regarded as an ordinary functor embeds in an exact sequence. More precisely,
we can consider the natural system $H^n$ $(n\geq 0)$ on $\C{nat}$ defined by the $n^\mathrm{th}$-cohomology of categories
$H^n(\C{C},D)$. The homomorphisms induced by a $1$-cell $(\varphi,\lambda)$ in $\C{nat}$ are $(\varphi,\lambda)_*=\lambda_*$
and $(\varphi,\lambda)^*=\varphi^*$. This is indeed a natural system because of the functorial properties of the cohomology of
categories. The exact sequence for $\pi$
\begin{equation}\label{esf}
H^2\st{+}\r\C{Pseudo}^{ab}_\simeq\st{\pi}\To\C{nat}\st{\theta}\r H^3
\end{equation}
is given by the obstruction operator $\theta$ in (\ref{teobs}) and
the action $+$ in Theorem \ref{2act}. Theorems \ref{obs1} and
\ref{2act} and Propositions \ref{1ldl}, \ref{factorhomo} and
\ref{2ldl} prove that the axioms of an exact sequence for a
functor are satisfied in this case.

In fact the theory developed in Section \ref{os} is much richer, because it also involves an obstruction theory for the realization of $2$-cells
through the $2$-functor $\pi$. However we have not developed general axioms for this sort of
exact sequence for a $2$-functor because in this paper we concentrate
only in the particular case of $\pi$.

The obstruction operator in (\ref{esf}) is an inner derivation in
the sense of \cite{ah} IV.7.3, the action $+$ is effective by
Theorem \ref{2act}, and any object in $\C{nat}$ is isomorphic to
an object in the image of $\pi$ by Proposition \ref{clasi1},
therefore we can obtain a weak linear extension from the exact
sequence (\ref{esf}) in the following canonical way: consider the
category $\C{CoH}$ whose objects are triples $(\C{C}, D, c)$ where
$\C{C}$ is a category, $D$ is a natural system on $\C{C}$ and
$c\in H^3(\C{C},D)$. A morphism
$(\varphi,\lambda)\colon(\C{C},D,c)\r(\C{D},E,d)$ in $\C{CoH}$ is
a morphism $(\varphi,\lambda)\colon(\C{C},D)\r(\C{D},E)$ in
$\C{nat}$ which is compatible with the cohomology classes, that is
$\lambda_*c=\varphi^*d$. By Theorem \ref{obs1} the functor $\pi$
extends to a functor
$$\bar{\pi}\colon\C{Pseudo}^{ab}_\simeq\To\C{CoH},\;\;\bar{\pi}(\C{A})=(\pi_0\C{A},\pi_1\C{A},\grupo{\C{A}}).$$
This functor fits into a weak linear extension
\begin{equation}\label{lec}
H^2\st{+}\rightarrowtail\C{Pseudo}^{ab}_\simeq\st{\bar{\pi}}\To\C{CoH}.
\end{equation}
Here $+$ is the same action as in (\ref{esf}), and $H^2$ is regarded as a natural system on $\C{CoH}$ through
the obvious forgetful functor $\C{CoH}\r\C{nat}$. 

\section{A $2$-category in terms of cocycles}\label{c2c}

In this section we describe an alternative and more algebraic
model of the $2$-category $\C{Pseudo}^{ab}_\simeq$ only in terms
of the Baues-Wirsching cochain complex as a $2$-functor. For this
we will use the exact sequence for functors (\ref{esf})
constructed in Section \ref{hcple}.

We define the $2$-category $\C{CoC}$ as follows: $0$-cells  $(\C{C},D,c)$ are triples where
$\C{C}$ is a category $D$ is a natural system on $\C{C}$ and $c\in\bar{F}^3(\C{C},D)$ is a $3$-cocycle, that is $\delta c=0$.
A $1$-cell $(\varphi,\lambda,e)\colon(\C{C},D,c)\r(\C{D},E,d)$ is given by a morphism $(\varphi,\lambda)\colon(\C{C},D)\r(\C{D},E)$
in $\C{nat}$ together with a $2$-cochain modulo coboundaries
$$e\in\frac{\bar{F}^2(\C{C},E\F(\varphi))}{\im\delta}$$ satisfying $\delta e=\lambda_*c-\varphi^*d$.
A $2$-cell $\alpha\colon(\varphi,\lambda,e)\rr(\psi,\zeta,b)$
between $1$-cells
$$(\varphi,\lambda,e),(\psi,\zeta,b)\colon(\C{C},D,c)\r(\C{D},E,d)$$
is just a $2$-cell $\alpha\colon(\varphi,\lambda)\rr(\psi,\zeta)$
in $\C{nat}$ such that
$$0=-(E\F(0^\vc_\varphi,\alpha))_*e-\alpha^\#d+(E\F(\alpha,0^\vc_\psi))_*b\in
H^2(\C{C},E\F(\alpha))\subset\frac{\bar{F}^2(\C{C},E\F(\alpha))}{\im\delta}.$$
One can check as in Proposition \ref{fmla2} that this linear
combination of $2$-cochains modulo coboundaries lies in the
cohomology group. Composition of $1$-cells
is determined by the formula
$$(\varphi,\lambda,e)(\psi,\zeta,b)=(\varphi\psi,(\lambda\F(\psi))\vc\zeta,(\lambda\F(\psi))_*b+\psi^*e).$$
Composition laws for $2$-cells are as in $\C{nat}$. It is not trivial at all to check that the horizontal or vertical composition
of two $2$-cells in $\C{CoC}$ is indeed in this $2$-category, however the necessary computations are similar to those in the proofs of
Propositions \ref{deriv2} and \ref{deriv3}, and for this reason we do not include them here.

There is an obvious $2$-functor
\begin{equation*}
\wp\colon\C{CoC}\To\C{nat},\;\;\wp(\C{C},D,c)=(\C{C},D),
\end{equation*}
which is faithful in dimension $2$, as $\pi$ is, see Proposition \ref{faith}. This functor also embeds in an exact sequence
\begin{equation}\label{esf2}
H^2\st{+}\r\C{CoC}\st{\wp}\To\C{nat}\st{\theta}\r H^3.
\end{equation}
The obstruction operator is defined as follows
$$\theta_{(\C{C},D,c),(\C{D},E,d)}(\varphi,\lambda)=\lambda_*\set{c}-\varphi^*\set{d}.$$
Here $\set{\cdot}$ denotes the cohomology class of the cocycle inside. Moreover, the action $+$ is defined as follows,
given a $1$-cell $(\varphi,\lambda,e)\colon(\C{C},D,c)\r(\C{D},E,d)$ in $\C{CoC}$ and $\set{s}\in H^2(\C{C},E\F(\varphi))$
$$(\varphi,\lambda,e)+\set{s}=(\varphi,\lambda,e+s).$$
Notice that this action is effective. The $2$-functor $\wp$ lifts
to $\C{CoH}$ through the forgetful functor $\C{CoH}\r\C{nat}$
giving rise to a weak linear extension of categories
\begin{equation}\label{leci}
H^2\st{+}\rightarrowtail\C{CoC}\st{\bar{\wp}}\To\C{CoH},
\end{equation}
where $\bar{\wp}(\C{C},D,c)=(\C{C},D,\set{c})$.

Proposition \ref{fmla1} proves that there is a $2$-functor
\begin{equation*}
\mathcal{C}\colon\C{Pseudo}^{ab}_\simeq\To\C{CoC},
\end{equation*}
such that the image of an abelian track category $\C{A}$ is
$\mathcal{C}(\C{A})=(\pi_0\C{A},\pi_1\C{A},c_\C{A})$, and given a
pseudofunctor $\varphi\colon\C{A}\pf\C{B}$ between abelian track
categories $\mathcal{C}(\varphi)=(\pi_0\varphi,
\pi_1\varphi,b_\varphi)$. Here we use the cochains defined in
Section \ref{gs} for a fixed global section, but notice that
Proposition \ref{fmla3} proves that up to $2$-natural isomorphism
this $2$-functor does not depend on that choice.

By Proposition \ref{suma} the $2$-functor $\mathcal{C}$
induces a map of exact sequences from (\ref{esf}) to (\ref{esf2}) in the sense of \cite{ah} 4.13
\begin{equation*}
\xymatrix{H^2\ar[r]^<(.25)+\ar@{=}[d]&\C{Pseudo}^{ab}_\simeq\ar[r]^<(.3)\pi\ar[d]^{\mathcal{C}}&\C{nat}\ar[r]^\theta\ar@{=}[d]&H^3\ar@{=}[d]\\
H^2\ar[r]^<(.35)+&\C{CoC}\ar[r]^\wp&\C{nat}\ar[r]^\theta&H^3}
\end{equation*}
hence by \cite{ah} 4.14 $\mathcal{C}$ is an equivalence of
categories regarded as an ordinary functor. Furthermore, the
definition of $2$-cells in $\C{CoC}$ together with (\ref{teobs2}),
Proposition \ref{faith} and Theorem \ref{obs2} prove that
$\mathcal{C}$ is fully faithful in dimension $2$, hence it is not
only an ordinary equivalence but an equivalence of enriched
categories.

\section{Proof of Theorem \ref{obs1}}\label{81}

The following result is a consequence of Proposition \ref{fmla1} (1).

\begin{prop}\label{es01}
If $\psi\colon\C{A}\pf\C{B}$ is a pseudofunctor between abelian
track categories then $\theta_{\C{A},\C{B}}(\pi\psi)=0$.
\end{prop}

In order to prove the converse we introduce the pull-back and push-forward of an abelian track category.

Let $\C{A}$ be an abelian track category and
$\varphi\colon\C{C}\r\pi_0\C{A}$ an ordinary functor, the
\emph{pull-back} of $\C{A}$ along $\varphi$ is the abelian track
category $\varphi^*\C{A}$ defined as follows. The category
$(\varphi^*\C{A})_0$ is the following pull-back in $\C{Cat}$
\begin{equation}\label{pulldiag}
\xymatrix{(\varphi^*\C{A})_0\ar[r]^<(.4){p'}\ar@{}[rd]|{\text{pull}}\ar[d]_{\varphi'}&\C{C}\ar[d]^\varphi\\\C{A}_0\ar[r]_p&\pi_0\C{A}}
\end{equation}
where $p$ is the quotient functor. Track sets in $\varphi^*\C{A}$ are given by the formula
\begin{equation}\label{pulltrack}
\homm{X,Y}_{\varphi^*\C{A}}(f,g)=\homm{\varphi'(X),\varphi'(Y)}_{\C{A}}(\varphi'(f),\varphi'(g)).
\end{equation}
This pull-back comes equipped with a track functor
$$\bar{\varphi}\colon\varphi^*\C{A}\To\C{A}$$
which is given by $\varphi'$ on objects and maps, and the equality
(\ref{pulltrack}) on tracks. This track functor determines
composition laws in $\varphi^*\C{A}$ involving tracks.

\begin{prop}\label{pull}
The pull back $\varphi^*\C{A}$ satisfies the equalities
$\pi_0\varphi^*\C{A}=\C{C}$ and
$\pi_1\varphi^*\C{A}=(\pi_1\C{A})\F(\varphi)$. Moreover, the track
functor $\bar{\varphi}$ induces $\pi_0\bar{\varphi}=\varphi$ and
$\pi_1\bar{\varphi}=0^\vc_{(\pi_1\C{A})\F(\varphi)}$.
\end{prop}

\begin{proof}
The natural projection $(\varphi^*\C{A})_0\r\C{C}$ is $p'$ in (\ref{pulldiag}), the equality $\pi_0\bar{\varphi}=\varphi$ follows
from the commutativity of (\ref{pulldiag}), and the identities $\pi_1\varphi^*\C{A}=(\pi_1\C{A})\F(\varphi)$ and
$\pi_1\bar{\varphi}=0^\vc_{(\pi_1\C{A})\F(\varphi)}$ are direct consequences of (\ref{pulltrack}).
\end{proof}

If $\C{A}$ is an abelian track category and
$\lambda\colon\pi_1\C{A}\rr D$ is a natural transformation between
natural systems on $\pi_0\C{A}$ the \emph{push-forward} of $\C{A}$
along $\lambda$ is an abelian track category $\lambda_*\C{A}$ with
the same objects and maps as $\C{A}$. A track $(\alpha|a)\colon
f\rr g$ in $\lambda_*\C{A}$ is represented by a track
$\alpha\colon f\rr g$ in $\C{A}$ together with an element $a\in
D_{\set{f}}$. Two of these pairs represent the same track
$$(\alpha|a)=(\beta|b)$$
if and only if
$$b=a+\lambda_{\set{f}}(\sigma^{-1}_f(\beta^\vi\vc\alpha)).$$
The vertical composition of tracks
$\bullet\st{(\beta|b)}\rr\bullet\st{(\alpha|a)}\rr\bullet$ in
$\lambda_*\C{A}$ is given by
$$(\alpha|a)\vc(\beta|b)=(\alpha\vc\beta|a+b).$$
Moreover, given a diagram in $\lambda_*\C{A}$ as follows
$$\xymatrix@C=35pt{\bullet\ar[r]^f&\bullet\ar@/^12pt/[r]_<(.25){\;}="a"\ar@/_12pt/[r]^<(.25){\,}="b"&
\bullet\ar[r]^g&\bullet\ar@{=>}^{(\alpha|a)}"a";"b"}$$
the two horizontal compositions in $\lambda_*\C{A}$ are given by
$$(\alpha|a)f=(\alpha f|\set{f}^*(a))\;\;\text{ and }\;\;g(\alpha|a)=(g\alpha|\set{g}_*(a)).$$
Notice that there is a track functor
$$\bar{\lambda}\colon\C{A}\To\lambda_*\C{A}$$
defined as the identity on objects and maps, and $\bar{\lambda}(\alpha)=(\alpha|0)$ on tracks.

\begin{prop}\label{push}
The push-forward $\lambda_*\C{A}$ satisfies $\pi_0\lambda_*\C{A}=\pi_0\C{A}$ and $\pi_1\lambda_*\C{A}=D$. Moreover,
the track functor $\bar{\lambda}$ induces $\pi_0\bar{\lambda}=1_{\pi_0\C{A}}$ and $\pi_1\bar{\lambda}=\lambda$.
\end{prop}

\begin{proof}
It is enough to notice that for any map $f\colon X\r Y$ in $\C{A}_0=(\lambda_*\C{A})_0$ the homomorphism
$$D_{\set{f}}\To\homm{X,Y}_{\lambda_*\C{A}}(f,f),\;\;a\mapsto (0^\vc_f|a),$$
is an isomorphism with inverse
$$\homm{X,Y}_{\lambda_*\C{A}}(f,f)\To D_{\set{f}},\;\;(\alpha|a)\mapsto \lambda_{\set{f}}(\sigma^{-1}_f(\alpha))+a.$$
\end{proof}

Now we are ready to prove Theorem \ref{obs1}.

\begin{proof}[Proof of Theorem \ref{obs1}]
The ``if'' part is Proposition \ref{es01}. Now suppose that we
have a $1$-cell $(\varphi,\lambda)\colon\pi\C{A}\r\pi\C{B}$ in
$\C{nat}$ such that
$\theta_{\C{A},\C{B}}(\varphi,\lambda)=\lambda_*\grupo{\C{A}}-\varphi^*\grupo{\C{B}}=0$.
By Propositions \ref{es01}, \ref{pull} and \ref{push} we have that
$\grupo{\lambda_*\C{A}}=\lambda_*\grupo{\C{A}}$ and
$\grupo{\varphi^*\C{B}}=\varphi^*\grupo{\C{B}}$, therefore
$\grupo{\lambda_*\C{A}}=\grupo{\varphi^*\C{A}}$. By \cite{stt}
2.4.2 and \cite{catc} 2.4 there exists a pseudofunctor
$\phi\colon\lambda_*\C{A}\pf\varphi^*\C{B}$ with
$\pi\phi=(1_{\pi_0\C{A}},0^\vc_{(\pi_1\C{B})\F(\varphi)})$, hence
again by Propositions \ref{pull} and \ref{push} we see that
$\psi=\bar{\varphi}\phi\bar{\lambda}$ satisfies the requirements
of Theorem \ref{obs1}.
\end{proof}

\section{Proof of Theorem \ref{obs2}}\label{84}

The ``if'' part follows from Proposition \ref{fmla1} (2).

Suppose now that a $2$-cell $\alpha$ in $\C{nat}$ is given with $\vartheta_{\varphi,\psi}(\alpha)=0$.
We choose a $1$-cochain
$$e\in\bar{F}^1(\pi_0\C{A},(\pi_1\C{B})\F(\alpha))$$ with
\begin{equation}\label{yqs}
\delta(e)=-((\pi_1\C{B})\F(0^\square_{\pi_0\varphi},\alpha))_*b_\varphi-\alpha^\#c_\C{B}
+((\pi_1\C{B})\F(\alpha,0^\square_{\pi_0\psi}))_*b_\psi.
\end{equation}
Notice that implicitly we have chosen a global section
$(t,\mu,\nu,\eta)$ as in Section \ref{gs}. It will remain fixed
throughout this proof.

We claim that there exists a unique pseudonatural transformation $\tau\colon\varphi\rr\psi$ such that
for all objects $X$ and morphisms $\set{f}\colon X\r Y$ in $\pi_0\C{A}$,
$\tau_X=t\alpha_X$ and $\tau_{t\set{f}}$ is the following composite track
$$
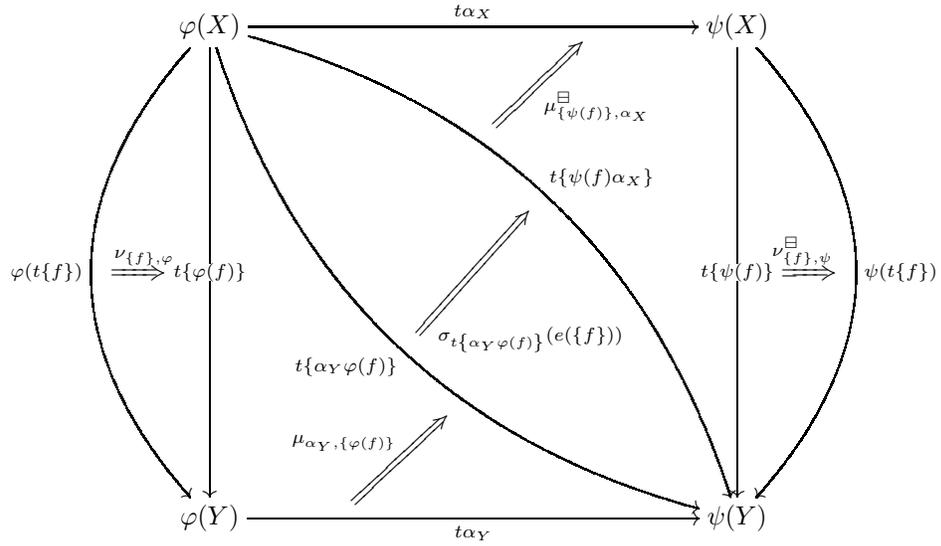
\xymatrix@C=170pt@R=170pt{\varphi(X)\ar[r]^{t\alpha_X}_<(.7){\;}="j"\ar@/_45pt/[d]_{\varphi (t\set{f})}^{\;}="c"
\ar[d]|{t\set{\varphi(f)}}_<(.478){\;\;\;\;\;}="d"
\ar@/^40pt/[rd]^{t\set{\psi(f)\alpha_X}}_{\;}="b"^<(.35){\;}="i"\ar@/_40pt/[rd]_{t\set{\alpha_Y\varphi(f)}}^{\;}="a"_<(.6){\;}="h"&
\psi(X)\ar@/^45pt/[d]^{\psi (t\set{f})}_{\;}="f"\ar[d]|{t\set{\psi(f)}}^<(0.478){\;\;\;\;\;}="e"\\
\varphi(Y)\ar[r]_{t\alpha_Y}^<(.2){\;}="g"&\psi(Y) \ar@{=>}"a";"b"_<(.1){\sigma_{t\set{\alpha_Y\varphi(f)}}(e(\set{f}))}
\ar@{=>}"c";"d"^{\nu_{\set{f},\varphi}}\ar@{=>}"e";"f"^{\nu_{\set{f},\psi}^\vi}\ar@{=>}"g";"h"^{\mu_{\alpha_Y,\set{\varphi(f)}}}
\ar@{=>}"i";"j"_{\mu_{\set{\psi(f)},\alpha_X}^\vi}}$$
that is
\begin{equation*}
\tau_{t\set{f}}=(\nu_{\set{f},\psi}^\vi t\alpha_X)\vc\mu_{\set{\psi(f)},\alpha_X}^\vi
\vc\sigma_{t\set{\alpha_Y\varphi(f)}}(e(\set{f}))\vc\mu_{\alpha_Y,\set{\varphi(f)}}\vc(t\alpha_Y\nu_{\set{f},\varphi}).
\end{equation*}
If such a $\tau$ exists then by (\ref{transnat}) for any map $f$
in $\C{A}$ if $\xi\colon f\rr t\set{f}$ is any track we have
$$\tau_f=(\psi(\xi)^\vi
t\alpha_X)\vc\tau_{t\set{f}}\vc(t\alpha_Y\varphi(\xi)),$$
In fact we can take this formula as a definition for
$\tau_f$. In order to check that $\tau_f$ does not depend on the
track $\xi$ chosen for its definition it is enough to see that the
equality holds when $\xi$ is a self-track of $t\set{f}$, or
equivalently that
$$\tau_{t\set{f}}=(\psi(\sigma_{t\set{f}}(a))^\vi t\alpha_X)\vc\tau_{t\set{f}}\vc(t\alpha_Y\varphi(\sigma_{t\set{f}}(a)))$$
for all $a\in(\pi_1\C{A})_{\set{f}}$. This follows from the next equalities, where we use
(\ref{lec1}), (\ref{lec2}) and (\ref{lec3}) in the first one, and the commutativity of (\ref{hdc}) in the second one
\begin{eqnarray*}
(\psi(\sigma_{t\set{f}}(a))^\vi t\alpha_X)\vc\tau_{t\set{f}}\vc(t\alpha_Y\varphi(\sigma_{t\set{f}}(a)))&=&
(\nu_{\set{f},\psi}^\vi t\alpha_X)\vc\mu_{\set{\psi(f)},\alpha_X}^\vi
\\&&\vc\sigma_{t\set{\alpha_Y\varphi(f)}}(-\alpha_X^*(\pi_1\psi)_*(a)\\&&+e(\set{f})+{\alpha_Y}_*(\pi_1\varphi)_*(a))\\&&\vc
\mu_{\alpha_Y,\set{\varphi(f)}}\vc(t\alpha_Y\nu_{\set{f},\varphi})\\
&=&\tau_{t\set{f}}.
\end{eqnarray*}

This can also be used to check that $\tau$ satisfies (\ref{transnat}). Moreover, (\ref{transnat3}) follows easily from (\ref{gs1}), \dots,
(\ref{gs4}), so it is only left to prove (\ref{transnat2}). For this given composable
maps $X\st{f}\r Y\st{g}\r Z$ in $\C{A}$ we choose fixed tracks $\xi_f\colon f\rr t\set{f}$ and
$\xi_g\colon g\rr t\set{g}$ for the definition of $\tau_f$ and $\tau_g$, respectively.
We will use the track $\mu_{\set{g},\set{f}}\vc(\xi_g\xi_f)\colon gf\rr t\set{gf}$ for the definition of $\tau_{gf}$.

\bigskip

$\begin{array}{ll}
(\psi_{g,f}t\alpha_X)\vc(\psi(g)\tau_f)\vc(\tau_g\varphi(f))\vc(t\alpha_Z\varphi_{g,f}^\vi)&=
\end{array}$

$
\begin{array}{ll}
=&(\psi_{g,f}t\alpha_X)\vc(\psi(g)\psi(\xi_f)^\vi
t\alpha_X)\vc(\psi(g)\nu_{\set{f},\psi}^\vi t\alpha_X)\vc(\psi(g)\mu_{\set{\psi(f)},\alpha_X}^\vi)
\\&\vc(\psi(g)\sigma_{t\set{\alpha_Y\varphi(f)}}(e(\set{f})))\vc
(\psi(g)\mu_{\alpha_Y,\set{\varphi(f)}})\vc(\psi(g)t\alpha_Y\nu_{\set{f},\varphi})\\&\vc(\psi(g)t\alpha_Y\varphi(\xi_f))
\vc(\psi(\xi_g)^\vi
t\alpha_Y\varphi(f))\vc(\nu_{\set{g},\psi}^\vi t\alpha_Y\varphi(f))\vc(\mu_{\set{\psi(g)},\alpha_Y}^\vi\varphi(f))
\\&\vc(\sigma_{t\set{\alpha_Z\varphi(f)}}(e(\set{g}))\varphi(f))\vc
(\mu_{\alpha_Z,\set{\varphi(g)}}\varphi(f))\vc(t\alpha_Z\nu_{\set{g},\varphi}\varphi(f))\\&\vc(t\alpha_Z\varphi(\xi_g)\varphi(f))
\vc(t\alpha_Z\varphi_{g,f}^\vi)
\end{array}$

$
\begin{array}{ll}
\st{\mathrm(a)}=&
(\psi_{g,f}t\alpha_X)\vc(\psi(\xi_g)^\vi\psi(\xi_f)^\vi
t\alpha_X)\vc(\nu_{\set{g},\psi}^\vi\nu_{\set{f},\psi}^\vi t\alpha_X)\vc(t\set{\psi(g)}\mu_{\set{\psi(f)},\alpha_X}^\vi)
\\&\vc\sigma_{t\set{\psi(g)}t\set{\alpha_Y\varphi(f)}}(\set{\psi(g)}_*e(\set{f}))\vc
(t\set{\psi(g)}\mu_{\alpha_Y,\set{\varphi(f)}})
\\&\vc(\mu_{\set{\psi(g)},\alpha_Y}^\vi t\set{\varphi(f)})\vc\sigma_{t\set{\alpha_Z\varphi(f)}t\set{\varphi(f)}}(\set{\varphi(f)}^*e(\set{g}))\\&\vc
(\mu_{\alpha_Z,\set{\varphi(g)}}t\set{\varphi(f)})
\vc(t\alpha_Z\nu_{\set{g},\varphi}\nu_{\set{f},\varphi})\vc(t\alpha_Z\varphi(\xi_g)\varphi(\xi_f))\vc(t\alpha_Z\varphi^\vi_{g,f})
\end{array}$

$
\begin{array}{ll}
\st{\mathrm{(b)}}=&
(\psi(\xi_g\xi_f)^\vi
t\alpha_X)\vc(\psi_{t\set{g},t\set{f}}t\alpha_X)\vc(\nu_{\set{g},\psi}^\vi\nu_{\set{f},\psi}^\vi t\alpha_X)\\&
\vc\sigma_{t\set{\psi(g)}t\set{\psi(f)}t\alpha_Y}(\set{\psi(g)}_*e(\set{f})+\set{\varphi(f)}^*e(\set{g}))\\&
\vc(t\set{\psi(g)}\mu_{\set{\psi(f)},\alpha_X}^\vi)\vc
(t\set{\psi(g)}\mu_{\alpha_Y,\set{\varphi(f)}})
\vc(\mu_{\set{\psi(g)},\alpha_Y}^\vi t\set{\varphi(f)})\\&\vc(\mu_{\alpha_Z,\set{\varphi(g)}}t\set{\varphi(f)})
\vc(t\alpha_Z\nu_{\set{g},\varphi}\nu_{\set{f},\varphi})\vc(t\alpha_Z\varphi^\vi_{t\set{g},t\set{f}})\vc(t\alpha_Z\varphi(\xi_g\xi_f))
\end{array}$

$
\begin{array}{ll}
\st{\mathrm{(c)}}=&(\psi(\xi_g\xi_f)^\vi
t\alpha_X)\vc(\psi_{t\set{g},t\set{f}}t\alpha_X)\vc(\nu_{\set{g},\psi}^\vi\nu_{\set{f},\psi}^\vi t\alpha_X)
\\&\vc\sigma_{t\set{\psi(g)}t\set{\psi(f)}t\alpha_Y}(\alpha_X^*b_\psi(\set{f},\set{g})-(\alpha^\#c_\C{B})(\set{f},\set{g})
+e(\set{gf}))\\&
\vc(t\set{\psi(g)}\mu_{\set{\psi(f)},\alpha_X}^\vi)\vc
(t\set{\psi(g)}\mu_{\alpha_Y,\set{\varphi(f)}})
\vc(\mu_{\set{\psi(g)},\alpha_Y}^\vi
t\set{\varphi(f)})\\&\vc(\mu_{\alpha_Z,\set{\varphi(g)}}t\set{\varphi(f)})
\vc(t\alpha_Z\nu_{\set{g},\varphi}\nu_{\set{f},\varphi})\vc(t\alpha_Z\varphi^\vi_{t\set{g},t\set{f}})
\\
&\vc\sigma_{t\alpha_Zt\set{\varphi(g)}t\set{\varphi(f)}\alpha_Y}(-{\alpha_Z}_*b_\varphi(\set{f},\set{g}))\vc(t\alpha_Z\varphi(\xi_g\xi_f))
\end{array}$

$\begin{array}{ll}
\st{\mathrm{(d)}}=&
(\psi(\xi_g\xi_f)^\vi
t\alpha_X)
\vc(\psi(\mu^\vi_{\set{g},\set{f}})t\alpha_X)\vc(\nu^\vi_{\set{gf},\psi}t\alpha_X)\vc(\mu_{\set{\psi(g)},\set{\psi(f)}}t\alpha_X)\\&
\vc\sigma_{t\set{\psi(g)}t\set{\psi(f)}t\alpha_Y}(-c_\C{B}(\set{\psi(g)},\set{\psi(f)},\alpha_X))
\vc(t\set{\psi(g)}\mu_{\set{\psi(f)},\alpha_X}^\vi)\\&\vc
\sigma_{t\set{\psi(g)}t\set{\alpha_Y\varphi(f)}}(e(\set{gf})+c_\C{B}(\set{\psi(g)},\alpha_Y,\set{\varphi(f)}))\\&
\vc(t\set{\psi(g)}\mu_{\alpha_Y,\set{\varphi(f)}})
\vc(\mu_{\set{\psi(g)},\alpha_Y}^\vi t\set{\varphi(f)})\\&\vc
\sigma_{t\set{\alpha_Z\varphi(g)}t\set{\varphi(f)}}(-c_\C{B}(\alpha_Z,\set{\varphi(g)},\set{\varphi(f)}))
\vc(\mu_{\alpha_Z,\set{\varphi(g)}}t\set{\varphi(f)})\\&
\vc(t\alpha_Z\mu^\vi_{\set{\varphi(g)},\set{\varphi(f)}})\vc(t\alpha_Z\nu_{\set{gf},\varphi})\vc(t\alpha_Z\varphi(\mu_{\set{g},\set{f}}))
\vc(t\alpha_Z\varphi(\xi_g\xi_f))
\end{array}$

$
\begin{array}{ll}
\st{\mathrm{(e)}}=&(\psi((\xi_g^\vi\xi_f^\vi)\vc\mu^\vi_{\set{g},\set{f}})t\alpha_X)\vc(\nu^\vi_{\set{gf},\psi}t\alpha_X)
\vc\mu^\vi_{\set{\psi(gf)},\alpha_X}\vc\mu_{\set{\psi(g)},\set{\psi(f)}\alpha_Y}\\&\vc\sigma_{t\set{\psi(g)}t\set{\alpha_Y\varphi(f)}}(e(\set{gf}))
\vc\mu^\vi_{\set{\psi(g)},\alpha_Y\set{\varphi(f)}}\vc\mu_{\set{\psi(g)}\alpha_Y,\set{\varphi(f)}}\\&
\vc\mu^\vi_{\alpha_Z\set{\varphi(g)},\set{\varphi(f)}}\vc\mu_{\alpha_Z,\set{\varphi(gf)}}
\vc(t\alpha_Z\nu_{\set{gf},\varphi})\vc(t\alpha_Z\varphi(\mu_{\set{g}\set{f}}\vc(\xi_g\xi_f)))
\end{array}$

$\begin{array}{ll}
\st{\mathrm{(f)}}=&\tau_{gf}
\end{array}$

\bigskip

In (a) we apply (\ref{lec2}); for (b) we use (\ref{ps4}) and (\ref{lec1}); in (c) we use (\ref{yqs}) and (\ref{lec1});
for (d) we apply (\ref{almo}), (\ref{d}), (\ref{lec1}) and (\ref{lec2}); for (e) we use (\ref{c}) and (\ref{lec1});
and finally we use (\ref{lec1}) for (f) as well as the naturality of $\alpha$.

\section{Proof of Theorem \ref{2act}}\label{p2a}

The key of the proof of Theorem \ref{2act} is the following
result.

\begin{prop}\label{suma}
Given a pseudofunctor $\varphi\colon\C{A}\pf\C{B}$ and a
$2$-cocycle
$$s\in\bar{F}^2(\pi_0\C{A},(\pi_1\C{B})\F(\pi_0\varphi))$$ there is
a new pseudofunctor $\varphi+s\colon\C{A}\pf\C{B}$ defined as
$\varphi$ on objects, maps and tracks such that for any object $X$
in $\C{A}$ we have $(\varphi+s)_X=\varphi_X$ and for maps
$\bullet\st{g}\r\bullet\st{f}\r\bullet$ in $\C{A}$
$$(\varphi+s)_{f,g}=\sigma_{\varphi(fg)}(s(\set{f},\set{g}))\vc\varphi_{f,g}.$$ Moreover, for any choice of a global section as
in Section \ref{gs} we have the following equality
$$b_{\varphi+s}=b_\varphi+s.$$
\end{prop}

The proof is tedious but straightforward.

\begin{proof}[Proof of Theorem \ref{2act}]
If $s$ is a cocycle representing $\omega$ then by Proposition \ref{suma} and (\ref{teobs2})
the homotopy class $[\varphi+s]$ satisfies the properties of $[\varphi]+\omega$ in the statement. The uniqueness
as well as the fact that this defines an effective and transitive action follow easily from
Propositions \ref{homotono} and \ref{deriv2} applied to the vertical composite $0^\vc_{\pi_0\varphi}=0^\vc_{\pi_0\varphi}\vc0^\vc_{\pi_0\varphi}$.
\end{proof}

\section{Global sections and cochains}\label{gs}

A \emph{global section}\index{global section} $(t,\mu, \nu, \eta)$
on the abelian track categories is a choice of maps and tracks
\begin{itemize}
\item $t\set{f}\in\set{f},$

\item $\mu_{\set{f},\set{g}}\colon t\set{f}t\set{g}\rr t\set{fg},$

\item $\nu_{\set{f},\varphi}\colon \varphi (t\set{f})\rr t\set{\varphi(f)},$

\item $\eta_{X,\alpha}\colon \alpha_X  \rr t\set{\alpha_X}$,
\end{itemize}
for all abelian track categories $\C{A}$, objects $X$ and (composable) morphisms
$\set{f},\set{g}$ in $\pi_0\C{A}$, pseudofunctors between
abelian track categories $\varphi\colon\C{A}\pf\C{B}$, and pseudonatural transformations $\alpha\colon\varphi\rr\psi$. These choices
must satisfy the following properties
for all objects $X$ and morphisms $\set{f}\colon X\r Y$ in
$\pi_0\C{A}$
\begin{eqnarray}
\label{gs1} t\set{1_X}&=&1_X,\\
\mu_{\set{f},\set{1_X}}&=&0^\vc_{t\set{f}},\\
\mu_{\set{1_Y},\set{f}}&=&0^\vc_{t\set{f}},\\
\label{gs4} \nu_{\set{1_X},\varphi}&=&\varphi_X.
\end{eqnarray}

Fixed a global section we define normalized cochains
$$c_\C{A}\in\bar{F}^3(\pi_0\C{A},\pi_1\C{A}),$$
$$b_{\varphi}\in\bar{F}^2(\pi_0\C{A},(\pi_1\C{B})\F(\pi_0\varphi)),$$
$$o_{\phi,\varphi}\in\bar{F}^1(\pi_0\C{A},(\pi_1\C{C})\F(\pi_0(\phi\varphi))),$$
$$e_\alpha\in\bar{F}^1(\pi_0\C{A},(\pi_1\C{B})\F(\pi_0\alpha)),$$
for all abelian track categories $\C{A}$, pseudofunctors
$\varphi\colon\C{A}\pf\C{B}$ and $\phi\colon\C{B}\pf\C{C}$, and
pseudonatural transformations $\alpha\colon\varphi\rr\psi$. The
cochain $c_\C{A}$ is the cocycle defined in \cite{ccglg}, it is
given by the formula
\begin{eqnarray}
\label{c}
c_\C{A}(\set{f},\set{g},\set{h})&=&\sigma_{t\set{fgh}}^{-1}(\mu_{\set{fg},\set{h}}\vc
(\mu_{\set{f},\set{g}}t\set{h})\\\nonumber&&\vc(t\set{f}\mu_{\set{g},\set{h}}^\vi)\vc\mu_{\set{f},\set{gh}}^\vi).
\end{eqnarray}
The formula for the cochain $b_\varphi$ is
\begin{eqnarray}
\label{d}
b_\varphi(\set{f},\set{g})&=&\sigma^{-1}_{t\set{\varphi(fg)}}(\mu_{\set{\varphi(f)},\set{\varphi(g)}}\vc
(\nu_{\set{f},\varphi}\nu_{\set{g},\varphi})\\\nonumber
&&\vc\varphi^\vi_{t\set{f},t\set{g}}\vc\varphi(\mu_{\set{f},\set{g}}^\vi)
\vc\nu^\vi_{\set{fg},\varphi}),
\end{eqnarray}
for $o_{\phi,\varphi}$ we have
$$o_{\phi,\varphi}(\set{f})=\sigma_{t\set{\phi\varphi(f)}}(\nu_{\set{f},\phi\varphi}\vc\phi(\nu_{\set{f},\varphi}^\vi)\vc
\nu_{\set{\varphi(f)},\phi}^\vi),$$ and finally if $\set{f}\colon
Y\r X$ is a morphism in $\pi_0\C{A}$
\begin{eqnarray*}
e_\alpha(\set{f})&=&\sigma^{-1}_{t\set{\alpha_X\varphi(f)}}(\mu_{\set{\psi(f)},\set{\alpha_Y}}\vc
(\nu_{\set{f},\psi}\eta_{Y,\alpha})\vc\alpha_{t\set{f}}\\&&\vc
(\eta_{X,\alpha}^\vi\nu_{\set{f},\varphi}^\vi)\vc\mu_{\set{\alpha_X},\set{\varphi(f)}}^\vi).
\end{eqnarray*}
The reader can check that these cochains are indeed normalized by
using (\ref{gs1}), \dots, (\ref{gs4}), (\ref{ps2}), (\ref{ps3})
and (\ref{ps6}).

\begin{prop}\label{fmla1}
The following equalities hold
\begin{enumerate}
\item $\delta(b_\varphi)=(\pi_1\varphi)_*c_\C{A}-(\pi_0\varphi)^*c_\C{B}$,

\item
$\delta(e_\alpha)=-((\pi_1\C{B})\F(0^\square_{\pi_0\varphi},\alpha))_*b_\varphi-\alpha^\#c_\C{B}
+((\pi_1\C{B})\F(\alpha,0^\square_{\pi_0\psi}))_*b_\psi$,

\item
$\delta(o_{\phi,\varphi})=b_{\phi\varphi}-((\pi_1\phi)\F(\pi_0\varphi))_*b_\varphi-(\pi_0\varphi)^*b_\phi$.
\end{enumerate}
\end{prop}

\begin{proof}
The proof of this proposition is a tedious but straightforward
computation which uses intensively the properties of linear track
extensions of categories, pseudofunctors and pseudonatural
transformations, see Section \ref{atc} and the Appendix. In fact
the proof of (2) consists essentially of reversing the long
computation carried out in the proof of Theorem \ref{obs2}, see
Section \ref{84}. So that the reader can get an idea of how to
carry out these calculations we include here a detailed proof of
(3).

Given composable morphisms
$\bullet\st{\set{f}}\r\bullet\st{\set{g}}\r\bullet$ in
$\pi_0\C{A}$

\begin{eqnarray*}
\delta(o)(\set{f},\set{g})&=
&\set{f}_*o(\set{g})-o(\set{fg})+\set{g}^*o(\set{f})
\\&=&\set{f}_*\sigma_{t\set{\phi\varphi(g)}}(\nu_{\set{g},\phi\varphi}\vc\phi(\nu_{\set{g},\varphi}^\vi)\vc
\nu_{\set{\varphi(g)},\phi}^\vi)\\&&-\sigma_{t\set{\phi\varphi(fg)}}(\nu_{\set{fg},\phi\varphi}
\vc\phi(\nu_{\set{fg},\varphi}^\vi)\vc
\nu_{\set{\varphi(fg)},\phi}^\vi)\\&&
+\set{g}^*\sigma_{t\set{\phi\varphi(f)}}(\nu_{\set{f},\phi\varphi}\vc\phi(\nu_{\set{f},\varphi}^\vi)\vc
\nu_{\set{\varphi(f)},\phi}^\vi)\\&\st{\mathrm{(a)}}=&
\sigma_{t\set{\phi\varphi(fg)}}(\mu_{\set{\phi\varphi(f)},\set{\phi\varphi(g)}}
\vc(t\set{\phi\varphi(f)}\nu_{\set{g},\phi\varphi})
\\&&\vc(t\set{\phi\varphi(f)}\phi(\nu_{\set{g},\varphi}^\vi))\vc
(t\set{\phi\varphi(f)}\nu_{\set{\varphi(g)},\phi}^\vi)\\&&\vc\mu_{\set{\phi\varphi(f)},\set{\phi\varphi(g)}}^\vi\vc
\nu_{\set{\varphi(fg)},\phi}\vc\phi(\nu_{\set{fg},\varphi})\vc\nu_{\set{fg},\phi\varphi}^\vi\\&&
\vc\mu_{\set{\phi\varphi(f)},\set{\phi\varphi(g)}}
(\nu_{\set{f},\phi\varphi}t\set{\phi\varphi(g)})\vc(\phi(\nu_{\set{f},\varphi}^\vi)t\set{\phi\varphi(g)})\\&&\vc
(\nu_{\set{\varphi(f)},\phi}^\vi t\set{\phi\varphi(g)})
\vc\mu_{\set{\phi\varphi(f)},\set{\phi\varphi(g)}}^\vi)\\&\st{\mathrm{(b)}}=&
\sigma_{t\set{\phi\varphi(fg)}}(\nu_{\set{\varphi(fg)},\phi}\vc\phi(\nu_{\set{fg},\varphi})
\vc\nu_{\set{fg},\phi\varphi}^\vi\\&&\vc\mu_{\set{\phi\varphi(f)},\set{\phi\varphi(g)}}
\vc(\nu_{\set{f},\phi\varphi}\nu_{\set{g},\phi\varphi})
\\&&\vc(\phi(\nu_{\set{f},\varphi}^\vi)\phi(\nu_{\set{g},\varphi}^\vi))\vc
(\nu_{\set{\varphi(f)},\phi}^\vi\nu_{\set{\varphi(g)},\phi}^\vi)
\\&&\vc\mu_{\set{\phi\varphi(f)},\set{\phi\varphi(g)}}^\vi)
\end{eqnarray*}
In (a) we use (\ref{lec1}) and (\ref{lec2}) and in (b) we apply
(\ref{lec1}) again.

On the other hand
\begin{small}
\begin{eqnarray*}
b_{\phi\varphi}(\set{f},\set{g})&&\\-((\pi_1\phi)\F(\pi_0\varphi))_*(b_\varphi)(\set{f},\set{g})
\\-(\pi_0\varphi)^*(b_\phi)(\set{f},\set{g})&\st{\mathrm{(c)}}=&
\sigma^{-1}_{t\set{\phi\varphi(fg)}}(\mu_{\set{\phi\varphi(f)},\set{\phi\varphi(g)}}\vc
(\nu_{\set{f},\phi\varphi}\nu_{\set{g},\phi\varphi})\\&&
\vc(\phi\varphi)^\vi_{t\set{f},t\set{g}}\vc\phi\varphi(\mu_{\set{f},\set{g}}^\vi)
\vc\nu^\vi_{\set{fg},\phi\varphi})\\&&
-\sigma^{-1}_{\phi(t\set{\varphi(fg)})}(\phi(\mu_{\set{\varphi(f)},\set{\varphi(g)}})\vc
\phi(\nu_{\set{f},\varphi}\nu_{\set{g},\varphi})\\&&\vc(\phi\varphi)^\vi_{t\set{f},t\set{g}}
\vc\phi\varphi(\mu_{\set{f},\set{g}}^\vi)
\vc\phi(\nu^\vi_{\set{fg},\varphi}))\\&&   
-\sigma^{-1}_{t\set{\phi\varphi(fg)}}(\mu_{\set{\phi\varphi(f)},\set{\phi\varphi(g)}}\vc
(\nu_{\set{\varphi(f)},\phi}\nu_{\set{\varphi(g)},\phi})\\&&
\vc\phi^\vi_{t\set{\varphi(f)},t\set{\varphi(g)}}\vc\phi(\mu_{\set{\varphi(f)},\set{\varphi(g)}}^\vi)
\vc\nu^\vi_{\set{\varphi(fg)},\phi})\\&\st{\mathrm{(d)}}=&
\sigma^{-1}_{t\set{\phi\varphi(fg)}}(\mu_{\set{\phi\varphi(f)},\set{\phi\varphi(g)}}\vc
(\nu_{\set{f},\phi\varphi}\nu_{\set{g},\phi\varphi})\\&&
\vc\phi^\vi_{\varphi(t\set{f}),\varphi(t\set{g})}\vc\phi(\varphi^\vi_{t\set{f},t\set{g}})\vc\phi\varphi(\mu_{\set{f},\set{g}}^\vi)
\\&&\vc\nu^\vi_{\set{fg},\phi\varphi}
\vc\nu_{\set{\varphi(fg)},\phi}\vc\phi(\nu_{\set{fg},\varphi})\vc\phi\varphi(\mu_{\set{f},\set{g}})
\\&&\vc\phi(\varphi_{t\set{f},t\set{g}})
\vc\phi(\nu^\vi_{\set{f},\varphi}\nu^\vi_{\set{g},\varphi})\\
&&\vc\phi_{t\set{\varphi(f)},t\set{\varphi(g)}} \vc
(\nu^\vi_{\set{\varphi(f)},\phi}\nu^\vi_{\set{\varphi(g)},\phi})\\&&
\vc\mu^\vi_{\set{\phi\varphi(f)},\set{\phi\varphi(g)}})\\&\st{\mathrm{(e)}}=&
\sigma^{-1}_{\phi\varphi(t\set{fg})}(\nu^\vi_{\set{fg},\phi\varphi}
\vc\nu_{\set{\varphi(fg)},\phi}\vc\phi(\nu_{\set{fg},\varphi}))\\&&
+\sigma^{-1}_{\phi(\varphi(t\set{f})\varphi(t\set{g}))}((\phi(\nu^\vi_{\set{f},\varphi})\phi(\nu^\vi_{\set{g},\varphi}))
\\&&\vc
(\nu^\vi_{\set{\varphi(f)},\phi}\nu^\vi_{\set{\varphi(g)},\phi})
\vc
(\nu_{\set{f},\phi\varphi}\nu_{\set{g},\phi\varphi}))\\&\st{\mathrm{(f)}}=&
\sigma_{t\set{\phi\varphi(fg)}}(\nu_{\set{\varphi(fg)},\phi}\vc\phi(\nu_{\set{fg},\varphi})
\vc\nu_{\set{fg},\phi\varphi}^\vi\\&&\vc\mu_{\set{\phi\varphi(f)},\set{\phi\varphi(g)}}
\vc(\nu_{\set{f},\phi\varphi}\nu_{\set{g},\phi\varphi})
\\&&\vc(\phi(\nu_{\set{f},\varphi}^\vi)\phi(\nu_{\set{g},\varphi}^\vi))\vc
(\nu_{\set{\varphi(f)},\phi}^\vi\nu_{\set{\varphi(g)},\phi}^\vi)
\\&&\vc\mu_{\set{\phi\varphi(f)},\set{\phi\varphi(g)}}^\vi)
\end{eqnarray*}
\end{small} In (c) we use (\ref{lec3}) and (\ref{ps5}); for (d) we apply
(\ref{lec1}); in (e) we use (\ref{ps4}) as well as (\ref{lec1});
and finally for (f) we use (\ref{lec1}) again. These two chains of
equalities establish (3).
\end{proof}

The following proposition follows from Proposition \ref{fmla1}
(1), (\ref{hdc2}), (\ref{homo13}) and the fact that $c_\C{A}$ and
$c_\C{B}$ are cocycles.

\begin{prop}\label{fmla2}
If $\varphi,\psi\colon\C{A}\pf\C{B}$ are pseudofunctors between abelian track categories and $\alpha\colon\pi\varphi\rr\pi\psi$ is a
$2$-cell in $\C{nat}$ then
$$\delta(((\pi_1\C{B})\F(0^\square_{\pi_0\varphi},\alpha))_*b_\varphi-((\pi_1\C{B})\F(\alpha,0^\square_{\pi_0\psi}))_*b_\psi)=
-\delta\alpha^\#(c_\C{B}).$$
\end{prop}

Let $(\bar{t},\bar{\mu},\bar{\nu},\bar{\eta})$ be another global
section. A \emph{transition}\index{transition} from
$(\bar{t},\bar{\mu},\bar{\nu},\bar{\eta})$ to $(t,\mu, \nu, \eta)$
is given by tracks
$$\rho_{\set{f}}\colon t\set{f}\rr\bar{t}\set{f}$$
indexed by the morphisms $\set{f}$ in $\pi_0\C{A}$ for all abelian
track categories $\C{A}$. Such a transition determines cochains
$$v_{\C{A}}\in \bar{F}^2(\pi_0\C{A},\pi_1\C{A}),$$
$$u_{\varphi}\in \bar{F}^2(\pi_0\C{A},(\pi_1\C{B})\F(\pi_0\varphi)),$$
where $\C{A}$ is any abelian track category and
$\varphi\colon\C{A}\pf\C{B}$ is an arbitrary pseudofunctor between
abelian track categories. These cochains are defined by the
following formulas
$$v_{\C{A}}(\set{f},\set{g})=\sigma^{-1}_{\bar{t}\set{fg}}(\bar{\mu}_{\set{f},\set{g}}\vc
(\rho_{\set{f}}\rho_{\set{g}})\vc\mu^\vi_{\set{f},\set{g}}\vc\rho_{\set{fg}}^\vi),$$
$$u_{\varphi}(\set{f})=\sigma^{-1}_{\bar{t}\set{f}}(\bar{\nu}_{\set{f},\varphi}\vc\varphi(\rho_{\set{f}})
\vc\nu_{\set{f},\varphi}^\vi\vc\rho_{\set{\varphi(f)}}^\vi).$$

Let us call $\bar{c}_\C{A}$ and $\bar{b}_\varphi$ the cochains
defined as above by the global section
$(\bar{t},\bar{\mu},\bar{\nu},\bar{\eta})$.

\begin{prop}\label{fmla3}
The following equalities hold
\begin{enumerate}
\item $\delta(v_{\C{A}})=c_\C{A}-\bar{c}_\C{A}$,

\item $\delta(u_{\varphi})=-b_\varphi+\bar{b}_\varphi-(\pi_0\varphi)^*v_{\C{B}}+(\pi_1\varphi)_*v_\C{A}$.
\end{enumerate}
\end{prop}

The proof of this proposition is tedious but straightforward, as
the proof of Proposition \ref{fmla1}. We leave it to the reader.

\begin{prop}\label{fmla4}
Given pseudofunctors $\varphi,\psi\colon\C{A}\pf\C{B}$ between
abelian track categories and a $2$-cell
$\alpha\colon\pi\varphi\rr\pi\psi$ in $\C{nat}$ the following
equality holds
\begin{eqnarray*}
\delta(((\pi_1\C{B})\F(0^\square_{\pi_0\varphi},\alpha))_*u_\varphi
-\alpha^\#v_\C{B}\\-((\pi_1\C{B})\F(\alpha,0^\square_{\pi_0\psi}))_*u_\psi)
&=&-((\pi_1\C{B})\F(0^\square_{\pi_0\varphi},\alpha))_*(b_\varphi-\bar{b}_\varphi)
\\&&-\alpha^\#(c_\C{B}-\bar{c}_\C{B})
\\&&+((\pi_1\C{B})\F(\alpha,0^\square_{\pi_0\psi}))_*(b_\psi-\bar{b}_\psi).
\end{eqnarray*}
\end{prop}

This proposition can be checked by using Proposition \ref{fmla3},
(\ref{hdc2}) and (\ref{homo13}).

\appendix

\section{Pseudofunctors and pseudonatural transformations}

A \emph{pseudofunctor}\index{pseudofunctor}
$\varphi\colon\C{A}\pf\C{B}$ between track categories is an
assignment of objects, maps and tracks together with additional
tracks
$$\varphi_{f,g}\colon\varphi(f)\varphi(g)\rr\varphi(fg)\;\;\text{ and }\;\;\varphi_X\colon\varphi(1_X)\rr1_{\varphi(X)}$$
for all objects $X$ and composable maps
$\bullet\st{g}\r\bullet\st{f}\r\bullet$ in $\C{A}$. These tracks
must satisfy the following conditions: the following three
composite tracks are identity tracks for all maps $f\colon X\r Y$
and composable maps $\bullet\st{h}\r\bullet\st{g}\r\bullet
\st{f}\r\bullet$ in $\C{A}$,
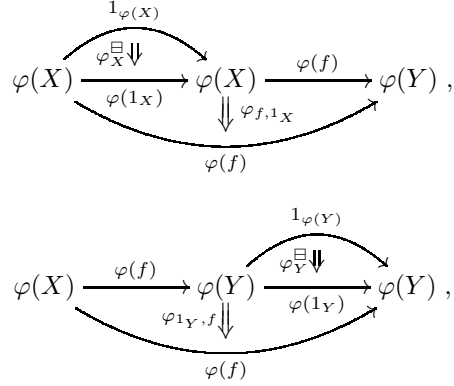
$$\xymatrix@C=40pt{\bullet\ar[r]^<(1){\;}="f"_{\varphi(h)}
\ar@/_40pt/[rrr]^<(.45){\;}="d"_{\varphi(fgh)}\ar@/^20pt/[rr]^{\varphi(gh)}_{\;}="e"^<(.745){\;}="h"
\ar@/^40pt/[rrr]_<(.5){\;}="g"^{\varphi(fgh)}
&\bullet\ar[r]|{\varphi(g)}_<(1){\;}="a"\ar@/_20pt/[rr]_{\varphi(fg)}^{\;}="b"_<(.18){\;}="c"&\bullet\ar[r]^{\varphi(f)}&
\bullet\ar@{=>}^{\;\varphi_{f,g}}"a";"b"\ar@{=>}_{\;\varphi_{fg,h}}"c";"d"\ar@{=>}^{\;\varphi_{g,h}^\vi}"e";"f"
\ar@{=>}^{\;\varphi_{f,gh}^\vi}"g";"h"},$$

$$\xymatrix@C=40pt{\varphi(X)\ar[r]^{\;}="d"_{\varphi(1_X)}_<(1){\;}="a"
\ar@/_25pt/[rr]_{\varphi(f)}^{\;}="b"\ar@/^20pt/[r]^{1_{\varphi(X)}}_{\;}="c"&\varphi(X)\ar[r]^{\varphi(f)}&
\varphi(Y)\ar@{=>}^{\;\varphi_{f,1_X}}"a";"b"
\ar@{=>}_{\;\varphi_X^\vi}"c";"d"},$$

$$\xymatrix@C=40pt{\varphi(X)\ar[r]^{\varphi(f)}_<(1){\;}="a"
\ar@/_25pt/[rr]_{\varphi(f)}^{\;}="b"&\varphi(Y)\ar[r]^{\;}="d"_{\varphi(1_Y)}\ar@/^20pt/[r]^{1_{\varphi(Y)}}_{\;}="c"&
\varphi(Y)\ar@{=>}_{\;\varphi_{1_Y,f}}"a";"b"
\ar@{=>}_{\;\varphi_Y^\vi}"c";"d"},$$
For any horizontally
composable tracks $\alpha$, $\beta$ in $\C{A}$ the following
composite track is $\varphi(\alpha\beta)$
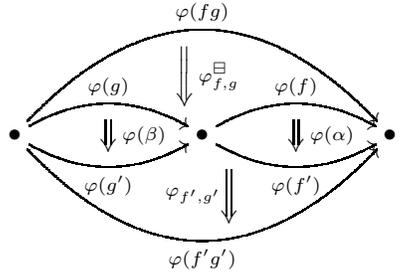
$$\xymatrix@C=60pt{\bullet\ar@/^12pt/[r]^{\varphi(g)}_{\;}="a"^<(.865){\;}="f"\ar@/_12pt/_{\varphi(g')}[r]^{\;}="b"
\ar@/^40pt/[rr]^{\varphi(fg)}_(.45){\;}="e"
\ar@/_40pt/[rr]_{\varphi(f'g')}^<(.562){\;}="h"
&\bullet\ar@{=>}^{\;\varphi(\beta)}"a";"b"
\ar@/^12pt/[r]^{\varphi(f)}_{\;}="c"\ar@/_12pt/[r]_{\varphi(f')}^{\;}="d"_<(.1){\;}="g"
&\bullet\ar@{=>}^{\;\varphi(\alpha)}"c";"d"
\ar@{=>}^{\;\varphi_{f,g}^\vi}"e";"f"
\ar@{=>}_{\;\varphi_{f',g'}}"g";"h"},$$  
and $\varphi$ preserves vertical composition of tracks as well as identity tracks.
These conditions can be summarized in the following six equations:
\begin{equation}\label{ps1}
\varphi_{fg,h}\vc(\varphi_{f,g}\varphi(h))=\varphi_{f,gh}\vc(\varphi(f)\varphi_{g,h}),
\end{equation}
\begin{equation}\label{ps2}
\varphi_{f,1_X}=\varphi(f)\varphi_X,
\end{equation}
\begin{equation}\label{ps3}
\varphi_{1_Y,f}=\varphi_Y\varphi(f),
\end{equation}
\begin{equation}\label{ps4}
\varphi_{f',g'}\vc(\varphi(\alpha)\varphi(\beta))=\varphi(\alpha\beta)\vc\varphi_{f,g},
\end{equation}
\begin{equation}\label{ps5}
\varphi(\gamma\vc\varepsilon)=\varphi(\gamma)\vc\varphi(\varepsilon),
\end{equation}
\begin{equation}\label{ps6}
\varphi(0^\vc_f)=0^\vc_{\varphi(f)}.
\end{equation}

Notice that a \emph{track functor}\index{track functor}
$\varphi\colon\C{A}\r\C{B}$ is exactly a pseudofunctor such that
$\varphi_{f,g}$ and $\varphi_X$ are always identity tracks.

Consider now two pseudofunctors $\varphi,\psi\colon\C{A}\pf\C{B}$.
A \emph{pseudonatural transformation}\index{pseudonatural
transformation} $\alpha\colon\varphi\rr\psi$ is given by maps
$\alpha_X\colon\varphi(X)\r\psi(X)$ for all objects $X$ in $\C{A}$
together with tracks
$$\xymatrix{\varphi(X)\ar[r]^{\alpha_X}_<(.5){\;}="b"\ar[d]_{\varphi(f)}^<(.5){\;}="a"&\psi(X)\ar[d]^{\psi(f)}\\
\varphi(Y)\ar[r]_{\alpha_Y}&\psi(Y)\ar@{=>}"a";"b"_{\;\alpha_f}}$$
for all maps $f\colon X\r Y$ in $\C{A}$. These tracks satisfy the following properties: for all tracks $\tau\colon f\rr g$
in $\C{A}$ the following composite track coincides with $\alpha_g$
$$\xymatrix@R=35pt@C=35pt{\varphi(X)\ar@/_30pt/[d]_{\varphi(g)}^{\;}="c"\ar[r]^{\alpha_X}_<(.5){\;}="b"
\ar[d]|<(.6){\varphi(f)}^<(.5){\;}="a"_{\;}="d"&
\psi(X)\ar[d]|<(.2){\psi(f)}^{\;}="e"\ar@/^30pt/[d]_{\;}="f"^{\psi(g)}\\
\varphi(Y)\ar[r]_{\alpha_Y}&\psi(Y)\ar@{=>}"a";"b"_{\;\alpha_f}\ar@{=>}"c";"d"^{\varphi(\tau)^\vi}\ar@{=>}"e";"f"_{\psi(\tau)}}$$
Given composable maps $X\st{f}\r Y\st{g}\r Z$ in $\C{A}$ the
following composite track coincides with $\alpha_{fg}$
$$\xymatrix{\varphi(X)\ar@/_35pt/[dd]_{\varphi(fg)}^{\;}="e"\ar[r]^{\alpha_X}_<(.5){\;}="b"\ar[d]_{\varphi(f)}_<(1){\;\;\;\;}="f"^<(.5){\;}="a"&
\psi(X)\ar@/^35pt/[dd]^{\psi(fg)}_{\;}="h"\ar[d]^{\psi(f)}^<(1){\;\;\;\;}="g"\\
\varphi(Y)\ar[r]^{\alpha_Y}_<(.5){\;}="d"\ar[d]_{\varphi(g)}^<(.5){\;}="c"&\psi(Y)\ar@{=>}"a";"b"_{\;\alpha_f}\ar[d]^{\psi(g)}\\
\varphi(Z)\ar[r]^{\alpha_Z}&\psi(Z)\ar@{=>}"c";"d"_{\;\alpha_g}\ar@{=>}"e";"f"_{\varphi_{g,f}^\vi}\ar@{=>}"g";"h"_{\psi_{g,f}}}$$
For any object $X$ in $\C{A}$ the following composite track is the
identity track $0^\vc_{\alpha_X}$.
$$\xymatrix@R=35pt@C=35pt{\varphi(X)\ar@/_30pt/[d]_{1_{\varphi(X)}}^{\;}="c"\ar[r]^{\alpha_X}_<(.5){\;}="b"
\ar[d]|<(.6){\varphi(1_X)}^<(.5){\;}="a"_{\;}="d"&
\psi(X)\ar[d]|<(.2){\psi(1_X)}^{\;}="e"\ar@/^30pt/[d]_{\;}="f"^{1_{\psi(X)}}\\
\varphi(X)\ar[r]_{\alpha_X}&\psi(X)\ar@{=>}"a";"b"_{\;\alpha_{1_X}}\ar@{=>}"c";"d"^{\varphi_X^\vi}\ar@{=>}"e";"f"_{\psi_X}}$$
These three conditions can be restated as the following equations
\begin{equation}\label{transnat}
(\psi(\tau)\alpha_X)\vc\alpha_f=\alpha_g\vc(\alpha_Y\varphi(\tau)),
\end{equation}
\begin{equation}\label{transnat2}
(\psi_{g,f}\alpha_X)\vc(\psi(g)\alpha_f)\vc(\alpha_g\varphi(f))=\alpha_{gf}\vc(\alpha_Z\varphi_{g,f}),
\end{equation}
\begin{equation}\label{transnat3}
\varphi_X=(\psi_X\alpha_X)\vc\alpha_{1_X}.
\end{equation}

If $\varphi$ and $\psi$ were track functors, $\alpha$ would be a
\emph{track natural transformation}\index{track natural
transformation} provided $\alpha_f$ were always an identity track.

Recall that a \emph{homotopy}\index{homotopy between
pseudofunctors} $\xi\colon\varphi\rr\psi$ between pseudofunctors
$\varphi,\psi\colon\C{A}\pf\C{B}$ has been defined as a
pseudonatural transformation with
$\xi_X=1_{\varphi(X)}=1_{\psi(X)}$ for all objects $X$ in $\C{A}$.
In order to regard homotopies as a special kind of pseudofunctors
we define the \emph{reduced cylinder track category}\index{reduced
cylinder track category} $\D$. It has a unique $0$-cell $*$, a
unique non-trivial morphism $\imath\colon*\r*$ with
$\iota^2=\iota$, and only one non-identity track $\ell\colon
1_*\rr \imath$ which also satisfies $\ell^2=\ell$. The reader can
check that this category can actually be drawn as a reduced
cylinder, that is, a cylinder over the circle with a shrunk
segment. There are exactly two pseudofunctors from the final track
category $j^0,j^1\colon*\pf\D$. One of them, let us say $j^0$, is
a track functor, and the other one is not. The unique
pseudofunctor in the opposite direction is denoted by
$q\colon\D\pf*$. Notice that all these pseudofunctors induce the
identity $\pi\D=\pi*$ when applying the functor $\pi$, that is,
they are weak equivalences.

\begin{prop}\label{carhomo}
Any homotopy $\xi\colon\varphi\rr\psi$ induces a
pseudofunctor $\bar{\xi}\colon\C{A}\times\D\pf\C{B}$ with
$\bar{\xi}j^0=\varphi$ and $\bar{\xi}j^1=\psi$. Conversely a
pseudofunctor $\zeta\colon\C{A}\times\D\pf\C{B}$ induces a
homotopy $\bar{\zeta}\colon\zeta j^0\rr\zeta j^1$ with
$\bar{\zeta}_f=\zeta(0^\vc_f,\ell)$. Moreover
$\bar{\bar{\xi}}=\xi$ and $\bar{\bar{\zeta}}=\zeta$.
\end{prop}

The homotopy $\bar{\zeta}$ is completely defined in the statement,
so all one has to do is to check that the axioms of a
pseudonatural transformation hold. We leave this task to the
reader. We also leave to the reader to check that there is a
unique choice for $\bar{\xi}$ satisfying the properties of the
statement.

The next proposition shows that pseudofunctors are flabby, in the
sense that one can change de definition on maps within the same
homotopy class in a quite general way.

\begin{prop}\label{prehom}
Let $\varphi\colon\C{A}\pf\C{B}$ be a pseudofunctor and
$\xi_f\colon\varphi(f)\rr\psi^f$ a collection of arbitrary tracks
indexed by the maps of $\C{A}$. There is a well-defined
pseudofunctor $\psi\colon\C{A}\pf\C{B}$ which coincides with
$\varphi$ on objects given on maps by $\psi(f)=\psi^f$ on tracks
$\alpha\colon f\rr g$ by
$\psi(\alpha)=\xi_g\varphi(\alpha)\xi_f^\vi$, and by the formulas
$\psi_{f,g}=\xi_{fg}\vc\varphi_{f,g}\vc(\xi^\vi_f\xi^\vi_g)$, and
$\psi_X=\varphi_X\xi_{1_X}^\vi$. Moreover, the tracks
$\xi_f$ define a homotopy $\xi\colon\varphi\rr\psi$. 
\end{prop}

A \emph{reduced pseudofunctor}\index{reduced!pseudofunctor}
$\varphi$ is a pseudofunctor that preserves identity maps
$\varphi(1_X)=1_{\varphi(X)}$ and $\varphi_X=0^\vc_{1_X}$.

\begin{prop}
Any pseudofunctor is homotopic to a reduced one.
\end{prop}

\begin{proof}
It is enough to apply Proposition \ref{prehom} to the tracks $\xi_f\colon\varphi(f)\rr\psi^f$ such that
$\psi^f=\varphi(f)$ and $\xi_f=0^\vc_{\varphi(f)}$ if $f$ is not an identity map, $\psi^{1_X}=1_{\varphi(X)}$
and $\xi_{1_X}=\varphi_X$.
\end{proof}


A track category $\C{A}$ has a \emph{strict zero
object}\index{strict zero object} $*$ if the morphism groupoids
$\homm{X,*}_\C{A}$ and $\homm{*,X}_\C{A}$ are trivial for all
objects $X$. The strict zero object is unique up to isomorphism in
$\C{A}$. In this case all morphism groupoids $\homm{X,Y}_\C{A}$
have a distinguished map $0_{X,Y}\colon X\r Y$ which is the unique
one that factors as $X\r*\r Y$.

\begin{prop}
Let $\C{A}$ be an abelian track category such that $\pi_0\C{A}$
has a zero object $*$ and $\pi_1\C{A}$ vanishes over all maps $X\r
*$ and $*\r X$, then there is an abelian track category $\C{B}$
with a strict zero object $*$ and a track functor
$\varphi\colon\C{A}\r\C{B}$ with $\pi_0\varphi=1_{\pi_0\C{A}}$ and
$\pi_1\varphi=0^\vc_{\pi_1\C{A}}$, in particular $\varphi$ is a
weak equivalence and $\C{A}$ is homotopy equivalent to $\C{B}$.
\end{prop}

\begin{proof}
The objects of $\C{B}$ are the same as in $\C{A}$. In order to
define maps in $\C{B}$ we  say that two maps $f,g\colon X\r Y$ in
$\C{A}$ are equivalent $f\sim g$ if $f=g$ or if there is a diagram
\begin{equation*}\tag{a}
\xymatrix{X\ar@/^10pt/[r]_{\;}="a"\ar@/_10pt/[r]^{\;}="b"&{*}\ar@/_10pt/[r]^{\;}="d"\ar@/^10pt/[r]_{\;}="c"&Y
\ar@{=>}"a";"b" \ar@{=>}"c";"d"}
\end{equation*}
in $\C{A}$ where the composition of the upper arrows is $f$ and
the composition of the lower ones is $g$. Notice that the tracks
in (a) are unique by the vanishing condition imposed to
$\pi_1\C{A}$. The reader can check that this is a natural
equivalence relation on the ordinary category $\C{A}_0$, therefore
the quotient category $\C{B}_0=\C{A}_0/\!\sim$ is defined. Let us
write $[f]$ for the equivalence class of a map $f$ in $\C{A}$.
Tracks $[f]\rr [g]$ in $\C{B}$ are defined as tracks $f\rr g$ in
$\C{A}$ between two arbitrarily chosen representatives. This
definition is consistent because, as we noticed before, diagram
(a) is unique when it exists, hence vertical composition of tracks
in $\C{B}$ can be defined in a unique reasonable way.
\end{proof}

Suppose that $\varphi\colon\C{A}\pf\C{B}$ is a pseudofunctor
between track categories with strict zero object $*$ such that
$\varphi$ preserves the zero object $\varphi(*)=*$. The
pseudofunctor $\varphi$ is \emph{normalized at zero
maps}\index{normalized at zero maps!pseudofunctor} if it preserves
zero maps $\varphi(0_{X,Y})=0_{\varphi(X),\varphi(Y)}$, and given
maps $f\colon Y\r Z$ and $g\colon W\r X$ the equalities
$\varphi_{f,0_{X,Y}}=0^\vc_{0_{X,Z}}$ and
$\varphi_{0_{X,Y},g}=0^\vc_{0_{W,X}}$ hold.

\begin{prop}
Any pseudofunctor $\varphi\colon\C{A}\pf\C{B}$ preserving strict
zero objects is homotopic to a pseudofunctor $\psi$ normalized at
zero maps.
\end{prop}

\begin{proof}
The pseudofunctor $\psi$ is obtained by applying Proposition
\ref{prehom} to the tracks $\xi_f\colon\varphi(f)\rr\psi^f$ such
that $\psi^f=\varphi(f)$ and $\xi_f=0^\vc_{\varphi(f)}$ if $f$ is
not a zero map, $\psi^{0_{X,Y}}=0_{\varphi(X),\varphi(Y)}$ and
$\xi_{0_{X,Y}}=\varphi_{0_{*,Y},0_{X,*}}^\vi$. The reader can
check that $\psi$ is indeed normalized by using (\ref{ps4}) and
the fact that $*$ is a strict zero object.
\end{proof}


\bibliographystyle{alpha}
\bibliography{Fernando}
\end{document}


\title{{[I]}\\\hspace{5pt}\\Triangulated track categories}
\author{H.-J. Baues}
\date{}

\blfootnote{\textit{$2000$ Mathematics Subject Classification:} 18D05, 18E30}%
\blfootnote{\textit{Key words and phrases:} triangulated categories, groupoid-enriched categories}%

\maketitle
\begin{abstract}
The concept of a triangulated homotopy category can be canonically
lifted to the level of a groupoid-enriched category. This way two
natural axioms on track triangles replace the four somewhat
obscure axioms of a triangulated category.
\end{abstract}

\section*{Introduction}

The axioms of a triangulated category are deduced from properties
of cofiber sequences (or fiber sequences) in a stable homotopy
category $\A$. Here $\A$ is an additive category with a
translation functor $t$ as for example the homotopy
category of chain complexes or the stable homotopy category of
spectra. Cofibers or fibers, however, are defined in terms of
homotopies and, as observed by Gabriel-Zisman \cite{GZ}, represent
``cone functors'' obtained by tracks. Here tracks are homotopy
classes of homotopies. Therefore it is natural to consider
properties of cofiber sequences in a {\em track category}, i.e. a
groupoid enriched category or a 2-category for which all 2-cells
are invertible. Gabriel-Zisman deal with the general unstable
fiber sequences in a track category. In this paper, however, we
consider stable fiber-cofiber sequences in a track category and we
lift the classical concept of triangulation of a category (as
introduced by Puppe and Verdier) to the new concept of
triangulation of a track category. A triangulated category is an
additive category with translation functor and a distinguished
class of {\em exact triangles}. Similarly a triangulated track
category is an additive track category with a translation track
functor and a distinguished class of {\em track triangles}
satisfying two natural axioms (TTr1) and (TTr2) which replace the
rather obscure axioms of a triangulated category. As a main result
we show that a triangulated track category always has a homotopy
category which is a triangulated category in the classical sense.
We also obtain this way a description of good morphisms between
exact triangles in the sense of \cite{Ne}. Further applications
are obtained in \cite{Ba2} where we compute an algebraic model for
the
category of principal maps between 2-stage Postnikov spectra. 

In \cite{ccctc} one finds the cohomological interpretation of the results on track categories in this paper.

The author is grateful to Fernando Muro and Teimuraz Pirashvili for discussions and hints concerning this paper.

\section{Triangulated categories}

 To fix notation we repeat the following well known concepts. Let $\A$
be an additive category. A {\em translation
functor}\index{translation functor} $t:\A\rightarrow \A$ is an
additive automorphism of $\A$. Hence
\begin{equation}\label{1.1}
t: \Hom_{\A}(X,Y)\cong \Hom_{\A} (tX,tY)
\end{equation}
is an isomorphism of abelian groups and $t(X\oplus Y)=
(tX)\oplus (tY)$. The inverse $t^{-1}$ of $t$ admits natural
isomorphisms $tt^{-1}X\cong X$ and $t^{-1}tX\cong X$ in $\A$. We
say that $t$ is {\em strict} if these isomorphisms are the
identity of $X$.

Now let $(\A,t)$ be an additive category with a translation
functor. A {\em triangle}\index{triangle}  in $\A$ is a diagram in
$\A$ of the form
\begin{equation}\label{1.2}
\xymatrix{
A \ar[r]^{u}& B \ar[r]^v& C \ar[r]^w& tA
}
\end{equation}
Each such triangle yields an infinite sequence of maps
$d_n : X_{n-1} \longrightarrow X_n$, $n\in \ZZ$, in $\A$ together
with isomorphisms $\tau_X : X_n \cong tX_{n-3}$ satisfying
$(td_{n-3})\tau_X = \tau_X d_n$ such that the following diagram
commutes.

\begin{equation}\tag{1}
\xymatrix{
... \ar[r] & X_{-3} \ar[r]& X_{-2} \ar[r]^{d_{-1}} \ar@{=}[d]& X_{-1}
\ar[r]^{d_0} \ar@{=}[d]& X_0 \ar[r]^{d_1}\ar@{=}[d] & X_1
\ar[r]\ar@{=}[d] & X_2 \ar[r]& ...\\
    &        &  A   \ar[r]^u  & B  \ar[r]^v    &      C \ar[r]^w& tA
&
}
\end{equation}
We define $X_n$ for $n<-2$ by $X_n = t^{-1}X_{n+3}$ and
for $n>1$ by $X_n = tX_{n-3}$. Then $tt^{-1}\cong 1$ yields the
isomorphism $X_n \cong tX_{n-3}$ for $n\leq 0$. If $X=(X_n, d_n)$
is a cochain complex, that is $d_nd_{n-1}=0$ for all $n$, then $X$
is a {\em candidate triangle}\index{triangle!candidate} in $\A$,
see \cite{Ne}.
A morphism between triangles is a commutative diagram\\

\begin{equation}\tag{2}
\xymatrix{
 A   \ar[r]  \ar[d]_{f}& B  \ar[r]   \ar[d]^{g} &      C
\ar[r]\ar[d]^{h} & tA   \ar[d]^{tf} \\
  A'   \ar[r]  & B'  \ar[r]    &      C' \ar[r]& tA'
}
\end{equation}
This is an isomorphism if $f,g,h$ are isomorphisms in $\A$. In the
case of candidate triangles such a morphism induces a chain map $
\alpha : X \longrightarrow X'$ between the associated cochain
complexes with $\alpha_{-2}=f$ , $\alpha_{-1}=g$, $\alpha_{0}=h$
and $\alpha_{n}=t\alpha_{n-3}$ for $n \in \ZZ$. Following
\cite{Pu} and \cite{Ve}, see also \cite{We} and \cite{GM}, we have
the following notion of a triangulated category.

\begin{Def}\label{1.3}
 Let $\A$ be an additive category with a translation functor $t$. Then
$\A$ is a {\em triangulated category}\index{triangulated category}
if $\A$ is equipped with a distinguished family $\mathcal{E}$ of
triangles (called the {\em exact}\index{triangle!exact} triangles)
which is subject to the following axioms:

\begin{itemize}
\item[(Tr0)] Any triangle isomorphic to an exact triangle is exact.
For any object $A$ in $\A$ the identity $1_A$ of $A$ yields the
exact triangle:

\begin{equation*}
\xymatrix{
 A   \ar[r]^{1} & A  \ar[r]   &  \ast  \ar[r] & tA .
}
\end{equation*}
where $\ast$ is the zero object of $\A$.

\item[(Tr1)] Any morphism in $A \longrightarrow B$ in $\A$ is part of
an exact triangle:

\begin{equation*}
\xymatrix{
 A   \ar[r]^{u} & B  \ar[r]^v   &  C  \ar[r]^w & tA .
}
\end{equation*}
Here the axiom (Tr3) below implies that the isomorphism
type of $C$ is well defined by $u$. We also write $C=C_u$ and we
call $v= i_u : B \longrightarrow C_u$ the {\em
inclusion}\index{inclusion} and $w=q_u : C_u \longrightarrow tA$
the {\em projection}\index{projection}.

\item[(Tr2)] If a triangle $ A \stackrel{u}{\longrightarrow} B
\stackrel{v}{\longrightarrow} C \stackrel{w}{\longrightarrow} tA $
is exact then the triangles $ B \stackrel{v}{\longrightarrow} C
\stackrel{w}{\longrightarrow} tA \stackrel{-tu}{\longrightarrow}
tB$ and $t^{-1}C \stackrel{-t^{-1}w}{\longrightarrow} A
\stackrel{u}{\longrightarrow} B \stackrel{v}{\longrightarrow} C$
are exact.

\item[(Tr3)] Any commutative diagram
\begin{equation*}
\xymatrix{
 A   \ar[r] \ar[d]& B  \ar[r]   \ar[d]& C  \ar[r] & tA \\
A' \ar[r]&B' \ar[r] &C'\ar[r]& tA'
}
\end{equation*}
where the rows are exact extends to a morphism of
triangles as defined in (\ref{2.1}).

\item[(Tr4)] {\em The octahedral axiom}\index{octahedral axiom}:
For a composite $ A \stackrel{u}{\longrightarrow} B
\stackrel{v}{\longrightarrow} C $ let $C_u$, $C_v$ and $C_{vu}$ be chosen as in
(Tr1). Then there exists an exact triangle
\begin{equation}\tag{\#}
\xymatrix{ C_{vu} \ar[r]^{\bar{u}}& C_v \ar[r]^{\bar{v}}
&tC_u\ar[r]^{-t\bar{w}}& tC_{vu} }
\end{equation}
with the following properties. The maps $\bar{u}$ and
$\bar{w}$ are maps for which the diagram
\begin{equation}\tag{\#\#}
\xymatrix{
B \ar[r]^{v} \ar[d]_{i_u}& C \ar@{=}[r]  \ar[d]^{i_{uv}} &  C
\ar[d]^{i_v}\\
C_u \ar[r]^{\bar{w}} \ar[d]_{q_u}& C_{vu}
\ar[r]^{\bar{u}}\ar[d]^{q_{vu}}&C_v\ar[d]^{q_v}\\
tA \ar@{=}[r]& tA \ar[r]^{tu}& tB
}
\end{equation}
commutes and $\bar{v}$ is the composite
\begin{equation*}
\xymatrix{
\bar{v}: C_v \ar[r]^{q_v}& tB \ar[r]^{ti_u} &tC_u
}
\end{equation*}
\end{itemize}
\end{Def}

Axiom (Tr4) is a reformulation of Verdier's octahedral
axiom as described in Weibel \cite{We}. This axiom does not say
that all maps $\bar{u}$, $\bar{w}$ satisfying ($\#\#$) yield an
exact triangle ($\#$). It is only required that there exist such
maps. Neeman \cite{Ne}, \cite{Ne1} raises the question to
characterize maps $\bar{u}$, $\bar{w}$ for which ($\#$) is an
exact triangle, compare section 9 below.

It is easy to see that (Tr0), (Tr2) and (Tr3) imply that for all objects $U$ in $\Ct{A}$ the induced sequences of abelian groups
$$\leftarrow\homa(X_0,U)\longleftarrow\homa(X_1,U)\longleftarrow\homa(X_2,U)\leftarrow$$
$$\rightarrow\homa(U,X_0)\longrightarrow\homa(U,X_1)\longrightarrow\homa(U,X_2)\rightarrow$$
are long exact sequences.

\section{Additive track categories}
\setcounter{equation}{0}

A {\em track category}\index{track category} is a category
enriched in groupoids; in particular, for all of its objects $X$,
$Y$ their hom-groupoid $\ll X,Y\rr$ is given, whose objects are
maps $f:X \rightarrow Y$ and whose morphisms, denoted $\alpha : f
\Rightarrow f'$, are called tracks. Equivalently, a track category
is a 2-category all of whose 2-cells are invertible. For a track
$\alpha : f \Rightarrow f'$ and maps $g: Y \rightarrow Y'$,
$e:X'\rightarrow X$, the resulting composite tracks will be
denoted by $g\alpha : gf \Rightarrow gf'$ and $\alpha e : fe
\Rightarrow f'e$. Moreover there is a vertical composition of
tracks, i.e. composition of morphisms in the groupoids $\ll
X,Y\rr$; for $\alpha  : f \Rightarrow f'$ and $\beta : f'
\Rightarrow f''$ it will be denoted by $\be\square\al :
f\Rightarrow f'$. An inverse of a track $\alpha$ with respect to
this composition will be denoted $\alpha^{\boxminus}$. The
identity track $f \Rightarrow f$ is denoted by
${0}_{f}^{\square}=0^{\square}$ and is called the trivial track. A
track category has the {\em homotopy category} which is an
ordinary category obtained by identifying {\em homotopic} maps $f
\cong f'$, i.e. maps $f$, $f'$ for which there exists a track $f
\Rightarrow f'$. If $\A$ is the homotopy category then we denote
the associated track category by

\begin{equation}\label{2.1}
\B=(\B_1 \rightrightarrows \B_0 \rightarrow \A)
\end{equation}
Here $\B_0$ is the underlying ordinary category of
1-cells. We have the quotient functor $\B_0\rightarrow \A$ which
is the identity on objects and which carries a 1-cell $f$ to its
homotopy class $\{ f \}$. Moreover $\B_1$ denotes the category of
2-cells with the same objects as in $\B_0$ but with morphisms from
$X$ to $Y$ being tracks $\al : f \Rightarrow f'$ with $f$, $f':X
\rightarrow Y$ in $\B_0$, composite of $\alpha$ and $\beta$ in the
diagram

\begin{equation*}
\xymatrix{
Z \ar@{}[r]|{\Downarrow \alpha}& \ar@/^1pc/[l]^{f'} \ar@/_1pc/[l]_{f} Y
\ar@{}[r]|{\Downarrow \beta}& X  \ar@/^1pc/[l]^{g'} \ar@/_1pc/[l]_{g}
}
\end{equation*}
being
\begin{equation}\label{2.2}
\alpha\beta = \alpha g' \square f\beta = f' \beta \square \alpha g
: fg \Rightarrow fg'.
\end{equation}

There are thus two functors $\B_1 \rightarrow \B_0$
which are identity on objects and which send a morphism $\alpha :
f \Rightarrow f'$ to $f$, resp. $f'$.

An object $X\times Y$ in $\B$ is a {\em strong
product}\index{strong product} if projections $p_1 :X\times Y
\rightarrow X$ and $p_2 :X\times Y \rightarrow Y$ are given which
induce for all objects $U$ in $\B$ an isomorphism of groupoids

\[
(p_{1*}, p_{2*}): \ll U, X \times Y \rr \rightarrow \ll U,X
\rr\times\ll U,Y\rr
\]

If this is an equivalence of groupoids then $X\times Y$
is a weak product.

We now assume given a track category $\B$ such that its
homotopy category is an additive category like $\A$ from section 1,
\[
\B_{\simeq}=\A
\]
and that moreover $\B$ has a strong zero object, that is, an
object * such that for every object $X$ of $\B$, $\ll X,*\rr $ and
$\ll *,X\rr $ are trivial groupoids with a single morphism. It
then follows that in each $\ll X,Y\rr$ there is a distinguished
map $0_{X,Y}$ obtained by composing the unique maps $X
\rightarrow *$ and $*  \rightarrow Y$. The identity track of this
map will be denoted just by $0^{\square}$. Note that $0_{X,Y}$ may
also admit non-identity self-tracks; one however has

\begin{equation}\label{2.3}
0_{Y,Z}\beta= 0^{\square}=\alpha  0_{X,Y}
\end{equation}
for any $\alpha : f \Rightarrow f',f,f':Y  \rightarrow Z,\beta :g
\Rightarrow g',g,g':X'  \rightarrow  Y$.

 The secondary analogue
of an additive category is an additive track category defined as
follows.

\begin{Def}\label{2.4}
A track category $\B$ is called {\em additive}\index{track
category! additive} if it has a strong zero object *, the homotopy
category $\A=\B_{\simeq}$ is additive and moreover $\B$ is a
linear track extension

\[
D\rightarrow\B_1 \rightrightarrows \B_0\rightarrow\A
\]
of $\A$ by a biadditive bifunctor
\[
D:\A^{op}\times \A \rightarrow \Ab .
\]
Explicitly, this means the following: a biadditive
bifunctor $D$ as above is given together with a system of
isomorphisms

\begin{equation}\label{2.5}
\sigma_f :D(X,Y)  \rightarrow \Aut_{\ll X,Y\rr } (f)
\end{equation}
for each 1-arrow $f: X  \rightarrow Y$ in $\B$, such
that for any $f:X  \rightarrow Y$, $g:Y \rightarrow Z$, $a\in
D(X,Y)$, $b \in D(Y,Z)$, $\alpha: f \Rightarrow f'$ one has

\begin{equation*}
\begin{aligned}
\sigma_{gf}(ga)&=g\sigma_f (a);\\
\sigma_{gf}(bf)&=\sigma_g (b)f;\\
\al\square \sigma_f (a)&= \sigma_{f'}(a)\square \al.
\end{aligned}
\end{equation*}
\end{Def}

\begin{Rem}\label{2.6}
Using 2.5 we can identify the bifunctor $D$ via the natural
equation
\[
D(X,Y)=\Aut(0_{X,Y}).
\]
Here $\Aut(0_{X,Y})$ is easily seen to be a well defined
bifunctor.
\end{Rem}

We use a translation functor $t$ as in (1.1) for the definition of
the {\em coefficient bifunctor}

\begin{equation}\label{2.7}
\begin{aligned}
& D:\A^{op}\times\A\rightarrow\Ab\\
& D(X,Y)=hom_{\A}(tX,Y)
\end{aligned}
\end{equation}
Then $t$ yields a natural isomorphism
\[
t_D : D \stackrel{\cong}{\longrightarrow} D(t^{op}\times t)=t^* D
\]
which on objects is given up to sign by the functor $t$,
that is
\[
t_{D}:D(X,Y)=hom_{\A}(tX,Y)\cong hom_{\A}(ttX,tY)=(t^*D)(X,Y)
\]
is defined by $t_{D}(\alpha)=-t(\alpha).$

\begin{Def}\label{2.8}
Let $\B$ be an additive track category with homotopy
category $\A$ and coefficient bifunctor $D(X,Y)=hom_{\A}(tX,Y)$
where $t$ is a strict translation functor on $\A$. Then a {\em
strict translation functor for $\B$} is a track functor (i.e. a
functor enriched in groupoids) as in the commutative
diagram

\begin{equation*}
\xymatrix{ D \ar[r]\ar[d]_{t_D} & \B_1  \ar[d]^{t^{(1)}}
\ar@<2pt>[r]
\ar@<-2pt>[r]& \B_0  \ar[r]\ar[d]^{t^{(0)}}& \A \ar[d]^{t}\\
D  \ar[r]& \B_1 \ar@<2pt>[r] \ar@<-2pt>[r]& \B_0  \ar[r]& \A }
\end{equation*}
such that $t^{(0)}$, $t^{(1)}$ are strict automorphisms of
categories inducing $t$ and $t_D$ respectively. We also write $t=t^{(0)}$
and $t=t^{(1)}$.
\end{Def}

\begin{exa}\label{2.9}
The category $Ch(\A)$ of chain complexes in an additive category
$\A$ is a track category with tracks given by homotopy classes of
homotopies. In fact $Ch(\A)$ is an additive track category with a
strict translation functor given by the shift functor.
\end{exa}

\begin{exa}\label{2.10}
The category of spectra {\em Spec} as introduced by Puppe
\cite{Pu} or Adams \cite{Ma} is a track category with tracks
defined by homotopy classes of homotopies. Moreover {\em Spec} is
an additive track category with strict translation functor given
by the shift functor.
\end{exa}

\begin{Rem}\label{2.11}
The main results below are proved in the presence of a strict
translation functor for $\B$. It is also possible to deal with with translation
functors for $\B$ which are not strict. Such (non strict) translation
functors are given by a pseudo functor $t: \B \rightarrow \B$ inducing
$t$ on $\A$ and $-t$ on $D$ as in (2.8). 
In \cite{ccctc} we study the more general case of non strict translation functors. It will be, however, more convenient for the
reader to consider first the case of strict translation functors as this is done in this paper.
\end{Rem}

\begin{exa}\label{2.12}
Let $\C$ be a cofibration category with a zero object, for
example given by a Quillen model category, see \cite{Ba}. Assume that
the suspension functor

\[
\Sigma : \C_{cf}/\simeq \rightarrow \C_{cf}/ \simeq
\]
is an equivalence of categories. Here $\C_{cf}$ is the full
subcategory of cofibrant and fibrant objects in $\C$. Then
$\C_{cf}$ is an additive track category with a (non strict)
translation functor $\Sigma$ as in (\ref{2.11}). The dual result
holds for fibration categories. Tracks are defined as in
\cite{Ba}(II.5).
\end{exa}

\section{Cone functors in additive track categories}
\setcounter{equation}{0}

In topology one defines the mapping cone $C_f$ of a map
$f:A\rightarrow B$ by the push out $B \cup_f CA$ of $CA \supset A
\rightarrow B$. Here $CA$ is the cone of $A$. Given a map $\alpha
: B \rightarrow U$ and a homotopy $\hat{\alpha}: \alpha f \simeq
0$ one gets a map $\alpha \cup \hat{\alpha}: C_f \rightarrow U$.
Therefore the set of homotopy classes of maps $C_f \rightarrow U$
denoted by $Cone_f (U)$ can be described only in terms of tracks
which are homotopy classes of homotopies. This leads to the
following definition of cone functors in track categories as
introduced by Gabriel-Zisman, see the remark (3.7) below.

Let $\B$ be an additive track category as in (\ref{2.4}) with
homotopy category $\A$. Let $\Ab$ be the category of abelian
groups. Given a morphism $f:X\rightarrow Y$ in $\B_0$ we define
the {\em cone functor}\index{cone functor}

\begin{equation}\label{3.1}
Cone_f : \A \longrightarrow \Ab
\end{equation}
and the {\em dual cone functor}\index{cone functor!dual}
\[
Cone^f : \A^{op} \longrightarrow \Ab
\]
as follows. For an object $U$ in $\A$ let $cone_f (U)$
be the set of all pairs $(\alpha , \hat{\alpha})$ where $\alpha :
Y \rightarrow U$ is a map in $\B_0$ and

\[
\hat{\alpha} : \alpha f \Rightarrow 0
\]
is a track in $\B$. Two such pairs $(\alpha ,
\hat{\alpha})$ and $(\beta , \hat{\beta})$ in $cone_f (U)$ are
equivalent if there exists a track $H: \alpha \Rightarrow \beta$
such that pasting of tracks in the diagram

\begin{equation*}
\xymatrix{
 && 0 && \\
X \ar[rr]^f \ar@/^2.5pc/[rrrr] \ar@/_2.5pc/[rrrr] && Y
\ar@/^1pc/[rr]^{\alpha} \ar@/_1pc/[rr]_{\beta} \ar@{}[u]|{\Uparrow \hat{\alpha}}
\ar@{}[d]|{\Downarrow \hat{\beta}} \ar@{}[rr]|{\Downarrow H}&& U\\
 && 0 &&
}
\end{equation*}
yields the trivial track $0^{\square}$. Now let
\begin{equation}\tag{1}
Cone_f(U)=cone_f (U)/ \sim
\end{equation}
be the set of equivalence classes $\{ \alpha
,\hat{\alpha}\}$. Given a map $v:U\rightarrow V$ in $\A$ we choose
$v_0$ in $\B_0$ representing $v$ and we define the induced map
\[
v_* : Cone_f (U) \longrightarrow Cone_f (V)
\]
by $v_* : \{ \alpha ,\hat{\alpha}\} = \{ v_0\alpha
,v_0\hat{\alpha}\}$. One readily checks that $v_*$ is well
defined. Below we show that $Cone_f (U)$ is an abelian group.

Similarly let $cone^f (U)$ be the set of all pairs $(\alpha
,\hat{\alpha})$ where $\alpha : U \rightarrow X$ is a map in $\B_0$ and
$\hat{\alpha} : f\alpha \Rightarrow 0$ is  a track in $\B$. Now two
such pairs $(\alpha ,\hat{\alpha})$, $(\beta ,\hat{\beta})$ are equivalent
if there exists a track $G:\alpha \Rightarrow\beta$ such that pasting
of tracks in the diagram

\begin{equation*}
\xymatrix{
&& 0 && \\
U \ar@/^1pc/[rr]^{\alpha} \ar@/_1pc/[rr]_{\beta}  \ar@{}[rr]|{\Uparrow
G}  \ar@/^2.5pc/[rrrr] \ar@/_2.5pc/[rrrr] && X  \ar[rr]^f
\ar@{}[u]|{\Uparrow \hat{\alpha}} \ar@{}[d]|{\Downarrow \hat{\beta}} && Y\\
 && 0 &&
}
\end{equation*}
yields the trivial track $0^{\square}$. As above let
\begin{equation}\tag{2}
Cone^f (U)= cone^f (U)/ \sim
\end{equation}
be the set of equivalence classes $\{ \alpha
,\hat{\alpha}\}$. Then a map $v:V \rightarrow U$ in $\A$
represented by $v_0$ in $\B_0$ yields the induced map
\[
v^* : Cone^f (U) \longrightarrow Cone^f (V)
\]
by $v^* \{ \alpha ,\hat{\alpha}\} = \{ \alpha
v_0,\hat{\alpha}v_0\}$.

A track $H:f' \Rightarrow f$ in $\B$ induces natural isomorphisms
\begin{equation}\tag{3}
\begin{aligned}
H_* : Cone_f &\cong Cone_{f'} \;\;\;\; , resp.\\
H_* : Cone^f &\cong Cone^{f'}
\end{aligned}
\end{equation}
which carries $\{ \alpha ,\hat{\alpha}\}$ to $\{ \alpha
,\hat{\alpha}\square\alpha H\}$, resp. to $\{ \alpha
,\hat{\alpha}\square H\alpha \}$. Hence up to natural isomorphism
cone functors depend only on the homotopy class of $f$.

\begin{Lem}\label{3.2}
Let $f: X \rightarrow Y$ be a map in $\B_0$ and let $U,W$ be objects in
$\A$ with direct sum $U\oplus W$ and projections $p_1 : U\oplus W
\rightarrow U$, $p_2 : U\oplus W \rightarrow W$ in $\A$. Then
\[
(p_{1*}, p_{2*}): Cone_f (U\oplus W) \rightarrow Cone_f (U) \times
Cone_f (W)
\]
is a bijection of sets.
\end{Lem}

\begin{proof}
Let $\{ \al ,\hat{\al}\} \in Cone_f (U) ,  \{ \beta ,\hat{\beta}\}
\in Cone_f (W)$. Then there exists $\gamma$ together with tracks
$Q_1 : p_1\ga \R \al$, $Q_2 : p_2\ga\R\be$ since $U\oplus W$ is a
direct sum in $\A$. Moreover there exists $H:\ga f \R 0$ since
$(\al ,\be)f=0$ in $\A$. Now $\hal\square Q_1 f$ and $p_1 H$
differ by $h_1\in D(X,U)$ and $\hbe\square Q_2 f$ and $p_2 H$
differ by $H_2\in D(X,W)$. Therefore we can alter $H$ by
$(H_1,H_2)\in D(X,U\oplus W)$ and get $\hga$ such that
$(p_1)_*\{\ga , \hga\}=\{\al ,\hal\}$ and $(p_2)_*\{\ga ,
\hga\}=\{ \beta ,\hat{\beta}\}$. Hence the map in (\ref{6.1}) is
surjective. Now let $\{\ga , \hga\}$, $\{\veps ,\hveps\}\in Cone_f
(U\oplus W)$ be given with
\begin{equation}\tag{1}
(p_i)_* \{\ga , \hga\}=(p_i)_*\{\ep ,\hep\}\;\;\;\mathrm{for}\;i=1,2.
\end{equation}

We have to show $\{\ga , \hga\}=\{\veps ,\hveps\}$, that is, there
is a track $H:\ga\R\veps$ with $\hveps\square (Hf)=\hga$. By (1)
we obtain tracks $H_i : p_i\ga\R p_i\veps$ with $(p_i
\hveps)\square(H_i f)=p_i \hga$. Hence $\ga$ and $\veps$ represent
the same homotopy class so that there exists a track
$G:\ga\R\veps$. Now $H_i$ and $p_iG$ differ by an element $D_i$
with $D_1\in D(X,U)$, $D_2\in D(X,W)$. We can alter $G$ by $D=(D_1
,D_2)\in D(X,U\oplus W)$ and get $H$ with $p_iH=H_i$. Moreover

\begin{equation}\tag{2}
p_i(\hveps\square Hf)=p_i\hga
\end{equation}

Now $\hveps\square Hd$ differs from $\hga$ by an element
$A\in D(Y, U\oplus W)$. Then (2) shows $A=0$. Hence $\hveps\square
Hf = \hga$ is satisfied.

 \end{proof}

\begin{cor}\label{3.3}
The sets $Cone_f (U)$ and $Cone^f (U)$ have the structure of abelian
groups and the functors (\ref{3.1}) are well defined additive functors
\end{cor}

\begin{proof}
Using the isomorphism (\ref{3.2}) we define addition by the composite
\[
Cone_f (U)\times Cone_f (U)\cong Cone_f (U\oplus U)
\stackrel{(1,1)_*}{\lra} Cone_f (U)
\]
where $(1,1):U\oplus U\ra U$ is the folding map. We
obtain for $Cone^f$ a result similar to (\ref{3.2}) and define the
addition in $Cone^f (U)$ in a similar way.

 \end{proof}

\begin{prop}\label{3.4}
Let $\B$ be an additive category with $\B_{\simeq}=\A$ and bifunctor
$D: \A^{op}\times \A \ra \Ab$. Then the cone functors are embedded in
natural exact sequences in $\Ab$

\begin{equation*}
\xymatrix{
D(Y,U) \ar[r]^{f^*}& D(X,U) \ar[r]^{i}& Cone_f (U) \ar[r]^{q}& hom_{\A}
(Y,U) \ar[r]^{f^*}& hom_{\A} (X,U)
}
\end{equation*}

\begin{equation*}
\xymatrix{
D(U,X) \ar[r]^{f_*}& D(U,Y) \ar[r]^{\bar{i}}& Cone^f (U)
\ar[r]^{\bar{q}}& hom_{\A} (U,X) \ar[r]^{f_*}& hom_{\A} (U,Y)
}
\end{equation*}
where $f:X\ra Y$ is a map in $\B_0$. Moreover $i,q$ and
$\bar{i}, \bar{q}$ are compatible with $H_*$ in (\ref{3.1})(3),
that is, $H_*i=i$, $qH_* =q$.
\end{prop}

\begin{proof}
We define $q\{\al ,\hal\}=\{\al\}$. Then $f^*q=0$ since $\hal :\al f \R
0$. On the other hand any element $\{\al\}\in hom_{\A} (Y,U)$ with
$f^*\{\al\}=0$ admits a track $\hal :\al f\R 0$ representing an element
$\{\al ,\hal\}\in Cone_f (U)$ with $q\{\al ,\hal\}=\{\al\}$. We define
$i(\delta)$ by the track $\sigma\delta :0\R 0$ in $\B_1$, see (\ref{2.5}),
so that we get the element $i(\delta)=(0,\sigma\delta)$ by the diagram

\begin{equation*}
\xymatrix{
 &&&& \\
 &&&& \\
X \ar[rr]_{f} \ar@/^1pc/[rrrr]^0 \ar@/^3.5pc/[rrrr]^0&& Y \ar[rr]_{0}
\ar@{}|{\Uparrow \sigma\delta}[uu]&&U
}
\end{equation*}
with $0f=0$ since $*$ is a strict zero object of $\B$.
Cleary $iq=0$. Moreover assume $q\{\al ,\hal\} =0$. Then there
exists a track $A:0\R \al$ and $\{\al ,\hal\}=\{ 0, \hal\square
Af\}$. Moreover by (\ref{2.5}) there is $\delta$ with $\hal\square
Af : 0\R 0$ being $\sigma (\delta)$. Next assume that $\delta =
f^* \delta '$ with $\delta ' \in D(Y,U)$. Then $\{ 0, \sigma\delta
\}=\{ 0,0^{\square}\}$ is the trivial element since $\sigma\delta$
coincides with

\begin{equation*}
\xymatrix{
 &&&& \\
 &&&& \\
X \ar[rr]  \ar@/^3.5pc/[rrrr]&& Y
\ar[rr]_{0}^{\Uparrow\sigma\delta '} \ar@{}|{\Uparrow
0^{\square}}[uu] \ar@/^2pc/[rr]_0&&U }
\end{equation*}

This shows $if^*=0$. Moreover let $\delta$ be given with $i\delta
= \{ 0,\sigma\delta\} =0$. Then there exixts a track $A: 0\R 0$
with $\sigma\delta =Af$. Here $A=\sigma\delta '$ by (\ref{2.5}).
Hence $\delta = f^*\delta '$.

In a similar way one proves exactness of the second sequence
in (\ref{3.4}) with

\begin{equation*}
\begin{aligned}
\bar{q}\{\al ,\hal\}=\{\al\} & \;\;\;\;\mathrm{for}\;\{\al ,\hal\}\in
Cone^f (U),\\
\bar{i}\delta =\{ 0,\sigma\delta\} & \;\;\;\;\mathrm{for}\;\delta\in
D(U,Y).
\end{aligned}
\end{equation*}

 \end{proof}

\begin{cor}\label{3.5}
Let $\B$ be an additive track category with $\B_{\simeq}=\A$ and with a
bifunctor $D(X,Y)= hom_{\A}(tX,Y)$ as in (\ref{2.7}). Let $f:X\ra Y$ be
a map in $\A$ and let $f(n):t^n X\ra t^nY$ be a map in $\B_0$
representing the homotopy class $t^n\{ f\}$. Then one has long exact sequences
in $\Ab$ as follows, $n\in \ZZ$.

\begin{small}
\begin{equation*}
hom_{\A}(t^{n+1}Y,U) \stackrel{f^*}{\ra} hom_{\A}(t^{n+1}X,U)
\stackrel{i}{\ra}Cone_{f(n)} (U)\stackrel{q}{\ra}  hom_{\A} (t^nY,U)
\stackrel{f^*}{\ra} hom_{\A} (t^nX,U)
\end{equation*}
\end{small}

\begin{small}
\begin{equation*}
hom_{\A}(U,t^{n+1}X) \stackrel{f^*}{\ra} hom_{\A}(U,t^{n+1}Y)
\stackrel{\bar{i}}{\ra}Cone^{f(n)} (U)\stackrel{\bar{q}}{\ra}  hom_{\A} (U,t^nX)
\stackrel{f_*}{\ra} hom_{\A} (U,t^nY)
\end{equation*}
\end{small}
\end{cor}

\begin{proof}
\[
D(X,Y) = hom_{\A}(tX,Y)\cong\homa(X,t^{-1}Y)
\]
and we can apply (\ref{3.4}).

 \end{proof}

\begin{prop}\label{3.6}
Let {\B} be an additive track category with $\B_{\simeq}=\A$ and with a
bifunctor $D(X,Y)=\homa (tX,Y)$ as in (\ref{2.7}). Let
\begin{equation*}
\xymatrix{
A\ar[r]^u & B \ar[r]^v & C
}
\end{equation*}
be maps in $\B_0$ and let be $X$ an object in $\A$. Then we have a
long exact sequence in $\Ab$.
\[
Cone_{vu}(t^{-1}X)\stackrel{v^*}{\ra}Cone_{u}(t^{-1}X)\stackrel{\delta}{\ra}Cone_{v}(X)\stackrel{u^*}{\ra}Cone_{vu}(X)\stackrel{v^*}{\ra}Cone_{u}(X)
\]
Here we set $v^*\{\al ,\hal\}=\{\al v,\hal\}$ and
$u^*\{\be ,\hbe\} =\{\be ,\hbe u\}$. Moreover $\delta$ is defined
by $\delta\{\ga ,\hga\} =it\{\ga\}$ where $\ga :B\ra t^{-1}X$ and
$t\{\ga\}\in\homa (tB,X)=D(B,X)$ and $i$ is given as in
(\ref{3.4}).
\end{prop}

We point out that the exact sequence (\ref{3.6}) is the
analogue of the exact triangle (\#) in the octahedral axiom (\ref{1.3}). The
operators $v^*,\delta,u^*$ in (\ref{3.6}), however, are canonically
given while in (\ref{1.3})(\#) the maps $\bar{u},\bar{w}$ are certain
choices.

\begin{proof}[Proof of (\ref{3.6})]
For $\{\be ,\hbe\}\in Cone_v (X)$ with $\hbe :\be v\R 0$ we get
$\hbe u:\be vu\R 0$ so that $v^*u^*\{\be ,\hbe\} =\{\be vu,\hbe
u\} =0$ in $Cone_u (X)$.

Now let $\{\al ,\hal\}\in Cone_{vu} (X)$ with $\hal : \al vu\R 0$
and assume $u^*\{\al ,\hal \}=\{\al v,\hal\} =0$ in $Cone_u (X)$.
Then there exists $z:\al v\R 0$ with $zu=\hal$. Hence $\{\al
,\hal\} =\{\al , zu\} = u^*\{\al ,z\}$. Hence we proved that $im
(u^*)= ker (v^*)$.

Next let $\{\ga ,\hga\}\in Cone_u (t^{-1}X)$ with $\hga :\ga
u \R 0$. Let $z:0\R 0$ be given by $\sigma t\{\ga\}$ where $0:B\ra *
\ra X$. Then $u^*\delta\{\ga ,\hga\} = u^*\{ 0,z\} =\{ 0,zu\} =\{ 0,
\sigma(t\{\ga u\} )\} =0$ since $\ga u =0$ in $\A$ by the track $\hga$. Now
let $\pbe\in Cone_{v} (X)$ with $u^*\pbe = \{\be ,\hbe u\}=0$ in
$Cone_{vu}(X)$. Then there exixts $z:\be\R 0$ with $\hbe u = zvu$. Hence we
get an element
\[
\sigma (y)= \hbe \square (zv)^{\boxminus} : 0\R 0
\]
with $y\in D(B,X)$ and $u^* y =0$. Thus we can find an
element $\bar{y}\in Cone_u (t^{-1} X)$ with $\delta\bar{y}=\pbe$.
This completes the proof that $im(\delta) = ker (u^*)$.

Finally let $\prho \in Cone_{vu} (t^{-1}X)$ with
$\hat{\rho} : \rho vu \R 0$ and $\rho : C \ra t^{-1}X$. Then
$\delta v^* \prho = \delta \{\rho v ,\hrho\} = it\{\rho v\} = \{
0, \sigma (t(\rho))v\} =0$ in $Cone_v (X)$.

Moreover let $\pga\in Cone_u (t^{-1}X)$ with $\hga : \ga u \R 0$
and $\delta\pga =0$ in $Cone_v (X)$. Then there exists $z\in
D(C,X)$ such that in the diagram
\begin{equation*}
\xymatrix{
 &&&& \\
X &\ar@{}|{\Downarrow \sigma z}[d]& C \ar@{}|{\Uparrow \sigma t
\{\ga\}}[u] \ar[ll]^{0} \ar@/^2pc/[ll]^0   && B \ar[ll]^{v}
\ar@/_3.5pc/[llll]_0\\
&&&&\\
}
\end{equation*}
we have $\sigma (z)v=\sigma t\{\ga\}$. Let $\rho : C \ra
t^{-1}X$ be a map in $\B_0$ which represents $t^{-1}z$. Then there
exixts $\hrho :\rho vu \R 0$ since $\ga u =0$ in $\A$. This
completes the proof that $im(v^*) = ker (\delta)$.

 \end{proof}
\begin{Rem}
Gabriel-Zisman consider cone functors $Cone^{f}$ in a track category $\B$ with zero
object * as follows. Given an object $U$ and a map $f:X\longrightarrow Y$ in $\B$ one
gets the morphism
\begin{equation*}
f^{U}: \ll U,X \rr\longrightarrow \ll U,Y\rr
\end{equation*}
between pointed groupoids induced by $f$. Then the groupoid $\Gamma(f^{U})$ is
defined by the objects $(\alpha,\hat{\alpha})$ and by the tracks $G:(\beta,\hat{\beta})
\Longrightarrow (\alpha,\hat{\alpha})$ satisfying the condition in (3.1)(2). Hence
\begin{equation*}
Cone^{f}(U)=\pi_{0}(\Gamma(f^{U}))
\end{equation*}
is the set of path components of $\Gamma(f^{U})$. Here
the groupoid $\Gamma(f^{U})$ coincides with the construction V.3.1
\cite{GZ} and axiom C in V.3.2 \cite{GZ} implies that $Cone^{f}$
is representable as in (4.1) below. Gabriel-Zisman use the axiom C
for the proof of the "exact fiber sequence", see V.4.2 \cite{GZ},
which in the context of this paper is easily achieved by the dual
of (3.5) and (4.3).
\end{Rem}

\section{Representable cone functors}
 \setcounter{equation}{0}

Let $\B$ be an additive track category with homotopy category $\A$
and with bifunctor $D(X,Y) =\homa (tX,Y)$ as in (\ref{2.7}). An
additive functor $F:\A\ra\Ab$ is {\em representable}\index{cone
functor!representable} if there exists an object $X$ in $\A$ and a
natural isomorphism of functors
\begin{equation}\label{4.1}
\chi : F \cong \homa (X,-)
\end{equation}
in $\Ab$. Then $X$ is well defined up to isomorphism in
$\A$ and is called a {\em representation} of $F$. We now consider
representations of cone functors.

If $C=C_f$ represents $Cone_f$ with $f:A\ra B$ in $\B_0$ then we
obtain the triangle
\begin{equation}\label{4.2}
\xymatrix{
A  \ar[r]^f& B  \ar[r]^{i_f}&  C_f \ar[r]^{q_f}& tA
}
\end{equation}
in $\A$ as follows. For $U=C_f$ we have the operator
\begin{equation*}
\xymatrix{
\homa (C_f,C_f) \ar[r]^{\chi_{f}^{-1}} & Cone_f (C_f) \ar[r]^q & \homa
(B,C_f)
}
\end{equation*}
by (\ref{3.4}) which carries the identity of $C_f$ to
$i_f$. Moreover for $U=tA$ we have the operator
\begin{equation*}
\xymatrix{
\homa (tA,tA) \ar[r]^{i} & Cone_f (tA) \ar[r]^{\chi_f} & \homa (C_f,
tA)
}
\end{equation*}
by (\ref{3.4}) which carries the identity of $tA$ to
$q_f$.

Of course $i_f$ and $q_f$ depend on the representation ($C_f
, \chi_f$)  of $Cone_f$. Using the exact sequence (\ref{3.4}) we get by
a representation the following diagram in $\Ab$

\begin{small}
\begin{equation}\label{4.3}
\xymatrix{
\homa (tB,U) \ar[r]^{f^*} \ar@{=}[d] & \homa (tA,U) \ar[r]^i \ar@{=}[d]
& Cone_f (U) \ar[r]^q  \ar[d]_{\cong}^{\chi_f}& \homa (B,U)
\ar[r]^{f^*} \ar@{=}[d] & \homa (A,U) \ar@{=}[d]\\
\homa(tB,U) \ar[r]^{f^*} & \homa (tA,U) \ar[r]^{q_f^*} & \homa (C_f,U)
\ar[r]^{{i_f}^*} & \homa (B,U) \ar[r]^{f^*} & \homa (A,U)
}
\end{equation}
\end{small}
Here we have $i^*_f \chi_f =q$ since
\begin{equation*}
\begin{aligned}
q\chi^{-1}_{f} (\al)& = q\chi^{-1}_{f} (\al {1}_{C_f}) = \al_* q
\chi^{-1}_{f} (1_{C_f})\\
&= \al_* i_f = i_f^* (\al).
\end{aligned}
\end{equation*}

Similarly we get $\chi_f i = q^*_f$ by
\begin{equation*}
\begin{aligned}
\chi_{f} i (\be)& = \chi_{f} i(\be {1}_{t_A}) = \be_* \chi_{f}i
(1_{tA})\\
&= \be_* q_f = q_f^* \be.
\end{aligned}
\end{equation*}
Hence diagram (\ref{4.3}) commutes and therefore the bottom of
(\ref{4.3}) is an exact sequence.Compare also (1.5).

\begin{prop}\label{4.4}
If $Cone_f$ is representable by $(C_f ,\chi_f)$ then there exixts a
diagram in $\B$

\begin{equation*}
\xymatrix{
 &&&&&& \\
A \ar[rr]^f  \ar@/^3pc/[rrrr]^0 && B \ar[rr]^{i_f} \ar@{}|{\Uparrow
\hat{i}_f}[u] \ar@/_3pc/[rrrr]&& C_f \ar[rr]^{q_f}  \ar@{}|{\Downarrow
\hat{q}_f }[d]&& tA\\
&&&&&&\\
}
\end{equation*}
with the properties
\[
\chi_f \{i_f , \hat{i}_f \} = 1_{C_f},
\]

\[
\sigma^{-1}(q_f \hat{i}_f \square (\hat{q}^{\boxminus}_{f})f)=
1_{tA}\in D (A, tA).
\]
\end{prop}

\begin{proof}
We have

\[
i (1_{tA})= \chi^{-1}_{f}q^*(1_{tA})= \chi^{-1}_{f}(q_f 1_{C_f})
= q_{f^*}\{i_f ,\hat{i}_f \}
\]
Hence by definition of $i$ and the equivalence relations
in $Cone_f$ there exists $\hat{q}_f$ such that the following
diagram represents the trivial track.

\begin{equation*}
\xymatrix{
 &&&&&& \\
A \ar[rr]^f \ar@/^3pc/[rrrr]^0  \ar@/_2pc/[rrrd] && B \ar@{}|{\Uparrow
\sigma 1_{tA}}[u]\ar[rr]^0 \ar[dr]^{i_f} \ar@{}|{\Downarrow \hat{i}_f
}[d] & \ar@{}|{\Uparrow \hat{q}_f }[d]& tA\\
&& & C_{f} \ar[ur]_{q_f}& &&\\
}
\end{equation*}
This yields the second formula in (\ref{4.4}).

 \end{proof}

The converse of (\ref{4.4}) is also true. We get the
following criterion for representability of cone functors.

\begin{theo}\label{4.5}
Let $f:A\ra B$ be a map in $\B_0$. Then $Cone_f$ is representable if
and only if there exists a diagram in $\B$

\begin{equation*}
\xymatrix{
 &&&&&& \\
A \ar[rr]^f \ar@/^3pc/[rrrr]^0  && B \ar@{}|{\Uparrow \hat{i}}[u]
\ar[rr]^{i_f} \ar@/_3pc/[rrrr]_0 && C_f \ar@{}|{\Uparrow \hat{q}^{\boxminus}
}[d] \ar[rr]^{q_f} && tA\\
&& && &&\\
}
\end{equation*}
which by pasting yields $\varepsilon 1_{tA}\in D (A,
tA)$, with $\veps=1$ or $\veps =-1$ and for which the bottom row
in (\ref{4.3}) is exact for all $U$.
\end{theo}

\begin{proof}
If $Cone_f$ is representable then a diagram as in (\ref{4.5}) exists by
(\ref{4.4}). Conversely given a diagram as in (\ref{4.5}) we define the
function

\[
\chi^{-1}_{f}:\homa (C_f ,U)\ra Cone_f (U)
\]
which carries $\al$ to $\al_* \{i_f,\hat{i}\}$. Now we
obtain a diagram as in (\ref{4.3})  satisfying $i^{*}_{f} \chi_f
=q$. Moreover $\chi^{-1}_{f}q^*_f =\veps i$ since the diagram in
(\ref{4.5}) represents $\veps 1_{tA}$, compare the proof of
(\ref{4.4}). Hence diagram (\ref{4.3})  commutes up to sign and
therefore $\chi_f$ is a bijection by the five lemma.

 \end{proof}

\section{Toda complexes}

\setcounter{equation}{0} As in (\ref{2.4}) let $\B$ an additive
track category with homotopy category $\A$ and bifunctor $D$. A
{\em Toda pair}\index{Toda!pair} $(H,G)$ is a diagram in $\B$ of
the form

\begin{equation}\label{5.1}
\xymatrix{
 &&&&&& \\
Z \ar[rr]^{\ga} \ar@/^3pc/[rrrr]^0  && Y \ar@{}|{\Uparrow H}[u]
\ar[rr]^{\be} \ar@/_3pc/[rrrr]_0 && X \ar@{}|{\Downarrow G }[d] \ar[rr]^{\al}
&& W\\
&& && &&\\
}
\end{equation}

We say that $(H,G)$ is {\em associated} to the triple
$(a,b,c)$ where $a=\{\al\}$, $b=\{\be\}$ and $c=\{\ga\}$ are homotopy classes
in $\A$. Moreover we say that $(H,G)$ represents $\xi\in D(Z,W)$ or
$[H|G]\in \Aut_{\ll z,w\rr}(\al\be\ga)$ if the equation
\begin{equation}\tag{1}
\sigma (\xi)=[H|G]=\al H \square G^{\boxminus}\ga
\end{equation}
holds where we use $\sigma$ in (\ref{2.4}). We point out
that the convention in (1) uses the positive direction of $H$ to
define $\xi$. Given a triple $(a,b,c)$ in $\A$ with $ab=0$ and
$bc=0$ the {\em Toda bracket}\index{Toda!bracket}
\begin{equation}\tag{2}
<a,b,c>\subset D(Z,W)
\end{equation}
consists of all elements $\xi$ represented by Toda pairs
associated to $(a,b,c)$. This is a coset of
\[
\ga^* D(Y,W)+\al_*D(Z,X)\subset D(Z,W).
\]

We now define a {\em Toda complex}\index{Toda!complex}
by a sequence of maps $d=d_n : X_{n-1}\ra X_n$, $n\in\ZZ$, in
$\B_0$ together with a sequence of tracks $H_n$ as in the diagram

\begin{equation}\label{5.2}
\xymatrix{
 &&&&&& \\
X_{n-3} \ar[rr]_{d_{n-2}} \ar@/^3pc/[rrrr]^0  && X_{n-2}
\ar@{}|{\Uparrow H_{n-1}}[u] \ar[rr]_{d_{n-1}} \ar@/_3pc/[rrrr]_0 && X_{n-1}
\ar@{}|{\Downarrow H_n }[d] \ar[rr]_{d_n} && X_n\\
&& && &&\\
}
\end{equation}
Hence $(H_{n-1},H_n)$ is a Toda pair representing $\al_n\in
D(X_{n-3},X_n)$ for $n\in\ZZ$.

Let $X=(X_n,d_n,H_n,n\in\ZZ)$ and $X'=(X'_n,d'_n,H'_n,n\in\ZZ)$ be
Toda complexes. A {\em morphism between Toda complexes} $f:X\ra Y$
is a sequence $f=(f_n,F_n,n\in\ZZ)$ of maps $f_n$ in $\B_0$ and
tracks $F_n$ in $\B_1$ as in the diagram

\begin{equation}\label{5.3}
\xymatrix{
&&&&\\
X_{n-2} \ar[rr] \ar@/^3pc/[rrrr]^0
\ar[dd]^{f_{n-2}}&\ar@{}|{\Downarrow F_{n-1} }[dd]  & X_{n-1} \ar@{}|{\Uparrow H_{n}}[u] \ar[rr]
\ar[dd]^{f_{n-1}} & \ar@{}|{\Downarrow F_n }[dd]& X_{n}  \ar[dd]^{f_n}\\
&& && \\
X'_{n-2} \ar[rr] \ar@/_3pc/[rrrr]_0 &&X'_{n-1}\ar[rr]
\ar@{}|{\Downarrow H'_n }[d]&&X'_n\\
&&&&
}
\end{equation}
such that this diagram represents the trivial track for
all $n$, that is
\[
f_n H_n= H'_n f_{n-2}\square d'_nF_{n-1}\square F_n d_{n-1}.
\]

Let $\Tod$ be the category of Toda complexes and morphisms as in
(\ref{5.3}). We have the {\em shift functor}\index{shift functor}

\begin{equation}\label{5.4}
sh:\Tod\ra\Tod
\end{equation}
which is a strict isomorphism of categories defined as
follows. For $X=(X_n,d_n,H_n,n\in\ZZ)$ in $\Tod$ let
$sh(X)=(X'_n,d'_n,H'_n,n\in\ZZ)$ be given by $X'_n=X_{n+1},
d'_n=d_{n+1},H'_n=H_{n+1}$.

We also define the category of rows denoted by $\row$ as
follows. A row in $\B_0$ is a sequence of maps $(X_n,d_n:X_{n-1}\ra
X_n,n\in\ZZ)$. A map $f:(X_n,d_n)\ra(X'_n,d'_n)$ between rows is a sequence
$f=(f_n,F_n,n\in\ZZ)$ of diagrams in $\B$

\begin{equation}\label{5.5}
\xymatrix{
X_{n-1} \ar[rr]^{d_{n}} \ar[dd]_{f_{n-1}}&  \ar@{}|{\Downarrow
F_{n}}[dd]& X_{n} \ar[dd]^{f_n}\\
&& && \\
X'_{n-1} \ar[rr]_{d'_{n}}&& X'_n
}
\end{equation}

Two such maps $f,f':X\ra X'$ are {\em homotopic} if there
exists  tracks $G_n:f_n\R f'_n$ such that for $n\in\ZZ$

\[
(d'_n G_{n-1})\square F_n =F'_n\square (G_nd_n).
\]

A {\em weak equivalence} in $\row$ is a map $f$ for
which all $f_n$ are isomorphisms in the homotopy category $\A$.

We have the shift functor also on the category $\row$ and the
faithful forgetful functor
\begin{equation}\label{5.6}
\Tod\lra\row
\end{equation}
carries $X=(X_n,d_n,H_n)$ to $(X_n,d_n)$. Then we say
that the Toda complex $X$ {\em extends} the row $(X_n,d_n)$.

\section{Triangulated track categories}

\setcounter{equation}{0}

Let $\B$ be an additive track category with a strict translation
functor $t:\B\ra\B$ as in (\ref{2.8}). Let $\A$ be the homotopy
category of $\B$ and $t:\A\ra\A$ the induced translation functor
on $\A$.

A {\em track triangle} is a diagram in $\B$

\begin{equation}\label{6.1}
\xymatrix{
&&&&&&&&\\
A \ar[rr]_f \ar@/^3pc/[rrrr]^0 && B \ar@{}|{\Uparrow H_{0}}[u]
\ar[rr]_{i_f} \ar@/_3pc/[rrrr]_0 && C \ar[rr]_{q_f} \ar@/^3pc/[rrrr]^0
\ar@{}|{\Downarrow H_1 }[d] && tA \ar@{}|{\Uparrow H_{2}}[u] \ar[rr]_{tf} &&
tB\\
&&&&&&&&\\
}
\end{equation}
for which the Toda pairs $(H_0,H_1)$, $(H_1,H_2)$ and
$(H_2,tH_0)$ represent $1_{tA},-1_{tB}$ and $1_{tC}$ respectively,
see (\ref{5.1}).  As in (\ref{1.2})(1) a track triangle can be
equivalently described as follows, compare the notation in section
5.

\begin{Def}\label{6.2}
A {\em track triangle}\index{track triangle} $X$ is a Toda complex
in $\B$ satisfying the following properties:

\begin{equation}\tag{1}
sh^3(X)=tX
\end{equation}
that is, for $X=(X_n,d_n,H_n,n\in\ZZ)$ we have
$X_{n+3}=tX_n,d_{n+3}=td_n,H_{n+3}=tH_n$. Moreover the Toda pairs
($H_{n-1},H_n$) represent for $n\in\ZZ$ the element

\begin{equation*}
\sigma[H_{n-1}|H_n]=(-1)^{n+1}1_{X_n}\in D(X_{n-3},X_n)=\homa
(tX_{n-3},X_n)
\end{equation*}
\begin{equation}\tag{2}
\ \ \ \ \ \ \ \ \ \ \ \ \ \ \ \ \ \ \ \ \ \ \ \ \ \ \ \ \ \ \ \ \ \ \ \ \ \ \ \ \ \ \ \ \ \ \ \ \ \ \ \ \;
=\homa (X_n,X_n)
\end{equation}
given by the identity $1_{X_n}$ of $X_n$.  A {\em map}
between track triangles $X,X'$ is map $f:X\ra X'$ between Toda
complexes such that

\begin{equation}\tag{3}
sh^3(f)=tf,
\end{equation}
that is, for $f=(f_n,F_n,n\in\ZZ)$ we have
$f_{n+3}=tf_n,\;F_{n+3}=tF_n$.
\end{Def}

\begin{Rem}
If $\B$ has a (non strict) translation functor we define a track
triangle $X$ as in (\ref{6.1}) but we replace the track $tH_0$ by the obvious track $(ti_f)(tf)\Rightarrow 0$ given by
$tH_0$ and the pseudofunctor track $t_{i_f,f}$. Similarly we define maps between such track triangles. For details see \cite{ccctc}.
It is easy to generalize the results and proofs in this paper to the non strict case.
\end{Rem}

The definition of a track triangle in $\B$ is motivated by the
representability result for cone functors in (\ref{4.5}) which
implies:

\begin{prop}\label{6.3}
Given a track triangle $X$ one gets for $n\in\ZZ$ representations of
cone functors
\[
\homa (X_n,U)\cong Cone_{d_{n-1}}(U)
\]
which carry $\al$ to $\al_* \{ d_n,H_n\}$. Moreover the exact
sequences (\ref{4.3}) for these representations coincide up to the
sign with the long exact sequence
\[
...\ra \homa (X_{n+1},U)\stackrel{d^*}{\ra} \homa
(X_n,U)\stackrel{d^*}{\ra}\homa (X_{n-1},U)
\]
induced by $X$ with $n\in \ZZ$.
\end{prop}

\begin{proof}
Using (4.5) we only have to show that the sequence is exact. We do this
for $n=1$ so that by the convention (1.2)(1) we have $X_1 =tA$. The
track $H_2$ implies that $d^*d^*=0$. Hence it remains to show
$ker(d^*)\subset image(d^*)$. Let $\al : tA \ra U$ be a map with $G: \al q_f \R
0$. We then construct a map $ \be : tB \ra U$ in $\A$ with $\be (tf)= \pm
t\al$. Let $\be$ be given by the Toda pair $(H_1,G)$, that is, $\sigma
(\be)= [H_1 |G]$ as in (\ref{5.1}). Then we get

\begin{equation*}
\begin{aligned}
\sigma ((tf)^*\be)&=f^* [H_1|G]= \pm \al_* [H_0 |H_1]= \pm
&\al_*\sigma (1_{tA})=\sigma (\pm\al).
\end{aligned}
\end{equation*}

This completes the proof \end{proof}

\begin{Def}\label{6.4}
Let $\B$ be an additive track category with a (strict) translation
functor $t:\B\ra\B$ as in (\ref{2.8}) or (\ref{2.11}). Then $\B$
is a {\em triangulated track category}\index{triangulated track
category} if $\B$ is equipped with a distinguished
family of track triangles which is subject to the following axioms:\\

\begin{itemize}
\item[(TTr1)] For each morphism $f:A\ra B$ in $\B_0$ there exists a
distinguished track triangle $X$ extending $f$, that is,
$d_{-1}=f$.





\item[(TTr2)] Let $X$ and $X'$ be distinguished track triangles
extending $f:A\ra B$ and $f':A'\ra B'$ respectively and let

\begin{equation*}
\xymatrix{
A \ar[rr]^f \ar[dd]_a &\ar@{}|{\Downarrow H }[dd]& B \ar[dd]^{b}\\
&&\\
A' \ar[rr]_{f'} && B'
}
\end{equation*}
be a diagram in $\B$. Then there exists a morphism
$h=(h_n,H_n,n\in\ZZ):X\ra X'$ between track triangles extending
$(a,b,H)$, that is $h_{-2}=a$, $h_{-1}=b$, $H_{-1}=H$.
\end{itemize}
\end{Def}

The category $\Ch (\A)$ of chain complexes in
(\ref{2.9}) and the category {\bf Spec} in (\ref{2.10}) are
examples of triangulated track categories with a strict
translation functor. Moreover a stable cofibration category as in
(\ref{2.12}) is an example of a triangulated track category with a
non strict translation functor, see (\ref{2.11}).

The author does not know any example of a triangulated category
which is not derived from a triangulated track category. We shall
use triangulated track categories in \cite{Ba2}.





We are ready to describe the main result of this paper.

\begin{theo}\label{6.5}
Let $(\Ct{B},t)$ be a triangulated track category and let $\Ct{A}$
be the homotopy category of $\Ct{B}$ and $t\colon\Ct{A}\r\Ct{A}$
be induced by $t\colon\Ct{B}\r\Ct{B}$. Let $\mathcal{E}$ be the
class of triangles in $(\Ct{A},t)$ which are induced by track
triangles in $(\Ct{B},t)$ via the quotient functor
$\Ct{B}_0\r\Ct{A}$. Then $(\Ct{A}, t,\mathcal{E})$ is a
triangulated category.
\end{theo}

More precisely we prove the following addendum.

\begin{adde}
Let $\Ct{B}$ be an additive track category with a strict translation functor $t\colon\Ct{B}\r\Ct{B}$ inducing $t$ on the homotopy category $\Ct{A}$
of $\Ct{B}$. Let $\mathcal{E}$ be the class of triangles in $(\Ct{A},t)$ which are induced by track triangles in $(\Ct{B},t)$. Then
$(\Ct{A},t,\mathcal{E})$ satisfies axioms (Tr0), (Tr2) and (Tr3). Moreover, if $(\Ct{B},t)$ is a triangulated track category, that is, $\Ct{B},t)$
satisfies (TTr1) and (TTr2), then $(\Ct{A},t,\mathcal{E})$ is a triangulated category, that is, (TTr1) obviously implies (Tr1) and (TTr1), (TTr2)
both imply (Tr4).
\end{adde}

The theorem has a generalization for an additive track category with a non-strict translation functor $t\colon\Ct{B}\r\Ct{B}$, see \cite{ccctc}.
In fact in \cite{ccctc} we show that $(\Ct{B},t)$ determines a characteristic cohomology class
$$\nabla=\langle\Ct{B},t\rangle\in H^3(\Ct{A},t)$$
in the translation cohomology of $(\Ct{A},t)$ defined in \cite{BM}. Moreover, the triangulated structure $\mathcal{E}$ of $(\Ct{A}. t)$ depends only on
this class, that is $\mathcal{E}=\mathcal{E}_\nabla$.

We point out that the family of exact triangles in $\A$ is well
defined only by the additive track category $\B$ and the
translation functor $t$ on $\B$, that is, the distinguished track
triangles in (\ref{6.4}) are not used in the definition of exact
triangles in $\A$. The next result shows that each exact triangle
in $\A$ is isomorphic to an exact triangle induced by a
distinguished track triangle. For this we use (TTr1).

\begin{prop}\label{6.6}
Let $X$ and $X'$ be track triangles both extending $f:A\ra B$. Then
there exists a commutative diagram in $\A$
\begin{equation*}
\xymatrix{
A \ar[rr]^f \ar@{=}[dd] && B \ar@{=}[dd] \ar[rr]&& C \ar[dd]^{h}
\ar[rr]&& tA \ar@{=}[dd]\\
&&&&&&\\
A \ar[rr]^f && B \ar[rr] && C' \ar[rr] && tA
}
\end{equation*}
where the rows are the exact triangles induced by $X$
and $X'$ and where $h$ is an isomorphism.
\end{prop}

\begin{proof}
We have a commutative diagram of natural transformations, see
(\ref{4.3}) and (\ref{4.4}),

\begin{equation*}
\xymatrix{
&\homa (C,U)\ar[rd]&\\
\homa (tA,U) \ar[r]^i  \ar[ru] \ar[dr]& Cone_f (U)
\ar@{}|{\wr\Arrowvert}[u] \ar@{}|{ \wr\Arrowvert}[d] \ar[r]^q & \homa (B,U)\\
 & \homa (C',U) \ar[ur]&
}
\end{equation*}

Now the isomorphism in the middle carries the identity of $C'$ to
$h$ and carries the identity of $C$ to the inverse of $h$.

 \end{proof}

\begin{proof}[Proof of (\ref{6.5})]
We show (Tr0), (Tr2) and (Tr4) in the following sections. We here only
consider (Tr1) and (Tr3). It is clear that (TTr1) implies (Tr1).
Moreover given track triangles $X$ and $X'$ extending $f$ and $f'$ we find by
(TTr1) distinguished track triangles $Y$ and $Y'$ extending $f$ and
$f'$ and we can choose $(a,b)$ as in (TTr2) inducing the square in
(Tr3). Then (TTr2) yields a map $Y\ra Y'$ between distinguished track
triangles extending $(a,b)$. Using this map and (\ref{6.6}) we get (Tr3).

 \end{proof}
\begin{Rem}
Grandis \cite{G} describes triangulated categories derived from an
abstract setting of homotopical algebra. This setting axiomatizes
the behaviour of continuous maps and homotopies between them and
does not use tracks as above which are homotopy classes of
homotopies.
\end{Rem}

\section{Induced Toda complexes and proof of (Tr0)}
\setcounter{equation}{0}

\begin{Lem}\label{7.1}
Let $X$ be a Toda complex extending the row ($X_n,d_n$) and let
$f:(X_n,d_n)\ra (X'_n,d'_n)$ be weak equivalence between rows. Then there
exists a unique Toda complex $X'$ extending $(X'_n,d'_n)$ such that $f:X\ra
X'$ is a morphism between Toda complexes.
\end{Lem}

We denote $X'=f_*X$ the {\em Toda complex induced by the
weak equivalence}\index{Toda!complex!induced} $f$.

\begin{proof}
Since $f_n$ is an isomorphism in $\A$ we can choose a track
$G_n:d'_nd'_{n-1}\R 0$. Then the tracks $(H_n,F_{n-1},F_n,G_n)$ represent an
element $\xi_n\in D(X_{n-2},X'_n)$ and hence we can alter $G_n$ by
$(f^*_{n-2})^{-1}(\xi)\in D(X'_{n-2},X'_n)$ so that we get $H'_n$. Then
$(H_n,F_{n-1},F_n,H'_n)$ represents the trivial track and hence $f:X\ra X'$ is
a morphism between Toda complexes with $X'=(X'_n,d'_n,H'_n,n\in \ZZ)$.

 \end{proof}

\begin{Lem}\label{7.2}
Let $X$ be a track triangle in $\B$ with row $(X_n,d_n)$ and let
$f:(X_n,d_n)\ra (X'_n,d'_n)$ be a weak equivalence of rows with $sh^3 f=tf$.
Then the induced Toda complex $f_*X$ is a track triangle.
\end{Lem}

\begin{proof}
We only show that the Toda pair $H'_0,H'_1$ represents $1_{tA'}$. Let
$F=F_{-1},\;G=F_0,\;H=F_1$ and $a=f_{-2},\;b=f_{-1},\;c=f_0,\;e=f_1$ be
given by $f=(f_n,H_n,n\in\ZZ)$. Moreover we write $d=d_i$ and $d=d'_i$.
Assume the Toda pair $(H'_0,H'_1)$ represents $\xi\in D(A',tA')$. Then
we get the equations:\\

\begin{equation*}
\begin{aligned}
\sigma (\al^*\xi)&=(dH'_0\square (H'_1d)^{\boxminus })a\\
&= dH'_0a\square (H'_1)^{\boxminus }da\\
&=d(cH_0\square G^{\boxminus }d\square dF^{\boxminus})\square
((H'_1)^{\boxminus}bd)\\
&=d(cH_0\square G^{\boxminus }d\square dF^{\boxminus})\square
(dG\square Hd\square
eH^{\boxminus}_{1})d\\
&=dcH_0\square dG^{\boxminus }d\square ddF^{\boxminus}\square
dGd\square Hdd\square
eH^{\boxminus}_{1})d\\
&=edH_0\square dG^{\boxminus }d\square dGd\square eH^{\boxminus}_{1} d\\
&= e(dH_0\square H^{\boxminus}_{1} d)\\
&= \sigma ((t\al)_* 1_{tA})\\
&= \sigma (1_{tA'}t\al)\\
&=\sigma (\al^* 1_{tA'})
\end{aligned}
\end{equation*}

Here $\al :A\ra A'$ is the isomorphism induced by $a$.
We know that $e$ induces $t\al$. The equations above show that
$\xi = 1_{tA'}$.

 \end{proof}

\begin{proof}[Proof of (Tr0) in (\ref{6.5})]
The identity $1_A$ yields a track triangle by considering

\begin{equation*}
\xymatrix{
 && && \\
A \ar@/^3pc/[rr] \ar[r]_{1} & A \ar@{}|{\Uparrow H_{0}}[u]
\ar@/_3pc/[rr]_0 \ar[r] & \ast \ar[r]  \ar@/^3pc/[rr] \ar@{}|{\Downarrow H_{1}}[d]
& tA \ar@{}|{\Uparrow H_{2}}[u] \ar[r]_1 & tA\\
&& &&
}
\end{equation*}

Here $H_0 =0,\; H_2 =0$ since $\ast$ is a strict zero object
and $H_1 =\sigma(-1_{tA})$. Therefore

\begin{equation*}
\xymatrix{
A  \ar[r]^{1} & A  \ar[r]& \ast \ar[r] &tA
}
\end{equation*}
is exact. Next let $\al : T\cong T'$ be an isomorphism
of triangles where $T$ is induced by the track triangles $X$. We
obtain for $T'$ a cochain complex $(X'_n,d'_n)$ in $\A$ as in
(\ref{2.1}) and we choose a row $(X'_n,d'_n)$ in $\B$ representing
this cochain complex. Then $\al$ shows that there is a weak
equivalence of rows $f:(X_n,d_n)\ra (X'_n,d'_n)$ with $sh^3f=tf$
inducing $\al$ so that the track triangle $f_* X$ is defined by
(\ref{7.2}). Hence $T'$ is induced by $f_* X$ and therefore $T'$
is exact.

 \end{proof}

\section{Negative Toda complexes and the proof of (Tr2)}
\setcounter{equation}{0}

Let $\B$ be an additive track category with $\B_{\simeq} =\A$. In this
section we use the existence of strong products $X\times Y$ in $\B$. The
case of strong coproducts is treated in a dual way. We can choose for
each object $X$ maps in $\B_0$
\[
\mu : X\times X \lra X
\]

\[
-1_{X} : X\lra X
\]
which induce $(1_X,1_X)$ and $-1_X$ in the additive
category $\A$. Moreover we choose tracks

\begin{equation*}
\begin{aligned}
\tau_X &: \mu_X (1_X,-1_X)\R 0\\
\tau'_X &: \mu_X(0_X,1_X)\R 1_X\\
\tau''_X &: \mu_X(1_X,0_X)\R 1_X
\end{aligned}
\end{equation*}

Here $0_X: X\ra\ast\ra X$ is given by the strict zero object
$\ast$ in $\B$. We now define for each morphism $f:X\ra Y$ in
$\B_0$ tracks $\Gamma_f$ and $\Gamma^{+}_{f}$ as in the diagrams
(compare \cite{Ba1}):

\begin{equation}\label{8.1}
\xymatrix{
X\times X \ar[rr]^{f\times f} \ar[dd]_{\mu_X} &
\ar@{}|{\stackrel{\Gamma^{+}_{f}}{\R}}[dd] & Y\times Y  \ar[dd]^{\mu_Y} && X \ar[rr]^{f}
\ar[dd]_{-1_X} & \ar@{}|{\stackrel{\Gamma_{f}}{\R}}[dd] & Y \ar[dd]^{-1_Y}
\\
&& && && \\
X \ar[rr]_{f} && Y && X \ar[rr]_{f} && Y
}
\end{equation}

Here $\Gamma^{+}_{f}$ is the unique track for which the
diagram

\begin{equation*}
\xymatrix{
& X \ar[rr]^f \ar[dd]_{i_{\veps}} \ar@/_3pc/[dddd]_{1}&& Y
\ar@/^3pc/[dddd]^{1} \ar[dd]^{i_{\veps}} &\\
&& && \\
& X\times X \ar[rr] \ar[dd] \ar@{}|{\Leftarrow}[l] &
\ar@{}|{\stackrel{\Gamma^{+}_{f}}{\R}}[dd]& Y\times Y \ar@{}|{\Rightarrow}[r] \ar[dd]& \\
&& && \\
& X \ar[rr]&& Y & \\
}
\end{equation*}
represents the trivial track $0^{\square}:f\R f$ for
$\veps =1$ and $\veps =2$ with $i_1 = (0,1)$ and $i_2 =(1,0)$.
Moreover $\Gamma_f$ is the unique track for which the following
diagram represents the trivial track $0^{\square} : 0\R 0$.

\begin{equation*}
\xymatrix{ & X \ar[rr]^f \ar[dd]_{(1,-1)} \ar@/_3.5pc/[dddd]&
\ar@{}|{\stackrel{(0^{\square}, \Gamma_f)}{\R}}[dd] & Y
\ar@/^3.5pc/[dddd]^{0}
\ar[dd]^{{(1,-1)}} & \\
&& && \\
& X\times X \ar[rr] \ar[dd] \ar@{}|{\Leftarrow}[l] &
\ar@{}|{\stackrel{\Gamma^{+}_{f}}{\R}}[dd]& Y\times Y \ar@{}|{\Rightarrow}[r] \ar[dd]& \\
&& && \\
& X \ar[rr]&& Y &
}
\end{equation*}

The uniqueness of $\Gamma^{+}_{f}$ and $\Gamma_{f}$
follows from the biadditivity of the coefficient functor $D$ in
the additive track category $\B$.

\begin{Lem}\label{8.2}
Pasting in the diagram

\begin{equation}\tag{1}
\xymatrix{
X \ar[rr]^f \ar[dd]_{-1} & \ar@{}|{\stackrel{\Gamma_{f}}{\R}}[dd]& Y
\ar[rr]^g \ar[dd] &\ar@{}|{\stackrel{\Gamma_{g}}{\R}}[dd] & Z
\ar[dd]^{-1}\\
&& && \\
X \ar[rr] && Y \ar[rr] && Z
}
\end{equation}
yields the track $\Gamma_{gf}$. Moreover for each track
$H:f\R g$ pasting in the diagram

\begin{equation}\tag{2}
\xymatrix{
&&\\
X \ar[rr]^g \ar[dd]_{-1} \ar@/^3pc/[rr]^f & \ar@{}|{\Uparrow
H^{\boxminus}}[u]& Y \ar[dd]^{-1} \\
\ar@{}|{\stackrel{\Gamma_{g}}{\R}}[rr]&&\\
X \ar[rr]^g \ar@/_3pc/[rr]_f & \ar@{}|{\Uparrow H}[d] & Y\\
&&
}
\end{equation}
yields $\Gamma_f$. 


\end{Lem}

Using $\Gamma$-tracks (\ref{8.1}) we can define the {\em
negative of a Toda complex}\index{Toda!complex!negative}

\begin{equation}\label{8.3}
neg(X)=X^- = (X^-_n, d_{n}^{-},H^-_n,n\in\ZZ)
\end{equation}
with $X^-_n =X_n$ and  $d^-_n=d_n$ for $n\equiv
0,1\;mod\;3$ and $d^-_n=(-1_{X_n})d_n$ for $n\equiv 2\;mod\;3$.
Moreover $H^-_n$  is obtained by $H^-_n=H_n$ for $n\equiv
1\;mod\;3$ and $H^-_n = (-1_{X_n})H_n$ for $n\equiv 2\;mod\;3$ and
for $n\equiv 0\;mod\;3$ the track $H^-_n$ is given by pasting the
diagram

\begin{equation}\label{8.4e}
\xymatrix{
&&&&\\
 && X_{n-1} \ar[rr]^d \ar[dd]^{-1} \ar@{}|{\Uparrow H_{n}}[u] &
\ar@{}|{\stackrel{\Gamma_{d}}{\R}}[dd]& X_n  \ar[dd]^{-1}\\
&& && \\
X_{n-2} \ar[rr]_{d^-_{n-1}} \ar[rruu]^{d} \ar@/^6pc/[rrrruu]^0 &&
X_{n-1} \ar[rr]_{d^-_{n}} && X_n
}
\end{equation}
Here we have $d^-_{n-1}=(-1)d_{n-1}$ and $d^-_n=d_n$.

\begin{prop}\label{8.4p}
Toda pairs $(H^-_{n-1}, H_n^- )$ in the negative Toda complex
$neg(X)=X^-$ satisfy the formula

\[
\sigma^{-1}[H^-_{n-1}|H^-_{n}]=-\sigma^{-1}[H_{n-1}|H_{n}]\in
D(X_{n-3},X_n)
\]
\end{prop}

\begin{proof}[Proof of (\ref{8.4p})]
We check three cases

\begin{equation*}
\begin{aligned}
& [ H^{-}_{1}|H^{-}_{2}] \\
                  & = [H_1|(-1)H_2]\\
                  &= (-1)d_2H_1\square (-1)H^{\boxminus}_{2}d_0\\
                  &= (-1)(d_2H_1\square H^{\boxminus}_{2}d_0)=
(-1)[H_1,H_2]
\end{aligned}
\end{equation*}

\begin{equation*}
\begin{aligned}
&[H^-_{2}|H^-_{3}]\\
&= d_3(-1)H_2\square\Gamma^{\boxminus}_{d}d_2d_1\square
(-1)H^{\boxminus}_{3}d_1\\
&= (-1)d_3H_2\square (-1)H^{\boxminus}_{3}d_1\\
&= (-1)(d_3H_2\square H^{\boxminus}_{3}d_1)=(-1)[H_2|H_3]
\end{aligned}
\end{equation*}

\begin{equation*}
\begin{aligned}
&[H^-_{3}|H^-_{4}]\\
&= d_4((-1)H_3\square\Gamma_{d_3}d_2)\square H^{\boxminus}_{4}(-1)d_2\\
&= (-1)d_4H_3\square(\Gamma_{d_4}d_3d_2\square
d_4\Gamma_{d_3}d_2)\square
H^{\boxminus}_{4}(-1)d_2 \\
&= (-1)d_4H_3\square\Gamma_{d_4d_3}d_2\square
H^{\boxminus}_{4}(-1)d_2
\;\;\mathrm{,see}\;(\ref{8.2})(1),\\
&= (-1)d_4H_3\square (\Gamma_{d_4d_3}\square H^{\boxminus}_{4}(-1))d_2\\
&= (-1)d_4H_3\square ((-1)H^{\boxminus}_{4} \square \Gamma_0)d_2
\;\;\mathrm{,see}\;(\ref{8.2})(1),\\
&= (-1)(d_4H_3\square H^{\boxminus}_{4}d_2)\square\Gamma_0 d_2\\
&= (-1)[H_3|H_4]
\end{aligned}
\end{equation*}
Here $\Gamma_0$ for $0:X\ra\ast\ra Z$ is trivial by
(\ref{8.2})(1).

 \end{proof}

The shift functor $sh$ in (\ref{5.4}) does not carry track
triangles $X$ to track triangles since $sh(X)$ does not satisfy
condition (\ref{6.2})(2) on Toda pairs. We get, however, by
(\ref{8.4p}) and (\ref{8.2})(3) the next result.



\begin{prop}\label{8.6}
If $X=(X_n,d_n,H_n,i\in\ZZ)$ is a track triangle then $(H_1^-,H_2^-,H_3^-)$ and $(H_{-1}^-,H_0^-,H_2^-)$
are track triangles in the sense of (6.1).
\end{prop}

\begin{proof}[Proof of (Tr2)]
Let $X$ be a track triangle in $\Ct{B}$ inducing $(u,v,w)$ in
(Tr2). Then by (8.6) we get track triangles $(H_1^-,H_2^-,H_3^-)$
and $(H_{-1}^-,H_0^-,H_2^-)$ inducing $(v,w,-tu)$ and
$(-t^{-1}w,u,v)$ respectively. Here we assume that $\Ct{B}$ has
strong products. If this is not satisfied we can choose by
\cite{BJP} a track category $\Ct{B}'$ with strong products and a
pseudofunctor $\lambda\colon\Ct{B}\r\Ct{B}'$ with homotopy inverse
$\rho\colon\Ct{B}'\r\Ct{B}$ inducing the identity on $\Ct{A}$ and
coefficients $D$. For $X$ in $\Ct{B}$ we get the Toda complex
$\lambda X$ in $\Ct{B}'$ and $neg(\lambda
X)=(X^-_n,d_n^-,H^-_n,n\in\ZZ)$ is defined as above. Then $(\rho
H_1^-,\rho H_2^-,\rho H_3^-)$ and $(\rho H_{-1}^-,\rho H_0^-,\rho
H_2^-)$ are track triangles in $\Ct{B}$ inducing $(v,w,-tu)$ and
$(-t^{-1}w,u,v)$ respectively.  \end{proof}

\begin{proof}[Proof of (6.6)]
Let $\Ct{B}$ be an additive track category with (strict)
translation functor $t$. Then track triangles in $\Ct{B}$ yield
exact triangles in $\Ct{A}$ which satisfy (Tr0) and (Tr2) by the
proof above. Moreover, by an argument as in the proof of (6.7) we
see that (Tr3) holds, compare \cite{ccctc} 7.5. If $\Ct{B}$ is a
triangulated track category then (TTr1) shows that
$(\Ct{A},t,\mathcal{E})$ satisfies (Tr1). Moreover
$(\Ct{A},t,\mathcal{E})$ satisfies (Tr4) by the proof in the next
section. Hence $(\Ct{A},t,\mathcal{E})$ is a triangulated
category.  \end{proof}

\section{Proof of (Tr4)}
\setcounter{equation}{0}

Using (\ref{6.6}) we see that (Tr4)($\#$) is isomorphic to a
sequence of the form (Tr4)($\#$) where, however, $C_u,C_v$ and
$C_{uv}$ are chosen by distinguished track triangles. Then (Tr0)
already proved above shows that it suffices to prove (Tr4) for
such distinguished choices. We consider the following diagram in
$\B$.

\begin{equation*}
\xymatrix{
& A \ar@{=}[rr] \ar[dd]^u && A \ar[rr]^u \ar[dd]^{vu} && B \ar[dd]^v & \\
&&&&&& \\
& B \ar[rr]^v \ar[dd]^d & \ar@{}|{\stackrel{F}{\R}}[dd] & C \ar@{=}[rr] \ar[dd]^{d'}& \ar@{}|{\stackrel{F'}{\R}}[dd] & C \ar[dd]^{d''} & \\
&&&&&& \\
 & C_u \ar[dd]^d \ar[rr]^{\bar{w}} & \ar@{}|{\stackrel{G}{\R}}[dd] & C_{vu}  \ar[rr]^{\bar{u}} \ar[dd]^{d'} & \ar@{}|{\stackrel{G'}{\R}}[dd]& C_v \ar[dd]^{d''} \ar@/^3pc/[dddd]^0 & \\
&&&&&& \\
 & tA \ar@{=}[rr] \ar[dd]^{tu} \ar@/_3pc/[dddd]_{0} && tA \ar[dd]^{tvu} \ar[rr]^{tu}&& tB \ar[dd]^{tv} \ar@{}|{\stackrel{H''_2}{\R}}[r] &\\
&&&&&&\\
& tB \ar@{}|{\stackrel{tH_0}{\Leftarrow}}[l] \ar[dd]^d \ar[rr]^{tv} & \ar@{}|{\stackrel{tF}{\R}}[dd]& tC \ar@{=}[rr] \ar[dd]^{d'} && tC \\
&&&&&&\\
& tC_u \ar[rr]^{t\bar{w}}&& tC_{vu} &&  & \\
}
\end{equation*}

The subdiagrams without tracks are commutative. Let
$\bar{C}_u = (H_0,H_1,H_2)$
be a distinguished track triangle extending the left column, let
$\bar{C}_{vu}=(H'_0,H'_1,H'_2)$ be a distinguished track triangle extending the column in the middle and let
$\bar{C}_v = (H''_0,H''_1,H''_2)$
be a distinguished track triangle extending the right hand column. Then (TTr2) shows that we can find tracks $F,G$ and $F',G'$ which determine maps between track triangles $\bar{C}_u \ra \bar{C}_{vu}$ and $\bar{C}_{vu} \ra \bar{C}_v$ respectively. This implies the equations

\begin{equation}\tag{1}
(t\bar{w})(tH_0) = (tH'_0)\square (tF)(tu)\\
\end{equation}

\begin{equation}\tag{2}  H'_2 = (H''_2\bar{u})\square (tv)G'
\end{equation}

Moreover we know $(H'_2,tH'_0)$ represents

\begin{equation}\tag{3}
\sigma (1_{tC_{vu}})=d'H'_2\square  (tH'_0)^{\boxminus}d'
\end{equation}
Therefore we get the equation

\begin{equation}\tag{4}
\sigma (1_{tC_{vu}})=d'H''_2\bar{u} \square d'(tv)G'\square
(tF)(tu)d'\square (t\bar{w})(tH_0)^{\boxminus}d'
\end{equation}

We now choose a distinguished track triangle for the map $\bar{u}:C_{vu}\ra C_v$ which exists by (TTr1). This yields $C_{\bar{u}}$ and the following diagram in the homotopy category $\A$ where the bottom row corresponds to (\ref{2.3})(\#).

\begin{equation}\tag{5}
\xymatrix{
 && && C_{\bar{u}} \ar[ddrr]^{-q_{\bar{u}}} \ar[dd]^{\varphi}&& \\
&& && && \\
C_{vu} \ar[rr]_{\bar{u}} && C_v \ar[rr]_{dd''} \ar[uurr]^{i_{\bar{u}}} && tC_u \ar[rr]_{t\bar{w}}&& tC_{vu}
}
\end{equation}

We shall construct a map $\varphi$ for which this diagram commutes. For this we use the representation $\chi: Cone_{\bar{u}}(U)\ra\homa (C_{\bar{u}},U)$ defining $i_{\bar{u}}$, $q_{\bar{u}}$ as in (\ref{8.4p}) with $U=tC_u$. Moreover we use the following diagram of tracks:\\

\begin{equation}\tag{6}
\xymatrix{
C_{vu}  \ar[dd]_{d'} \ar[rr]^{\bar{u}}&& C_{v} \ar[dd]_{d''} \ar[ddrr]^{0}_{\stackrel{H''_2}{\R}}&& \\
\ar@{}|{\stackrel{G'}{\R}}[rr] && && \\
tA \ar[rr]^{tu} \ar[ddrr]_{0}^{\stackrel{(tH_0)^{\boxminus}}{\R}}&& tB \ar[rr]^{tv} \ar[dd]^d&& tC \ar[dd]^{d'}\\
 &&\ar@{}|{\stackrel{tF}{\R}}[rr]  && \\
 && tC_u \ar[rr]_{t\bar{w}}&& tC_{vu}
}
\end{equation}

According to (4) this diagram represents $\sigma(1_{tC_{vu}})$. Now let
\begin{equation}\tag{7}
 \varphi = \chi \{ \al,\hal\} \;\;\textrm{with}\;\;\al =dd'' \;\;and \;\;\hat{\al}=(tH_0)d'\square d(G')^{\boxminus}
\end{equation}
Then we have $\varphi i_{\bar{u}} =dd'$ and we get
$(t\bar{w})\varphi = -\bar{q}_v$ since using $R=(d'H''_2\square
(tF)d'')^{\boxminus}$ yields the equation

 \begin{equation}\tag{8}
 \chi ((t\bar{w})\varphi)= (t\bar{w})_* \{\al ,\hal\}=\{(t\bar{w})\al , (t\bar{w})\hal\}=\{0,\sigma (-1_{tC_{vu}})\}=\chi(-q_{\bar{u}})
\end{equation}
Hence the proof of the commutativity of (5) is complete.

Finally we observe that there is a commutative diagram

\begin{equation}\tag{9}
\xymatrix{
\homa (tC_{vu},X) \ar[r]^{(t\bar{w})^*}& \homa (tC_{u},X) \ar[r]^{(dd'')^*}& \homa (C_v,X) \ar[r]^{\bar{u}^*}& \homa (C_{vu},X)\\
&&&\\
Cone_{vu}(t^{-1}X) \ar[r]^{v^*}\ar[uu]_{t\chi}& Cone_{u}(t^{-1}X) \ar[r]^{\delta} \ar[uu]_{t\chi}& Cone_{v}(X) \ar[uu]_{\chi} \ar[r]^{u^*}& Cone_{vu}(X)\ar[uu]_{\chi}
}
\end{equation}
Here the bottom row is the exact sequence in (\ref{3.6}). The
vertical arrows are isomorphisms as in (\ref{4.4}) by the track
triangles $\bar{C}_{v}$, $\bar{C}_{u}$ and $\bar{C}_{vu}$
respectively. One readily checks by the properties of maps between
track triangles that the diagram commutes. Hence the top row is
exact by (\ref{3.6}) and therefore $\varphi$ in (5) is an
isomorphism. Now (Tr0) shows that the bottom row of (5) is an
exact triengle. This completes the proof of (Tr4).

\section{The dual of a triangulated track category}
\setcounter{equation}{0}

We have seen in (1.5) that an exact triangle
$A\stackrel{f}{\ra}B\ra C\ra tA$ in a triangulated category $\A$ induces for each
object $U$ of $\A$ the long exact sequence

\begin{equation}\label{10.1}
\homa (U,t^{-1}A)\ra \homa (U,t^{-1}B)\ra\homa (U,t^{-1}C)\ra\homa
(U,A)\stackrel{f_*}{\ra}\homa (U,B)
\end{equation}
which resembles the exact sequence for the dual cone
functor $Cone^f$ in (\ref{3.4}). In fact we get:

\begin{prop}\label{10.2}
Let $\B$ be a triangulated track category. Then for each $f:A\ra B$ in
$\B_0$ the dual cone functor $Cone^f$ has a canonical representation

\[
\chi : Cone^f (U)\cong\homa (U,t^{-1}C_f)
\]
where $A\stackrel{f}{\ra}B\ra C_f\ra tA$ is induced by a
track triangle extending $f$. Moreover $\chi$ yields an
isomorphism of exact sequences (\ref{3.4}) and (\ref{10.1}).
\end{prop}

\begin{proof}
We have the tracks
\begin{equation*}
\xymatrix{
&&& &\\
t^{-1}A \ar[r]^{t^{-1}f} \ar@/_3pc/[rr]& t^{-1}B \ar[r]^u
\ar@/^3pc/[rr] \ar@{}|{\Downarrow H_{-3}}[d]&  t^{-1}C_f \ar[r]^w \ar@{}|{\Uparrow
H_{-2}}[u]  \ar@/_3pc/[rr]& A \ar[r]^f \ar@{}|{\Downarrow H_{-1}}[d]& B
\\
&&&&\\
}
\end{equation*}
which define the element

\[
\{ w,H_{-1}\} \in Cone^f (t^{-1}C_f).
\]

Now we define the inverse of $\chi$ in (\ref{10.1}) by

\[
\chi^{-1}(\al)=\al^*\{ w,H_{-1}\}.
\]

We claim that the following diagram commutes up to sign
\begin{equation*}
\xymatrix{
&& Cone^f (U) \ar[drr]^{\bar{q}}&& \\
\homa (U,t^{-1}B) \ar[urr]^{\bar{i}} \ar[drr]_{u_*}&&&& \homa (U,A)\\
 && \homa (U,t^{-1}C_f)\ar[uu]_{\chi^{-1}} \ar[urr]_{w_*}&&
}
\end{equation*}
where $\bar{i}$ and $\bar{q}$ are defined as in
(\ref{3.4}). In fact $\bar{q}\chi^{-1}=w_*$ since

\[
\bar{q}\chi^{-1} (\al )= \bar{q} \{ w\al ,H_{n-1}\al\} = w\al = w_*
(\al ).
\]

Moreover, for $\delta\in\homa (U,t^{-1}B)$ we have to show
that $\bar{i}\delta =\chi^{-1}u_*\delta = \chi^{-1} (u\delta )$ where

\begin{equation*}
\begin{aligned}
\bar{i}\delta &= \{ 0, \sigma\delta\} = \delta^*\{ 0,\sigma
1_{t^{-1}B}\}\\
\chi^{-1} (u\delta) &=\delta^* \chi^{-1} (u)= \delta^*\{ wu, H_{-1}u\}
\end{aligned}
\end{equation*}
Hence it suffices to prove that in $Cone^f (t^{-1}B)$ we
have $\{ 0,\sigma 1_{t^{-1}B}\}= -\{ wu,H_{-1}u\}$. But  this
follows since the pair $(H_{-2},H_{-1})$ represents $\sigma 1_B$.
Hence we proved $\chi^{-1}u_x =-\bar{i}$. Using the exact
sequences (\ref{3.4}) and (\ref{10.1}) we see that $\chi^{-1}$ is
an isomorphism and the proposition (\ref{10.2}) is proved.

 \end{proof}

Proposition (\ref{10.1}) corresponds to the following result:

\begin{theo}\label{10.3}
The categorical dual $\B^{op}$ of a triangulated track category $\B$ is
again a triangulated track category with translation functor given by
the inverse $t^{-1}$ of the strict translation functor $t$ for $\B$. The
distinguished track triangles in $\B^{op}$ are the duals $X^{op}$ of
distinguished track triangles $X$ in $\B$ with $X^{op}_n =X_{-n}$.
\end{theo}

\section{The natural notion of good morphism}
\setcounter{equation}{0}

This section should be compared with section 3 in \cite{Ne} which has
the same title. Let $\B$ be a triangulated track category with homotopy
category $\A$ and quotient functor $p:\B_0 \ra\A$. Then $p$ induces the
functor

\begin{equation}\label{11.1}
p: {\bf TTr}_{\B} \ra {\bf Tri}_{\A}
\end{equation}
Here ${\bf Tri}_{\A}$ is the category of triangles in $\A$ and
morphisms as in (\ref{2.1}). Moreover ${\bf TTr}_{\B}$ is the
category of distinguished track triangles in $\B$ and maps as in
(\ref{6.2}).

A {\em special} isomorphism in ${\bf Tri}_{\A}$ is an
isomorphism of the form

\begin{equation}\label{11.2}
\xymatrix{
A \ar[rr] \ar@{=}[dd]&& B \ar[rr] \ar@{=}[dd]&& C \ar[rr] \ar[dd]&& tA
\ar@{=}[dd]\\
&&&&&&\\
A \ar[rr]&& B \ar[rr]&& C'\ar[rr] && tA
}
\end{equation}

For example (\ref{6.6}) yields such special isomorphisms. We use
the functor $p$ in (\ref{11.1}) and special isomorphisms for the
following definition:

\begin{Def}\label{11.3}
A {\em good triangle}\index{triangle!good} is a triangle $T$ in
${\bf Tri}_{\A}$ together with an equivalence class of special
isomorphisms

\[
T\cong p(X)
\]
where $X$ is an object in  ${\bf TTr}_{\B}$. Two such isomorphisms
$T\cong p(X)$, $T\cong p(Y)$ are equivalent if there exists a map
$X\ra Y$ in  ${\bf TTr}_{\B}$ inducing the composite $p(X)\cong
T\cong p(Y)$. A {\em good morphism}\index{good morphism} between
good triangles $T,T'$ is a morphism $f:T\ra T'$ in ${\bf
Tri}_{\A}$ for which there exists a commutative diagram in ${\bf
Tri}_{\A}$

\begin{equation*}
\xymatrix{
p(X) \ar[rr]^{p(\varphi)} && p(X') \\
&&\\
T \ar[uu]^{\al } \ar[rr]_{f} && T'\ar[uu]_{\al '}
}
\end{equation*}
where $\al$ and $\al '$ are representing the equivalence
classes of isomorphisms for $T$ and $T'$ respectively and where
$p(\varphi)$ is induced by a map $\varphi : X\ra X'$ in ${\bf
TTr}_{\B}$. Let ${\bf GTr}$ be the category of good triangles and
good morphisms. We have the faithful forgetful functor

\[
F: {\bf GTr} \longrightarrow {\bf Tri}_{\A}.
\]
\end{Def}

\begin{theo}\label{11.4}
For a triangulated track category $\B$ the homotopy category $\A$
of $\B$ together with the functor $F: {\bf GTr}\longrightarrow
{\bf Tri}_{\A}$ satisfy the axioms (GTR1)...(GTR10) of a good
triangulated category defined in \cite{Ne}(3.4).
\end{theo}

We do not want to recall the 10 axioms of Neeman and leave the
proof of the theorem to the reader. The author does not know any
examples of triangulated categories or good triangulated
categories which are not given by triangulated track categories.
Therefore it seems reasonable to replace the complicated set of
axioms (Tr$n$) for $n=0,...,4$ or (GTR$n$) for $n=0,...,10$ by two
very natural axioms (TTr1) and (TTr2) of a triangulated track
category. Many of the numerous results on triangulated categories
in the literature have an  analogue in a triangulated track
category  yielding a natural and interesting extension of such
results.


\title{{[II]}\\\hspace{5pt}\\Principal maps in triangulated categories
and maps between $2$-stage spectra}
\author{H.-J. Baues}
\date{}
\blfootnote{\textit{$2000$ Mathematics Subject Classification:} 18D05, 18E30}%
\blfootnote{\textit{Key words and phrases:} triangulated
categories, groupoid-enriched categories, homotopy pairs,
principal maps}

\maketitle

\begin{abstract}
Using the two natural axioms on track triangles in a triangulated
track category we show that the mapping cone construction is
functorial, that is, principal maps between mapping cones in a
triangulated category can be identified with the homotopy classes
of homotopy pairs. This is applied to obtain an algebraic model of
the category of principal maps between $2$-stage spectra.
\end{abstract}

\section*{Introduction}

Consider the following commutative diagram in a triangulated category with
translation functor \emph{t}

\begin{equation*}
\xymatrix{
A \ar[d]_f \ar[r]^d   & B \ar[d]^g \ar[r] & C \ar[d]^h \ar[r] & t(A) \ar[d]^{t(f)} \\
A^{'}\ar[r]^{d'}         & B^{'}\ar[r]       & C^{'}\ar[r]       & t(A^{'})          }
\end{equation*}
Here the rows are exact (i.e. distinguished) triangles and $C$ is
called a \emph{cone} of $d$. Then one of the axioms of a
triangulated category says: \emph{Given $f$ and $g$ with $gd=d'f$,
then there exists $h$ such that the diagram commutes.} Conversely
we define that $a$ map $h:C\longrightarrow C'$ is \emph{principal}
with respect to the cone structure of $C$ and $C'$ if there exists
$(f,g)$ such that the diagram commutes.

It is an old mystery in triangulated categories how pairs $(f,g)$
and principal maps $h$ are related. For this compare Gelfand-Manin
\cite{GM} IV.7 who remark: \emph{"This nonfunctoriality of a cone
is the first symptom that something is going wrong in the axioms
of a triangulated category. Unfortunately, at the moment, we don't
have a more satisfactory version."} We think that triangulated
track categories as introduced in \cite{Ba1} are in this respect
more satisfactory since we show  that in a triangulated track
category the category of principal maps is equivalent to the
homotopy category of the category of homotopy pairs. This leads
(in the case of the stable homotopy category of spectra) to the
algebraic computation of all principal maps $C\longrightarrow C'$
where $A$, $B$, $A'$ and $B'$ are finite direct sums of
Eilenberg-Mac Lane spectra over a prime field $\mathbb{F}$. In
this case we call $C$ and $C'$ $2$-\emph{stage spectra}, see (1.5)
and section 6.

\section{{Principal maps}}
\setcounter{equation}{0}

Let {\bf A} be a triangulated category with translation functor $t: {\bf A}\longrightarrow {\bf A}$ which
is an additive equivalence. Then a distinguished family of triangles

\begin{equation}\label{1.1}
\xymatrix{
X \ar @{} [r] |{=}  & ( A \ar[r]^{d} & B \ar[r]^{d_{0}} & C \ar[r]^{d_{1}} & tA )} ,
\end{equation}
termed \emph{exact triangles}\index{triangle!exact}, is
given satisfying axioms (Tr0)...(Tr4), see \cite{Ne1}, \cite{We}.
We also write $C=C_{d}$. A morphism between exact triangles
$X\longrightarrow X'$ is a commutative diagram in ${\bf A}$ of the
form

\begin{equation}\label{1.2}
\xymatrix{
X \ar @{} [r] |{=}  & ( A \ar[d]^f \ar[r]   & B \ar[d]^g \ar[r] & C \ar[d]^h \ar[r] & tA ) \ar[d]^{t(f)} \\
X' \ar @{} [r] |{=} & ( A' \ar[r]            & B' \ar[r]        & C' \ar[r]        & tA' )          }
\end{equation}
We say that the morphism $(f,g,h):X\longrightarrow X'$
extends the map $h:C\longrightarrow C'$ and we say that a map
$h:C\longrightarrow C'$ in $\A$ is
\emph{principal}\index{principal map} with respect to $X$ and $X'$
if there exists a morphism $X\longrightarrow X'$ extending $h$.
Let {\bf ET} be the category of exact triangles (1.1) and
morphisms (1.2). Then we have the forgetful functor

\begin{equation}\label{1.3}
{\bf ET}\longrightarrow\A
\end{equation}
which carries $X$ to $C$ and $(f,g,h)$ to $h$. The image
category of this functor is the category {\bf Prin} of principal
maps. Objects in {\bf Prin} are exact triangles $X$, $X'$ and
morphisms $X\longrightarrow X'$ in {\bf Prin} are the maps
$C\longrightarrow C'$ in $\A$ which are principal, i.e. in the
image of the functor (1.3). We have the faithful functor ${\bf
Prin}\longrightarrow \A$.

\begin{Lem}\label{eq:1.4}
A map $h:C\longrightarrow C'$ is principal with respect to $X$ and
$X'$ if there exists a commutative diagram in $\A$  as in $(1)$ or
$(2)$ respectively.
\begin{equation}
\xymatrix{
B \ar[r] \ar[d] \ar @{}[rd]|{(1)} & C \ar @{}[r] \ar[d]^{h} & C \ar @{}[rd]|{(2)} \ar[d]_{h} \ar[r] & tA
\ar[d]\\
B'\ar[r]                         & C' \ar @{}[r]           & C' \ar[r]                             & tA'}
\end{equation}
\end{Lem}

Axiom $(Tr3)$ says that given $(f,g)$ as above with $gd=df$ there
exists $(f,g,h):X\longrightarrow X'$ extending $(f,g)$. Here,
however, the map $h$ is not well defined by the pairs $(f,g)$.
Therefore the axiom yields a misterious relationship between the
pairs $(f,g)$ and the principal maps $h:C\longrightarrow C'$ in a
triangulated category. In this paper we clarify this relationship
in the context of a triangulated track category \cite{Ba1}. This
leads in section 6 below to a computation as in the following
example.

\begin{exa}\label{1.5}
Let $\A={\bf Spec}$ be the stable homotopy category of spectra. Then  $\A$  is a triangulated category. For a prime
field
$\mathbb{F}$ we have the full subcategory

\begin{equation*}
2\textrm{- }{\bf stage} \subset  {\bf Prin}.
\end{equation*}

The objects of $2$-{\bf stage} are triangles $X$ for
which $A$ and $B$ are finite direct sums of Eilenberg-Mac Lane
spectra over $\mathbb{F}$. Hence for an object $X=C_{d}$ in
$2$-{\bf stage} the map $d:A\longrightarrow B$ in $\A$ is
algebraicly determined by a map between free modules over the
Steenrod algebra. Also for objects $X\longrightarrow X'$ in
$2$-{\bf stage} the map $(f,g):d\longrightarrow d'$ is easily
expressed algebraicly. It is, however, an old problem to compute
the principal maps $C_{d}\longrightarrow C_{d'}$ in the category
$2$-{\bf stage}, see \cite{KM}. As an application of the main
result in this paper we obtain an algebraic description of the
category $2$-{\bf stage} in terms of the pair algebra
$\mathcal{B}$ of secondary cohomology operations computed in
\cite{Ba3}, see section 6.
\end{exa}

\section{{The category of homotopy pairs}}

\setcounter{equation}{0}

Let $\B$  be a triangulated track category (see \cite{Ba1}) with
strict translation functor $t:\B\longrightarrow\B$ and homotopy
category $\A$. Then we have a linear extension of categories

\begin{equation}\label{2.1}
D \longrightarrow \B_{{\bf 1}} \rightrightarrows \B_{{\bf 0}} \longrightarrow \A
\end{equation}
with coefficient $D(A,B)=Hom_{\A}(tA,B)$. Here $\B =(\B_{{\bf
1}}\rightrightarrows \B_{{\bf 0}})$ is a track category (i.e. a
groupoid enriched category) with hom-groupoids $\ll A,B \rr$ for
all objects $A$, $B$ in $\B$. The objects of the groupoid $\ll A,B
\rr$ are maps $f:A\longrightarrow B$ (or 1-cells) and morphisms,
denoted $F:f\Longrightarrow f'$, are called tracks (or 2-cells).
The 1-cells are the morphisms of the category $\B_{{\bf 0}}$. Two
1-cells $f$, $f'$ are \emph{homotopic} if there exists a track
$f\Longrightarrow f'$. The homotopy classes of 1-cells in
$\B_{{\bf 0}}$ are the morphisms in $\A$ and $\B_{{\bf
0}}\longrightarrow \A$ is the quotient functor. The track category
$\B$  has a strong zero object * with $\ll A,* \rr=0=\ll *,A \rr$
which yields the zero maps $0:A\longrightarrow *\longrightarrow B$
in $\B_{{\bf 0}}$. Moreover, we assume that $\B$ has strong
coproducts $A\vee B$ (which yield the direct sum $A\oplus B$ in
\A) with $\ll A\vee B,C \rr=\ll A,C \rr\times \ll B,C \rr$.

A \emph{track-triangle}\index{track triangle}
$X=(H_{0},H_{1},H_{2})$ in $\B$ is a diagram of the form
\begin{equation} \label{2.2}
\xymatrix{
&  &&  & \\
A \ar[r]_{d}\ar@/^2.5pc/[rr]^{0} & B \ar@{}[u]|{\ \ \Uparrow H_{0}} \ar[r]\ar@/_2.5pc/[rr]_{0} & C
 \ar@{}[d]|{\ \ \Downarrow H_{1}} \ar[r]\ar@/^2.5pc/[rr]^{0} & tA \ar@{}[u]|{\ \ \Uparrow H_{2}} \ar[r]_{td} & tB \\
&&  &&
}
\end{equation}
satisfying certain properties \cite{Ba1}. In a
triangulated track category $\B$ one has a distinguished class of
track triangles subject to two axioms (TTr1) and (TTr2)'. We show
in \cite{Ba1} that the two axioms imply that $\A$ is a
triangulated category in which exact triangles are
induced by track triangles via the quotient functor $\B_{{\bf 0}}\longrightarrow \A$. 

\begin{Def}\label{2.3}
Let $\B$ be a track category. Then the \emph{category} {\bf HP}
\emph{of homotopy pairs in}\index{homotopy pair} {\bf B} is
defined as follows. Objects are the morphisms $d:A\longrightarrow
B$ in $\B_{{\bf 0}}$ and morphisms $F:d\longrightarrow d'$ in {\bf
HP} are diagrams

\begin{equation*}
\xymatrix{
A \ar[r]^{f} \ar[d]_{d} \ar@ {}[rd]|{\overset{F}{\Rightarrow}} & A' \ar[d]^{d'} \\
B \ar[r]_{g} & B'
}
\end{equation*}
in $\B$. Composition in {\bf HP} of the morphisms
\begin{equation*}
\xymatrix{
A \ar[r]^{f} \ar[d]_{d} \ar@ {}[rd]|{\overset{F}{\Rightarrow}} & A' \ar[d] \ar[r]^{f'} \ar[d] \ar@ {}
[rd]|{\overset{F'}{\Rightarrow}} & A'' \ar[d]^{d''} \\
B \ar[r]_{g} & B'\ar[r]_{g'} & B''
}
\end{equation*}
is defined by $ff'$, $gg'$ and
\begin{equation*}
F' \ \bar{\square}\ F \ = \ F' f \ \square \ g' F.
\end{equation*}
Here  $F \;\square\;G$ denotes the composition of tracks. Let
$F^{\,\small{\boxminus}}$ be the inverse of the track $F$.
\end{Def}

\begin{prop}
Let $\B$ be a triangulated track category. Then there is a functor
\begin{equation*}
\psi :{\bf HP}\longrightarrow {\bf Prin}
\end{equation*}
from the category of homotopy pairs in $\B$  to the category of
principal maps in $\A$  such that this functor is full and
representative.
\end{prop}

This result shows that there is a natural equivalence relation
$\sim$ on {\bf HP} such that quotient category ${\bf HP}/\sim$ is
equivalent to the category {\bf Prin}. We shall compute the
equivalence relation in section 4. The functor $\psi$ in (2.4)
corresponds to the construction of principal maps between mapping
cones in \cite{Ba2}.

\begin{adden}
If $\B$  has strong coproducts then there is a functor $\psi$ in
such a way that $\psi$ carries coproducts in {\bf HP} to
coproducts in $\A$.
\end{adden}

\begin{proof}[Proof of (2.4)] For each object $d:
A\longrightarrow B$ we choose a distinguished track triangle
$X=X_{d}$ as in (2.2) extending $d$, this corresponds to axiom
(TTr1) in $\B$. We define the functor $\psi$ on the object $d$ by
$\psi (d)= C_{d}$ where $C=C_{d}$ is given by $X$ as in (2.2),
that is $\psi (d)$ is the triangle induced by $X_{d}$ in \A. Let
$(f,g,F):d\longrightarrow d'$ be a morphism in {\bf HP}. Then we
can choose by the axiom (TTr2)' in $\B$  a map
$\varphi:X_{d}\longrightarrow X_{d'}$ between track triangles
extending $(f,g,F)$. The map $\varphi$ defines a map $h:
C_{d}\longrightarrow C_{d'}$ and a track $G$ as in the diagram
\begin{equation}\tag{1}
\xymatrix{
&&0&&\\
A \ar[rr] \ar@/^2.5pc/[rrrr] \ar[d]_{f} \ar @{} \ar @{}[rrd]|{\Downarrow F} && B \ar @{}[u]|{\ \ \Uparrow H_{0}} \ar[rr]
\ar[d]^{g} \ar @{}[rrd]|{\Downarrow G} && C \ar[d]^{h}\\
A' \ar[rr] \ar@/_2.5pc/[rrrr] && B' \ar[rr] \ar @{}[d]|{\ \ \Downarrow H'_{0}} && C'\\
&& 0 && \\
}
\end{equation}
which represents the trivial track $0\Longrightarrow 0$. The
homotopy class \{h\} of $h$ in $\A$ corresponds via the representation

\begin{equation}\tag{2}
\{h\}\in hom_{\A}(C,C')\cong Cone_{d}(C')
\end{equation}
in \cite{Ba1} (4.4) to the element $\{hd,hH_{0}\}$ in
$Cone_{d}(C')$. By the equivalence relation defining the cone
functor in \cite{Ba1} (3.1) (1)  we see that for any other map
$\varphi ':X_{d}\longrightarrow X_{d'}$ extending $(f,g,F)$ we
have

\begin{equation*}
\{hd\,,\,hH_{0}\}\,=\,\{d'g\,,\,hH_{0} \,\square\,
G^{\,\boxminus}d\}\,\overset{(*)}{=}\,\{d'g\,,\,H'_{0}f
\,\square\, d'F\}
\end{equation*}
\begin{equation}\tag{3}
\ \ \ \ \ \ \ \ \,\,\overset{(*)}{=}\,\{d'g\,,\,h'H_{0}\,\square\,
G'^{\,\boxminus}d\}\,=\,\{h'd\,,\,h'H_{0}\}
\end{equation}
Here $(h',G')$ is given by $\varphi '$ and in (*) we use the fact
that diagram (1) represents the trivial track. This shows that the
homotopy class $\{h\}$ of $h$ does not depend on the choice of
$\varphi$. We now define the functor $\psi$ on morphisms by $\psi
(f,g,F)=\{h\}$. It is clear that $\psi$ is a well defined functor.
By \cite{Ba1} (6.6) we see that $\psi$ is representative, that is,
for each exact triangle (1.1) in {\bf Prin} there is  an object
$d$ in {\bf HP} with $\psi (d)\cong C$ in {\bf Prin}. Finally we
check that the functor $\psi$ is full. For this let $h:
\psi(d)\longrightarrow \psi(d')$ be a map in {\bf Prin}. Then we
can choose $(g,G)$ as in (1) such that
$(g,h,G):d_{0}\longrightarrow d'_{0}$ is a map in {\bf HP} with
$d_{0}=d:B\longrightarrow C$. We can find a map $\varphi
':X_{d_{0}}\longrightarrow X_{d'_{0}}$ extending $(g,h,G)$ and we
define $(f,F)$ in such a way that the strict translation functor
$t$ on $\B$ yields $(tf,tg,tF)$ as part of $\varphi_{0}$. Then we
get $\psi (f,g,F)=\{h\}$. This completes the proof of (2.4).
\end{proof}

For the proof of the addendum we observe that the coproduct of track
triangles is a track triangle.

\section{{Cogroups and $\Gamma$-tracks}}
\setcounter{equation}{0}

Let {\bf gr} be the category of finitely generated free groups $<E>$.
\begin{Def}
Let $C$ be a category with coproducts $X\vee Y$. Then a {\em cogroup} $X$ in {\bf C}
is given by a coproduct preserving functor
\begin{equation*}
{\bf gr}\longrightarrow {\bf C}
\end{equation*}
which carries $\mathbb {Z}$ to $X$. If {\bf B} is a track category with strong coproducts
$X\vee Y$ then a {\em pseudo cogroup} $X$ in {\bf B} is a pseudo functor
\begin{equation*}
{\bf gr}\longrightarrow {\bf B}
\end{equation*}
which preserves coproducts. Hence we have for $\alpha : <E>\longrightarrow <E'>$ in {\bf gr} a map
$\alpha _{*}:\vee_{E}X\longrightarrow \vee_{E'}X$ and for $\beta :<E>\longrightarrow <E''>$ we have a track
\begin{equation*}
\xymatrix {
\vee_{E}X \ar[r]^{\alpha_{*}} \ar@/_2.5pc/[rr] & \vee_{E'}X \ar [r]^{\beta_{*}} \ar @{}[d] |{\Downarrow \tau} & \vee_{E''}X \\
& (\,\beta \alpha)_{*}&
}
\end{equation*}
We call $\tau$ a {\em linearity track} for $X$.

For $-1:\mathbb{Z}\rightarrow \mathbb{Z}$ we write $-1_{X}=(-1)_{*}:X\rightarrow X$ and for
$i_{1}+i_{2}:\mathbb{Z}\rightarrow \mathbb{Z}\vee \mathbb{Z}$ in {\bf gr} we write
$(i_{1}+i_{2})_{*}=i_{1}+ i_{2}:X\rightarrow X\vee X$. This is the {\em comultiplication} of $X$.
\end{Def}

Now let {\bf B} be an additive track category with strong
coproducts and a strong zero object $*$. Then we obtain as in
\cite{Ba1} (8.1) by use of the linearity tracks above well defined
$\Gamma$-Tracks for each $f:X\longrightarrow Y$ in ${\bf B}_{0}$.
\begin{equation}
\xymatrix{
X\vee X \ar[r]^{f\vee f}\ar @{}[rd]|{\overset{\Gamma}{\Longrightarrow}} & Y\vee Y\\
X \ar[u]^{i_{1}+i_{2}} \ar[r]_{f}& Y\ar[u]_{i_{1}+i_{2}}}
\end{equation}

\begin{equation}\tag{1}
\xymatrix{
X \ar[r]^{f}\ar @{}[rd]|{\overset{\Gamma}{\Longrightarrow}} & Y\\
X \ar[u]^{-1} \ar[r]_{f}& X\ar[u]_{-1}}
\end{equation}
We use these tracks to define for
$-i_{1}+i_{2}:X\longrightarrow X\vee X$ also the $\Gamma$-track
\begin{equation}\tag{2}
\Gamma : (-i_{1}+i_{2})f\Longrightarrow -i_{1}f+i_{2}f
\end{equation}
given by pasting in the diagram
\begin{equation}\tag{3}
\xymatrix{
X\vee X \ar[r]^{f\vee f}\ar @{}[rdd]|{\Gamma\vee 0^{\square}_{f}} & Y\vee Y\\
&&\\
X\vee X \ar @{}[rdd]|{\overset{\Gamma}{\Longrightarrow}}\ar[uu]^{-1\vee 1} \ar[r]^{f\vee f}& Y\vee Y\ar[uu]_{-1\vee 1}\\
&&\\
X\ar[uu]^{i_{1}+i_{2}}\ar[r]_{f} & Y \ar[uu]_{i_{1}+i_{2}}
}
\end{equation}
Here $0^{\square}_{f}: f\Longrightarrow f$ is the
trivial track. In \cite{BP} we prove the following result.
\begin{Lem}
Let ${\bf B}$ be an additive track category with strong
coproducts. Then each object $X$ in ${\bf B}$ admits the structure
of a pseudo cogroup.
\end{Lem}

In some of the proofs below we assume that objects X in ${\bf B}$
are actually cogroups so that we omit this way the linearity
tracks for X. We leave it to the reader to generalize the argument
in case the linearity tracks are non trivial.
\begin{exa}
Let ${\bf S_{P_{F}}}$ be the category of simplicial group spectra
$FX$ where $X$ is a simplicial spectrum, see \cite{Ka} (5.2). Then
${\bf S_{P_{F}}}$ has a cylinder functor and the homotopy category
${\bf S_{P_{F}}}/\simeq \ \,= {\bf Spec}$ is the stable homotopy
category of the spectra. Moreover the associated track category
for ${\bf S_{P_{F}}}$ is a triangulated track category ${\bf B}$
with a strict translation functor $t$ in which each object is a
cogroup in ${\bf B_{0}}$. Hence the results of this paper hold in
particular in stable homotopy theory.
\end{exa}

\section{{The homotopy relation for homotopy pairs}}
\setcounter{equation}{0}

Let ${\bf B}$ be an additive track category with strong coproducts and a strong zero object. For each object $X$ in $\mathbf{B}$ we choose the structure of a
pseudo cogroup as in Lemma 3.3.
\begin{Def}
For a morphism $d:A\longrightarrow B$ in ${\bf B_{0}}$ we define the {\em cylinder object} $I(d)$ by the diagram in ${\bf B}$
\begin{equation*}
\xymatrix {
A \ar[r]^{\partial_{2}\ \ \ \ \ \,} \ar@/_2.5pc/[rr]_{0} & A\vee B\vee A \ar[r]^{\ \ \ \partial_{1}}\ar @{}[d]|{\ \ \ \Downarrow \partial_{12}}
& B\vee B \\
&&
}
\end{equation*}
Here we set $\partial_{1}=(i_{1}d, -i_{1}+i_{2},i_{2}d)$ and $\partial_{2}= i_{1}+i_{2}d-i_{3}$. Moreover the
track $\partial_{12}$ is given by $\Gamma$ in(3.2)(2) and by linearity tracks, that is, for $\partial_{1}\partial_{2}=
i_{1}d+(-i_{1}+i_{2})d-i_{2}d$ and $\delta=i_{1}d+d(-i_{1}+i_{2})-i_{2}d=(d\vee d)(i_{1}+(-i_{1}+i_{2})-i_{2})$ let
$\partial_{1}\partial_{2}$ be the composite
\begin{equation*}
\xymatrix{
\partial_{1}\partial_{2}\ar@{=>}[rr]^{i_{1}d+\Gamma-i_{2}d}&&\delta \ar@{=>}[rr]^{(d\vee d)\tau}&& 0.
}
\end{equation*}
Here $\tau$ is the appropriate linearity track.
\end{Def}
\begin{Def}
Let $d:A\longrightarrow B$ and $d':A'\longrightarrow B'$ be objects in the category of homotopy pairs. Then two
morphisms
\begin{equation*}
(d,g,F)\ ,\ (f',g',F'):\ d\longrightarrow d'
\end{equation*}
in {\bf HP} are {\em homotopic} if and only if there exists a map $a:B\longrightarrow A$ and there exist
tracks
\begin{equation*}
\alpha: (g,g') (-i_{1}+i_{2})\Longrightarrow d'a
\end{equation*}
\begin{equation*}
\beta: (f,a,f')\partial_{2}\Longrightarrow 0
\end{equation*}
such that pasting in the following diagram yields the trivial track $0\Longrightarrow 0$.
\begin{equation*}
\xymatrix{
& A\ar[dd]_{\partial_{2}} \ar@/^2.0pc/[rrdd]^0 \ar@/_4.7pc/[dddd]_0 \ar@{}[rrdd]|{\overset{\beta}{\Longrightarrow}}&&\\
&&&\\
&A\vee B\vee A \ar@{=>}[l]\ar[rr]^{(f,a,f')}\ar[dd]_{\partial_{1}}\ar@{}[rrdd]|{\overset{(F,\alpha,F')}{\Longrightarrow}}
&& A'\ar[dd]^{d'}\\
&&&\\
&B\vee B\ar[rr]_{(g,g')}&&B'
}
\end{equation*}
Here the left hand side of the diagram is given by the cylinder
$I(d)$ in (4.1).
\end{Def}
\begin{Rem}
A Toda complex $X=(X_{n},d_{n},H_{n},n\in \mathbb {Z})$ is a
sequence of maps $d_{n}:X_{n+1}\longrightarrow X_{n}$ and tracks
$H_{n}: d_{n+1}d_{n}\Longrightarrow 0$ in {\bf B}. This is a {\em
secondary chain complex} if all Toda pairs $(H_{n+1},H_{n})$
represent the zero element, see \cite{Ba1} (5.1). Let {\bf SCh} be
the category of secondary chain complexes with morphisms defined
as maps between Toda complexes,see \cite{Ba1} (5.3). We consider
an object $d:A \longrightarrow B$ in {\bf HP} as an object $X$ in
{\bf SCh} concentrated in degrees 0 and 1, that is, $X_{n}=*$ for
$n\neq 0,1$. This way we have the full inclusion of categories
${\bf HP}\subset{\bf SCh}$. In \cite{BG} we shall describe a
cylinder $I(X)$ in the category {\bf SCh} which generalizes the
cylinder $I(d)$ in (4.1). The cylinder $I(X)$ yields a homotopy
relation $\simeq$ on {\bf SCh} such that one has the full
inclusion of homotopy categories ${\bf HP}/\simeq\ \subset{\bf
SCh}/\simeq$.
\end{Rem}

\begin{theo}
Let {\bf B} be a triangulated track category with strict coproducts $X\vee Y$. Then the functor $\psi$ in (2.4) induces
an equivalence of categories
\begin{equation*}
\psi:{\bf HP}/\simeq\ \overset{\sim}{\longrightarrow} {\bf Prin}.
\end{equation*}
\end{theo}

Hence the homotopy category of the homotopy pairs in
{\bf B} is equivalent to the category of principal maps in
$\A=\B_{0}/\simeq$. Results of this kind under connectivity and
dimensional restrictions are described in \cite{Ba2}. We prove the
theorem in the next section.

\section{{Proof of the equivalence (4.4)}}
\setcounter{equation}{0}

Let $d:A\longrightarrow B$ be a map of ${\bf B_{0}}$ and let $C=C_{d}$ be given by a distinguished track triangle
extending $d$, see (2.2). We consider the following diagram in ${\bf B}$
\begin{equation*}
\xymatrix{
A\ar[r]^{d} \ar[d]\ar@{}[rd]|{\overset {H^{\boxminus}_{0}}{\Longrightarrow}}& B\ar@{=}[r]\ar[d]^{d_{0}}\ar@{}[rd]|{\overset{\Gamma}{
\Longrightarrow}}&B\ar[d]^{\partial}&\\
{*}\ar[r]_{0}\ar[d]&C\ar[r]_{-i_{1}+i_{2}\ \ \ }\ar[d]^{j}&C\vee C\ar[r]^{\ \ (\varphi,\varphi')}\ar[d]^{i}&X\\
tA\ar[r]^{v}&C_{d_{0}}\ar[r]^{w}&C_{\partial}\ar[ru]^{u}&
}
\end{equation*}
Here $\partial=i_{1}d_{0}+i_{2}d_{0}$ and $\Gamma$ is
defined by (3.2)(2). Moreover $H_{0}$ is given by (2.2). The
tracks in the diagram describe homotopy pairs in {\bf HP} and $wv$
is a map in ${\bf B_{0}}$ representing
$\psi(\Gamma\bar{\square}H_{0}^{\boxminus})$ where $\psi$ is the
functor in (2.4) and $w=\psi(\Gamma)$ and
$v=\psi(H_{0}^{\boxminus})$. Here $v$ is an isomorphism in $\A$ by
(2.2).
\begin{Lem}
Let $\varphi,\varphi':C\longrightarrow X$ in ${\bf B_{0}}$. Then $\varphi\simeq\varphi'$ if and only if there exists
$u:C_{\partial} \longrightarrow X$ in ${\bf B_{0}}$ with $ui=(\varphi,\varphi')$ and $uw=0$ in $\A$.
\end{Lem}

\begin{proof}
We have $\varphi\simeq\varphi'$ if and only if
$(\varphi,\varphi')(-i_{1}+i_{2})\simeq 0$. Hence
$\varphi\simeq\varphi'$ implies that
$(\varphi,\varphi')\partial=0$ in $\A$ and hence there is $u'$
with $u'i=(\varphi,\varphi')$ in $\A$. Then $u'wj=0$ so that $u'w$
factors in $\A$ through the cofiber $tB$ of $j$. This cofiber is
isomorphic to the cofiber of $i$. Therefore we can alter $u'$ by a
map $tB\longrightarrow X$ in such a way that we get $u$ with
$uw=0$. Conversely assume $u$ is given. Then $uv=0$ and hence
$(\varphi,\varphi')(-i_{1}+i_{2})=uwj=0$ in $\A$.
\end{proof}

Next we consider the following diagram in $\B$ where we use the cylinder object $I(d)$ in (4.1)
\begin{equation*}
\xymatrix{
A\ar[r]^{\partial_{2}\ \ \ \ \ } \ar[d]\ar@{}[rd]|{\overset {\partial^{\boxminus}_{12}}{\Longrightarrow}}& A\vee B\vee A
\ar[r]^{\ \ \ \ \ (0,1,0)}\ar[d]^{\partial_{1}}\ar@{}[rd]|{\overset{V}{\Longrightarrow}}& B\ar[d]^{\partial} \\
{*} \ar[r] \ar[d]& B\vee B \ar[r]_{d_{0}\vee d_{0}} \ar[d]& C\vee C \ar[d] \\
tA\ar[r]^{w'}& C_{\partial_{1}}\ar[r]^{v'}& C_{\partial}
}
\end{equation*}
Here $V$ is given by $(i_{1}H_{0}, 0^{\square},
i_{2}H_{0})$. We consider $V$ and $\partial^{\boxminus}_{12}$ as
morphisms in {\bf HP}. Let $v'$ and $w'$ be maps which represent
$\psi(V)$ and $\psi(\partial_{12}^{\boxminus})$ respectively. We
observe $v'$ is an isomorphism in $\A$.
\begin{Lem}
We have $v'w'=wv$ in $\A$.
\end{Lem}

\begin{proof}
We have $wv=\psi(\Gamma\,\bar{\square}\,H_{0}^{\boxminus})$ and
$v'w'=\psi(V\,\bar{\square}\,\partial_{12}^{\boxminus})$.
Therefore it suffices to show that in ${\bf HP}$ we have the
equation
\begin{equation}\tag{1}
\Gamma\,\bar{\square}\,H_{0}^{\boxminus}=V\,\bar{\square}\,\partial_{12}^{\boxminus}.
\end{equation}
Let $p_{1}=(1,0):C\vee C\longrightarrow C$ and $p_{2}=(0,1):C\vee C\longrightarrow C$ be the projections.
Then (1) holds if and only if (2) holds.
\begin{equation}\tag{2}
p_{i}(\Gamma\,\bar{\square}\,H_{0}^{\boxminus})=p_{i}(V\,\bar{\square}\,\partial_{12}^{\boxminus})
\ \ \ \ \textrm{for } i=1,2.
\end{equation}
This follows from the bilinearity of the coefficients
$D$ in (2.1). For convenience we now assume that additive tracks
$\tau$ in (3.1) are trivial tracks. Then we get
$p_{2}V=(0^{\square}, H_{0}),\ p_{2}\partial_{12}^{\boxminus}=
0^{\square}$ and $p_{2}\Gamma=0^{\square}$. Hence (2) for $i=2$ is
equivalent to
\begin{equation}\tag{3}
H_{0}^{\boxminus}=d_{0}d-H_{0}.
\end{equation}
This holds since for any track $H:f\Rightarrow 0$ one has
$0^{\square}=H-H=H\,\square\,(f-H)$.

Next we consider the case $i=1$ in (2). Then we have
$p_{1}V=(H_{0},0^{\square}),\
p_{1}\partial_{12}^{\boxminus}=d+\bar{\Gamma}$ with
$\bar{\Gamma}:d(-1)\Rightarrow (-1)d,\
p_{1}\Gamma=\bar{\Gamma}_{0}^{\boxminus}:(-1)d_{0}\Rightarrow
d_{0}(-1)$. Therefore (2) for $i=1$ is equivalent to
\begin{equation}\tag{4}
\bar{\Gamma}_{0}^{\boxminus}\,\bar{\square}\,H_{0}^{\boxminus}=(H_{0}^{\boxminus},0^{\square})\,\bar{\square}\,(d+\bar{\Gamma}).
\end{equation}
In fact, we have
\begin{equation}\tag{5}
(H_{0},0^{\square})\,\bar{\square}\,(d+\bar{\Gamma})=(H_{0},0^{\square})(i_{1}+i_{2}d)\,\square\,d_{0}
(d+\bar{\Gamma})
\end{equation}
\begin{equation*}
\ \ \ \ \ \ \ \ \ \ \ \ \ \ \ \ \ \ \ \ \ \ \ \ \ =
(H_{0}+d_{0}(-1)d)\,\square\,(d_{0}d+d_{0}\bar{\Gamma})
\end{equation*}
\begin{equation*}
\ \ \ \ \ \ \ \ \ \ \ \ \ \ \ \ \ \ \ =H_{0}+d_{0}\bar{\Gamma}:0\Longrightarrow d_{0}(-1)d
\end{equation*}
On the other hand we get:
\begin{equation}\tag{6}
\bar{\Gamma}_{0}^{\boxminus}\,\bar{\square}\,H_{0}^{\boxminus}=
\bar{\Gamma}_{0}^{\boxminus}d\,\square\,(-1)H_{0}^{\boxminus}
\end{equation}
\begin{equation*}
\ \ \ \ \ \ \ \ \ \ \ \ \ \ \ \ \ \ \ \ \ \ \ \ \ \ \ \ \ \ \ \ \
=d_{0}\bar{\Gamma}\,\square\,(-H_{0}^{\boxminus})\ \ \ \ \ \ \ \ \
\ \ \textrm{, see (7),}
\end{equation*}
\begin{equation*}
\ \ \ \ \ \ \ \ \ \ \ \ \ \ \ \ \ \ \ \ \ \ \ \ \ \ \ \ \ \ \ \
=d_{0}\bar{\Gamma}\,\square\,(H_{0}-d_{0}d)\ \ \ \ \ \ \textrm{,
see (3),}
\end{equation*}
\begin{equation*}
\ \ \ \ \ \,=H_{0}+d_{0}\bar{\Gamma}.
\end{equation*}
Here we use the fact that by the rules for $\Gamma$-tracks the following diagram represents the trivial track.
\begin{equation}\tag{7}
\xymatrix{
&&0&&\\
{\bullet}\ar@{}[rrd]|{\overset{\bar{\Gamma}}{\Longrightarrow}}\ar[rr]^{d}\ar[d]_{-1}\ar@/^2.5pc/[rrrr]&&
{\bullet} \ar@{}[rrd]|{\overset{\bar{\Gamma}_{0}}{\Longrightarrow}}\ar@{}[u]|{\ \ \Uparrow H_{0}}\ar[rr]^{d_{0}}\ar[d]^{-1}&&{\bullet}\ar[d]^{-1} \\
{\bullet}\ar@/_2.5pc/[rrrr]\ar[rr]_{d}&&{\bullet} \ar[rr]_{d_{0}}\ar@{}[d]|{\ \ \Downarrow H_{0}}&&{\bullet} \\
&&0&&}
\end{equation}
This completes the proof of (5.2).
\end{proof}

\begin{proof}[Proof of (4.4)] Let $d:A\longrightarrow B$ and
$d':A'\longrightarrow B'$ be objects in {\bf HP} and let
$F=(f,g,F)$ and $F'=(f',g',F')$ be morphisms $d\longrightarrow d'$
in {\bf HP}. We have to show that $F$ is homotopic to $F'$ as in
(4.2) if and only if $\psi(F)=\psi(F')$. By (4.2) we obtain the
following diagram in {\bf B}.
\begin{equation*}
\xymatrix{
&&0&&&\\
A \ar@/^2.5pc/[rrrr]\ar[rr]^{\partial_{2}}\ar[d]\ar@{}[rrd]|{\overset{\partial_{12}^{\boxminus}}{\Longrightarrow}} &&
A\vee B\vee A \ar@{}[rrd]|{\overset{(F,\alpha,F')}{\Longrightarrow}}\ar[rr]^{(f,a,f')}\ar[d]^{\partial_{1}}\ar@{}[u]|{\ \ \Uparrow \beta} && A'\ar[d]^{d'}\ar@/^2.5pc/[dd]&\\
{*}\ar[rr]\ar[d]&&B\vee B \ar[d]\ar[rr]_{(g,g')}&&B'\ar[d]\ar@{}[r]|{\overset{H'_{0}}{\Rightarrow}}&0\\
tA\ar[rr]^{w'}&&C_{\partial_{1}}\ar[rr]^{W}&&C_{d'}&\\
&&C_{d}\vee C_{d}\ar[u]_{i}\ar[rru]_{(\psi,\psi')}&&&
}
\end{equation*}
Here $w'$, $W$ and $i$ are maps representing
$\psi(\partial_{12}^{\boxminus})$, $\psi(F,\alpha,F')$ and $d\vee
d\subset\partial_{1}$ respectively and $\psi$ and $\psi'$ are maps
representing $\psi(F)$, resp. $\psi(F')$. Since $\psi$ is a
functor we have $Wi=(\psi,\psi')$ in $\A$.

Now assume $F\simeq F'$ in ${\bf HP}$ so that $a, \alpha, \beta$
are given. Then we can choose $W$ as above and we observe that
$Ww'$ in $\A$ is represented by pasting tracks
$\partial_{12}^{\boxminus}$, $(F, \alpha, F')$ and $H'_{0}$.
Compare the proof of (2.4). This pasting yields the same element
as the pasting of $\partial_{12}^{\boxminus}$, $(F, \alpha, F')$
and $\beta$ which by definition of $F\simeq F'$ in (4.2)
represents 0. Hence $Ww'=0$ and therefore (5.1) and (5.2) shows
$\psi\simeq \psi'$ since $v$ and $v'$ are isomorphisms in $\A$.

On the other hand assume now $\psi\simeq \psi'$. Then by (5.1) and
(5.2) there exists $W$ with $Wi=(\psi,\psi')$ in $\A$ and $Ww'=0$
in $\A$. Since the functor $\psi$ in (2.4) is full we find $a$,
$\alpha$ as above with $\psi(F,\alpha,F')=W$ and $\psi(F)=\psi$,
$\psi(F')=\psi'$ in $\A$. Moreover there exists
$\beta':(f,a,f')\partial_{2}\Rightarrow 0$ since $Ww'=0$ induces
$t((f,a,f')\partial_{2})=0$. Now pasting
$\partial_{12}^{\boxminus}$, $(F,\alpha,F')$ and $\beta$ yields an
element $\xi\in Hom_{\A}(tA,B')$. This element maps to $0$ in
$Hom_{\A}(tA,C_{d'})$ since $Ww'=0$. Hence there exists $\xi'\in
Hom(tA,A')$ which maps to $\xi$, that is, $d'\xi'=\xi$. Now we
alter $\beta'$ by $\xi'$ and get $\beta$ such that pasting of
$\partial_{12}^{\boxminus}$, $(F,\alpha,F')$ and $\beta$
represents the zero element. This shows $F\simeq F'$ in {\bf HP}.

\end{proof}

\section{{The category of $2$-stage spectra}}
\setcounter{equation}{0}

Let $k$ be a commutative ring with unit and let {\bf Mod} be the category of $k$-modules and $k$-linear maps. This
is a symmetric monoidal category via the tensor product $A\otimes B$ over $k$. A {\em pair} of
modules is a morphism
\begin{equation*}
X=\left(X_{1}\overset{\partial}{\longrightarrow}X_{0}\right)
\end{equation*}
in {\bf Mod}. A {\em morphism} $f: X\longrightarrow Y$ of pairs is a commutative diagram
\begin{equation*}
\xymatrix{
X_{1}\ar[r]^{f_{1}}\ar[d]_{\partial}&Y_{1}\ar[d]^{\partial}\\
X_{0}\ar[r]^{f_{0}}&Y_{0}
}
\end{equation*}
Clearly a pair in {\bf Mod} coincides with a chain complex concentrated in degrees 0 and 1. For two pairs $X$
and $Y$ the tensor product of the chain complexes corresponding to them is concentrated in degrees 0, 1 and 2 and is given by
\begin{equation*}
X_{1}\otimes Y_{1}\overset{\partial_{1}}{\longrightarrow} X_{1}\otimes Y_{0}\oplus X_{0}\otimes Y_{1}\overset{\partial_{0}}{\longrightarrow}
X_{0}\otimes Y_{0}
\end{equation*}
with $\partial_{0}=(\partial \otimes 1,1\otimes \partial)$ and $\partial_{1}=(-1\otimes \partial, \partial\otimes 1)$.
Truncating $X\otimes Y$ we get the pair
\begin{equation*}
X\bar{\otimes}Y=\left((X\bar{\otimes}Y)_{1}=coker(\partial_{1})\overset{\partial}{\rightarrow}X_{0}\otimes Y_{0}=(X\bar{\otimes}Y)_{0}\right)
\end{equation*}
with $\partial$ induced by $\partial_{0}$.

 We define the tensor
product $A\otimes B$ of two graded modules in the usual way, i.e.
by
\begin{equation*}
(A \otimes B)^{n}=\bigoplus_{i+j=n}A^{i}\otimes B^{j}.
\end{equation*}

A {\em (graded) pair module} is a sequence
$X^{n}=(\partial:X_{1}^{n}\rightarrow X_{0}^{n})$ with $n\in
\mathbb{Z}$ of pairs in {\bf Mod}.
 We identify such a graded pair module $X$ with the underlying morphism $\partial$ of degree 0 between graded modules
\begin{equation*}
X=\left(X_{1}\overset{\partial}{\longrightarrow}X_{0}\right).
\end{equation*}
Now the tensor product $X\bar{\otimes}Y$ of graded pair modules $X$, $Y$ is defined by
\begin{equation*}
(X \bar{\otimes} Y)^{n}=\bigoplus_{i+j=n}X^{i}\bar{\otimes} Y^{j}.
\end{equation*}
This defines a monoidal structure on the category of graded pair modules. Morphisms in this category are of
degree 0.
For two morphisms $f$, $g : X\longrightarrow Y$ between graded pairs modules, a {\em homotopy} $H:f\Longrightarrow g$
is a morphism $H:X_{0}\longrightarrow Y_{1}$ of degree 0 as in the diagram
\begin{equation*}
\xymatrix{
X_{1}\ar[d]_{\partial}\ar@<0.5ex> [r]^{f_{1}}\ar@<-0.5ex> [r]_{g_{1}}& Y_{1}\ar[d]^{\partial}\\
X_{0}\ar@{->}[ru]|{H}\ar@<0.5ex> [r]^{f_{0}}\ar@<-0.5ex> [r]_{g_{0}}&Y_{0}
}
\end{equation*}
satisfying $f_{0}-g_{0}=\partial H$ and $f_{1}-g_{1}=H\partial$.
A {\em pair algebra} $B$ is a monoid in the monoidal category of graded pair modules, with multiplication
\begin{equation*}
\mu : B\bar{\otimes}B\longrightarrow B.
\end{equation*}

We assume that B is concentrated in nonnegative degrees, that is
$B^{n}=0$ for $n<0$. A {\em left B-Module} is a graded pair module
$M$ together with a left action
\begin{equation*}
\mu:B\bar{\otimes}M\longrightarrow M
\end{equation*}
of the monoid $B$ on $M$.
More explicitly pair algebras and modules over them can be described as follows.
\begin{Def}
A {\em pair algebra} $B$ is a graded pair
\begin{equation*}
\partial:B_{1}\longrightarrow B_{0}
\end{equation*}
in {\bf Mod} with $B^{n}_{1}=B^{n}_{0}=0$ for $n<0$ such that $B_{0}$ is a graded algebra in {\bf Mod}, $B_{1}$ is a graded
$B_{0}$-$B_{0}$-bimodule, and $\partial$ is a bimodule homomorphism. Moreover for $x$, $y$ $\in$ $B_{1}$ the equality
\begin{equation*}
\partial(x)y=x\partial(y)
\end{equation*}
holds in $B_{1}$.
\end{Def}
\begin{Def}
A {\em (left) module} over a pair algebra $B$ is a graded pair $M=(\partial:M_{1}\longrightarrow M_{0})$ in {\bf Mod}
such that $M_{1}$ and $M_{0}$ are left $B_{0}$-modules and $\partial$ is $B_{0}$-linear. Moreover a $B_{0}$-linear map
\begin{equation*}
\bar{\mu}:B_{1}\otimes_{B_{0}}M_{0}\longrightarrow M_{1}
\end{equation*}
is given fitting in the commutative diagram
\begin{equation*}
\xymatrix{
B_{1}\otimes_{B_{0}}M_{1}\ar[d]_{\mu}\ar[r]^{1\otimes\partial}&B_{1}\otimes_{B_{0}}M_{0}\ar[d]^{\mu}\ar[dl]^{\bar{\mu}}\\
M_{1}\ar[r]_{\partial}&M_{0}
}
\end{equation*}
where $\mu(b\otimes m)=\partial(b)m$ for $b \in B_{1}$ and $m \in M_{1} \cup M_{0}$.
\end{Def}
For an indeterminate element $x$ of degree $n=|x|$ let $B[x]$ denote the $B$-module with $B[x]_{i}$ consisting
of expressions $bx$ with $b \in B_{i}$, $i=$0,1, with $bx$ having degree $|b| + n$, and structure maps given by
$\partial(bx)=\partial(b)x$, $\mu(b'\otimes bx)= (b'b)x$ and $\bar{\mu}(b' \otimes bx)=(b'b)x$.
A {\em free B}-module is direct sum of several copies of modules of the form $B[x]$, with $x\in I$ for some set $I$
of indeterminates of possibly different degrees. It will be denoted
\begin{equation*}
B[I]=\bigoplus_{x\in I}B[x].
\end{equation*}

Let $B$-{\bf Mod} be the category of left modules over the pair algebra $B$. Morphisms $f=(f_{0},f_{1}):M\longrightarrow N$
are pair morphisms which are $B$-equivariant, that is, $f_{0}$ and $f_{1}$ are $B_{0}$-equivariant and compatible with $\bar{\mu}$
above, i.e. the diagram
\begin{equation*}
\xymatrix{
B_{1}\otimes_{B_{0}}M_{0}\ar[d]_{1\otimes f_{0}}\ar[r]^{\ \ \ \ \ \bar{\mu}}&M_{1}\ar[d]^{f_{1}}\\
B_{1}\otimes_{B_{0}}N_{0}\ar[r]^{\ \ \ \ \ \bar{\mu}}&N_{1}
}
\end{equation*}
commutes.

For two such maps $f$, $g : M\longrightarrow N$ a track $H: f \Longrightarrow g$ is a degree zero map
\begin{equation*}
H:M_{0}\longrightarrow N_{1}
\end{equation*}
satisfying $f_{0}-g_{0}=\partial H$ and $f_{1}-g_{1}=H\partial$
such that $H$ is $B_{0}$-equivariant. For tracks
$H:f\Longrightarrow g$, $K:g\Longrightarrow h$ their composition
$K\square H: f\Longrightarrow h$ is $K+H$.

 The category $B$-{\bf Mod} with these tracks is a well defined
additive track category with strong direct sums $M\oplus N$. Let
\begin{equation}
B\textrm{-}{\bf mod}\subset B\textrm{-}{\bf Mod}
\end{equation}
be the full subcategory consisting of finitely generated free $B$-modules. The category of homotopy pairs
{\bf HP}($B$-{\bf mod}) is defined dually to (2.3). Morphisms $d'\longrightarrow d$ are diagrams
\begin{equation}
\xymatrix{
A\ar@{}[rd]|{\overset {F}{\Leftarrow}}& A'\ar[l]_{f}\\
B\ar[u]^{d}& B'\ar[l]^{g}\ar[u]_{d'}
}
\end{equation}
in $B$-{\bf mod}. Two such morphisms $(f, g, F)$, $(f', g', F):d'\longrightarrow d$ are {\em homotopic} if there
exists $(a, \alpha, \beta)$ as in the diagram
\begin{equation*}
\xymatrix{
A\ar@{}[rrd]|{\overset{\beta}{\Leftarrow}}&&\\
A\oplus B\oplus A \ar[u]^{\partial_{2}}\ar@{}[rrd]|{\overset{(F,\alpha,F')}{\Longleftarrow}}&&A'\ar@/_1.5pc/[llu]\ar[ll]_{\ \ \ \ \ \ (f,a,f')}\\
B\oplus B\ar[u]^{\partial_{1}}&&B'\ar[ll]^{\ \ \ \ \ (g,g')}\ar[u]_{d'}
}
\end{equation*}
with $\partial_{1}=(dp_{1},-p_{1}+p_{2},dp_{2})$ and $\partial_{2}=p_{1}+dp_{2}-p_{3}$. Here $p_{1}$, $p_{2}$,
$p_{3}$ denote the projections. This defines a natural equivalence relation $\simeq$ on {\bf HP}($B$-{\bf mod}).
\begin{theo}
Let $B$ be the pair algebra of secondary cohomology operations
computed in \cite{Ba3}. Then one has an equivalence of categories
\begin{equation*}
\left({\bf HP}(B\textrm{- }{\bf mod})/\simeq\right)^{op}\overset{\sim}{\longrightarrow} 2\textrm{- }{\bf stage}.
\end{equation*}
Here $2$-{\bf stage} is the category of $2$-stage spectra and
principal maps in (1.5).

\end{theo}
\begin{proof} It is proved in \cite{Ba3} (5.5.6) that the track category
of finite direct sums of Eilenberg-Mac Lane spectra is equivalent
to the track category ($B$-{\bf mod}$)^{op}$. This equivalence is
compatible with strong coproducts and hence induces the
equivalence in the theorem by use of theorem (4.4). Since
composition in $B$-{\bf mod} is bilinear $\Gamma$-tracks are
trivial in the notion of homotopy in (6.4).
\end{proof}


\title{{[IV]}\\\hspace{5pt}\\Cohomologically triangulated categories}

\author{H.-J. Baues and F. Muro}%

\blfootnote{\textit{$2000$ Mathematics Subject Classification:} 18E30, 18D05}%
\blfootnote{\textit{Key words and phrases:} triangulated categories, pretriangulated categories, groupoid-enriched categories}%
\blfootnote{The second author was partially supported by the MEC
grant MTM2004-01865 and postdoctoral fellowship EX2004-0616}

\maketitle

\begin{abstract}
We define cohomologically triangulated categories as triples
$(\C{A},t,\nabla)$ given by an additive category $\C{A}$, an
additive equivalence $t\colon\C{A}\r\C{A}$ and a cohomology class
$\nabla$ in the translation cohomology $H^3(\C{A},t)$. Here $\nabla$
is triangulated if it is represented by a good translation track category
$(\C{B},s)$ satisfying two natural axioms (TTr1) and (TTr2). It is
proved that in this case $(\C{A},t)$ is a triangulated category in
which exact triangles depend only on the class $\nabla$.
\end{abstract}


\section*{Introduction}

A stable homotopy theory (like the homotopy category of spectra or
chain complexes) yields a groupoid-enriched category $\C{B}$ in
which $2$-cells are tracks, i. e. homotopy classes of homotopies.
The suspension determines a pseudofunctor $s\colon\C{B}\pf\C{B}$
such that 
by the classification result in
\cite{ptc} one obtains a characteristic 
cohomology class
$$\nabla=\grupo{\C{B},[s]}\in H^3(\C{A},t).$$
Here $\C{A}$ is the homotopy category of $\C{B}$ and $t$ is induced
by $s$. Extending the main result in \cite{ttc} we show that if the
pair $(\C{B},s)$ satisfies two natural axioms (TTr1) and (TTr2) on
the so-called track triangles then $(\C{A},t)$ is a triangulated
 category in which exact triangles are induced by the track
triangles. Moreover, we show that the exact triangles are determined
by the cohomology class $\nabla$ and do not depend on the choice of
$(\C{B},s)$ representing $\nabla$. Therefore we call the triple
$(\C{A},t,\nabla)$ a cohomologically triangulated category and we
call $\nabla$ a triangulated cohomology class. In \cite{ccctc2} we
characterize triangulated cohomology classes in a purely
cohomological way.

The authors are very grateful to Bernhard Keller for discussions on
an earlier version of this paper.

\subsection*{Notation}

In this paper the arrow $\r$ is used for morphisms in ordinary
categories, $1$-cells in $2$-categories and functors. We use $\rr$
for $2$-cells in $2$-categories and natural transformations, and
$\pf$ for pseudofunctors. Identity morphisms or $1$-cells are
denoted by $1$, with the object as a subscript when necessary. For
the identity $2$-cells and natural transformations we use the
symbol $0^\vc$. Horizontal composition in $2$-categories is
denoted by juxtaposition, as the composition law in ordinary
categories, while we use the symbol $\vc$ for vertical
composition. The vertical inverse of an invertible $2$-cell or
natural transformation $\alpha$ is $\alpha^\vi$.

\section{Triangulated categories}\label{tridef}


Let $\C{A}$ be an additive category with zero object $*$ and let
$t\colon\C{A}\r\C{A}$ be an additive equivalence, hence
$$t\colon\hom_\C{A}(X,Y)\cong\hom_\C{A}(tX,tY)$$
is an isomorphism of abelian groups and $t(X\oplus
Y)=(tX)\oplus(tY)$. A \emph{small candidate
triangle}\index{triangle!small candidate} in $\C{A}$ is a diagram
in $\C{A}$ of the form
$$A\st{f}\To B\st{i}\To C\st{q}\To tA$$
such that $if=0$, $qi=0$ and $(tf)q=0$. A morphism $k$ between small
candidate triangles is a commutative diagram
$$\xymatrix{A\ar[r]^f\ar[d]^{k_0}&B\ar[r]^i\ar[d]^{k_1}&C\ar[r]^q\ar[d]^{k_2}&tA\ar[d]^{tk_0}\\
\bar{A}\ar[r]^{\bar{f}}&\bar{B}\ar[r]^{\bar{i}}&\bar{C}\ar[r]^{\bar{q}}&t\bar{A}}$$
It is an isomorphism if the vertical arrows are isomorphisms in
$\C{A}$. The category of small candidate triangles will be denoted
by $\C{cand}(\C{A},t)$.

\begin{defn}\label{pt}
We say that $\C{A}$ is a \emph{pretriangulated
category}\index{pretriangulated category} if it is equipped with a
distinguished family $\mathcal{E}$ of candidate triangles (called
\emph{exact}\index{triangle!exact} triangles) which is subject to
the following axioms:
\begin{itemize}
\item[(Tr0)] Any candidate triangle isomorphic to an exact triangle is exact. For any object $A$ in $\C{A}$ the following candidate triangle is exact
$$A\st{1_A}\To A\To *\To tA.$$

\item[(Tr1)] Any morphism $f\colon A\r B$ in $\C{A}$ is part of an exact triangle
$$A\st{f}\To B\st{i}\To C\st{q}\To tA.$$

\item[(Tr2)] A candidate triangle $A\st{f}\r B\st{i}\r C\st{q}\r tA$ is exact if and only if
$B\st{i}\r C\st{q}\r tA\st{-tf}\r tB$ is exact.

\item[(Tr3)] Any commutative diagram
$$\xymatrix{A\ar[r]^f\ar[d]_{k_0}&B\ar[r]^i\ar[d]^{k_1}&C\ar[r]^q&tA\ar[d]^{tk_0}\\\bar{A}\ar[r]_{\bar{f}}&\bar{B}\ar[r]_{\bar{i}}&\bar{C}\ar[r]_{\bar{q}}&t\bar{A}}$$
whose rows are exact triangles can be extended to a morphism $k$ between candidate triangles.
\end{itemize}

We say that $\C{A}$ is a \emph{triangulated
category}\index{triangulated category} if in addition it
satisfies:

\begin{itemize}
\item[(Tr4)] The \emph{octahedral axiom}\index{octahedral axiom}. Given two morphisms $A\st{f}\r B\st{g}\r C$ in $\C{A}$, if the rows of the commutative diagram
$$\xymatrix{A\ar[r]^f\ar@{=}[d]&B\ar[r]^{i_f}\ar[d]^{g}&C_f\ar[r]^{q_f}&tA\ar@{=}[d]\\
A\ar[d]^f\ar[r]^{gf}&C\ar@{=}[d]\ar[r]^{i_{gf}}&C_{gf}\ar[r]^{q_{gf}}&tA\ar[d]^{tf}\\
B\ar[r]^g&C\ar[r]^{i_g}&C_g\ar[r]^{q_g}&tB}$$ are exact triangles
then this diagram admits a commutative extension
$$\xymatrix{A\ar[r]^f\ar@{=}[d]&B\ar[r]^{i_f}\ar[d]^{g}&C_f\ar[r]^{q_f}\ar[d]^{\bar{g}}&tA\ar@{=}[d]\\
A\ar[d]^f\ar[r]^{gf}&C\ar@{=}[d]\ar[r]^{i_{gf}}&C_{gf}\ar[d]^{\bar{f}}\ar[r]^{q_{gf}}&tA\ar[d]^{tf}\\
B\ar[r]^g&C\ar[r]^{i_g}&C_g\ar[r]^{q_g}&tB}$$
such that
$$C_f\st{\bar{g}}\To C_{gf}\st{\bar{f}}\To C_g\st{t(i_f)q_g}\To tC_f$$
is an exact triangle.
\end{itemize}

\end{defn}

\section{Triangulated track categories and cohomologically triangulated categories}\label{track}

A \emph{track category}\index{track category} $\C{A}$ is a
category enriched in groupoids. Hence for objects $X, Y$ in
$\C{A}$ one has the hom-groupoid $\homm{X,Y}_{\C{A}}$ in $\C{A}$.
The objects $f\colon X\r Y$ of this groupoid are called
\emph{maps}, and the morphisms $\alpha\colon f\rr g$ are called
\emph{tracks}. Let $\homm{X,Y}_\C{A}(f,g)$ be the set of all
tracks $f\rr g$. The symbol for the identity track in $f$ is
$0^\vc_f$. Equivalently, a track category is a $2$-category all of
whose $2$-cells are invertible. Ordinary categories can be
regarded as track categories with only the identity tracks.

A \emph{pseudofunctor}\index{pseudofunctor} between track
categories $\varphi\colon\C{A}\pf\C{B}$ is an assignment of
objects, maps and tracks which preserves composition and
identities only up to certain given tracks
$$\varphi_{f,g}\colon\varphi(f)\varphi(g)\rr\varphi(fg)\text{ and
} \varphi_X\colon\varphi(1_X)\rr1_{\varphi(X)}.$$ These tracks
must satisfy well-known coherence and naturality properties, see
for example  the Appendix of \cite{ptc} for details.  A
\emph{pseudonatural transformation}\index{pseudonatural
transformation} $\alpha\colon\varphi\rr\psi$ between two
pseudofunctors $\varphi,\psi\colon\C{A}\pf\C{B}$ is a collection
of maps $\alpha_X\colon\varphi(X)\r\psi(X)$ in $\C{B}$ indexed by
the objects of $\C{A}$ such that the usual squares only commute up
to given tracks
$\alpha_f\colon\alpha_Y\varphi(f)\rr\psi(f)\alpha_X$ indexed by
the maps $f\colon X\r Y$ in $\C{A}$. These tracks must satisfy
certain axioms, see also the Appendix of \cite{ptc}.

A \emph{homotopy}\index{homotopy between pseudofunctors}
$\xi\colon\varphi\rr\psi$ between two pseudofunctors
$\varphi,\psi\colon\C{A}\pf\C{B}$ is a pseudonatural
transformation $\xi$ such that the maps
$\xi_X\colon\varphi(X)\r\psi(X)$
are identities, i. e. $\varphi(X)=\psi(X)$ and 
$\xi_X=1_{\varphi(X)}=1_{\psi(X)}$. Two pseudofunctors $\varphi$
and $\psi$ are \emph{homotopic} $\varphi\simeq\psi$ if there
exists a homotopy between them. Since homotopies are invertible
pseudonatural transformations the homotopy relation is a natural
equivalence relation in the category of track categories and
pseudofunctors. See \cite{ptc} 3 for details.

In the rest of this section $(\C{A},t)$ will be a translation
category with $\C{A}$ and $t$ additive and $t$ an equivalence. We
will consider the following $\C{A}$-bimodule (i. e. functor
$\C{A}^\op\times\C{A}\r\C{Ab}$ where $\C{Ab}$ is the category of
abelian groups)
\begin{equation}\label{homt}
\hom^t=\hom_\C{A}(t,-)
\end{equation}
and
the natural transformation
\begin{equation}\label{tbar}
\bar{t}\colon\hom^t=\hom_\C{A}(t,-)\rr\hom_\C{A}(t^2,t)=\hom^t(t,t)
\end{equation} given by $(-1)t$.
Recall that in this case the translation cohomology of
$(t,\bar{t})$ defined in \cite{ptc} 5 is denoted by
$$H^*(\C{A},t)=H^*(t,\bar{t}).$$

\begin{defn}\label{trantra}
A \emph{translation track category}\index{translation track
category} $(\C{B},[s])$ over $(\C{A},t)$ is a linear track
extension in the sense of \cite{ccglg}
\begin{equation}
\hom^t\st{\sigma}\r\C{B}_1\rightrightarrows\C{B}_0\st{p}\rightarrow\C{A}
\end{equation}
consisting of a track category $\C{B}$, a functor
$p\colon\C{B}_0\r\C{A}$ from the underlying ordinary category
$\C{B}_0$ of $\C{B}$ such that $p$ is the identity on objects and
$p(f)=p(g)$ if and only if there exists a track $f\rr g$, and
isomorphisms
$$\sigma_f\colon\hom^t(X,Y)\cong\aut_{\homm{X,Y}_\C{B}}(f)=\homm{X,Y}_\C{B}(f,f)$$
for all maps $f\colon X\r Y$ in $\C{B}$. These isomorphisms
satisfy the following properties, given a track $\alpha\colon f\rr
g$ in $\C{B}$,
\begin{equation}\label{lec1}
\alpha\vc\sigma_f(x)=\sigma_g(x)\vc\alpha;
\end{equation}
and given composable maps
$\bullet\st{h}\r\bullet\st{g}\r\bullet\st{f}\r\bullet$ in $\C{B}$
\begin{equation}\label{lec2}
f\sigma_g(x)=\sigma_{fg}(p(f)_*x) \text{ and }
\sigma_g(x)h=\sigma_{gh}(p(h)^*x).
\end{equation}
Also as part of the structure of a translation track category
there is a homotopy class of pseudofunctors
$[s]\colon\C{B}\pf\C{B}$ such that for some representative $s$
(and hence for anyone) one has that for all maps $f\colon X\r Y$
in $\C{B}$ and $a\in\hom^t(X,Y)$ the equalities $p(s(f))=t(p(f))$
and $s(\sigma_f(a))=\sigma_{s(f)}(\bar{t}(a))$ hold.

Two translation track categories $(\C{B},[s])$ and $(\C{C},[r])$ over $(\C{A},t)$
are equivalent if there exists a pseudofunctor
$\varphi\colon\C{B}\pf\C{C}$ such that for all maps $f\colon X\r
Y$ in $\C{B}$ and $a\in\hom^t(X,Y)$ the equalities
$p(\varphi(f))=p(f)$ and
$\sigma_{\varphi(f)}(a)=\varphi(\sigma_f(a))$ are satisfied and $\varphi s$ is homotopic to $r\varphi$.
\end{defn}

The next result is the main theorem in \cite{ptc}.

\begin{thm}\label{clastrans}
Equivalence classes of translation track categories over $(\C{A},t)$ are
in bijective correspondence with the translation cohomology classes in
$H^3(\C{A},t)$. This bijection is determined by the characteristic
translation cohomology class of a translation track category
$(\C{B},[s])$ over $(\C{A},t)$
$$\grupo{\C{B},[s]}\in H^3(\C{A},t)$$
defined in the proof of \cite{ptc} 6.3.
\end{thm}

Hence for each class $\nabla\in H^3(\C{A},t)$ we can choose an associated
translation track category $(\C{B},[s])$ with
$$\grupo{\C{B},[s]}=\nabla.$$
By \cite{ptc} A.13 and A.14 we can always assume that the zero
object $*$ of $\C{A}$ is a strict zero object in $\C{B}$, i. e.
$\homm{*,X}_\C{B}$ and $\homm{X,*}_\C{B}$ are the trivial groupoids
for any object $X$ in $\C{B}$, in particular the zero map $0\colon
X\r Y$ is defined for any two objects $X, Y$ in $\C{B}$ as the
composite $X\r*\r Y$. We can also assume that the pseudofunctor $s$
representing the homotopy class $[s]$ (which preserves $*$, since
$s$ coincides on objects with $t$) is normalized at zero maps, i. e.
$s(0)=0$, $s_{0,f}=0_0^\vc$ and $s_{f,0}=0_0^\vc$ for any map $f$ in
$\C{B}$. In this situation we say tat $(\C{B},s)$ is a \emph{good
translation track category}\index{translation track category!good}.

A \emph{track triangle}\index{track triangle} in a good translation
track category $(\C{B},s)$, compare \cite{ttc}, is a diagram in
$\C{B}$
\begin{equation}\label{dtt}
\xymatrix{A\ar@/^30pt/[rr]^0_{\;}="a"\ar[r]_{f}&B\ar@/_30pt/[rr]_0^{\;}="c"
\ar[r]^<(-0.22){\;}="b"|{i}&
C\ar[r]_<(-0.17){\;}="d"|{q}\ar@/^30pt/[rr]^0_{\;}="f"&sA\ar[r]_{s(f)}^<(-0.23){\;}="e"&
sB\ar@{=>}"b";"a"_{H_0}\ar@{=>}"d";"c"^{H_{1}}\ar@{=>}"e";"f"_{H_{2}}}
\end{equation}
such that the following equalities hold
\begin{eqnarray*}
(qH_0)\vc(H_1^\vi f)&=&\sigma_{0}(1_{tA}),\\
(s(f)H_1)\vc(H_2^\vi i)&=&\sigma_{0}(-1_{tB}),\\
(s(i)H_2)\vc((s(H_0)\vc s_{i,f})^\vi q)&=&\sigma_{0}(1_{tC}).
\end{eqnarray*}
A morphism of track triangles is a diagram in
$\C{B}$
\begin{equation}\label{morctra}
\xymatrix{A\ar[r]^f_<(.6){\;}="a"\ar[d]_{k_0}^<(.6){\;}="b"&B\ar[r]^i_<(.6){\;}="c"\ar[d]|{k_1}^<(.6){\;}="d"&
C\ar[r]^q_<(.6){\;}="e"\ar[d]|{k_2}^<(.6){\;}="f"&sA\ar[d]^{sk_0}\\
\bar{A}\ar[r]_{\bar{f}}&\bar{B}\ar[r]_{\bar{i}}&\bar{C}\ar[r]_{\bar{q}}&s\bar{A}
\ar@{=>}"a";"b"^{K_0}\ar@{=>}"c";"d"^{K_1}\ar@{=>}"e";"f"^{K_2}}
\end{equation}
such that the next three equalities are satisfied
\begin{eqnarray*}
k_2H_0&=&(\bar{H}_0k_0)\vc(\bar{i}K_0)\vc(K_1f),\\
s(k_0)H_1&=&(\bar{H}_1k_1)\vc(\bar{q}K_1)\vc(K_2i),\\
s(k_1)H_2&=&(\bar{H}_2k_2)\vc(s(\bar{f})K_2)\vc((s_{\bar{f},k_0}^\vi\vc s(K_0)\vc s_{k_1,f})q).\\
\end{eqnarray*}
Here $\bar{H}_i$ are the structure tracks of the lower track
triangle in (\ref{morctra}). Another morphism
\begin{equation}\label{morctra2}
\xymatrix{A\ar[r]^f_<(.6){\;}="a"\ar[d]_{\bar{k}_0}^<(.6){\;}="b"&B\ar[r]^i_<(.6){\;}="c"\ar[d]|{\bar{k}_1}^<(.6){\;}="d"&
C\ar[r]^q_<(.6){\;}="e"\ar[d]|{\bar{k}_2}^<(.6){\;}="f"&sA\ar[d]^{s\bar{k}_0}\\
\bar{A}\ar[r]_{\bar{f}}&\bar{B}\ar[r]_{\bar{i}}&\bar{C}\ar[r]_{\bar{q}}&s\bar{A}
\ar@{=>}"a";"b"^{\bar{K}_0}\ar@{=>}"c";"d"^{\bar{K}_1}\ar@{=>}"e";"f"^{\bar{K}_2}}
\end{equation}
is \emph{homotopic} to the morphism in (\ref{morctra}) if there
exist tracks $L_i\colon k_i\rr \bar{k}_i$ $(i=0,1,2)$ such that
the following equations are satisfied
\begin{eqnarray*}
(\bar{f}L_0)\vc K_0&=&\bar{K}_0(L_1f),\\
(\bar{i}L_1)\vc K_1&=&\bar{K}_1(L_2i),\\
(\bar{q}L_2)\vc K_2&=&\bar{K}_2(s(L_0)q).\\
\end{eqnarray*}
Notice that if the tracks $L_i$ are given first then the diagram (\ref{morctra2}) defined by the previous three equations automatically defines a morphism of
track triangles which, in addition, is homotopic to (\ref{morctra}).

The category of track triangles in $(\C{B},s)$ and
homotopy classes of morphisms will be denoted by
$\C{candt}(\C{B},s)$. There is an obvious forgetful functor
\begin{equation}\label{barro}
\varrho=\varrho_{(\C{B},s)}\colon\C{candt}(\C{B},s)\To\C{cand}(\C{A},t)
\end{equation} sending a morphism in the source category represented by (\ref{morctra}) to
the following morphism of small candidate triangles in the sense
of Section \ref{tridef}
$$\xymatrix{A\ar[r]^{p(f)}\ar[d]_{p(k_0)}&B\ar[r]^{p(i)}\ar[d]^{p(k_1)}&C\ar[r]^{p(q)}\ar[d]^{p(k_2)}&
tA\ar[d]^{tp(k_0)}\\\bar{A}\ar[r]_{p(\bar{f})}&\bar{B}\ar[r]_{p(\bar{i})}&\bar{C}\ar[r]_{p(\bar{q})}&t\bar{A}}$$

\begin{prop}\label{sontodos}
Let $(\C{B},s)$ and $(\C{C},r)$ be good translation track categories
representing the same cohomology class in $H^3(\C{A},t)$. Then a
small candidate triangle in $(\C{A},t)$ is in the image of
$\varrho_{(\C{B},s)}$ if and only if it is in the image of
$\varrho_{(\C{C},r)}$.
\end{prop}

\begin{proof}
Let $\varphi\colon\C{B}\pf\C{C}$ be a pseudofunctor realizing the
equivalence between $(\C{B},s)$ and $(\C{C},r)$ as in Definition
\ref{trantra} and let $\xi\colon r\varphi \rr \varphi s$ be a
homotopy. Recall that the pseudofunctors $s$ and $r$ coincide on
objects (they both actually coincide on objects with $t$) and that
$\varphi$ is the identity on objects. In particular by \cite{ptc}
A.14 we can suppose without loss of generality that $\varphi$ is
normalized with respect to zero maps. The pseudofunctor $\varphi$
and the homotopy $\xi$ induce a functor
$$(\varphi,\xi)_*\colon\C{candt}(\C{B},s)\To\C{candt}(\C{C},r)$$
sending a track triangle
\begin{equation*}
\xymatrix{A\ar@/^30pt/[rr]^0_{\;}="a"\ar[r]_{f}&B\ar@/_30pt/[rr]_0^{\;}="c"
\ar[r]^<(-0.22){\;}="b"|{i}&
C\ar[r]_<(-0.17){\;}="d"|{q}\ar@/^30pt/[rr]^0_{\;}="f"&sA\ar[r]_{s(f)}^<(-0.23){\;}="e"&
sB\ar@{=>}"b";"a"_{H_0}\ar@{=>}"d";"c"^{H_{1}}\ar@{=>}"e";"f"_{H_{2}}}
\end{equation*}
to
\begin{equation*}
\xymatrix@C=50pt{A\ar@/^40pt/[rr]|{\varphi(if)}_{\;}="a"^{\;}="i"\ar@/^80pt/[rr]^0_{\;}="j"\ar[r]_{\varphi(f)}&
B\ar@/_40pt/[rr]|{\varphi(qi)}^{\;}="c"_{\;}="k"\ar@/_80pt/[rr]_0^{\;}="l"
\ar[r]^<(-0.12){\;}="b"|{\varphi(i)}&
C\ar[r]_<(-0.1){\;}="d"|{\varphi(q)}\ar@/^40pt/[rr]|{\varphi(s(f)q)}_{\;}="f"^{\;}="m"\ar@/^80pt/[rr]^0_{\;}="n"&r
A\ar@/^20pt/[r]|{\varphi s(f)}_{\;}="h"\ar[r]_{r
\varphi(f)}^<(-0.14){\;}="e"^{\;}="g"& r
B\ar@{=>}"b";"a"_{\varphi_{i,f}}\ar@{=>}"d";"c"^{\varphi_{q,i}}\ar@{=>}"e";"f"^{\varphi_{s(f),q}}\ar@{=>}"g";"h"^{\xi_f}
\ar@{=>}"i";"j"_{\varphi(H_0)}\ar@{=>}"k";"l"^{\varphi(H_{1})}\ar@{=>}"m";"n"_{\varphi(H_{2})}}
\end{equation*}
The functor $(\varphi,\xi)_*$ sends the homotopy class of a track
triangle morphism in $(\C{B},s)$ given by the diagram
\begin{equation*}
\xymatrix{A\ar[r]^f_<(.6){\;}="a"\ar[d]_{k_0}^<(.6){\;}="b"&B\ar[r]^i_<(.6){\;}="c"\ar[d]|{k_1}^<(.6){\;}="d"&
C\ar[r]^q_<(.6){\;}="e"\ar[d]|{k_2}^<(.6){\;}="f"&sA\ar[d]^{sk_0}\\
\bar{A}\ar[r]_{\bar{f}}&\bar{B}\ar[r]_{\bar{i}}&\bar{C}\ar[r]_{\bar{q}}&s\bar{A}
\ar@{=>}"a";"b"^{K_0}\ar@{=>}"c";"d"^{K_1}\ar@{=>}"e";"f"^{K_2}}
\end{equation*}
to the homotopy class of track triangle morphisms in $(\C{C},r)$
given by
\begin{equation*}
\xymatrix@R=90pt@C=90pt{A\ar@/^20pt/[rd]|{\varphi(k_1f)}^{\;}="g"_{\;}="h"\ar@/_20pt/[rd]|{\varphi(\bar{f}k_0)}^{\;}="i"_{\;}="j"
\ar[r]^{\varphi(f)}_<(1){\;}="a"\ar[d]|{\varphi(k_0)}^<(1){\;}="b"&
B\ar@/^20pt/[rd]|{\varphi(k_2i)}^{\;}="k"_{\;}="l"\ar@/_20pt/[rd]|{\varphi(\bar{i}k_1)}^{\;}="m"_{\;}="n"
\ar[r]^{\varphi(i)}_<(1){\;}="c"\ar[d]|{\varphi(k_1)}^<(1){\;}="d"&
C\ar@/^20pt/[rd]|{\varphi(s(k_0)q)}^{\;}="o"_{\;}="p"\ar@/_20pt/[rd]|{\varphi(\bar{q}k_2)}^{\;}="q"_{\;}="r"
\ar[r]^{\varphi(q)}_<(1){\;}="e"\ar[d]|{\varphi(k_2)}^<(1){\;}="f"&rA\ar[d]|<(.3){\varphi
s(k_0)}^{\;}="t"
\ar@/^30pt/[d]|<(.7){r\varphi (k_0)}_{\;}="s"\\
\bar{A}\ar[r]_{\varphi(\bar{f})}&\bar{B}\ar[r]_{\varphi(\bar{i})}&\bar{C}\ar[r]_{\varphi(\bar{q})}&r\bar{A}
\ar@{=>}"a";"g"_{\varphi_{k_1,f}}\ar@{=>}"h";"i"_{\varphi(K_0)}\ar@{=>}"j";"b"^{\varphi^\vi_{\bar{f},k_0}}
\ar@{=>}"c";"k"_{\varphi_{k_2,i}}\ar@{=>}"l";"m"_{\varphi(K_1)}\ar@{=>}"n";"d"^{\varphi^\vi_{\bar{i},k_1}}
\ar@{=>}"e";"o"_{\varphi_{s(k_0),q}}\ar@{=>}"p";"q"_{\varphi(K_2)}\ar@{=>}"r";"f"^{\varphi^\vi_{\bar{q},k_2}}
\ar@{=>}"s";"t"^{\xi_{k_0}}}
\end{equation*}
The following diagram of functors is commutative
$$\xymatrix{\C{candt}(\C{B},s)\ar[rd]^{\varrho_{(\C{B},s)}}\ar[dd]_{(\varphi,\xi)_*}&\\&\C{cand}(\C{A},t)\\
\C{candt}(\C{C},r)\ar[ru]_{\varrho_{(\C{C},r)}}}$$ therefore the
small candidate triangles in $(\C{A},t)$ in the image of
$\varrho_{(\C{B},s)}$ are also in the image of
$\varrho_{(\C{C},r)}$. The other inclusion follows by symmetry.
\end{proof}

\begin{defn}\label{snabla}
Given a cohomology class $\nabla\in H^3(\C{A},t)$ we say that a
small candidate triangle $$A\st{f}\To B\st{i}\To C\st{q}\To tA$$ in
$(\C{A},t)$ is a \emph{small $\nabla$-triangle}\index{triangle!small
$\nabla$-} if it is in the image of $\varrho_{(\C{B},s)}$ for some
(and hence for any) good translation track category $(\C{B},s)$
representing $\nabla$. This definition is consistent by Proposition
\ref{sontodos}.
\end{defn}

The following definition generalizes the concept of a triangulated
track category as introduced in \cite{ttc} to the case of a
``non-strict'' translation functor.

\begin{defn}\label{wttc}
A \emph{pretriangulated track category}\index{triangulated track
category} is a good translation track category $(\C{B},s)$ together
with a distinguished family of track triangles which is subject to
the following axiom:
\begin{itemize}
\item[(TTr1)] For each map
$f\colon A\r B$ in $\C{B}$ there exists a distinguished track
triangle
$$\xymatrix{A\ar@/^30pt/[rr]^0_{\;}="a"\ar[r]_{f}&B\ar@/_30pt/[rr]_0^{\;}="c"
\ar[r]^<(-0.22){\;}="b"|{i}&
C\ar[r]_<(-0.17){\;}="d"|{q}\ar@/^30pt/[rr]^0_{\;}="f"&sA\ar[r]_{s(f)}^<(-0.23){\;}="e"&
sB\ar@{=>}"b";"a"_{H_0}\ar@{=>}"d";"c"^{H_{1}}\ar@{=>}"e";"f"_{H_{2}}}.$$
\end{itemize}

It is a \emph{triangulated track category} if in addition the
following axiom is satisfied:
\begin{itemize}
\item[(TTr2)] Given distinguished track triangles $$\xymatrix{A\ar@/^30pt/[rr]^0_{\;}="a"\ar[r]_{f}&B\ar@/_30pt/[rr]_0^{\;}="c"
\ar[r]^<(-0.22){\;}="b"|{i}&
C\ar[r]_<(-0.17){\;}="d"|{q}\ar@/^30pt/[rr]^0_{\;}="f"&sA\ar[r]_{s(f)}^<(-0.23){\;}="e"&
sB\ar@{=>}"b";"a"_{H_0}\ar@{=>}"d";"c"^{H_{1}}\ar@{=>}"e";"f"_{H_{2}}}$$
and
$$\xymatrix{\bar{A}\ar@/^30pt/[rr]^0_{\;}="a"\ar[r]_{\bar{f}}&\bar{B}\ar@/_30pt/[rr]_0^{\;}="c"
\ar[r]^<(-0.22){\;}="b"|{\bar{i}}&
\bar{C}\ar[r]_<(-0.17){\;}="d"|{\bar{q}}\ar@/^30pt/[rr]^0_{\;}="f"&s\bar{A}\ar[r]_{s(\bar{f})}^<(-0.23){\;}="e"&
s\bar{B}\ar@{=>}"b";"a"_{\bar{H}_0}\ar@{=>}"d";"c"^{\bar{H}_{1}}\ar@{=>}"e";"f"_{\bar{H}_{2}}},$$
any diagram
$$\xymatrix{A\ar[r]^f_(.7){\;}="a"\ar[d]_{k_0}^(.7){\;}="b"&B\ar[d]^{k_1}\\\bar{A}\ar[r]_{\bar{f}}&\bar{B}\ar@{=>}"a";"b"^{K_0}}$$ in $\C{B}$ extends to a track triangle
morphism
$$\xymatrix{A\ar[r]^f_<(.6){\;}="a"\ar[d]_{k_0}^<(.6){\;}="b"&B\ar[r]^i_<(.6){\;}="c"\ar[d]|{k_1}^<(.6){\;}="d"&
C\ar[r]^q_<(.6){\;}="e"\ar[d]|{k_2}^<(.6){\;}="f"&sA\ar[d]^{sk_0}\\
\bar{A}\ar[r]_{\bar{f}}&\bar{B}\ar[r]_{\bar{i}}&\bar{C}\ar[r]_{\bar{q}}&s\bar{A}
\ar@{=>}"a";"b"^{K_0}\ar@{=>}"c";"d"^{K_1}\ar@{=>}"e";"f"^{K_2}}$$
\end{itemize}
\end{defn}



\begin{thm}\label{iff}
Let $(\C{B},s)$ and $(\C{C},r)$ be good translation track categories
representing the same cohomology class in $H^3(\C{A},t)$. Then
$(\C{B},s)$ admits the structure of a (pre)triangulated track
category if and only if $(\C{C},r)$ does.
\end{thm}



The proof of this theorem will be given in Section \ref{hp}.

\begin{defn}\label{cohpt}
Let $\C{A}$ be an additive category, $t\colon\C{A}\r\C{A}$ an
additive equivalence and $\nabla\in H^3(\C{A},t)$ a translation
cohomology class. We say that $(\C{A},t,\nabla)$ is a
\emph{cohomologically (pre)triangulated
category}\index{cohomologically pretriangulated
category}\index{cohomologically triangulated category} if for one
(and hence for any) good translation track category $(\C{B},s)$
representing $\nabla$ the pair $(\C{B},s)$ admits a
(pre)triangulated track category structure. In this case $\nabla$ is
called a \emph{(pre)triangulated cohomology
class.}\index{triangulated cohomology class} \index{pretriangulated
cohomology class} This definition is consistent by Theorem
\ref{iff}.
\end{defn}

Generalizing the main result in \cite{ttc} we show:

\begin{thm}\label{elt}
If $(\C{A},t,\nabla)$ is a cohomologically (pre)triangulated
category and $\mathcal{E}_\nabla$ is the class of small
$\nabla$-triangles then $(\C{A},t,\mathcal{E}_\nabla)$ is a
(pre)triangulated category.
\end{thm}

\begin{proof}
Axiom (Tr1) follows from (TTr1) and (Tr3) is a consequence of
Proposition \ref{basic1}. Moreover, (Tr0) and (Tr2) will be proved
in \cite{ccctc2} 5.8 by using a purely cohomological
characterization of small $\nabla$-triangles. We prove now that the
octahedral axiom holds provided $(\C{A},t,\nabla)$ is
cohomologically triangulated.
Let  $(\C{B},s)$ be a good translation track category representing
$\nabla$. The pair $(\C{B},s)$ admits a triangulated track category
structure by Definition \ref{cohpt}. Let us fix such a structure.

Consider morphisms $A\st{p(f)}\r B\st{p(g)}\r C$ in $\C{A}$.
Indeed one can always suppose that this diagram comes from a
diagram $A\st{f}\r B\st{g}\r C$ in $\C{B}$. By (TTr1) and (TTr2)
there is a diagram of track triangle morphisms
$$\xymatrix{A\ar[r]^f\ar@{=}[d]&B\ar[r]^{i_f}_(.7){\;\;}="a"\ar[d]|g^(.7){\;\;}="b"&C_f\ar[r]^{q_f}_(.7){\;\;}="c"\ar[d]|{\bar{g}}^(.7){\;\;}="d"&sA\ar@{=}[d]\\
A\ar[r]|{gf}\ar[d]_f&C\ar@{=}[d]^(.7){\;\;}="f"\ar[r]|{i_{gf}}_(.7){\;\;}="e"&C_{gf}\ar[d]|{\bar{f}}^(.7){\;\;}="h"\ar[r]|{q_{gf}}_(.7){\;\;}="g"&sA\ar[d]^{sf}\\
B\ar[r]_{g}&C\ar[r]_{i_{g}}&C_g\ar[r]_{q_g}&sB\ar@{=>}"a";"b"^{F}\ar@{=>}"c";"d"^{G}\ar@{=>}"e";"f"^{F'}\ar@{=>}"g";"h"^{G'}}$$
where the rows are distinguished track triangles, and now by
Proposition \ref{octa1}
$$C_f\st{p(\bar{g})}\To C_{gf}\st{p(\bar{f})}\To C_g\st{(tp(i_f))p(q_g)}\To tC_f$$
is a small $\nabla$-triangle. This means that axiom (Tr4) holds
for the diagram
$$\xymatrix{A\ar[r]^{p(f)}\ar@{=}[d]&B\ar[r]^{p(i_f)}\ar[d]^{g}&C_f\ar[r]^{p(q_f)}&tA\ar@{=}[d]\\
A\ar[d]^{p(f)}\ar[r]^{p(gf)}&C\ar@{=}[d]\ar[r]^{p(i_{gf})}&C_{gf}\ar[r]^{p(q_{gf})}&tA\ar[d]^{tp(f)}\\
B\ar[r]^{p(g)}&C\ar[r]^{p(i_g)}&C_g\ar[r]^{p(q_g)}&tB}$$ but it also
must hold for any other diagram
$$\xymatrix{A\ar[r]^{p(f)}\ar@{=}[d]&B\ar[r]\ar[d]^{g}&\bar{C}_f\ar[r]&tA\ar@{=}[d]\\
A\ar[d]^{p(f)}\ar[r]^{p(gf)}&C\ar@{=}[d]\ar[r]&\bar{C}_{gf}\ar[r]&tA\ar[d]^{tp(f)}\\
B\ar[r]^{p(g)}&C\ar[r]&\bar{C}_g\ar[r]&tB}$$ whose rows are small
$\nabla$-triangles, however this is a consequence of the fact that
the extension of a morphism $h\colon X\r Y$ in $\C{A}$ to a small
$\nabla$-triangle given by (Tr1) is unique up to an isomorphism of
small $\nabla$-triangles given by the identity on $X$ and $Y$ since
the functor $\varsigma$ in (\ref{varsigma}) restricted to small
$\nabla$-triangles fits into a linear extension of categories
defined in Proposition \ref{varsiwle}, and therefore it reflects
isomorphisms, see \cite{ah} IV.4.11.
\end{proof}

Notice that by Proposition \ref{sontodos} the (pre)triangulated structure determined on $\C{A}$ by
this theorem only depends on the translation cohomology class
$\nabla\in H^3(\C{A},t)$.

\section{Obstruction theory for track triangles}

Track triangles, introduced in Section \ref{track}, are the basic
elements for the definition of a (pre)triangulated track category.
In this section develop the obstruction theory for the realization
of small candidate triangle morphisms by track triangle morphisms
through the functor $\varrho$ in (\ref{barro}). This obstruction
theory will play an important role in the purely cohomological
characterization of cohomologically triangulated categories carried
out in \cite{ccctc2}.

In the statement of the following proposition we use the
$\C{cand}(\C{A},t)$-bimodules
$$H^*(sh^3,[-,-]^*\bar{t})$$
defined by the translation cohomology groups of the triangle
category, see \cite{ccctc2} 4. Actually we only need them in
dimensions $1$ and $2$. The cochain complex
$\tilde{G}^*([T,\bar{T}]^*\bar{t})$ computing the cohomology groups
$H^*(sh^3,[T,\bar{T}]^*\bar{t})$ for two small candidate triangles
$$T=\set{A\st{f}\To B\st{i}\To C\st{q}\To tA},$$
$$\bar{T}=\set{\bar{A}\st{\bar{f}}\To \bar{B}\st{\bar{i}}\To \bar{C}\st{\bar{q}}\To t\bar{A}},$$
is given in low dimensions by
$$\xymatrix{\tilde{G}^0([T,\bar{T}]^*\bar{t})\ar@{}[r]|<(.1)=\ar[d]^\delta&\hom_{\C{A}}(tA,\bar{A})\oplus\hom_{\C{A}}(tB,\bar{B})\oplus\hom_{\C{A}}(tC,\bar{C})
\ar[d]^{\begin{tiny}\left(
          \begin{array}{ccc}
            \bar{f}_* & -f^* & 0 \\
            0 & \bar{i}_* & -i^* \\
            -q^*\bar{t}_{(A,\bar{A})} & 0 & \bar{q}_* \\
          \end{array}
        \right)
        \end{tiny}
}\\
\tilde{G}^1([T,\bar{T}]^*\bar{t})\ar@{}[r]|<(.1)=\ar[d]^\delta&\hom_{\C{A}}(tA,\bar{B})\oplus\hom_{\C{A}}(tB,\bar{C})\oplus\hom_{\C{A}}(tC,t\bar{A})
\ar[d]^{\begin{tiny}\left(
          \begin{array}{ccc}
            \bar{i}_* & f^* & 0 \\
            0 & \bar{q}_* & i^* \\
            q^*\bar{t}_{(A,\bar{B})} & 0 & (t\bar{f})_* \\
          \end{array}
        \right)
        \end{tiny}
}\\
\tilde{G}^2([T,\bar{T}]^*\bar{t})\ar@{}[r]|<(.1)=\ar[d]^\delta&\hom_{\C{A}}(tA,\bar{C})\oplus\hom_{\C{A}}(tB,t\bar{A})\oplus\hom_{\C{A}}(tC,t\bar{B})
\ar[d]^{\begin{tiny}\left(
          \begin{array}{ccc}
            \bar{q}_* & -f^* & 0 \\
            0 & (t\bar{f})_* & -i^* \\
            -q^*\bar{t}_{(A,\bar{C})} & 0 & (t\bar{i})_* \\
          \end{array}
        \right)
        \end{tiny}
}\\
\tilde{G}^3([T,\bar{T}]^*\bar{t})\ar@{}[r]|<(.1)=&\hom_{\C{A}}(tA,t\bar{A})\oplus\hom_{\C{A}}(tB,t\bar{B})\oplus\hom_{\C{A}}(tC,t\bar{C})}$$
See \cite{ccctc2} 4.14. Here $(-)^*$ and $(-)_*$ denote morphisms
induced by the $\C{A}$-bimodule $\hom_\C{A}(t,-)$. The functoriality
on $T$ and $\bar{T}$ is clear.

\begin{prop}\label{exa}
The functor (\ref{barro}) fits into an exact sequence
$$H^1(sh^3,[-,-]^*\bar{t})\st{+}\To \C{candt}(\C{B},s)\st{\varrho}\To\C{cand}(\C{A},t)\st{\vartheta}\To H^2(sh^3,[-,-]^*\bar{t})$$
in the sense of \cite{ah} IV.4.10. Moreover, the action $+$ is
effective.
\end{prop}

An exact sequence for a functor encodes the usual properties of
obstruction theory in homotopy theory. The abstract definition will
not be recalled here, we refer the reader to \cite{ah} IV.4.10.
However the verification of the axioms for the case of Proposition
\ref{exa} constitutes the statement of the following five lemmas,
which prove Proposition \ref{exa}.

In the proof of the following results we will use a fixed global
section $(\ell,\mu,\nu,\eta)$ in the sense of \cite{ptc} 14. We are
actually interested in $\ell$, $\mu$ and part of $\nu$ only. This is
given by a choice of a map $\ell f\colon A\r B$ in $\C{B}$ for any
$f\colon A\r B$ in $\C{A}$ with $p(\ell f)=f$ and $\ell(1_X)=1_X$,
tracks $\mu_{f,g}\colon(\ell f)(\ell g)\rr\ell(fg)$ with
$\mu_{f,1}=0^\vc_{\ell f}=\mu_{1,f}$, and tracks $\nu_{f,s}\colon
s(\ell f)\rr\ell(tf)$ with $\nu_{1_X,s}=s_X$. Since $\C{B}$ has a
strict zero object $*$ which is preserved by $s$ and $s$ is
normalized at zero maps we can also impose here $\ell(0)=0$,
$\mu_{f,0}=0^\vc_{0}=\mu_{0,f}$ and $\nu_{0,s}=0^\vc_0$.  
We will also need a collection $\gamma$ of tracks
$\gamma_f\colon\ell p(f)\rr f$ indexed by all maps $f$ in $\C{B}$
such that $\gamma_{\ell g}=0^\vc_{\ell g}$ for any morphism $g$ in
$\C{A}$.

\begin{lem}\label{A}
Let $T$ and $\bar{T}$ be track triangles given by
$$\xymatrix{A\ar@/^30pt/[rr]^0_{\;}="a"\ar[r]_{f}&B\ar@/_30pt/[rr]_0^{\;}="c"
\ar[r]^<(-0.22){\;}="b"|{i}&
C\ar[r]_<(-0.17){\;}="d"|{q}\ar@/^30pt/[rr]^0_{\;}="f"&sA\ar[r]_{s(f)}^<(-0.23){\;}="e"&
sB\ar@{=>}"b";"a"_{H_0}\ar@{=>}"d";"c"^{H_{1}}\ar@{=>}"e";"f"_{H_{2}}}$$
and
$$\xymatrix{\bar{A}\ar@/^30pt/[rr]^0_{\;}="a"\ar[r]_{\bar{f}}&\bar{B}\ar@/_30pt/[rr]_0^{\;}="c"
\ar[r]^<(-0.22){\;}="b"|{\bar{i}}&
\bar{C}\ar[r]_<(-0.17){\;}="d"|{\bar{q}}\ar@/^30pt/[rr]^0_{\;}="f"&s\bar{A}\ar[r]_{s(\bar{f})}^<(-0.23){\;}="e"&
s\bar{B}\ar@{=>}"b";"a"_{\bar{H}_0}\ar@{=>}"d";"c"^{\bar{H}_{1}}\ar@{=>}"e";"f"_{\bar{H}_{2}}}$$
For any morphism $k\colon\varrho T\r\varrho\bar{T}$ of small
candidate triangles as in the following diagram
\begin{equation*}
\xymatrix{A\ar[r]^{p(f)}\ar[d]^{k_0}&B\ar[r]^{p(i)}\ar[d]^{k_1}&C\ar[r]^{p(q)}\ar[d]^{k_2}&tA\ar[d]^{tk_0}\\
\bar{A}\ar[r]^{p(\bar{f})}&\bar{B}\ar[r]^{p(\bar{i})}&\bar{C}\ar[r]^{p(\bar{q})}&t\bar{A}}
\end{equation*}
there is a well-defined cohomology class
$$\vartheta_{T,\bar{T}}(k)\in H^2(sh^3,[\varrho T,\varrho\bar{T}]^*\bar{t})$$
which vanishes if and only if $k$ is the image by $\varrho$ of a
track triangle morphism $T\r\bar{T}$.
\end{lem}

\begin{proof}
The element $\vartheta_{T,\bar{T}}(k)$ measures the failure in
$$\xymatrix@C=100pt@R=100pt{
A\ar[r]^f_{\;}="a"\ar@/_20pt/[r]|{\ell
p(f)}^{\;}="b"\ar[rd]|{\ell(k_1 p(f))}^{\;}="d"_{\;}="e"\ar[d]_{\ell
k_0}^<(.7){\;}="f"& B\ar[r]^i_{\;}="g"\ar@/_20pt/[r]|{\ell
p(i)}^{\;}="h"\ar[rd]|{\ell(k_2 p(i))}^{\;}="j"_{\;}="k"\ar[d]|{\ell
k_1}_<(.2){\;}="c"^<(.7){\;}="l"&
C\ar[r]^q_{\;}="m"\ar@/_20pt/[r]|{\ell q}^{\;}="n"\ar[rd]|{\ell
(p(\bar{q})k_2)}^{\;}="p"_{\;}="q"\ar[d]|{\ell
k_2}_<(.2){\;}="i"^<(.7){\;}="r"&
sA\ar[d]^{s\ell k_0}_{\;}="s"\ar@/_25pt/[d]|{\ell t k_0}^{\;}="t"_<(.2){\;}="o"\\
\bar{A}\ar[r]_{\bar{f}}^{\;}="v"\ar@/^20pt/[r]|{\ell
p(\bar{f})}_{\;}="u"&
\bar{B}\ar[r]_{\bar{i}}^{\;}="x"\ar@/^20pt/[r]|{\ell
p(\bar{i})}_{\;}="w"&
\bar{C}\ar[r]_{\bar{q}}^{\;}="z"\ar@/^20pt/[r]|{\ell
p(\bar{q})}_{\;}="y"& s\bar{A} \ar@{=>}"a";"b"^{\gamma_f^\vi}
\ar@{=>}"c";"d"_{\mu_{k_1,p(f)}}
\ar@{=>}"e";"f"_{\mu^\vi_{p(\bar{f}),k_0}\;}
\ar@{=>}"g";"h"^{\gamma_i^\vi} \ar@{=>}"i";"j"_{\mu_{k_2,p(i)}}
\ar@{=>}"k";"l"_{\mu^\vi_{p(\bar{i}),k_1}\;}
\ar@{=>}"m";"n"^{\gamma_q^\vi} \ar@{=>}"o";"p"_{\mu_{tk_0,p(q)}}
\ar@{=>}"q";"r"_{\mu^\vi_{p(\bar{q}),k_2}\;}
\ar@{=>}"s";"t"_{\nu_{k_0,s}} \ar@{=>}"u";"v"^{\gamma_{\bar{f}}}
\ar@{=>}"w";"x"^{\gamma_{\bar{i}}}
\ar@{=>}"y";"z"^{\gamma_{\bar{q}}} }$$ to define a morphism of track
triangles, i. e. it is represented by
 the cocycle
$(b_{-2},b_{-1},b_0)\in\tilde{G}^2([\varrho
T,\varrho\bar{T}]^*\bar{t})$ with
\begin{eqnarray*}
b_{-2}&=&\sigma_0^{-1}((\bar{H}_0(\ell k_0))\vc(\bar{i}M_0)\vc(M_1f)\vc((\ell k_2) H_0^\vi)),\\
b_{-1}&=&\sigma_0^{-1}((\bar{H}_1(\ell k_1))\vc(\bar{q}M_1)\vc(M_2i)\vc(s(\ell k_0) H_1^\vi)),\\
b_0&=&\sigma_0^{-1}((\bar{H}_2(\ell
k_2))\vc(s(\bar{f})M_2)\vc((s^\vi_{\bar{f},\ell k_0}\vc s(M_0)\vc
s_{\ell k_1,f})q)\vc(s(\ell k_1) H_2^\vi)),
\end{eqnarray*}
where
\begin{eqnarray*}
M_0&=&(\gamma_{\bar{f}}(\ell k_0))\vc\mu^\vi_{p(\bar{f}),k_0}\vc\mu_{k_1,p(f)}\vc((\ell k_1)\gamma_f^\vi),\\
M_1&=&(\gamma_{\bar{i}}(\ell k_1))\vc\mu^\vi_{p(\bar{i}),k_1}\vc\mu_{k_2,p(i)}\vc((\ell k_2)\gamma_i^\vi),\\
M_2&=&(\gamma_{\bar{q}}(\ell
k_2))\vc\mu^\vi_{p(\bar{q}),k_2}\vc\mu_{tk_0,p(q)}\vc(\nu_{k_0,s}\gamma_q^\vi).
\end{eqnarray*}

If the cocycle $(b_{-2},b_{-1},b_0)$ is the coboundary of
$(c_{-2},c_{-1},c_0)\in\tilde{G}^1([\varrho
T,\varrho\bar{T}]^*\bar{t})$ then $k$ is the image under $\varrho$
of the homotopy class of morphisms of track triangles $T\r \bar{T}$
represented by the diagram
$$\xymatrix{A\ar[r]^f_<(.6){\;}="a"\ar[d]_{\ell k_0}^<(.6){\;}="b"&B\ar[r]^i_<(.6){\;}="c"\ar[d]|{\ell k_1}^<(.6){\;}="d"&
C\ar[r]^q_<(.6){\;}="e"\ar[d]|{\ell k_2}^<(.6){\;}="f"&sA\ar[d]^{s\ell k_0}\\
\bar{A}\ar[r]_{\bar{f}}&\bar{B}\ar[r]_{\bar{i}}&\bar{C}\ar[r]_{\bar{q}}&s\bar{A}
\ar@{=>}"a";"b"^{K_0}\ar@{=>}"c";"d"^{K_1}\ar@{=>}"e";"f"^{K_2}}$$
with
\begin{eqnarray*}
K_0&=&\sigma_{\bar{f}(\ell k_0)}(-c_{-2})\vc M_0,\\
K_1&=&\sigma_{\bar{i}(\ell k_1)}(-c_{-1})\vc M_1,\\
K_2&=&\sigma_{\bar{q}(\ell k_2)}(-c_0)\vc M_2.
\end{eqnarray*}

On the other hand if $k$ is the image by $\varrho$ of the track
triangle morphism
$$\xymatrix{A\ar[r]^f_<(.6){\;}="a"\ar[d]_{\bar{k}_0}^<(.6){\;}="b"&B\ar[r]^i_<(.6){\;}="c"\ar[d]|{\bar{k}_1}^<(.6){\;}="d"&
C\ar[r]^q_<(.6){\;}="e"\ar[d]|{\bar{k}_2}^<(.6){\;}="f"&sA\ar[d]^{s\bar{k}_0}\\
\bar{A}\ar[r]_{\bar{f}}&\bar{B}\ar[r]_{\bar{i}}&\bar{C}\ar[r]_{\bar{q}}&s\bar{A}
\ar@{=>}"a";"b"^{\bar{K}_0}\ar@{=>}"c";"d"^{\bar{K}_1}\ar@{=>}"e";"f"^{\bar{K}_2}}$$
then the cochain $(c_{-2},c_{-1},c_0)\in\tilde{G}^1([\varrho
T,\varrho \bar{T}]^*\bar{t})$ defined by
\begin{eqnarray*}
-c_{-2}&=&\sigma_{\bar{f}\bar{k}_0}^{-1}(\bar{K_0}\vc(\gamma_{\bar{k}_1}f)\vc
M_0^\vi\vc(\bar{f}\gamma^\vi_{\bar{k}_0})),\\
-c_{-1}&=&\sigma_{\bar{i}\bar{k}_1}^{-1}(\bar{K_1}\vc
(\gamma_{\bar{k}_2}i)\vc M_1^\vi\vc(\bar{i}\gamma_{\bar{k}_1}^\vi)),\\
-c_0&=&\sigma_{\bar{q}\bar{k}_2}^{-1}(\bar{K_2}\vc
(s(\gamma_{\bar{k}_0})q)\vc
M_2^\vi\vc(\bar{q}\gamma^\vi_{\bar{k}_2})),
\end{eqnarray*}
satisfies $\delta(c_{-2},c_{-1},c_0)=(b_{-2},b_{-1},b_0)$ so
$\vartheta_{T,\bar{T}}(k)=0$.

The reader can check that the cohomology class
$\vartheta_{T,\bar{T}}(k)$ does not depend on the choices made for
its definition.
\end{proof}

\begin{lem}\label{B}
In the conditions of the statement of Lemma \ref{A} if another track
triangle
$$\xymatrix{\bar{\bar{A}}\ar@/^30pt/[rr]^0_{\;}="a"\ar[r]_{\bar{\bar{f}}}&\bar{\bar{B}}\ar@/_30pt/[rr]_0^{\;}="c"
\ar[r]^<(-0.22){\;}="b"|{\bar{\bar{i}}}&
\bar{\bar{C}}\ar[r]_<(-0.17){\;}="d"|{\bar{\bar{q}}}\ar@/^30pt/[rr]^0_{\;}="f"&s\bar{\bar{A}}\ar[r]_{s(\bar{\bar{f}})}^<(-0.23){\;}="e"&
s\bar{\bar{B}}\ar@{=>}"b";"a"_{\bar{\bar{H}}_0}\ar@{=>}"d";"c"^{\bar{\bar{H}}_{1}}\ar@{=>}"e";"f"_{\bar{\bar{H}}_{2}}}$$
 and another morphism $\bar{k}\colon \varrho \bar{T}\r\varrho\bar{\bar{T}}$ of small candidate
 triangles
\begin{equation*}
\xymatrix{\bar{A}\ar[r]^{p(\bar{f})}\ar[d]^{\bar{k}_0}&\bar{B}\ar[r]^{p(\bar{i})}\ar[d]^{\bar{k}_1}&\bar{C}\ar[r]^{p(\bar{q})}\ar[d]^{\bar{k}_2}&t\bar{A}\ar[d]^{t\bar{k}_0}\\
\bar{\bar{A}}\ar[r]^{p(\bar{\bar{f}})}&\bar{\bar{B}}\ar[r]^{p(\bar{\bar{i}})}&\bar{\bar{C}}\ar[r]^{p(\bar{\bar{q}})}&t\bar{\bar{A}}}
\end{equation*}
are given, then the following equality holds
$$\vartheta_{T,\bar{\bar{T}}}(\bar{k}k)=\bar{k}_*\vartheta_{T,\bar{T}}(k)+k^*\vartheta_{\bar{T},\bar{\bar{T}}}(\bar{k}).$$
\end{lem}

\begin{proof}
Following the definition of $\vartheta$ in the proof of Lemma
\ref{A} the reader can check that
\begin{equation*}\tag{a}
\bar{k}_*\vartheta_{T,\bar{T}}(k)+k^*\vartheta_{\bar{T},\bar{\bar{T}}}(\bar{k})
\end{equation*}
measures the failure in
$$\xymatrix@C=100pt@R=100pt{
A\ar[r]^f_{\;}="a"\ar@/_20pt/[r]|{\ell
p(f)}^{\;}="b"\ar[rd]|{\ell(k_1 p(f))}^{\;}="d"_{\;}="e"\ar[d]_{\ell
k_0}^<(.7){\;}="f"& B\ar[r]^i_{\;}="g"\ar@/_20pt/[r]|{\ell
p(i)}^{\;}="h"\ar[rd]|{\ell(k_2 p(i))}^{\;}="j"_{\;}="k"\ar[d]|{\ell
k_1}_<(.2){\;}="c"^<(.7){\;}="l"&
C\ar[r]^q_{\;}="m"\ar@/_20pt/[r]|{\ell q}^{\;}="n"\ar[rd]|{\ell
(p(\bar{q})k_2)}^{\;}="p"_{\;}="q"\ar[d]|{\ell
k_2}_<(.2){\;}="i"^<(.7){\;}="r"&
sA\ar[d]^{s\ell k_0}_{\;}="s"\ar@/_25pt/[d]|{\ell t k_0}^{\;}="t"_<(.2){\;}="o"\\
\bar{A}\ar[r]|{\ell
p(\bar{f})}_<(.7){\;}="aa"\ar[rd]|{\ell(\bar{k}_1p(\bar{f}))}^{\;}="bb"_{\;}="cc"\ar[d]_{\ell\bar{k}_0}^<(.7){\;}="dd"&
\bar{B}\ar[r]|{\ell
p(\bar{i})}_<(.7){\;}="ee"\ar[rd]|{\ell(\bar{k}_2p(\bar{i}))}^{\;}="ff"_{\;}="gg"\ar[d]|{\ell\bar{k}_1}^<(.7){\;}="hh"&
\bar{C}\ar[r]|{\ell
p(\bar{q})}_<(.7){\;}="ii"\ar[rd]|{\ell(\bar{\bar{q}}\bar{k}_2)}^{\;}="jj"_{\;}="kk"\ar[d]|{\ell\bar{k}_2}^<(.7){\;}="ll"&
s\bar{A}\ar[d]^{s\ell\bar{k}_0}_{\;}="mm"\ar@/_25pt/[d]|{\ell
t\bar{k_0}}^{\;}="nn"
\\
\bar{\bar{A}}\ar[r]_{\bar{\bar{f}}}^{\;}="v"\ar@/^20pt/[r]|{\ell
p(\bar{\bar{f}})}_{\;}="u"&
\bar{\bar{B}}\ar[r]_{\bar{\bar{i}}}^{\;}="x"\ar@/^20pt/[r]|{\ell
p(\bar{\bar{i}})}_{\;}="w"&
\bar{\bar{C}}\ar[r]_{\bar{\bar{q}}}^{\;}="z"\ar@/^20pt/[r]|{\ell
p(\bar{\bar{q}})}_{\;}="y"& s\bar{\bar{A}}
\ar@{=>}"a";"b"^{\gamma_f^\vi} \ar@{=>}"c";"d"_{\mu_{k_1,p(f)}\;}
\ar@{=>}"e";"f"_{\mu^\vi_{p(\bar{f}),k_0}\;}
\ar@{=>}"g";"h"^{\gamma_i^\vi} \ar@{=>}"i";"j"_{\mu_{k_2,p(i)}\;}
\ar@{=>}"k";"l"_{\mu^\vi_{p(\bar{i}),k_1}\;}
\ar@{=>}"m";"n"^{\gamma_q^\vi} \ar@{=>}"o";"p"_{\mu_{tk_0,p(q)}\;}
\ar@{=>}"q";"r"_{\mu^\vi_{p(\bar{q}),k_2}\;}
\ar@{=>}"s";"t"_{\nu_{k_0,s}}
\ar@{=>}"u";"v"^{\gamma_{\bar{\bar{f}}}}
\ar@{=>}"w";"x"^{\gamma_{\bar{\bar{i}}}}
\ar@{=>}"y";"z"^{\gamma_{\bar{\bar{q}}}}
\ar@{=>}"aa";"bb"_{\mu_{\bar{k}_1,p(\bar{f})}\;}
\ar@{=>}"cc";"dd"_{\mu^\vi_{p(\bar{\bar{f}}),\bar{k}_0}\;}
\ar@{=>}"ee";"ff"_{\mu_{\bar{k}_2,p(\bar{i})}\;}
\ar@{=>}"gg";"hh"_{\mu^\vi_{p(\bar{\bar{i}}),\bar{k}_1}\;}
\ar@{=>}"ii";"jj"_{\mu_{t\bar{k}_0,p(\bar{q})}\;}
\ar@{=>}"kk";"ll"_{\mu^\vi_{p(\bar{\bar{q}}),\bar{k}_2}\;}
\ar@{=>}"mm";"nn"_{\nu_{\bar{k}_0,s}} }$$ to define a track triangle
morphism $T\r\bar{\bar{T}}$, i. e. (a) is represented by
$(\bar{b}_{-2},\bar{b}_{-1},\bar{b}_0)\in\tilde{G}^2([\varrho
T,\varrho\bar{\bar{T}}]^*\bar{t})$ defined by
\begin{eqnarray*}
\bar{b}_{-2}&=&\sigma_0^{-1}((\bar{\bar{H}}_0(\ell \bar{k}_0)(\ell k_0))\vc(\bar{\bar{i}}\bar{M}_0)\vc(\bar{M}_1 f)\vc((\ell\bar{k}_2)(\ell k_2)H^\vi_0)),\\
\bar{b}_{-1}&=&\sigma_0^{-1}((\bar{\bar{H}}_1(\ell \bar{k}_1)(\ell k_1))\vc(\bar{\bar{q}}\bar{M}_1)\vc(\bar{M}_2 i)\vc(s((\ell\bar{k}_0)(\ell k_0))H^\vi_1)),\\
\bar{b}_0&=&\sigma_0^{-1}((\bar{\bar{H}}_2(\ell \bar{k}_2)(\ell
k_2))\vc(s(\bar{\bar{f}})\bar{M}_2)
\vc((s_{\bar{\bar{f}},(\ell\bar{k}_0)(\ell\bar{k}_0)}^\vi\vc
s(\bar{M}_0)\vc s_{(\ell\bar{k}_1)(\ell k_1),f})
q)\\&&\vc(s((\ell\bar{k}_1)(\ell k_1))H^\vi_2)),
\end{eqnarray*}
where
\begin{eqnarray*}
\bar{M}_0&=&(((\gamma_{\bar{\bar{f}}}(\ell
\bar{k}_0))\vc\mu^\vi_{p(\bar{\bar{f}}),\bar{k}_0}\vc\mu_{\bar{k}_1,p(\bar{f})})(\ell
k_0))
\vc((\ell\bar{k}_1)(\mu^\vi_{p(\bar{f}),k_0}\vc\mu_{k_1,p(f)}\vc((\ell k_1)\gamma^\vi_f))),\\
\bar{M}_1&=&(((\gamma_{\bar{\bar{i}}}(\ell
\bar{k}_1))\vc\mu^\vi_{p(\bar{\bar{i}}),\bar{k}_1}\vc\mu_{\bar{k}_2,p(\bar{i})})(\ell
k_1))
\vc((\ell\bar{k}_2)(\mu^\vi_{p(\bar{i}),k_1}\vc\mu_{k_2,p(i)}\vc((\ell k_2)\gamma^\vi_i))),\\
\bar{M}_2&=&(((\gamma_{\bar{\bar{q}}}(\ell
\bar{k}_2))\vc\mu^\vi_{p(\bar{\bar{q}}),\bar{k}_2}\vc\mu_{t\bar{k}_0,p(\bar{q})}\vc(\nu_{\bar{k}_0,s}(\ell
p(\bar{q}))))(\ell k_2))
\\&&
\vc(s(\ell\bar{k}_0)(\mu^\vi_{p(\bar{q}),k_2}\vc\mu_{tk_0,p(q)}\vc(\nu_{
k_0,s}\gamma^\vi_q))).
\end{eqnarray*}

On the other hand $\vartheta_{T,\bar{\bar{T}}}(\bar{k}k)$ measures
the failure in
$$\xymatrix@C=100pt@R=100pt{
A\ar[r]^f_{\;}="a"\ar@/_20pt/[r]|{\ell
p(f)}^{\;}="b"\ar[rd]|{\ell(\bar{k}_1k_1
p(f))}^{\;}="d"_{\;}="e"\ar[d]|{\ell (\bar{k}_0k_0)}^<(.7){\;}="f"&
B\ar[r]^i_{\;}="g"\ar@/_20pt/[r]|{\ell
p(i)}^{\;}="h"\ar[rd]|{\ell(\bar{k}_2k_2
p(i))}^{\;}="j"_{\;}="k"\ar[d]|{\ell
(\bar{k}_1k_1)}_<(.2){\;}="c"^<(.7){\;}="l"&
C\ar[r]^q_{\;}="m"\ar@/_20pt/[r]|{\ell
q}^{\;}="n"\ar@/_10pt/[rd]|{\ell
(p(\bar{\bar{q}})\bar{k}_2k_2)}^{\;}="p"_{\;}="q"\ar[d]|{\ell
(\bar{k}_2k_2)}_<(.2){\;}="i"^<(.7){\;}="r"&
sA\ar[d]|<(.7){s\ell (\bar{k}_0k_0)}_{\;}="s"\ar@/_38pt/[d]|<(.65){\ell t (\bar{k}_0k_0)}^{\;}="t"_<(.2){\;}="o"\\
\bar{\bar{A}}\ar[r]_{\bar{\bar{f}}}^{\;}="v"\ar@/^20pt/[r]|{\ell
p(\bar{\bar{f}})}_{\;}="u"&
\bar{\bar{B}}\ar[r]_{\bar{\bar{i}}}^{\;}="x"\ar@/^20pt/[r]|{\ell
p(\bar{\bar{i}})}_{\;}="w"&
\bar{\bar{C}}\ar[r]_{\bar{\bar{q}}}^{\;}="z"\ar@/^20pt/[r]|{\ell
p(\bar{\bar{q}})}_{\;}="y"& s\bar{\bar{A}}
\ar@{=>}"a";"b"^{\gamma_f^\vi}
\ar@{=>}"c";"d"_{\mu_{\bar{k}_1k_1,p(f)}\;}
\ar@{=>}"e";"f"^{\;\mu^\vi_{p(\bar{\bar{f}}),\bar{k}_0k_0}}
\ar@{=>}"g";"h"^{\gamma_i^\vi}
\ar@{=>}"i";"j"_{\mu_{\bar{k}_2k_2,p(i)}\;}
\ar@{=>}"k";"l"^{\;\mu^\vi_{p(\bar{\bar{i}}),\bar{k}_1k_1}}
\ar@{=>}"m";"n"^{\gamma_q^\vi}
\ar@{=>}"o";"p"_<(.2){\mu_{t(\bar{k}_0k_0),p(q)}\!}
\ar@{=>}"q";"r"^{\;\mu^\vi_{p(\bar{\bar{q}}),\bar{k}_2k_2}}
\ar@{=>}"s";"t"_{\nu_{\bar{k}_0k_0,s}}
\ar@{=>}"u";"v"^{\gamma_{\bar{\bar{f}}}}
\ar@{=>}"w";"x"^{\gamma_{\bar{\bar{i}}}}
\ar@{=>}"y";"z"^{\gamma_{\bar{\bar{q}}}} }$$ to define a track
triangle morphism $T\r\bar{\bar{T}}$, i. e.
$\vartheta_{T,\bar{\bar{T}}}(\bar{k}k)$ is represented by
$(b_{-2},b_{-1},b_0)\in\tilde{G}^2([\varrho
T,\varrho\bar{\bar{T}}]^*\bar{t})$ with
\begin{eqnarray*}
b_{-2}&=&\sigma_0^{-1}((\bar{\bar{H}}_0(\ell (\bar{k}_0k_0)))\vc(\bar{\bar{i}}M_0)\vc(M_1f)\vc((\ell(\bar{k}_2 k_2)) H_0^\vi)),\\
b_{-1}&=&\sigma_0^{-1}((\bar{\bar{H}}_1(\ell(\bar{k}_1 k_1)))\vc(\bar{\bar{q}}M_1)\vc(M_2i)\vc(s(\ell(\bar{k}_0 k_0)) H_1^\vi)),\\
b_0&=&\sigma_0^{-1}((\bar{\bar{H}}_2(\ell(\bar{k}_2
k_2)))\vc(s(\bar{\bar{f}})M_2)\vc((s^\vi_{\bar{\bar{f}},\ell
(\bar{k}_0k_0)}\vc s(M_0)\vc s_{\ell(\bar{k}_1 k_1),f})q)\\&&
\vc(s(\ell(\bar{k}_1 k_1)) H_0^\vi))),
\end{eqnarray*}
where
\begin{eqnarray*}
M_0&=&(\gamma_{\bar{\bar{f}}}(\ell (\bar{k}_0k_0)))\vc\mu^\vi_{p(\bar{\bar{f}}),\bar{k}_0k_0}\vc\mu_{\bar{k}_1k_1,p(f)}\vc((\ell (\bar{k}_1k_1))\gamma_f^\vi),\\
M_1&=&(\gamma_{\bar{\bar{i}}}(\ell (\bar{k}_1k_1)))\vc\mu^\vi_{p(\bar{\bar{i}}),\bar{k}_1k_1}\vc\mu_{\bar{k}_2k_2,p(i)}\vc((\ell (\bar{k}_2k_2))\gamma_i^\vi),\\
M_2&=&(\gamma_{\bar{\bar{q}}}(\ell
(\bar{k}_2k_2)))\vc\mu^\vi_{p(\bar{\bar{q}}),\bar{k}_2k_2}\vc\mu_{t(\bar{k}_0k_0),p(q)}\vc(\nu_{\bar{k}_0k_0,s}\gamma_q^\vi).
\end{eqnarray*}

Now a laborious but straightforward computation shows that the
cochain $(c_{-2},c_{-1},c_0)\in\tilde{G}^1([\varrho
T,\varrho\bar{\bar{T}}]^*\bar{t})$ given by
\begin{eqnarray*}
c_{-2}&=&\sigma^{-1}_{\bar{\bar{f}}\ell(\bar{k}_0k_0)}(M_0\vc(\mu_{\bar{k}_1,k_1}f)\vc\bar{M}_0^\vi\vc(\bar{\bar{f}}\mu_{\bar{k}_0,k_0}^\vi))\\
c_{-1}&=&\sigma^{-1}_{\bar{\bar{i}}\ell(\bar{k}_1k_1)}(M_1\vc(\mu_{\bar{k}_2,k_2}i)\vc\bar{M}_1^\vi\vc(\bar{\bar{i}}\mu_{\bar{k}_1,k_1}^\vi))\\
c_0&=&\sigma^{-1}_{\bar{\bar{q}}\ell(\bar{k}_2k_2)}(M_2\vc((s(\mu_{\bar{k}_0,k_0})\vc
s_{\ell\bar{k}_0,\ell k_0})q)
\vc\bar{M}_2^\vi\vc(\bar{\bar{f}}\mu_{\bar{k}_2,k_2}^\vi))
\end{eqnarray*}
satisfies
$$\delta(c_{-2},c_{-1},c_0)=(b_{-2},b_{-1},b_0)-(\bar{b}_{-2},\bar{b}_{-1},\bar{b}_0),$$
therefore the lemma follows.
\end{proof}

\begin{lem}\label{C}
For any candidate triangle $T$ given by
$$\xymatrix{A\ar@/^30pt/[rr]^0_{\;}="a"\ar[r]_{f}&B\ar@/_30pt/[rr]_0^{\;}="c"
\ar[r]^<(-0.22){\;}="b"|{i}&
C\ar[r]_<(-0.17){\;}="d"|{q}\ar@/^30pt/[rr]^0_{\;}="f"&sA\ar[r]_{s(f)}^<(-0.23){\;}="e"&
sB\ar@{=>}"b";"a"_{H_0}\ar@{=>}"d";"c"^{H_{1}}\ar@{=>}"e";"f"_{H_{2}}}$$
and any $c\in H^2(sh^3,[\varrho T,\varrho T]^*\bar{t})$ there is a candidate
triangle $\bar{T}$ as follows
$$\xymatrix{A\ar@/^30pt/[rr]^0_{\;}="a"\ar[r]_{f}&B\ar@/_30pt/[rr]_0^{\;}="c"
\ar[r]^<(-0.22){\;}="b"|{i}&
C\ar[r]_<(-0.17){\;}="d"|{q}\ar@/^30pt/[rr]^0_{\;}="f"&sA\ar[r]_{s(f)}^<(-0.23){\;}="e"&
sB\ar@{=>}"b";"a"_{\bar{H}_0}\ar@{=>}"d";"c"^{\bar{H}_{1}}\ar@{=>}"e";"f"_{\bar{H}_{2}}}$$
such that $\theta_{T,\bar{T}}(1_{\varrho T})=c$.
\end{lem}

\begin{proof}
If $c$
 represented by the
cocycle $(c_{-2},c_{-1},c_0)\in\tilde{G}^2([\varrho T,\varrho T]^*\bar{t})$
then it is enough to define
\begin{eqnarray*}
\bar{H}_0&=&\sigma_0(c_{-2})\vc H_0,\\
\bar{H}_1&=&\sigma_0(c_{-1})\vc H_1,\\
\bar{H}_2&=&\sigma_0(c_0)\vc H_2.
\end{eqnarray*}
\end{proof}

\begin{lem}\label{D}
Given track triangles $T$ and $\bar{T}$ defined by
$$\xymatrix{A\ar@/^30pt/[rr]^0_{\;}="a"\ar[r]_{f}&B\ar@/_30pt/[rr]_0^{\;}="c"
\ar[r]^<(-0.22){\;}="b"|{i}&
C\ar[r]_<(-0.17){\;}="d"|{q}\ar@/^30pt/[rr]^0_{\;}="f"&sA\ar[r]_{s(f)}^<(-0.23){\;}="e"&
sB\ar@{=>}"b";"a"_{H_0}\ar@{=>}"d";"c"^{H_{1}}\ar@{=>}"e";"f"_{H_{2}}}$$
and
$$\xymatrix{\bar{A}\ar@/^30pt/[rr]^0_{\;}="a"\ar[r]_{\bar{f}}&\bar{B}\ar@/_30pt/[rr]_0^{\;}="c"
\ar[r]^<(-0.22){\;}="b"|{\bar{i}}&
\bar{C}\ar[r]_<(-0.17){\;}="d"|{\bar{q}}\ar@/^30pt/[rr]^0_{\;}="f"&s\bar{A}\ar[r]_{s(\bar{f})}^<(-0.23){\;}="e"&
s\bar{B}\ar@{=>}"b";"a"_{\bar{H}_0}\ar@{=>}"d";"c"^{\bar{H}_{1}}\ar@{=>}"e";"f"_{\bar{H}_{2}}}$$
the abelian group $H^1(sh^3,[\varrho T,\varrho\bar{T}]^*\bar{t})$
acts effectively on the set of homotopy classes of track triangle
morphisms $T\r\bar{T}$ in such a way that given $c\in
H^1(sh^3,[\varrho T,\varrho\bar{T}]^*\bar{t})$ and a homotopy class
$k\colon T\r\bar{T}$ then $\varrho(k+c)=\varrho(k)$ and if
$\bar{k}\colon T\r\bar{T}$ is another homotopy class with
$\varrho(k)=\varrho(\bar{k})$ then there exists (a unique) $c\in
H^1(sh^3,[\varrho T,\varrho\bar{T}]^*\bar{t})$ such that
$\bar{k}=k+c$.
\end{lem}

\begin{proof}
The action $+$ is defined as follows: let $k\colon T\r \bar{T}$ be a
homotopy class of morphisms of track triangles represented by
\begin{equation*}\tag{a}
\xymatrix{A\ar[r]^f_<(.6){\;}="a"\ar[d]_{k_0}^<(.6){\;}="b"&B\ar[r]^i_<(.6){\;}="c"\ar[d]|{k_1}^<(.6){\;}="d"&
C\ar[r]^q_<(.6){\;}="e"\ar[d]|{k_2}^<(.6){\;}="f"&sA\ar[d]^{sk_0}\\
\bar{A}\ar[r]_{\bar{f}}&\bar{B}\ar[r]_{\bar{i}}&\bar{C}\ar[r]_{\bar{q}}&s\bar{A}
\ar@{=>}"a";"b"^{K_0}\ar@{=>}"c";"d"^{K_1}\ar@{=>}"e";"f"^{K_2}}
\end{equation*}
and let $(c_{-2},c_{-1},c_0)\in\tilde{G}^1([\varrho
T,\varrho\bar{T}]^*\bar{t})$ be a cocycle representing the
cohomology class $c\in H^1(sh^3,[\varrho
T,\varrho\bar{T}]^*\bar{t})$, then $k+c\colon T\r \bar{T}$ is
represented by the diagram
\begin{equation*}\tag{b}
\xymatrix{A\ar[r]^f_<(.6){\;}="a"\ar[d]_{k_0}^<(.6){\;}="b"&B\ar[r]^i_<(.6){\;}="c"\ar[d]|{k_1}^<(.6){\;}="d"&
C\ar[r]^q_<(.6){\;}="e"\ar[d]|{k_2}^<(.6){\;}="f"&sA\ar[d]^{sk_0}\\
\bar{A}\ar[r]_{\bar{f}}&\bar{B}\ar[r]_{\bar{i}}&\bar{C}\ar[r]_{\bar{q}}&s\bar{A}
\ar@{=>}"a";"b"^{\bar{K}_0}\ar@{=>}"c";"d"^{\bar{K}_1}\ar@{=>}"e";"f"^{\bar{K}_2}}
\end{equation*}
with
\begin{eqnarray*}
{\bar{K_0}}&=&\sigma_{\bar{f}{k_0}}(c_{-2})\vc K_0,\\
{\bar{K_1}}&=&\sigma_{\bar{i}{k_1}}(c_{-1})\vc K_1,\\
{\bar{K_2}}&=&\sigma_{\bar{q}{k_2}}(c_0)\vc K_2.
\end{eqnarray*}
Actually, the homotopy class of (b) only depends on $c$, because if
$(\bar{c}_{-2},\bar{c}_{-1},\bar{c}_0)-
(c_{-2},c_{-1},c_0)\in\tilde{G}^1([\varrho
T,\varrho\bar{T}]^*\bar{t})$ is the coboundary of
$(e_{-2},e_{-1},e_0)\in\tilde{G}^0([\varrho
T,\varrho\bar{T}]^*\bar{t})$ and
\begin{equation*}\tag{c}
\xymatrix{A\ar[r]^f_<(.6){\;}="a"\ar[d]_{k_0}^<(.6){\;}="b"&B\ar[r]^i_<(.6){\;}="c"\ar[d]|{k_1}^<(.6){\;}="d"&
C\ar[r]^q_<(.6){\;}="e"\ar[d]|{k_2}^<(.6){\;}="f"&sA\ar[d]^{sk_0}\\
\bar{A}\ar[r]_{\bar{f}}&\bar{B}\ar[r]_{\bar{i}}&\bar{C}\ar[r]_{\bar{q}}&s\bar{A}
\ar@{=>}"a";"b"^{\bar{\bar{K}}_0}\ar@{=>}"c";"d"^{\bar{\bar{K}}_1}\ar@{=>}"e";"f"^{\bar{\bar{K}}_2}}
\end{equation*}
is given by
\begin{eqnarray*}
\bar{\bar{K_0}}&=&\sigma_{\bar{f}{k_0}}(\bar{c}_{-2})\vc K_0,\\
\bar{\bar{K_1}}&=&\sigma_{\bar{i}{k_1}}(\bar{c}_{-1})\vc K_1,\\
\bar{\bar{K_2}}&=&\sigma_{\bar{q}{k_2}}(\bar{c}_0)\vc K_2,
\end{eqnarray*}
then the tracks $\sigma_{k_0}(e_{-2})\colon k_0\rr k_0$,
$\sigma_{k_1}(e_{-1})\colon k_1\rr k_1$ and $\sigma_{k_2}(e_0)\colon
k_2\rr k_2$ prove that (b) and (c) are homotopic.

Moreover, if (a) is homotopic to
\begin{equation*}\tag{d}
\xymatrix{A\ar[r]^f_<(.6){\;}="a"\ar[d]_{\bar{k}_0}^<(.6){\;}="b"&B\ar[r]^i_<(.6){\;}="c"\ar[d]|{\bar{k}_1}^<(.6){\;}="d"&
C\ar[r]^q_<(.6){\;}="e"\ar[d]|{\bar{k}_2}^<(.6){\;}="f"&sA\ar[d]^{s\bar{k}_0}\\
\bar{A}\ar[r]_{\bar{f}}&\bar{B}\ar[r]_{\bar{i}}&\bar{C}\ar[r]_{\bar{q}}&s\bar{A}
\ar@{=>}"a";"b"^{\tilde{K}_0}\ar@{=>}"c";"d"^{\tilde{K}_1}\ar@{=>}"e";"f"^{\tilde{K}_2}}
\end{equation*}
by the tracks $L_i\colon k_i\rr\bar{k}_i$ $(i=0,1,2)$ then (b) is
homotopic to
\begin{equation*}
\xymatrix{A\ar[r]^f_<(.6){\;}="a"\ar[d]_{\bar{k}_0}^<(.6){\;}="b"&B\ar[r]^i_<(.6){\;}="c"\ar[d]|{\bar{k}_1}^<(.6){\;}="d"&
C\ar[r]^q_<(.6){\;}="e"\ar[d]|{\bar{k}_2}^<(.6){\;}="f"&sA\ar[d]^{s\bar{k}_0}\\
\bar{A}\ar[r]_{\bar{f}}&\bar{B}\ar[r]_{\bar{i}}&\bar{C}\ar[r]_{\bar{q}}&s\bar{A}
\ar@{=>}"a";"b"^{\tilde{\tilde{K}}_0}\ar@{=>}"c";"d"^{\tilde{\tilde{K}}_1}\ar@{=>}"e";"f"^{\tilde{\tilde{K}}_2}}
\end{equation*}
with
\begin{eqnarray*}
\tilde{\tilde{K_0}}&=&\sigma_{\bar{f}{k_0}}({c}_{-2})\vc \tilde{K}_0,\\
\tilde{\tilde{K_1}}&=&\sigma_{\bar{i}{k_1}}({c}_{-1})\vc \tilde{K}_1,\\
\tilde{\tilde{K_2}}&=&\sigma_{\bar{q}{k_2}}({c}_0)\vc \tilde{K}_2,
\end{eqnarray*}
also by the tracks $L_i$. Therefore $k+c$ does not depend on the
choice of a representative for the homotopy class $k\colon
T\r\bar{T}$.

If (a) and (b) where homotopic the homotopy would be necessarily
given by tracks $\sigma_{k_0}(b_{-2})\colon k_0\rr k_0$,
$\sigma_{k_1}(b_{-1})\colon k_1\rr k_1$ and $\sigma_{k_2}(b_0)\colon
k_2\rr k_2$. In this case it is easy to check that
$\delta(b_{-2},b_{-1},b_0)=(c_{-2},c_{-1},c_0)$, so $c=0$. This
proves the effectiveness.

Finally suppose that we have two homotopy classes of morphisms
$k,\bar{k}\colon T\r\bar{T}$ defined by (a) and (d), respectively,
such that $p(k_i)=p(\bar{k}_i)$ $(i=0,1,2)$. Then if we define
\begin{eqnarray*}
c_{-2}&=&\sigma_{\bar{f}k_0}((\bar{f}(\gamma_{k_0}\vc\gamma_{\bar{k}_0}^\vi))\vc\tilde{K}_0\vc((\gamma_{\bar{k}_1}\vc\gamma_{k_1}^\vi)f)\vc K_0^\vi),\\
c_{-1}&=&\sigma_{\bar{i}k_1}((\bar{i}(\gamma_{k_1}\vc\gamma_{\bar{k}_1}^\vi))\vc\tilde{K}_1\vc((\gamma_{\bar{k}_2}\vc\gamma_{k_2}^\vi)i)\vc K_1^\vi),\\
c_0&=&\sigma_{\bar{q}k_2}((\bar{q}(\gamma_{k_2}\vc\gamma_{\bar{k}_2}^\vi))\vc\tilde{K}_2\vc(s(\gamma_{\bar{k}_0}\vc\gamma_{k_0}^\vi)q)\vc
K_2^\vi),
\end{eqnarray*}
it is not had to check that
$(c_{-2},c_{-1},c_0)\in\tilde{G}^2([\varrho
T,\varrho\bar{T}]^*\bar{t})$ is a cocycle and its cohomology class
$c$ satisfies $\bar{k}=k+c$, since the track triangle morphism (b)
defined by (a) and $(c_{-2},c_{-1},c_0)$ is homotopic to (d) by
$\gamma_{\bar{k}_i}\vc\gamma_{k_i}^\vi\colon k_i\rr\bar{k}_i$
$(i=0,1,2)$.
\end{proof}

\begin{lem}\label{E}
Given homotopy classes of track triangle morphisms
$$T\st{k}\To\bar{T}\st{\bar{k}}\To\bar{\bar{T}}$$ and cohomology
classes $c\in H^1(sh^3,[\varrho T,\varrho\bar{T}]^*\bar{t})$ and
$\bar{c}\in H^1(sh^3,[\varrho
\bar{T},\varrho\bar{\bar{T}}]^*\bar{t})$ then the action defined in
Lemma \ref{D} satisfies the following linear distributivity law
$$(\bar{k}+\bar{c})(k+c)=\bar{k}k+(\varrho(k)^*\bar{c}+\varrho(\bar{k})_*c).$$
\end{lem}

This lemma follows easily from the definition of the action $+$ in
the proof of Lemma \ref{D}.




\section{Small $\nabla$-triangles and homotopy pairs}\label{hp}

Consider an additive category $\C{A}$, an additive equivalence
$t\colon\C{A}\r\C{A}$, a translation cohomology class $\nabla\in
H^3(\C{A},t)$, and a good translation track category $(\C{B},s)$
representing $\nabla$. These data will remain fixed throughout this
section. Here we define the category of homotopy pairs in $\C{B}$.
We establish relations between this category, the category of track
triangles on $(\C{B},s)$, and the category of small
$\nabla$-triangles. These relations will play an important role in
the proof of the main theorems of \cite{ccctc2}. We also
characterize (pre)triangulated categories in terms of homotopy
pairs. By using this characterization we prove Theorem \ref{iff} at
the end of this section.

We are going to use the concept of (weak) linear extension of
categories in \cite{ah} IV.3.2 that we now recall. At the same time
we are going to define the weaker concept of (weak) linear action on a
functor, which is however different from the linear actions
considered in \cite{ah} IV.2.8. A \emph{linear
extension}\index{linear extension of categories} (resp.
\emph{action}\index{linear action})
\begin{equation}\label{wle}
D\st{+}{\To}\C{D}\st{p}{\To}\C{C}
\end{equation}
is given by a full functor $p$ which is the identity on objects, a
$\C{C}$-bimodule $D$ (i. e. a functor
$\C{C}^\op\times\C{C}\r\C{Ab}$), and (effective) right actions of
$D(pX,pY)$ on the morphism sets $\C{D}(X,Y)$ for all objects $X$ and
$Y$ in $\C{D}$ such that the orbits coincide with (resp. are
contained in) the fibers of the map
$p\colon\C{D}(X,Y)\twoheadrightarrow\C{C}(pX,pY)$. These actions
must satisfy the following \emph{linear distributivity law}: given
morphisms
$$X\st{f}\To Y\st{g}\To Z$$
in $\C{D}$ and elements $a\in D(pX,pY)$ and $b\in D(pY,pZ)$ we have
$$(g+b)(f+a)=gf+(p(f)^*b+p(g)_*a).$$

In a \emph{weak linear extension}\index{linear extension of
categories!weak} (resp. \emph{action}\index{linear action!weak}) we
do not require $p$ to be the identity on objects but only that any
object in $\C{C}$ is isomorphic to an object in the image of $p$. A
weak linear extension can always be replaced by a honest linear
extension up to equivalence of categories.
This definition
of weak linear extension is equivalent to the original definition in
\cite{ah} IV.3.3.

Let $\C{T}$ be any category. The category $\C{Pair}(\C{T})$ of pairs
on $\C{T}$ is defined as follows. Objects are morphisms in $\C{T}$
and a morphism $(h,k)\colon f\r g$ is given by a commutative square
in $\C{T}$
$$\xymatrix{\bullet\ar[r]^h\ar[d]_f&\bullet\ar[d]^g\\\bullet\ar[r]^k&\bullet}$$

If is $\C{T}$ additive and $D$ is a biadditive
$\C{T}$-bimodule we define $D^\#$ as the
$\C{Pair}(\C{T})$-bimodule given by
$$D^\#(f,g)=D(A,D)/(g_*D(A,C)+f^*D(B,D))$$
for $f\colon A\r B$ and $g\colon C\r D$ in $\C{T}$.

The category $\C{Hopair}(\C{B})$ of \emph{homotopy
pairs}\index{homotopy pair} in the track category $\C{B}$ is defined
as follows. Objects are maps $f\colon A\r B$ in $\C{B}$. A morphism
$[k_0,k_1,K]\colon f\r \bar{f}$ is represented by a diagram
\begin{equation}\label{homor}
\xymatrix{A\ar[d]_f^(.7){\;}="a"\ar[r]^{k_0}_<(.7){\;}="b"&\bar{A}\ar[d]^{\bar{f}}\\B\ar[r]_{k_1}&\bar{B}\ar@{=>}"a";"b"_{K}}
\end{equation}
in $\C{B}$. Moreover, $[k_0,k_1,K]=[\bar{k}_0,\bar{k}_1,\bar{K}]$
if and only if there are tracks $L_0\colon k_0\rr\bar{k}_0$ and
$L_1\colon k_1\rr\bar{k}_1$ such that $(\bar{f}L_0)\vc
K=\bar{K}\vc(L_1f)$. The obvious forgetful functor
$\tilde{p}_\C{B}\colon\C{Hopair}(\C{B})\r\C{Pair}(\C{A})$ fits
into a weak linear extension
\begin{equation}\label{lexho}
(\hom^t)^\#\st{+}\To\C{Hopair}(\C{B})\st{\tilde{p}_\C{B}}\To\C{Pair}(\C{A}),
\end{equation}
compare \cite{tbhp}. The action is defined as follows, given a
morphism $[k_0,k_1,K]$ represented by a diagram as above and
$\tilde{g}\in(\hom^t)^\#(p(f),p(\bar{f}))$ represented by $g\colon
tA\r \bar{B}$ then $[k_0,k_1,K]+\tilde{g}=[k_0,k_1,
K\vc\sigma_{k_1 f}(g)]$.

The category $\C{Pair}(\C{A})$ is obviously additive and
$(\hom^t)^\#$ is a biadditive $\C{Pair}(\C{A})$-bimodule. By
\cite{qcsr} 6.2 the category $\C{Hopair}(\C{B})$ is also additive in
a unique compatible way with the structure of weak linear extension
in (\ref{lexho}). One can easily check that the full subcategory of
$\C{Hopair}(\C{B})$ given by the maps $\iota X\colon *\r X$ is
isomorphic to $\C{A}$. We denote
$$\iota\colon\C{A}\To\C{Hopair}(\C{B})$$ to the full additive
inclusion. The \emph{cone functor}\index{cone functor} associated to
a map $f\colon A\r B$ in $\C{B}$ is the additive functor
\begin{equation}\label{conenew}
\cone_f=\hom_{\C{Hopair}(\C{B})}(f,\iota)\colon\C{A}\To\Ab.
\end{equation}
Cone functors are introduced in \cite{ttc} in a different but equivalent way. There $\cone_f(X)$ is defined as the quotient set
$$\cone_f(X)=\set{(\alpha,\hat{\alpha})\,;\,\alpha\colon B\r X\text{ in }\C{B}\text{ and }\hat{\alpha}\colon\alpha f\rr 0}/\simeq$$
where $(\alpha,\hat{\alpha})\simeq(\beta,\hat{\beta})$ if there
exists $H\colon\alpha\rr\beta$ such that
$\hat{\alpha}=\hat{\beta}\vc(Hf)$. The equivalence class represented
by $(\alpha,\hat{\alpha})$ is denoted by
$\set{\alpha,\hat{\alpha}}$. the correspondence with the definition
given here is determined by the identification
$\set{\alpha,\hat{\alpha}}=[0,\alpha,\hat{\alpha}]\colon f\r\iota
X$.

Recall that $\C{A}$-modules are additive functors $\C{A}\r\Ab$.
The category of $\C{A}$-modules and natural transformations between
them is denoted by $\mo(\C{A})$. The category $\C{A}^\op$ is fully
included in $\mo(\C{A})$ through the functor
$X\mapsto\hom_\C{A}(X,-)$. The advantage of the new approach to cone
functors is that we immediately obtain the following result.

\begin{prop}\label{conefun}
Cone functors define an additive functor
$$\cone\colon\C{Hopair}(\C{B})^\op\To\mo(\C{A}).$$
\end{prop}

\begin{proof}
The functor $\cone$ corresponds to
$$\hom_{\C{Hopair}(\C{B})}(-,\iota)\colon\C{Hopair}(\C{B})^\op\times\C{A}\To\C{Ab}$$
by the exponential law. With the notation in \cite{ttc} a morphism
of homotopy pairs $[k_0,k_1,K]\colon f\r\bar{f}$
 induces the $\C{A}$-module morphism
$$[k_0,k_1,K]^*\colon\cone_{\bar{f}}\To\cone_f$$
 defined by $\set{\alpha,\hat{\alpha}}\mapsto\set{\alpha
k_1,(\hat{\alpha}k_0)\vc(\alpha K)}$.
\end{proof}

Cone
functors are embedded in exact sequences of $\C{A}$-modules
\begin{equation}\label{exaa}
\hom_\C{A}(tB,-)\st{(tp(f))^*}\To\hom_\C{A}(tA,-)\st{\tilde{q}_f}\To\cone_f\st{\tilde{i}_f}\To\hom_\C{A}(B,-)\st{p(f)^*}\To
\hom_\C{A}(A,-).
\end{equation}
Here $\tilde{q}_f(g)=\set{0,\sigma_0(g)}$ and
$\tilde{i}_f\set{\alpha,\hat{\alpha}}=p(\alpha)$, see \cite{ttc}
3.4. In fact the exactness of this sequence is an immediate
consequence of the axioms of a weak linear extension for
(\ref{lexho}) and the compatibility with the additive structure on
$\C{Hopair}(\C{B})$, see \cite{qcsr} 6.2.

\begin{rem}\label{esnat}
This exact sequence is natural in $f$ in the following sense, given
a morphism of homotopy pairs $[k_0,k_1,K]\colon f\r\bar{f}$
represented by a diagram as in (\ref{homor}) the following diagram
of $\C{A}$-modules commutes
$$\xymatrix{\hom_\C{A}(tB,-)\ar[r]^{(tp(f))^*}&\hom_\C{A}(tA,-)\ar[r]^<(.3){\tilde{q}_f}&
\cone_f\ar[r]^<(.2){\tilde{i}_f}&\hom_\C{A}(B,-)\ar[r]^{p(f)^*}&
\hom_\C{A}(A,-)\\
\hom_\C{A}(t\bar{B},-)\ar[u]_{(tp(k_1))^*}\ar[r]_{(tp(\bar{f}))^*}&\hom_\C{A}(t\bar{A},-)\ar[u]_{(tp(k_0))^*}\ar[r]_<(.3){\tilde{q}_{\bar{f}}}&
\cone_{\bar{f}}\ar[r]_<(.2){\tilde{i}_{\bar{f}}}\ar[u]_{[k_0,k_1,K]^*}&\hom_\C{A}(\bar{B},-)\ar[u]_{p(k_1)^*}\ar[r]_{p(\bar{f})^*}&
\hom_\C{A}(\bar{A},-)\ar[u]_{p(k_0)^*}}$$
\end{rem}

\begin{rem}\label{defino}
Let $\C{Hopair}^{rep}(\C{B})\subset\C{Hopair}(\C{B})$ be the full
subcategory given by those maps in $\C{B}$ such that the associated
cone functor is representable, for any such a map $f$ in $\C{B}$ one
can choose an $\C{A}$-module isomorphism
$$\chi_f\colon\hom_\C{A}(C_f,-)\cong\cone_f.$$
A particular choice determines a functor
\begin{equation}\label{elcho}
C\colon\C{Hopair}^{rep}(\C{B})\To\C{A}
\end{equation}
such that the restriction of $\cone$ in Proposition \ref{conefun} to
$\C{Hopair}^{rep}(\C{B})^{op}$ coincides up to natural equivalence
with
$$\C{Hopair}^{rep}(\C{B})^{op}\st{C^{op}}\To\C{A}^{op}\subset\mo(\C{A}).$$
The natural equivalence $\chi$ is given by the representations
above. By using $C$ we can define a functor
\begin{equation}\label{lazeta}
\zeta\colon\C{Hopair}^{rep}(\C{B})\To\C{cand}(\C{A},t)
\end{equation}
as follows. Given an object $f\colon A\r B$ in
$\C{Hopair}^{rep}(\C{B})$ the small candidate triangle $\zeta(f)$ is
\begin{equation}\label{elcho2}
A\st{p(f)}\To B\st{i_f}\To C_f\st{q_f}\To tA,
\end{equation}
where $i_f$ and $q_f$ are determined by the $\C{A}$-module morphisms
$\tilde{i}_f$ and $\tilde{q}_f$ in (\ref{exaa}), respectively, the
isomorphism $\chi_f$ and Yoneda's lemma. The functor $\zeta$ sends a
morphism $[k_0,k_1,K]$ in $\C{Hopair}^{rep}(\C{B})$ represented by
diagram (\ref{homor}) to the following morphism of small candidate
triangles
$$\xymatrix@C=40pt{A\ar[r]^{p(f)}\ar[d]^{p(k_0)}&B\ar[r]^{i_f}\ar[d]^{p(k_1)}&C_f\ar[r]^{q_f}\ar[d]^{C[k_0,k_1,K]}&tA\ar[d]^{tp(k_0)}\\
\bar{A}\ar[r]_{p(\bar{f})}&\bar{B}\ar[r]_{i_{\bar{f}}}&C_{\bar{f}}\ar[r]_{q_{\bar{f}}}&t\bar{A}}$$
This diagram commutes by Remark \ref{esnat}. The functor $\zeta$
does not depend on the choice of the isomorphisms $\chi_f$, up to
natural equivalence.
\end{rem}

The following proposition shows that track triangles determine
representations of cone functors.

\begin{prop}\label{repre}
Given a track triangle in $(\C{B},s)$
$$\xymatrix{A\ar@/^30pt/[rr]^0_{\;}="a"\ar[r]_{f}&B\ar@/_30pt/[rr]_0^{\;}="c"
\ar[r]^<(-0.22){\;}="b"|{i}&
C\ar[r]_<(-0.17){\;}="d"|{q}\ar@/^30pt/[rr]^0_{\;}="f"&sA\ar[r]_{s(f)}^<(-0.23){\;}="e"&
sB\ar@{=>}"b";"a"_{H_0}\ar@{=>}"d";"c"^{H_{1}}\ar@{=>}"e";"f"_{H_{2}}}$$
the $\C{A}$-module morphism
$$\chi_{\set{i,H_0}}\colon\hom_\C{A}(C,-)\To\cone_f$$ represented by
$\set{i,H_0}\in\cone_f(C)$ is an isomorphism.
\end{prop}

\begin{proof}
The inverse sends $\set{\alpha,\hat{\alpha}}\in\cone_f(X)$ to the
unique morphism $l\colon C\r X$ in $\C{A}$ such that
$\sigma_0(\bar{t}l)=((s(\hat{\alpha})\vc s_{\alpha,f})q)\vc(\alpha
H_2^\vi)$. Compare \cite{ttc} 4.
\end{proof}

\begin{rem}\label{conmuta}
In the conditions of the statement of Proposition \ref{repre} the
following diagram of $\C{A}$-modules commutes
$$\xymatrix{&\hom_\C{A}(C,-)\ar[rd]^{p(i)^*}\ar[dd]^{\chi_{\set{i,H_0}}}_\cong&\\\hom_\C{A}(tA,-)\ar[ru]^{p(q)^*}\ar[rd]_{\tilde{q}_f}&&
\hom_\C{A}(B,-)\\&\cone_f\ar[ru]_{\tilde{i}_f}}$$
\end{rem}

By Proposition \ref{repre} there is an obvious forgetful functor
\begin{equation}\label{barrito}
\bar{\varrho}=\bar{\varrho}_{(\C{B},s)}\colon
\C{candt}(\C{B},s)\To\C{Hopair}^{rep}(\C{B})
\end{equation}
sending
the homotopy class of a track triangle morphism as in
(\ref{morctra}) to $[k_0,k_1,K_0]$.

There is also an obvious functor
\begin{equation}\label{varsigma}
\varsigma\colon\C{cand}(\C{A},t)\To\C{Pair}(\C{A})
\end{equation}
sending a small candidate triangle
$$A\st{f}\To B\st{i}\To C\st{q}\To tA$$
to $f\colon A\r B$ regarded as an object of $\C{Pair}(\C{A})$.

\begin{prop}\label{noc}
In the following diagram of functors
$$\xymatrix@C=40pt{&\C{Hopair}^{rep}(\C{B})\ar[dd]^\zeta\ar[rd]^{\tilde{p}}&\\\C{candt}(\C{B},s)\ar[ru]^{\bar{\varrho}}\ar[rd]_\varrho&&\C{Pair}(\C{A})\\&\C{cand}(\C{A},t)\ar[ru]_\varsigma&}$$
the triangle in the left commutes up to natural equivalence and the
triangle in the right is strictly commutative.
\end{prop}

\begin{proof}
The statement about the triangle in the right-hand side is obvious.
Let
\begin{equation*}\tag{a}
\xymatrix{A\ar@/^30pt/[rr]^0_{\;}="a"\ar[r]_{f}&B\ar@/_30pt/[rr]_0^{\;}="c"
\ar[r]^<(-0.22){\;}="b"|{i}&
D\ar[r]_<(-0.17){\;}="d"|{q}\ar@/^30pt/[rr]^0_{\;}="f"&sA\ar[r]_{s(f)}^<(-0.23){\;}="e"&
sB\ar@{=>}"b";"a"_{H_0}\ar@{=>}"d";"c"^{H_{1}}\ar@{=>}"e";"f"_{H_{2}}}
\end{equation*}
be any track triangle. Possibly the representation
$$\chi_f\colon\hom_\C{A}(C_f,-)\cong\cone_f$$
 chosen for the definition of $C$ in (\ref{elcho}) and $\zeta$ in (\ref{lazeta}) is different from
$$\chi_{\set{i,H_0}}\colon\hom_\C{A}(D,-)\cong\cone_f$$ in Proposition \ref{repre}.
However if the morphism $l\colon D\r C_f$ in $\C{A}$ represents
$$\hom_\C{A}(C_f,-)\st{\chi_f}\cong\cone_{f}\st{\chi^{-1}_{\set{{i},{H}_0}}}\cong\hom_\C{A}(D,-)$$
then the following isomorphism of small candidate triangles
$$\xymatrix{A\ar@{=}[d]\ar[r]^{p(f)}&B\ar@{=}[d]\ar[r]^{p(i)}&D\ar[d]_\cong^l
\ar[r]^{p(q)}&tA\ar@{=}[d]\\
{A}\ar[r]_{p({f})}&{B}\ar[r]_{i_{f}}&C_{{f}}\ar[r]_{q_{{f}}}&t{A}}$$
is natural in the track triangle (a). This can be checked by using
the definition of $\zeta$  and the next lemma.
\end{proof}


\begin{lem}\label{pepito}
Given a morphism between track triangles
$$\xymatrix{A\ar[r]^f_<(.6){\;}="a"\ar[d]_{k_0}^<(.6){\;}="b"&B\ar[r]^i_<(.6){\;}="c"\ar[d]|{k_1}^<(.6){\;}="d"&
C\ar[r]^q_<(.6){\;}="e"\ar[d]|{k_2}^<(.6){\;}="f"&sA\ar[d]^{sk_0}\\
\bar{A}\ar[r]_{\bar{f}}&\bar{B}\ar[r]_{\bar{i}}&\bar{C}\ar[r]_{\bar{q}}&s\bar{A}
\ar@{=>}"a";"b"^{K_0}\ar@{=>}"c";"d"^{K_1}\ar@{=>}"e";"f"^{K_2}}$$
the composite $\C{A}$-module morphism
$$\hom_\C{A}(\bar{C},-)\st{\chi_{\set{\bar{i},\bar{H}_0}}}\cong\cone_{\bar{f}}\st{[k_0,k_1,K_0]^*}\To\cone_f\st{\chi^{-1}_{\set{i,H_0}}}\cong\hom_\C{A}(C,-)$$
is represented by $p(k_2)\colon C\r\bar{C}$.
\end{lem}

\begin{proof}
It is enough to check that
$[k_0,k_1,K_0]^*\set{\bar{i},\bar{H}_0}=\set{\bar{i}k_1,\bar{H}_0\vc(\bar{i}K_0)}$
coincides with $\chi_{\set{i,H_0}}(p(k_2))=\set{k_2i,k_2H_0}$ in
$\cone_f(\bar{C})$. For this one uses $K_1\colon k_2 i\rr
\bar{i}k_1$ and the definition of track triangle morphism.
\end{proof}

The following result allows us to prove (Tr3) for cohomologically
pretriangulated categories. Let $\C{cand}(\nabla)$ be the category
of small $\nabla$-triangles.

\begin{prop}\label{basic1}
The functor $\varsigma$ is full over the full subcategory
$\C{cand}(\nabla)$ of small $\nabla$-triangles.
\end{prop}

\begin{proof}

By Proposition \ref{repre}, Remark \ref{conmuta} and Yoneda's lemma
any small $\nabla$-triangle is isomorphic to an object in the image
of the functor $\zeta$ in (\ref{lazeta}). By Proposition \ref{noc}
the composite functor $\varsigma\zeta$ coincides with $\tilde{p}$,
which is known to be full, hence this proposition follows.
\end{proof}

Suppose now that $(\C{A},t,\nabla)$ is cohomologically
pretriangulated. By (TTr1) and Proposition \ref{repre} 
all cone functors are representable. By Remark \ref{defino} a choice
of representations determines a functor
$$\zeta\colon\C{Hopair}(\C{B})\r\C{cand}(\C{A},t)$$
which is well defined up to natural equivalence.

\begin{prop}\label{mama}
If $(\C{A},t,\nabla)$ is cohomologically pretriangulated then the
functor $\zeta$ takes values in small
$\nabla$-triangles.
\end{prop}

\begin{proof}
By Proposition \ref{repre}, Remark \ref{conmuta} and Yoneda's lemma
for any map $f$ in $\C{B}$ the small candidate triangle $\zeta(f)$
is isomorphic to a small $\nabla$-triangle, now the proposition
follows from Theorem \ref{elt} (pretriangulated case) and (Tr0).
\end{proof}

In the next proposition we show that the functor $\bar{\varrho}$ in
(\ref{barrito}) is compatible with the actions of
$H^1(sh^3,[-,-]^*\bar{t})$ and $(\hom^t)^\#$ in Lemma \ref{D} and
(\ref{lexho}), respectively. For this we consider the natural
transformation of $\C{cand}(\C{A},t)$-bimodules
\begin{equation}
\varpi\colon
H^1(sh^3,[-,-]^*\bar{t})\rr(\hom^t)^\#(\varsigma,\varsigma)
\end{equation}
sending a cohomology class represented by a $1$-cocycle
$(c_{-2},c_{-1},c_0)\in\tilde{G}^1([T,\bar{T}]^*\bar{t})$ to the element
represented by $c_{-2}$ in $(\hom^t)^\#(\varsigma T,\varsigma
\bar{T})$.

\begin{prop}\label{equi1}
For any morphism $f\colon T\r\bar{T}$ in $\C{candt}(\C{B},s)$ and
any $\alpha\in H^1(sh^3,[\varrho T,\varrho\bar{T}]^*\bar{t})$ we
have that
$$\bar{\varrho}(f+\alpha)=\bar{\varrho}(f)+\varpi_{(\varrho T,\varrho\bar{T})}(\alpha).$$
\end{prop}
This result follows immediately from the definition of the linear
actions.

There is a similar result for the functor $\zeta$ in (\ref{lazeta}).
In order to state it we construct first a weak linear action on the
functor $\varsigma$ in (\ref{varsigma}). 
For this we consider
the $\C{Pair}(\C{A})$-bimodule
$$\hom_{\mo(\C{A}^\op)}(\ker t,\coker).$$
Here $\hom_{\mo(\C{A}^\op)}$ denotes morphism abelian groups in the
category of $\C{A}^\op$-modules, and
$$\ker t,\coker\colon\C{Pair}(\C{A})\To\mo(\C{A}^\op)$$
send a morphism $f\colon A\r B$ in $\C{A}$ to the kernel of $tf$ and
the cokernel of $f$ in the abelian category of $\C{A}^\op$-modules through the
Yoneda full inclusion $\C{A}\subset\mo(\C{A}^\op)$, respectively.

\begin{prop}\label{varsiwle}
There is a weak linear action
$$\hom_{\mo(\C{A}^\op)}(\ker t,\coker)\st{+}\To\C{cand}(\C{A},t)\st{\varsigma}\To\C{Pair}(\C{A}).$$
Moreover, if $(\C{A},t,\nabla)$ is a cohomologically pretriangulated
category then this weak linear action restricts to a weak linear
extension on the full subcategory of small $\nabla$-triangles
$$\hom_{\mo(\C{A}^\op)}(\ker t,\coker)\st{+}\To\C{cand}(\nabla)\st{\varsigma}\To\C{Pair}(\C{A}).$$
\end{prop}

\begin{proof}
The second part of this proposition actually holds for any
pretriangulated category if we change $\C{cand}(\nabla)$ by the
category of exact triangles. For any small candidate triangle
$$A\st{f}\To B\st{i}\To C\st{q}\To tA$$
we denote by
$$\hat{i}\colon\coker f\To\hom_\C{A}(-,C)$$ and
$$\hat{q}\colon\hom_\C{A}(-,C)\To\ker tf$$
to the $\C{A}^\op$-module morphisms induced by the equalities $if=0$
and $(tf)q=0$, respectively. Given a morphism $k$ between small
candidate triangles
$$\xymatrix{A\ar[r]^f\ar[d]^{k_0}&B\ar[r]^i\ar[d]^{k_1}&C\ar[r]^q\ar[d]^{k_2}&tA\ar[d]^{tk_0}\\
\bar{A}\ar[r]^{\bar{f}}&\bar{B}\ar[r]^{\bar{i}}&\bar{C}\ar[r]^{\bar{q}}&t\bar{A}}$$
and $\alpha\in\hom_{\mo(\C{A}^\op)}(\ker tf,\coker\bar{f})$ the
morphism $k+\alpha$ is given by the diagram
$$\xymatrix{A\ar[r]^f\ar[d]^{k_0}&B\ar[r]^i\ar[d]^{k_1}&C\ar[r]^q\ar[d]^{k_2+\hat{\bar{i}}\alpha\hat{q}}&tA\ar[d]^{tk_0}\\
\bar{A}\ar[r]^{\bar{f}}&\bar{B}\ar[r]^{\bar{i}}&\bar{C}\ar[r]^{\bar{q}}&t\bar{A}}$$
where we identify the $\C{A}^\op$-module morphism
$$\hat{\bar{i}}\alpha\hat{q}\colon\hom_\C{A}(-,C)\To\hom_\C{A}(-,\bar{C})$$
with the corresponding morphism $C\r\bar{C}$ in $\C{A}$ by Yoneda's
lemma.

The proof of the second part is an easy exercise of pretriangulated
category theory.
\end{proof}

Let
\begin{equation}
\Xi\colon(\hom^t)^\#\rr\hom_{\mo(\C{A}^\op)}(\ker t,\coker)
\end{equation}
be the natural transformation of $\C{Pair}(\C{A})$-bimodules defined
as follows: given $f\colon A\r B$ and
$\bar{f}\colon\bar{A}\r\bar{B}$ in $\C{A}$ and an element
$\alpha\in(\hom^t)^\#(f,\bar{f})$ represented by
$\tilde{\alpha}\colon tA\r \bar{B}$ in $\C{A}$ then
$\Xi(\alpha)\in\hom_{\mo(\C{A}^\op)}(\ker tf,\coker\bar{f})$ is
represented by the composite
$$\ker tf\hookrightarrow\hom_\C{A}(-,tA)\st{\tilde{\alpha}}\To\hom_\C{A}(-,\bar{B})\twoheadrightarrow\coker\bar{f},$$
where the first and the last arrows are the inclusion and the
projection of the kernel and the cokernel, respectively.

\begin{prop}\label{equi2}
Given a morphism $[k_0,k_1,K]\colon f\r\bar{f}$ in
$\C{Hopair}^{rep}(\C{B})$ and
$\alpha\in(\hom^t)^\#(\tilde{p}(f),\tilde{p}(\bar{f}))$ then
$$\zeta([k_0,k_1,K]+\alpha)=\zeta([k_0,k_1,K])+\Xi_{(\tilde{p}(f),\tilde{p}(\bar{f}))}(\alpha).$$
\end{prop}

\begin{proof}
By using the compatibility of the additive structure on
$\C{Hopair}(\C{B})$ given by \cite{qcsr} 6.2 with the linear
extension of categories (\ref{lexho}) and our new definition of cone
functors in (\ref{conenew}) one readily checks that if
$\tilde{\alpha}\colon tA\r\bar{B}$ represents $\alpha$ then
$$([k_0,k_1,K]+\alpha)^*=[k_0,k_1,K]^*+\tilde{q}_f\tilde{\alpha}^*\tilde{i}_{\bar{f}}\colon\cone_{\bar{f}}\To\cone_f.$$
Here we use the $\C{A}^\op$-module morphisms in the exact sequence
(\ref{exaa}). By definition of $C$ in (\ref{elcho}) this implies
that
$$C([k_0,k_1,K]+\alpha)=C[k_0,k_1,K]+i_{\bar{f}}\tilde{\alpha}q_f$$
where $i_{\bar{f}}$ and $q_f$ are given by the definition of $\zeta$
on objects in (\ref{elcho2}). This last equation yields the formula
in the statement.
\end{proof}

We record the following further property of the natural
transformations $\varpi$ and $\Xi$ which is going to be relevant
in \cite{ccctc2}. 

\begin{lem}\label{exacti}
The composition of the two following $\C{cand}(\C{A},t)$-bimodule
morphisms is trivial
$$H^1(sh^3,[-,-]^*\bar{t})\st{\varpi}\To(\hom^t)^\#(\varsigma,\varsigma)\st{\Xi(\varsigma^\op\times\varsigma)}\To\hom_{\mo(\C{A})}(\ker t\varsigma,\coker\varsigma).$$
Moreover, if $(\C{A},t,\nabla)$ is cohomologically pretriangulated
and we restrict to the full subcategory $\C{cand}(\nabla)$ of small
$\nabla$-triangles then the corresponding sequence of
$\C{cand}(\nabla)$-bimodule morphisms is exact.
\end{lem}

\begin{proof}
The first part of the statement is an easy exercise based on the
definition of the cochain complex $\tilde{G}^*$ in Section 3, and we
leave it to the reader.

Suppose now that $(\C{A},t,\nabla)$ is cohomologically pretriangulated and that we have two small $\nabla$-triangles
$$T=\set{A\st{f}\To B\st{i}\To C\st{q}\To tA},$$
$$\bar{T}=\set{\bar{A}\st{\bar{f}}\To \bar{B}\st{\bar{i}}\To \bar{C}\st{\bar{q}}\To t\bar{A}},$$
and an element $\alpha\in(\hom^t)^\#(\varsigma T,\varsigma \bar{T})$
represented by $\tilde{\alpha}\colon tA\r\bar{B}$ such that
$\Xi(\alpha)=0$, i. e. the composite
\begin{equation*}\tag{a}
\ker tf\hookrightarrow\hom_\C{A}(-,tA)\st{\tilde{\alpha}}\To\hom_\C{A}(-,\bar{B})\twoheadrightarrow\coker\bar{f}
\end{equation*}
vanishes. Since $T$ and $\bar{T}$ are small $\nabla$-triangles by
Theorem \ref{elt} (pretriangulated case) and well-known properties
of pretriangulated categories the sequences of $\C{A}^\op$-modules
$$\hom_\C{A}(-,C)\st{q}\To\hom_\C{A}(-,tA)\st{tf}\To\hom_\C{A}(-,tB),$$
\begin{equation*}\tag{b}
\hom_\C{A}(-,\bar{A})\st{\bar{f}}\To\hom_\C{A}(-,\bar{B})\st{\bar{i}}\To\hom_\C{A}(-,\bar{C}),
\end{equation*}
are exact, therefore the vanishing of (a) is equivalent to the
identity $\bar{i}\tilde{\alpha}q=0$. Since (b) is exact and
$\bar{i}\tilde{\alpha}q=0$ then by Yoneda's lemma there exists a
commutative diagram in $\C{A}$
$$\xymatrix{C\ar[r]^q\ar[d]_{k_0}&tA\ar[r]^{tf}\ar[d]^{\tilde{\alpha}}&tB\ar[r]^{ti}&tC\ar[d]^{tk_0}\\\bar{A}\ar[r]_{\bar{f}}&\bar{B}\ar[r]_{\bar{i}}&
\bar{C}\ar[r]_{\bar{q}}&t\bar{A}}$$ Now by (Tr0) and (Tr2) the upper
row is an exact triangle in the pretriangulated structure for
$\C{A}$ given by Theorem \ref{elt}, hence by (Tr3) there is a
commutative extension
$$\xymatrix{C\ar[r]^q\ar[d]_{k_0}&tA\ar[r]^{tf}\ar[d]^{\tilde{\alpha}}&tB\ar[r]^{ti}\ar[d]^{k_2}&tC\ar[d]^{tk_0}\\\bar{A}\ar[r]_{\bar{f}}&\bar{B}\ar[r]_{\bar{i}}&
\bar{C}\ar[r]_{\bar{q}}&t\bar{A}}$$ By using the commutativity of
this diagram the reader can easily check that
$(\tilde{\alpha},-k_2,tk_0)\in\tilde{G}^1([T,\bar{T}]^*\bar{t})$ is
a $1$-cocycle whose cohomology class $c$ satisfies
$\varpi(c)=\alpha$.

\end{proof}




The following lemmas state translations of axioms (TTr1) and (TTr2).

\begin{lem}\label{unacar}
The pair $(\C{B},s)$ admits a pretriangulated track category
structure if and only if all cone functors are representable and the
functor
$$\zeta\colon\C{Hopair}(\C{B})\To\C{cand}(\C{A},t)$$
defined in (\ref{lazeta}) takes values in small $\nabla$-triangles.
\end{lem}

\begin{proof}
In order to prove the ``if'' part we only need to construct a class
of distinguished track triangles satisfying (TTr1). For any map $f$
in $\C{B}$ the small candidate triangle $\zeta(f)$ is a small
$\nabla$-triangle so it is the image by $\varrho$ of a track
triangle
\begin{equation*}
\xymatrix{A\ar@/^30pt/[rr]^0_{\;}="a"\ar[r]_{\tilde{f}}&B\ar@/_30pt/[rr]_0^{\;}="c"
\ar[r]^<(-0.22){\;}="b"|{i}&
C\ar[r]_<(-0.17){\;}="d"|{q}\ar@/^30pt/[rr]^0_{\;}="f"&sA\ar[r]_{s(\tilde{f})}^<(-0.23){\;}="e"&
sB\ar@{=>}"b";"a"_{H_0}\ar@{=>}"d";"c"^{H_{1}}\ar@{=>}"e";"f"_{H_{2}}}
\end{equation*}
In particular $p(\tilde{f})=f$. We choose a track $\gamma\colon
f\rr\tilde{f}$ and set the track triangle
\begin{equation*}
\xymatrix@C=70pt{ {A}\ar@/^20pt/[r]|{
{\tilde{f}}}_{\;}="h"\ar@/^40pt/[rr]^0_{\;\;\;\;\;}="a"\ar[r]_{f}^{\;}="g"&
{B}\ar@/_30pt/[rr]_0^{\;}="c" \ar[r]^<(-.09){\;}="b"|{ {i}}&
 {C}\ar[r]_<(-0.075){\;}="d"|{ {q}}\ar@/^40pt/[rr]^0_{\;}="f"&s {A}
\ar@/^20pt/[r]|{s
(\tilde{f})}_{\;}="j"\ar[r]_{s(f)}^<(-0.1){\;\;\;\;\;}="e"^{\;}="i"&
s {B}\ar@{=>}"b";"a"_{ {H}_0}\ar@{=>}"d";"c"^{
{H}_{1}}\ar@{=>}"e";"f"^{ {H}_{2}}
\ar@{=>}"i";"j"_{s(\gamma)}\ar@{=>}"g";"h"_{\gamma}}
\end{equation*}
to be distinguished. Now (TTr1) is clearly satisfied.

The ``only if'' part follows from Propositions \ref{repre} and
\ref{mama}.
\end{proof}

\begin{lem}\label{otracar}
The pair $(\C{B},s)$ admits a triangulated track category structure
if and only if there is a full subcategory $\C{dist}(\C{B},s)$ of
$\C{candt}(\C{B},s)$ such that the composite functor
\begin{equation*}\tag{a}
\C{dist}(\C{B},s)\subset\C{candt}(\C{B},s)\st{\bar{\varrho}_{(\C{B},s)}}\To\C{Hopair}(\C{B})\st{\tilde{p}_\C{B}}\To\C{Pair}(\C{A})
\end{equation*}
is surjective on objects 
and the restriction of $\bar{\varrho}_{(\C{B},s)}$ to
$\C{dist}(\C{B},s)$ is full.
\end{lem}

\begin{proof}
If $(\C{B},s)$ is a triangulated track category and
$\C{dist}(\C{B},s)$ is given by distinguished track triangles then
(a) is surjective on objects by (TTr1), since $\tilde{p}$ is already
known to be surjective on objects, and $\bar{\varrho}_{(\C{B},s)}$
restricted to $\C{dist}(\C{B},s)$ is full by (TTr2).

Conversely if (a) is surjective on objects we can define
distinguished track triangles as follows. We choose for any morphism
$f\colon A\r B$ in $\C{A}$ a track triangle in $\C{dist}(\C{B},s)$
\begin{equation*}\tag{b}
\xymatrix{A\ar@/^30pt/[rr]^0_{\;}="a"\ar[r]_{\ell
f}&B\ar@/_30pt/[rr]_0^{\;}="c" \ar[r]^<(-0.22){\;}="b"|{i_f}&
C\ar[r]_<(-0.17){\;}="d"|{q_f}\ar@/^30pt/[rr]^0_{\;}="f"&sA\ar[r]_{s(\ell
f)}^<(-0.23){\;}="e"&
sB\ar@{=>}"b";"a"_{H_0^f}\ar@{=>}"d";"c"^{H_{1}^f}\ar@{=>}"e";"f"_{H_{2}^f}}
\end{equation*}
with $p(\ell f)=f$. Moreover, for any map $g\colon A\r B$ in $\C{B}$
we choose a track $\gamma_g\colon\ell p(g)\rr g$ with $\gamma_{\ell
f}=0^\vc_{\ell f}$. Now we define the distinguished track triangles
as the track triangles
\begin{equation*}\tag{c}
\xymatrix@C=70pt{ {A}\ar@/^20pt/[r]|{ {\ell
p(g)}}_{\;}="h"\ar@/^40pt/[rr]^0_{\;\;\;\;\;}="a"\ar[r]_{g}^{\;}="g"&
{B}\ar@/_30pt/[rr]_0^{\;}="c" \ar[r]^<(-0.08){\;}="b"|{ {i_{p(g)}}}&
 {C}\ar[r]_<(-0.07){\;}="d"|{ {q_{p(g)}}}\ar@/^40pt/[rr]^0_{\;}="f"&s {A}
\ar@/^20pt/[r]|{s \ell p(g)}_{\;}="j"\ar[r]_{s
g}^<(-0.1){\;\;\;\;\;}="e"^{\;}="i"& s {B}\ar@{=>}"b";"a"_{
{H}_0^{p(g)}}\ar@{=>}"d";"c"^{ {H}_{1}^{p(g)}}\ar@{=>}"e";"f"^{
{H}_{2}^{p(g)}}
\ar@{=>}"i";"j"_{s\gamma_g^\vi}\ar@{=>}"g";"h"_{\gamma_g^\vi}}
\end{equation*}
for all maps $g$ in $\C{B}$. The axiom (TTr1) is now clearly
satisfied.
Moreover, (TTr2) also holds since (c) is isomorphic to (b) for
$f=p(g)$. The isomorphism is given by diagram (c) itself.
\end{proof}

Now we are ready to prove Theorem \ref{iff}.

\begin{proof}[Proof of Theorem \ref{iff}]

As in the proof of Proposition \ref{sontodos} we choose a
pseudofunctor $\varphi\colon\C{B}\pf\C{C}$ in the conditions of
Definition \ref{trantra} realizing the equivalence between $(\C{B},s)$ and $(\C{C},r)$ and a homotopy $\xi\colon r\varphi \rr
\varphi s$. We recall that the pseudofunctors $s$ and $r$ coincide
on objects (they both actually coincide with $t$) and that
$\varphi$ is the identity on objects, in particular by \cite{ptc}
A.14 we can suppose that $\varphi$ is normalized with respect to
zero maps. In the proof of Proposition \ref{sontodos} we defined a
functor
$$(\varphi,\xi)_*\colon\C{candt}(\C{B},s)\To\C{candt}(\C{C},r)$$
induced by $\varphi$ and $\xi$. The pseudofunctor $\varphi$ also
induces a functor between the categories of homotopy pairs
$$\varphi_*\colon\C{Hopair}(\C{B})\To\C{Hopair}(\C{C}).$$
This functor sends a map $f\colon A\r B$ in $\C{B}$ regarded as an object in
$\C{Hopair}(\C{B})$ to $\varphi(f)\colon A\r B$. The image of the
morphism represented by the diagram
\begin{equation*}
\xymatrix{A\ar[d]_f^(.7){\;}="a"\ar[r]^{k_0}_<(.7){\;}="b"&\bar{A}\ar[d]^{\bar{f}}\\B\ar[r]_{k_1}&\bar{B}\ar@{=>}"a";"b"_{K}}
\end{equation*}
in $\C{B}$ is given by
\begin{equation*}
\xymatrix@R=90pt@C=90pt{A\ar@/_20pt/[rd]|{\varphi(k_1f)}^{\;}="c"_{\;}="e"\ar@/^20pt/[rd]|{\varphi(\bar{f}k_0)}_{\;}="d"^{\;}="f"\ar[d]_{\varphi(f)}^(1){\;}="a"\ar[r]^{\varphi(k_0)}_<(1){\;}="b"&
\bar{A}\ar[d]^{\varphi(\bar{f})}\\B\ar[r]_{\varphi(k_1)}&\bar{B}\ar@{=>}"c";"d"_{\varphi(K)}\ar@{=>}"a";"e"_{\varphi_{k_1,f}}\ar@{=>}"f";"b"^{\varphi^\vi_{\bar{f},k_0}}}
\end{equation*}
The following diagram of functors is commutative
\begin{equation*}\tag{a}
\xymatrix@C=40pt{\C{candt}(\C{B},s)\ar[r]^{\bar{\varrho}_{(\C{B},s)}}\ar[d]_{(\varphi,\xi)_*}&\C{Hopair}(\C{B})\ar[d]^{\varphi_*}\\
\C{candt}(\C{C},r)\ar[r]_{\bar{\varrho}_{(\C{C},r)}}&\C{Hopair}(\C{C})}
\end{equation*}
Moreover, $\varphi_*$ induces a map of weak linear extensions
\begin{equation*}\tag{b}
\xymatrix{(\hom^t)^\#\ar[r]^+\ar@{=}[d]&\C{Hopair}(\C{B})\ar[r]^{\tilde{p}_{\C{B}}}\ar[d]^{\varphi_*}&\C{Pair}(\C{A})\ar@{=}[d]\\
(\hom^t)^\#\ar[r]^+&\C{Hopair}(\C{C})\ar[r]^{\tilde{p}_{\C{C}}}&\C{Pair}(\C{A})}
\end{equation*}
In particular $\varphi_*$ is an equivalence of categories.

The equivalence $\varphi_*$ in diagram (b) determines isomorphisms
between the cone functors defined by $\C{B}$ and by $\C{C}$ since
cone functors are restriction of $\hom$-functors in the category of
homotopy pairs, see (\ref{conenew}). In particular $\varphi_*$
restricts to the full subcategories given by maps with representable
associated cone functor and the functors $\zeta$ defined in Remark
\ref{defino} fit into a diagram commuting up to equivalence
$$\xymatrix{\C{Hopair}^{rep}(\C{B})\ar[rd]^\zeta\ar[dd]_\sim^{\varphi_*}&\\&\C{cand}(\C{A},t)\\
\C{Hopair}^{rep}(\C{C})\ar[ru]_\zeta&}$$ This and Lemma \ref{unacar}
prove that $(\C{B},s)$ admits a pretriangulated track category
structure provided $(\C{C},r)$ does. The converse follows by
symmetry.

Suppose now that $(\C{B},s)$ admits a triangulated track category
structure. Then by Lemma \ref{otracar} there is a full subcategory
$\C{dist}(\C{B},s)$ of $\C{candt}(\C{B},s)$ such that, over
$\C{dist}(\C{B},s)$, the functor
$\tilde{p}_{\C{B}}\bar{\varrho}_{(\C{B},s)}$ is surjective on
objects and $\bar{\varrho}_{(\C{B},s)}$ is full. If we define
$\C{dist}(\C{C},r)$ as the full subcategory of $\C{candt}(\C{C},r)$
given by the images by $(\varphi,\xi)_*$ of the objects in
$\C{dist}(\C{B},s)$ then the commutativity of (a) and (b) and the
fact that $\varphi_*$ is an equivalence imply that, over
$\C{dist}(\C{C},r)$, the functor
$\tilde{p}_{\C{C}}\bar{\varrho}_{(\C{C},r)}$ is surjective on
objects and $\bar{\varrho}_{(\C{C},r)}$ is full, therefore again by
Lemma \ref{unacar} we have that $(\C{C},r)$ admits a triangulated
track category structure. The converse follows also by symmetry.
\end{proof}

\section{The octahedral axiom}\label{toa}

This section contains the results which allow us to prove the
octahedral axiom for cohomologically  triangulated categories. 
In this section $\C{A}$ is always an additive category and
$t\colon\C{A}\r\C{A}$ an additive equivalence.

\begin{prop}\label{octa1}
Let $(\C{B},s)$ be a good translation track category representing
$\nabla\in H^*(\C{A},t)$. Suppose in addition that
$(\C{A},t,\nabla)$ is cohomologically pretriangulated. Consider two
composable morphisms of track triangles as follows
$$\xymatrix{A\ar[r]^f\ar@{=}[d]&B\ar[r]^{i_f}_(.7){\;\;}="a"\ar[d]|g^(.7){\;\;}="b"&C_f\ar[r]^{q_f}_(.7){\;\;}="c"\ar[d]|{\bar{g}}^(.7){\;\;}="d"&sA\ar@{=}[d]\\
A\ar[r]|{gf}\ar[d]_f&C\ar@{=}[d]^(.7){\;\;}="f"\ar[r]|{i_{gf}}_(.7){\;\;}="e"&C_{gf}\ar[d]|{\bar{f}}^(.7){\;\;}="h"\ar[r]|{q_{gf}}_(.7){\;\;}="g"&sA\ar[d]^{sf}\\
B\ar[r]_{g}&C\ar[r]_{i_{g}}&C_g\ar[r]_{q_g}&sB\ar@{=>}"a";"b"^{F}\ar@{=>}"c";"d"^{G}\ar@{=>}"e";"f"^{F'}\ar@{=>}"g";"h"^{G'}}$$
then $$C_f\st{p(\bar{g})}\To C_{gf}\st{p(\bar{f})}\To C_g\st{(tp(i_f))p(q_g)}\To
tC_f$$ is a small $\nabla$-triangle.
\end{prop}

For the proof of this proposition we need the exact sequence in
the following lemma.

\begin{lem}\label{ifni}
Let $(\C{B},s)$ be a good translation track category over 
$(\C{A},t)$ 
and let $t^{-1}$ be an inverse
equivalence for $t$ with $\mu\colon tt^{-1}\cong 1_\C{A}$ a natural
isomorphism. Given morphisms $A\st{f}\r B\st{g}\r C$ in $\C{B}$
there is an exact sequence of $\C{A}$-modules
$$\xymatrix@C=17pt{\cone_{gf}t^{-1}\ar[rr]^<(.4){[1_A,g,0^\vc_{gf}]^* t^{-1}}&&\cone_f t^{-1}\ar[r]^<(.25){\partial}&\cone_g
\ar[rr]^{[f,1_C,0^\vc_{gf}]^*}&&\cone_{gf}
\ar[rr]^{[1_A,g,0^\vc_{gf}]^*}&&\cone_f }$$
where $\partial$ is the composite
$$\partial\colon\cone_ft^{-1}\st{\tilde{i}_ft^{-1}}\r\hom_\C{A}(B,t^{-1})\st{t}\cong\hom_\C{A}(tB,tt^{-1})\st{\mu_*}\cong\hom_\C{A}(tB,-)\st{\tilde{q}_g}\r\cone_g.$$
\end{lem}

Compare \cite{ttc} 3.6.

\begin{proof}[Proof of Proposition \ref{octa1}]
This proof is essentially and adaptation of \cite{ttc} 9 to the more general case considered in this paper.
Consider the following diagram
\begin{equation*}\tag{a}
\rotatebox{90}{\xymatrix@R=20pt@C=15pt{
\cone_{gf}t^{-1}\ar[r]^{[1_A,g,0^\vc_{gf}]^*t^{-1}}\ar[dd]_\cong^{\chi^{-1}_{\set{i_{gf},H^{gf}_0}}t^{-1}}&
\cone_ft^{-1}\ar[rr]^\partial\ar[dd]_\cong^{\chi^{-1}_{\set{i_f,H_0^f}}t^{-1}}\ar[rd]^<(.5){\tilde{i}_f
t^{-1}}&&
\cone_g\ar[r]^{[f,1_C,0^\vc_{gf}]^*}\ar[dddddd]^\cong_{\chi^{-1}_{\set{i_g,H^g_0}}}&
\cone_{gf}\ar[r]^{[1_A,g,0^\vc_{gf}]}\ar[dddddd]^\cong_{\chi^{-1}_{\set{i_{gf},H^{gf}_0}}}&
\cone_f\ar[dddddd]_{\chi^{-1}_{\set{i_f,H^f_0}}}^\cong\\
\ar@{}[r]|<(.7){\spadesuit}&\ar@{}[r]|{\clubsuit}&\hom_\C{A}(B,t^{-1})\ar[dd]_\cong^t&&&\\
\hom_\C{A}(C_{gf},t^{-1})\ar[r]^{p(\bar{g})^*}\ar[dd]_\cong^t&
\hom_\C{A}(C_f,t^{-1})\ar[ru]_{p(i_f)^*}\ar[dd]_\cong^t&&&&\\
&&\hom_\C{A}(tB,tt^{-1})\ar[dd]_\cong^{\mu_*}&&&\\
\hom_\C{A}(tC_{gf},tt^{-1})\ar[r]^<(.15){(tp(\bar{g}))^*}\ar[dd]_\cong^{\mu_*}&
\hom_\C{A}(tC_f,tt^{-1})\ar[ru]_{(tp(i_f))^*}\ar[dd]_\cong^{\mu_*}&\ar@{}[r]|{\clubsuit}&\ar@{}[r]|{\spadesuit}&\ar@{}[r]|{\spadesuit}&\\
&&\hom_\C{A}(tB,-)\ar[rd]^{p(q_g)^*}\ar@/_5pt/[ruuuuu]^ {\tilde{q}_g}&&&\\
\hom_\C{A}(tC_{gf},-)\ar[r]^{(tp(\bar{g}))^*}&
\hom_\C{A}(tC_f,-)\ar[ru]^{(tp(i_f))^*}\ar[rr]_{((tp(i_f))p(q_g))^*}&&
\hom_\C{A}(C_g,-)\ar[r]^{p(\bar{f})^*}&
\hom_\C{A}(C_{gf},-)\ar[r]^{p(\bar{g})^*}&
\hom_\C{A}(C_f,-)
}}
\end{equation*}
This diagram is commutative. The subdiagrams labeled with a
$\spadesuit$ commute by Lemma \ref{pepito}, and the subdiagrams
labeled with a $\clubsuit$ by Remark \ref{conmuta}. The
commutativity of the other subdiagrams is trivial.

Consider now the diagram
\begin{equation*}\tag{b}
\rotatebox{90}{ \xymatrix@R=20pt@C=15pt{
\cone_{gf}t^{-1}\ar[r]^{[1_A,g,0^\vc_{gf}]^*t^{-1}}\ar@{<-}[dd]_\cong^{\chi_{\set{i_{gf},H^{gf}_0}}t^{-1}}&
\cone_ft^{-1}\ar[rr]^\partial\ar@{<-}[dd]_\cong^{\chi_{\set{i_f,H_0^f}}t^{-1}}\ar[rd]^<(.5){\tilde{i}_f
t^{-1}}&&
\cone_g\ar[r]^{[f,1_C,0^\vc_{gf}]^*}\ar@{<-}[dddddd]^<(.2){[i_f,i_{gf},F^\vi]^*}&
\cone_{gf}\ar[r]^{[1_A,g,0^\vc_{gf}]}\ar@{<-}[dddddd]^\cong_{\chi_{\set{i_{gf},H^{gf}_0}}}&
\cone_f\ar@{<-}[dddddd]_{\chi_{\set{i_f,H^f_0}}}^\cong\\
&&\hom_\C{A}(B,t^{-1})\ar[dd]_\cong^t&&&\\
\hom_\C{A}(C_{gf},t^{-1})\ar[r]^{p(\bar{g})^*}\ar@{<-}[dd]_\cong^{t^{-1}}&
\hom_\C{A}(C_f,t^{-1})\ar[ru]_{p(i_f)^*}\ar@{<-}[dd]_\cong^{t^{-1}}&&&&\\
&&\hom_\C{A}(tB,tt^{-1})\ar[dd]_\cong^{\mu_*}&&&\\
\hom_\C{A}(tC_{gf},tt^{-1})\ar[r]^<(.15){(tp(\bar{g}))^*}\ar@{<-}[dd]_\cong^{\mu_*^{-1}}&
\hom_\C{A}(tC_f,tt^{-1})\ar[ru]_{(tp(i_f))^*}\ar@{<-}[dd]_\cong^{\mu_*^{-1}}&\ar@{}[r]|{\clubsuit}&\ar@{}[r]|{\spadesuit}&&\\
&&\hom_\C{A}(tB,-)\ar@/_10pt/[ruuuuu]^{\tilde{q}_g}&&&\\
\hom_\C{A}(tC_{gf},-)\ar[r]^{(tp(\bar{g}))^*}&
\hom_\C{A}(tC_f,-)\ar[ru]^{(tp(i_f))^*}\ar[rr]_{\tilde{q}_{\bar{g}}}&&
\cone_{\bar{g}}\ar[r]^<(.15){\tilde{i}_{\bar{g}}}&
\hom_\C{A}(C_{gf},-)\ar[r]^{p(\bar{g})^*}& \hom_\C{A}(C_f,-) }}
\end{equation*}
We claim that this diagram also commutes. Actually we only have to check the commutativity of $\clubsuit$ and $\spadesuit$.
The commutativity of $\clubsuit$ follows from the next equalities:
\begin{eqnarray*}
[f,1_C,0^\vc_{gf}]^*\tilde{q}_{\bar{g}}(h)&=&[f,1_C,0^\vc_{gf}]^*\set{0,\sigma_0(h)}\\
&=&\set{0i_{gf},\sigma_0(h)\vc(0F^\vi)}\\
&=&\set{0,\sigma_0(h)i_f}\\
&=&\set{0,\sigma_0(h(tp(i_f)))}\\
&=&\tilde{q}_g(h(tp(i_f)))\\
&=&\tilde{q}_g(tp(i_f))^*(h).
\end{eqnarray*}

Given now $\set{\alpha,\hat{\alpha}}\in\cone_{\bar{g}}(X)$ we have
that
$[f,1_C,0^\vc_{gf}]^*[i_f,i_{gf},F^\vi]^*\set{\alpha,\hat{\alpha}}$
is represented by
$$\xymatrix{&C_f\ar@/^30pt/[rr]^0_{\;}="b"\ar[r]|{\bar{g}}^<(1){\;\;\;\;\;}="a"&C_{gf}\ar[r]_\alpha&X\\
A\ar[r]_f&B\ar[u]^{i_f}_<(.7){\;}="d"\ar[r]_g^<(.7){\;}="c"&C\ar[u]_{i_{gf}}&\ar@{=>}"a";"b"_{\tilde{\alpha}}\ar@{=>}"c";"d"_{F^\vi}}$$
This composite track coincides clearly with
$$\xymatrix{&C_f\ar@/^30pt/[rr]^0_{\;}="b"\ar[r]|{\bar{g}}^<(1){\;\;\;\;\;}="a"&C_{gf}\ar[r]_\alpha&X\\
A\ar@/^15pt/[ru]^0_{\;}="e"\ar[r]_f^{\;}="f"&B\ar[u]|{i_f}_<(.7){\;}="d"\ar[r]_g^<(.7){\;}="c"&C\ar[u]_{i_{gf}}&\ar@{=>}"a";"b"_{\tilde{\alpha}}\ar@{=>}"c";"d"_{F^\vi}
\ar@{=>}"f";"e"_{H^f_0}}$$
This other one is obviously the same as
$$\xymatrix{&C_f\ar[r]^{\bar{g}}^<(1){\;\;\;\;\;}="a"&C_{gf}\ar[r]^\alpha&X\\
A\ar@/^15pt/[ru]^0_{\;}="e"\ar[r]_f^{\;}="f"&B\ar[u]|{i_f}_<(.7){\;}="d"\ar[r]_g^<(.7){\;}="c"&C\ar[u]_{i_{gf}}&\ar@{=>}"c";"d"_{F^\vi}
\ar@{=>}"f";"e"_{H^f_0}}$$
The conditions defining track triangle morphisms show that this composite track coincides with
$$\xymatrix{A\ar[r]_{gf}^(1.06){\;\;\;\;}="b"\ar@/^30pt/[rr]^0_{\;}="a"&C\ar[r]_{i_gf}&C_{gf}\ar[r]^\alpha&X\ar@{=>}"b";"a"_{H^{gf}_0}}$$
and this composite track represents
$\chi_{\set{i_{gf},H^{gf}_0}}\tilde{i}_{\bar{g}}\set{\alpha,\hat{\alpha}}\in\cone_{gf}(X)$,
so $\spadesuit$ commutes.

The commutativity of (b), Lemma \ref{ifni}, (\ref{exaa}) and the
five lemma show that $[i_f,i_{gf},F^\vi]^*$ is an isomorphism.
Moreover, since $(\C{A},t,\nabla)$ is cohomologically
pretriangulated there is a track triangle
$$\xymatrix{C_f\ar@/^30pt/[rr]^0_{\;}="a"\ar[r]_{\bar{g}}&C_{gf}\ar@/_30pt/[rr]_0^{\;}="c"
\ar[r]^<(-.35){\;}="b"|{i_{\bar{g}}}&
C_{\bar{g}}\ar[r]_<(-0.25){\;}="d"|{q_{\bar{g}}}\ar@/^30pt/[rr]^0_{\;}="f"&sC_f\ar[r]_{s(\bar{g})}^<(-0.25){\;}="e"&
sC_{gf}\ar@{=>}"b";"a"_{H_0^{\bar{g}}}\ar@{=>}"d";"c"^{H_{1}^{\bar{g}}}\ar@{=>}"e";"f"_{H_{2}^{\bar{g}}}}.$$
and a commutative diagram
\begin{equation*}\tag{c}
\rotatebox{90}{\xymatrix{&&\hom_\C{A}(C_{\bar{g}},-)\ar[rd]^{p(i_{\bar{g}})^*}\ar[dd]^{\chi_{\set{i_{\bar{g}},H_0^{\bar{g}}}}}_\cong&&\\
\hom_\C{A}(tC_{gf},-)\ar[r]^{(tp(\bar{g}))^*}&\hom_\C{A}(tC_f,-)\ar[ru]^{p(q_{\bar{g}})^*}\ar[rd]_{\tilde{q}_{\bar{g}}}&&
\hom_\C{A}(C_{gf},-)\ar[r]^{p(\bar{g})^*}&\hom_\C{A}(C_f,-),\\&&\cone_{\bar{g}}\ar[ru]_{\tilde{i}_{\bar{g}}}&}}
\end{equation*}
see Remark \ref{conmuta}. Therefore composing (a), (b) and (c) and applying Yoneda's lemma we obtain an isomorphism of small candidate triangles
$$\xymatrix@C=50pt{C_f\ar[r]^{p(\bar{g})}\ar@{=}[d]&C_{gf}\ar[r]^{p(\bar{f})}\ar@{=}[d]&C_g\ar[r]^{(tp(i_f))p(q_g)}\ar[d]_\cong^l&tC_f\ar@{=}[d]\\
C_f\ar[r]_{p(\bar{g})}&C_{gf}\ar[r]_{p(i_{\bar{g}})}&C_{\bar{g}}\ar[r]_{p(q_{\bar{g}})}&tC_f}$$
where $l\colon C_g\r C_{\bar{g}}$ represents
$$\chi^{-1}_{\set{i_{g},H^{g}_0}}[i_f,i_{gf},F^\vi]^*\chi_{\set{i_{\bar{g}},H^{\bar{g}}_0}}\colon
\hom_\C{A}(C_{\bar{g}},-)\cong\hom_\C{A}(C_g,-).$$ Since the lower
candidate triangle is a small $\nabla$-triangle, the upper one is
also a small $\nabla$-triangle by Theorem \ref{elt} (pretriangulated
case) and (Tr0), hence the proof is finished.
\end{proof}

\bibliographystyle{alpha}
\bibliography{Fernando}


\title{{[V]}\\\hspace{5pt}\\Triangulated cohomology classes}

\author{H.-J. Baues and F. Muro}%

\blfootnote{\textit{$2000$ Mathematics Subject Classification:} 18E30, 18G60}%
\blfootnote{\textit{Key words and phrases:} triangulated categories, pretriangulated categories, cohomology of categories}%
\blfootnote{The second author was partially supported by the MEC
grant MTM2004-01865 and postdoctoral fellowship EX2004-0616}

\maketitle

\begin{abstract}
A cohomologically triangulated category is an additive category
$\C{A}$ together with a translation functor $t$ and a cohomology
class $\nabla\in H^3(\C{A},t)$. We show that $\nabla$ satisfies
certain purely cohomological conditions if and only if $\nabla$ is
represented by a triangulated track category.  This yields a purely
cohomologically characterization of triangulated cohomology classes.
\end{abstract}


\section*{Introduction}

In \cite{ccctc1} we define a cohomologically triangulated category
as a triple $(\C{A},t,\nabla)$ given by an additive category
$\C{A}$, an additive equivalence $t\colon\C{A}\r\C{A}$ and a
cohomology class $$\nabla\in H^3(\C{A},t)$$ in the translation
cohomology of the pair $(\C{A},t)$ in the sense of \cite{ptc} such
that $\nabla$ is represented by a good translation track category
satisfying two axioms (TTr1) and (TTr2). For this we use the main
result of \cite{ptc}, which says that $3$-dimensional translation
cohomology classes classify translation track
categories.

The aim of this paper is to characterize cohomologically
triangulated categories in a purely cohomological way, i. e. without
invoking the fact that translation cohomology classes are
represented by translation track categories. This is achieved in
Theorem \ref{mamona} as a consequence of Theorems \ref{ñi} and
\ref{vanish}. There we obtain vanishing conditions on some new
cohomology classes associated to $\nabla$ which imply that
$(\C{A},t,\nabla)$ is a cohomologically triangulated category. These
new cohomology classes are obtained from $\nabla$ by using spectral
sequences constructed in Section \ref{ss}, see Section \ref{tmss}.

With this new characterization one can forget about the techniques
of track categories in the previous papers since the axioms (TTr1)
and (TTr2) are completely encoded in the translation cohomology
class $\nabla$. The original definition in terms of track categories
is the most convenient one to check that stable homotopy categories
give rise to cohomologically triangulated categories. However the
cohomological approach will be more convenient for a further study
of cohomologically triangulated categories, like their
localizations, categories of diagrams over them, $t$-structures,
monoidal structures, etc.

\subsection*{Notation}

In this paper the arrow $\r$ is used for morphisms in ordinary
categories, $1$-cells in $2$-categories and functors. We use $\rr$
for $2$-cells in $2$-categories and natural transformations, and
$\pf$ for pseudofunctors. Identity morphisms or $1$-cells are
denoted by $1$, with the object as a subscript when necessary. For
the identity $2$-cells and natural transformations we use the
symbol $0^\vc$. Horizontal composition in $2$-categories is
denoted by juxtaposition, as the composition law in ordinary
categories, while we use the symbol $\vc$ for vertical
composition. The vertical inverse of an invertible $2$-cell or
natural transformation $\alpha$ is $\alpha^\vi$.

\section{Translation categories}\label{tcs}

A \emph{translation category}\index{translation category}
$(\C{A},t)$ is a category $\C{A}$ with a zero object $*$ together
with an endofunctor $t\colon\C{A}\r\C{A}$ preserving the zero
object $*$. Notice that (pre)triangulated categories 
yield examples of translation categories.

It is well known that the category $\C{Cat}$ of all small categories
with zero object is a $2$-category with zero-preserving functors as
$1$-cells and natural transformations as $2$-cells. Small
translation categories also form a $2$-category $\C{TCat}$. A
morphism (or $1$-cell) of translation categories
$(\varphi,\alpha)\colon(\C{A},t)\r(\C{B},s)$ is a functor
$\varphi\colon\C{A}\r\C{B}$ preserving zero objects together with a
natural isomorphism $\alpha\colon s\varphi\rr\varphi t$, i. e. the
diagram
$$\xymatrix{\C{A}\ar[r]^t\ar[d]_\varphi&\C{A}\ar[d]^\varphi\\\C{B}\ar[r]_s&\C{B}}$$
commutes up to the natural isomorphism $\alpha$. Composition of morphisms is defined by
$$(\varphi,\alpha)(\psi,\beta)=(\varphi\psi,(\varphi\beta)\vc(\alpha\psi)).$$
A $2$-cell
$\gamma\colon(\varphi,\alpha)\rr(\psi,\beta)$ between two
morphisms $(\varphi,\alpha),(\psi,\beta)
\colon(\C{A},t)\r(\C{B},s)$ in $\C{TCat}$ 
is a natural transformation
$\gamma\colon\varphi\rr\psi$ such that
$$\beta\vc(s\gamma)=(\gamma t)\vc\alpha.$$
Vertical and horizontal composition of $2$-cells in $\C{TCat}$ are
given by the corresponding composition laws between natural
transformations in $\C{Cat}$.

Ordinary categories $\C{C}$ with zero object can be regarded as
translation categories by using identity functors
$\C{C}=(\C{C},1_{\C{C}})$.
This defines an inclusion of the $2$-category $\C{Cat}$ into
$\C{TCat}$. The $2$-category $\C{TCat}$ has products
$(\C{A},t)\times(\C{B},s)=(\C{A}\times\C{B},t\times s)$.

Morphisms between translation categories and $2$-cells between
them are also defined for non-small translation categories exactly
in the same way. 

We shall need the cohomology of translation categories, in particular, the cohomology of a triangulated category. For this we recall first the concept of
cohomology of categories, see \cite{csc}.

\section{Cohomology of categories}\label{cc}

Given a category $\C{C}$ a
\emph{$\C{C}$-bimodule}\index{$\mathbf{C}$-bimodule} $D$ is a
functor $D\colon\C{C}^\op\times\C{C}\r\C{Ab}$ where $\C{Ab}$ is
the category of abelian groups. We write $f^*=D(f,1)$ and
$g_*=D(1,g)$ for the induced homomorphisms. The cohomology
$H^*(\C{C},D)$ of $\C{C}$ with coefficients in $D$ is defined by
the cochain complex $F^*(\C{C},D)$ concentrated in dimensions
$\geq 0$. This complex is given in dimension $n$ by the following
product indexed by all sequences of morphisms of length $n$ in
$\C{C}$,
\begin{equation}\label{producto}
F^n(\C{C},D)=\prod_{X_0\st{\sigma_1}\l\cdots\st{\sigma_n}\l X_n}D(X_n,X_0).
\end{equation}
Here we understand that a sequence of length $0$ is the identity morphism of an object $X$ in
$\C{C}$.
The coordinate of $c\in F^n(\C{C},D)$ in
$X_0\st{\sigma_1}\l\cdots\st{\sigma_n}\l X_n$ will be
denoted by $c(\sigma_1,\dots,\sigma_n)\in D(X_n,X_0)$. The differential $\delta$
is defined as
\begin{eqnarray*}
  \delta(c)(\sigma_1,\dots,\sigma_{n+1}) &=& {\sigma_1}_*c(\sigma_2,\dots,\sigma_{n+1}) \\
   && +\sum_{i=1}^{n}(-1)^ic(\sigma_1,\dots,\sigma_i\sigma_{i+1},\dots,\sigma_{n+1}) \\
   && +(-1)^{n+1}\sigma_{n+1}^*c(\sigma_1,\dots,\sigma_n)
\end{eqnarray*}
over an $n$-cochain $c$ for $n\geq 1$, and
$\delta(c)(\sigma)=\sigma_*c(X)-\sigma^*c(Y)$ for $n=0$ and
$\sigma\colon X\r Y$.

In order $H^*(\C{C},D)$ to be a set we assume that $\C{C}$ is
small, one can use this notation accordingly if
$\C{C}$ is small with respect to some fixed universe.

This is a particular case of the Baues-Wirsching cohomology of $\C{C}$, see \cite{csc}. It also coincides with the Hochschild-Mitchell cohomology of
the free ringoid $\Z\C{C}$ generated by the category $\C{C}$, see \cite{rso}.

A functor $\varphi\colon\C{D}\r\C{C}$ and a natural transformation
$\lambda\colon D\rr E$ induce the following cochain homomorphisms
where $X_n$ is the source of $\sigma_n$ and $X_0$ is the target of
$\sigma_1$,
$$\varphi^*\colon F^*(\C{C},D)\To F^*(\C{D},D(\varphi,\varphi)),\;\; \varphi^*(c)(\sigma_1,\dots,\sigma_n)=
c(\varphi(\sigma_1),\dots,\varphi(\sigma_n));$$
\begin{equation}\label{lsa}
\lambda_*\colon F^*(\C{C},D)\To F^*(\C{C},E),\;\;
\lambda_*(c)(\sigma_1,\dots,\sigma_n)=
\lambda_{X_n,X_0}(c(\sigma_1,\dots,\sigma_n));
\end{equation}
and the corresponding homomorphisms in cohomology,
$$\varphi^*\colon H^*(\C{C},D)\To H^*(\C{D},D(\varphi,\varphi)),$$
$$\lambda_*\colon H^*(\C{C},D)\To H^*(\C{C},E).$$

A cochain $c\in F^n(\C{C},D)$ is said to be
\emph{normalized}\index{normalized cochain} (at identities) if
$c(\sigma_1,\dots,\sigma_n)=0$ provided $\sigma_i$ is an identity
morphism for some $1\leq i\leq n$. The inclusion of the subcomplex
$\bar{F}^*(\C{C},D)\subset F^*(\C{C},D)$ of normalized cochains
induces an isomorphism in cohomology, see \cite{ccglg}. Moreover,
the induced cochain homomorphisms above restrict to normalized
cochains.

In this paper we are mainly interested in categories $\C{C}$ with a
zero object $*$. In this case morphism sets $\C{C}(X,Y)$ in $\C{C}$
are pointed by the zero morphisms $0_{X,Y}=0\colon X\r*\r Y$.
Suppose that the $\C{C}$ bimodule $D$ vanishes whenever the zero
object is in one side $D(X,*)=0=D(*,X)$. In this case we say that
$D$ is \emph{$*$-normalized}\index{normalized
bimodule@$*$-normalized bimodule}. A cochain $c\in F^n(\C{C},D)$ is
said to be \emph{$*$-normalized}\index{normalized
cochaio@$*$-normalized cochain} (i. e. normalized at identities and
zero morphisms) if $c(\sigma_1,\dots,\sigma_n)=0$ provided
$\sigma_i$ is an identity or a zero morphism for some
$1\leq i\leq n$. 
If $D$ is $*$-normalized the inclusion of the subcomplex
$\tilde{F}^*(\C{C},D)\subset F^*(\C{C},D)$ of $*$-normalized
cochains also induces an isomorphism in cohomology, see
\cite{ccglg}. Moreover, the induced cochain homomorphisms
$\varphi^*$ and $\lambda_*$ above restrict to $*$-normalized
cochains when $\varphi$ preserves the zero object and $E$ is
$*$-normalized.

A natural transformation $\alpha\colon\varphi\rr\psi$
between functors $\varphi,\psi\colon\C{D}\r\C{C}$ induces a cochain homomorphism
$$\alpha^*\colon F^*(\C{C},D)\To F^*(\C{D},D(\varphi,\psi))$$ defined as $$\alpha^*(c)(\sigma_1,\cdots,\sigma_n)=
{\alpha_{X_0}}_* (c(\varphi(\sigma_1),\dots,\varphi(\sigma_n))).$$
This is $F^*(\alpha,0^\vc_{D(\varphi,\psi)})$ in \cite{fcc} 4.1.
Here $X_0$ is the target of $\sigma_1$. The induced homomorphism
in cohomology is denoted in the same way
$$\alpha^*\colon H^*(\C{C},D)\To H^*(\C{D},D(\varphi,\psi)).$$
Notice that the cochain homomorphism defined by an identity natural
transformation $0^\vc_\varphi\colon\varphi\rr\varphi$ coincides with
the cochain homomorphism defined by the corresponding functor, i. e.
$(0^\square_\varphi)^*=\varphi^*.$ The natural transformation
$\alpha\colon\varphi\rr\psi$ also induces cochain homotopies given
by degree $-1$ homomorphisms
$$\alpha_\#,\alpha^\#\colon F^*(\C{C},D)\To F^*(\C{D},D(\varphi,\psi)),$$
with
\begin{equation}\label{homo11}
\delta\alpha_\#+\alpha_\#\delta=-(D(0^\square_\varphi,\alpha))_*\varphi^*+\alpha^*,
\end{equation}
\begin{equation}\label{homo12}
\delta\alpha^\#+\alpha^\#\delta=-(D(\alpha,0^\square_\psi))_*\psi^*+\alpha^*.
\end{equation}
They are defined by the following formulas
$$\alpha_\#(c)(\sigma_1,\dots,\sigma_n)=
\sum\limits_{i=0}^n(-1)^i{\alpha_{X_0}}_*(
c(\varphi(\sigma_1),\dots,\varphi(\sigma_i),1_{\varphi(X_i)},\varphi(\sigma_{i+1}),\dots,\varphi(\sigma_n))),$$
\begin{equation*}\label{almo}
\alpha^\#(c)(\sigma_1,\dots,\sigma_n)=\sum_{i=0}^n(-1)^i
c(\psi(\sigma_1),\dots,\psi(\sigma_i),\alpha_{X_i},\varphi(\sigma_{i+1}),\dots,\varphi(\sigma_n)).
\end{equation*}
Here $X_i$ is the source of $\sigma_i$ and/or the target of
$\sigma_{i+1}$. These homotopies are part of the $2$-functorial
properties of the Baues-Wirsching complex $F^*(\C{C},D)$ studied
in \cite{fcc}. More precisely, $\alpha_\#$ and $\alpha^\#$ are
$h_{(0^\vc_\varphi,\alpha)}$ and $h_{(\alpha,0^\vc_\psi)}$ in
\cite{fcc} 4.1, respectively. Notice that these induced cochain
homomorphisms and homotopies restrict always to the subcomplex of
normalized cochains, and also to the subcomplex of $*$-normalized
cochains if the functors $\varphi$ and $\psi$ preserve the zero
object. Moreover, as one can readily notice on the normalized
cochain complex the formulas $\alpha_\#=0$, $(0^\vc_\varphi)^\#=0$
and $\alpha^*=(D(0^\square_\varphi,\alpha))_*\varphi^*$ hold, in
particular by (\ref{homo11}) and (\ref{homo12}) we have that
\begin{equation}\label{homo13}
\delta\alpha^\#+\alpha^\#\delta=-(D(\alpha,0^\square_\psi))_*\psi^*+(D(0^\square_\varphi,\alpha))_*\varphi^*.
\end{equation}

\section{Translation cohomology}\label{tc}

Let $(\C{A},t)$ be a translation category in the sense of Section
\ref{tcs}, $D$ a $*$-normalized $\C{A}$-bimodule and
$\bar{t}\colon D\rr D(t,t)$ a natural transformation. Then the
pair $(t,\bar{t})$ is termed a
\emph{translation}\index{translation} on $(\C{A},D)$, see
\cite{ptc} 6.
The translation cohomology $H^*(t,\bar{t})$ of
$(t,\bar{t})$, introduced in \cite{ptc}, is the cohomology of the cochain complex $\tilde{F}(t,\bar{t})$ which is the homotopy fiber
of the cochain homomorphism
$$\bar{t}_*-t^*\colon\tilde{F}^*(\C{A},D)\To\tilde{F}^*(\C{A},D(t,t)).$$
We will sometimes refer to the translation cohomology $H^*(t,\bar{t})$ as the translation cohomology of the translation category $(\C{A},t)$
with coefficients in the natural transformation $\bar{t}$.
In dimension $n$
this homotopy fiber is
$$\tilde{F}^n(t,\bar{t})=\tilde{F}^n(\C{A},D)\oplus\tilde{F}^{n-1}(\C{A},D(t,t))$$
and the differential is the matrix
$$\left(%
\begin{array}{cc}
  \delta & 0 \\
  \bar{t}_*-t^* & -\delta \\
\end{array}%
\right).$$ In particular there is a long exact sequence $(n\in\Z)$
\begin{equation}\label{exacta}
\cdots\r H^n(t,\bar{t})\st{j}\r H^n(\C{A},D)\st{\bar{t}_*-t^*}\To
H^n(\C{A},D(t,t))\st{\partial}\r H^{n+1}(t,\bar{t})\r\cdots.
\end{equation}

Given 
a morphism $(\varphi,\alpha)\colon(\C{B},s)\r(\C{A},t)$ in $\C{TCat}$
there is an induced translation $(s,(\varphi,\alpha)^*\bar{t})$  on $(\C{B},D(\varphi,\varphi))$ with
\begin{equation}\label{indt}
(\varphi,\alpha)^*\bar{t}=(D(\alpha^\vi,\alpha))\vc(\bar{t}(\varphi^\op\times\varphi)),
\end{equation}
and an induced cochain homomorphism
$$(\varphi,\alpha)^*\colon\tilde{F}^*(t,\bar{t})\To\tilde{F}^*(s,(\varphi,\alpha)^*\bar{t})$$
given by the matrix
$$(\varphi,\alpha)^*=\left(\begin{array}{cc}
\varphi^*&0\\
  (D(\alpha^\vi,0^\vc_{\varphi s}))_*\alpha^\# & (D(\alpha^\vi,\alpha))_*\varphi^*
\end{array}%
\right).$$
Here we use (\ref{homo13}). The induced homomorphisms in translation cohomology will be denoted in the same way
$$(\varphi,\alpha)^*\colon H^*(t,\bar{t})\To H^*(s,(\varphi,\alpha)^*\bar{t}).$$
Similarly two natural transformations $\lambda\colon D\rr E$ and $\bar{s}\colon E\rr E(t,t)$ between $*$-normalized $\C{A}$-bimodules on $\C{A}$
such that $(\lambda(t^\op\times t))\vc\bar{t}=\bar{s}\vc\lambda$ give rise to a translation $(t,\bar{s})$ on $(\C{A},E)$
together with a cochain homomorphism
$$\lambda_*\colon\tilde{F}^*(t,\bar{t})\To\tilde{F}^*(t,\bar{s})$$
given by the matrix
\begin{equation}\label{lsa2}
\lambda_*=\left(\begin{array}{cc}
\lambda_*&0\\
0 & (\lambda(t^\op\times t))_*
\end{array}%
\right).
\end{equation}
The corresponding homomorphism in translation cohomology will also be denoted by
$$\lambda_*\colon H^*(t,\bar{t})\To H^*(t,\bar{s}).$$
These homomorphisms can be used to describe the simple version of
the functorial properties of translation cohomology, see
\cite{ptc} 6.2.

Two morphisms
$(\varphi,\alpha),(\psi,\beta)\colon(\C{B},s)\r(\C{A},t)$ in $\C{TCat}$ also induce
a translation $(s,[(\varphi,\alpha),(\psi,\beta)]^*\bar{t})$ on $(\C{B},D(\varphi,\psi))$ given by
\begin{equation}\label{indt2}
[(\varphi,\alpha),(\psi,\beta)]^*\bar{t}=(D(\alpha^\vi,\beta))\vc(\bar{t}(\varphi^\op\times\psi)).
\end{equation}
Notice that
\begin{equation}\label{diag}
[(\varphi,\alpha),(\varphi,\alpha)]^*\bar{t}=(\varphi,\alpha)^*\bar{t}. 
\end{equation}

In \cite{ptc} we consider cochains which are only normalized with
respect to identities since we are not considering categories with
zero objects there. More precisely, in the cited paper we define
the translation cochain complex $\bar{F}(t,\bar{t})$ as the
homotopy fiber of
$$\bar{t}_*-t^*\colon\bar{F}^*(\C{A},D)\To\bar{F}^*(\C{A},D(t,t)).$$
Clearly $\tilde{F}^*(t,\bar{t})\subset\bar{F}^*(t,\bar{t})$ and this
inclusion induces isomorphisms in all cohomology groups. This
follows from the exact sequence (\ref{exacta}), the five lemma, and
the fact that the inclusion of
the subcomplexes of normalized and $*$-normalized cochains
in ordinary cohomology of categories 
induce isomorphisms in cohomology, see Section \ref{cc}.

\begin{rem}\label{fun}
As we pointed out in \cite{ptc} 6.6 the translation cochain
complex is not a functor itself since, in general, composition of
induced cochain homomorphisms is only associative up to cochain
homotopy, compare with the proof of \cite{ptc} 6.4.
However it is
immediate to notice by using (\ref{lsa}) and (\ref{lsa2}) that
once we fix a translation category $(\C{A},t)$ then
$\bar{F}(t,-)$ is a functor on the category
whose objects are natural transformations of $\C{A}$-bimodules
$\bar{t}\colon D\rr D(t,t)$, and morphisms $\lambda\colon
\bar{t}\r \bar{s}$ with $\bar{s}\colon E\rr E(t,t)$ are natural
transformations $\lambda\colon D\rr E$ with
 $(\lambda(t^\op\times t))\vc\bar{t}=\bar{s}\vc\lambda$. The $*$-normalized translation cochain complex $\tilde{F}(t,-)$ is also a functor if we restrict
to $*$-normalized $\C{A}$-bimodules.
\end{rem}

\section{The cohomology of the triangle category}\label{cotri}\label{ct}\label{ttc}


The \emph{triangle translation category}\index{triangle
translation category} $(\V,sh^3)$ is a universal model of
candidate triangles in triangulated categories. The set of objects
of $\V$ is $\Z\cup\set{*}$, $*$ is the zero object of $\V$ and the
unique morphisms apart from identities and zeroes are $(n\in\Z)$
$$\imath=\imath_n\colon n-1\To n,$$
in particular $\imath^2=0$. The zero-object-preserving functor
(\emph{shift functor}\index{shift functor})
$$sh\colon\V\To\V$$ is defined by $sh(n)=n+1$ and $sh(\imath_n)=\imath_{n+1}$. This is
actually an isomorphism of categories
and $sh^3$ denotes its third power.

Let $\C{A}$ be a category with zero object. A functor
$X\colon\V\r\C{A}$ preserving the zero object coincides with a
\emph{cochain complex} in $\C{A}$, i. e. a diagram $(n\in\Z)$
$$\cdots\r X_{n-1}\st{d}\To X_n\st{d}\To X_{n+1}\r\cdots$$
with $d^2=0$. We write $X_n=X(n)$ for the values of $X$ on objects
in $\V$ and $d=d_n=d_n^X=X(\imath_n)$ for the morphisms induced by
the non-trivial morphisms in the category $\V$.



Let $(\C{A},t)$ be a translation category. A \emph{candidate
triangle}\index{triangle!candidate} $(X,\alpha)$ in $(\C{A},t)$ is
a cochain complex $X$ in $\C{A}$
$$\cdots\r X_{n-1}\st{d}\To X_n\st{d}\To X_{n+1}\r\cdots$$
together with a degree $3$ cochain isomorphism $\alpha\colon tX\cong
X$, i. e. a commutative diagram
$$\xymatrix{\cdots\ar[r]&tX_n\ar[r]^{td}\ar[d]^{\alpha_n}_\cong&tX_{n+1}\ar[r]\ar[d]^{\alpha_{n+1}}_\cong&\cdots\\
\cdots\ar[r]&X_{n+3}\ar[r]^d&X_{n+4}\ar[r]&\cdots}$$ Morphisms of
candidate triangles $f\colon (X,\alpha)\r (Y,\beta)$ are
cochain homomorphisms $f\colon X\r Y$ such that
$\beta(tf)=f\alpha$.

Let $\C{Cand}(\C{A},t)$ be the category of candidate triangles. If
$\C{A}$ is additive and $t$ is an additive equivalence 
then $\C{Cand}(\C{A},t)$ is equivalent to  the category of small
candidate triangles $\C{cand}(\C{A},t)$ as defined in \cite{ccctc1}
1, compare \cite{ttc}. Recall that an object in this category is a
diagram in $\C{A}$ of the form
$$A\st{f}\To B\st{i}\To C\st{q}\To tA$$
such that $if=0$, $qi=0$ and $(tf)q=0$. A morphism $k$ between small
candidate triangles is a commutative diagram
$$\xymatrix{A\ar[r]^f\ar[d]^{k_0}&B\ar[r]^i\ar[d]^{k_1}&C\ar[r]^q\ar[d]^{k_2}&tA\ar[d]^{tk_0}\\
\bar{A}\ar[r]^{\bar{f}}&\bar{B}\ar[r]^{\bar{i}}&\bar{C}\ar[r]^{\bar{q}}&t\bar{A}}$$
There
is an equivalence of categories,
\begin{equation}\label{sb}
\C{Cand}(\C{A},t)\To\C{cand}(\C{A},t)
\end{equation}
sending a candidate triangle $(X,\alpha)$ to the small candidate triangle
$$X_{-2}\st{d_{-1}}\To X_{-1}\st{d_0}\To X_0\st{\alpha^{-1}_{-2}d_1}\To tX_{-2}.$$
In order to construct the inverse equivalence
\begin{equation}\label{sb'}
\C{cand}(\C{A},t)\To\C{Cand}(\C{A},t)
\end{equation}
we need to choose an inverse equivalence
$t^{-1}\colon\C{A}\r\C{A}$ for $t$ and a natural isomorphism
$\mu\colon tt^{-1}\cong 1_\C{A}$. With these choices the
equivalence (\ref{sb'})
sends a small candidate triangle $A\st{f}\r B\st{i}\r C\st{q}\r
tA$ to the candidate triangle $(X,\alpha)$ with $X_{3n}=t^nC$,
$X_{3n+1}=t^{n+1}A$, $X_{3n+2}=t^{n+1}B$, $d_{3n}=t^ni$,
$d_{3n+1}=t^nq$ and $d_{3n+2}=t^{n+1}f$ for all $n\in\Z$,
$\alpha_n=1_{tX_n}$ if $n\geq -2$, and $\alpha_n=\mu_{X_{n+3}}$
for $n<-2$.

Notice that a candidate triangle $(X,\alpha)$ as above is the same
as a morphism of translation categories
$(X,\alpha)\colon(\V,sh^3)\r(\C{A},t)$ in the sense of Section
\ref{tcs}. Moreover, the category
$\C{Cand}(\C{A},t)$ of candidate triangles on $(\C{A},t)$
coincides with the category given by 
morphisms from $(\V,sh^3)$  to $(\C{A},t)$ and $2$-cells between
them in the $2$-category $\C{TCat}$ defined in Section \ref{tcs},
$$\C{Cand}(\C{A},t)=\C{TCat}((\V,sh^3),(\C{A},t)).$$




The translation cohomology of the triangle translation category
$(\V,sh^3)$ is relevant in this
paper. Below we identify its translation cochain
complex. For this we begin by considering the cohomology of the
category $\V$.

Let $D$ be a $*$-normalized $\V$-bimodule. The reader can readily
notice that there is an isomorphism
\begin{equation}\label{ident}
\tilde{F}^n(\V,D)\cong\prod_{m\in\Z}D(m,m+n),
\end{equation}
which we use as an identification, where the $m^\mathrm{th}$
factor of the product corresponds to the sequence
$$m+n\st{\imath_{m+n}}\longleftarrow\cdots\st{\imath_{m+1}}\longleftarrow m.$$
If we write $c(m)$ for the $m^\mathrm{th}$ coordinate of an
$n$-cochain $c$ in the product (\ref{ident}) then
$$(\delta (c))(m)={\imath_{m+n+1}}_*c(m)+(-1)^{n+1}\imath^*_{m+1}c(m+1).$$ The identification in (\ref{ident}) is clearly natural in
$D$, in particular, given a natural transformation $\bar{t}\colon
D\rr D(sh^3,sh^3)$ there is an identification
\begin{equation}\label{ident2}
\tilde{F}^n(sh^3,\bar{t})\cong\prod_{m\in\Z}D(m,m+n)\oplus\prod_{m\in\Z}D(m+3,m+n+2)
\end{equation}
for $n>0$ and
\begin{equation}\label{ident3}
\tilde{F}^0(sh^3,\bar{t})\cong\prod_{m\in\Z}D(m,m).
\end{equation}
The translation differential is given by
\begin{eqnarray*}(\delta(c,b))(m)&=&({\imath_{m+n+1}}_*c(m)+(-1)^{n+1}\imath^*_{m+1}c(m+1),\\
&&\bar{t}_{(m,m+n)}c(m)-c(m+3)-{\imath_{m+n}}_*b(m)-(-1)^{n}\imath^*_{m+1}b(m+1))
\end{eqnarray*}
for an $n$-cochain $(c,b)$ with $n>0$, and
$$(\delta c)(m)=({\imath_{m+1}}_*c(m)-\imath^*_{m+1}c(m+1),\bar{t}_{(m,m)}c(m)-c(m+3))$$
 for a $0$-cochain $c$.
These identifications in the translation cochain complex are
natural in $\bar{t}$, see Remark \ref{fun}.

\begin{lem}\label{sur}
If $\bar{t}$ is a natural isomorphism then the cochain
homomorphism
$\bar{t}_*-{sh^3}^*\colon\tilde{F}^*(\V,D)\r\tilde{F}^*(\V,D(sh^3,sh^3))$
is surjective.
\end{lem}

\begin{proof}
Actually there are levelwise sections
$\omega\colon\tilde{F}^n(\V,D(sh^3,sh^3))\r\tilde{F}^n(\V,D)$
$(n\in\Z)$ defined as follows $(m\in\Z,0\leq i\leq 2)$
$$\omega(a)(3m-i)=\left\{\begin{array}{ll}-\sum_{j=1}^m\bar{t}^{m-j}(a(3(j-1)-i)),&\text{if $m>0$;}\\0,&\text{if $m=0$;}\\
\sum_{j=m+1}^0\bar{t}^{m-j}(a(3(j-1)-i)),&\text{if
$m<0$.}\end{array}\right.$$
The reader can easily check that $(\bar{t}_*-{sh^3}^*)\omega=1$.
Notice however that $\omega$ does not define a cochain
homomorphism. 
\end{proof}

Consider now the cochain complex $\tilde{G}^*(\bar{t})$ defined by
$$\tilde{G}^n(\bar{t})=D(-2,n-2)\oplus D(-1,n-1)\oplus D(0,n),$$
with differential
$${\delta}=\left(\begin{array}{ccc}{\imath_{n-1}}_*&(-1)^{n+1}\imath^*_{-1}&0\\
0&{\imath_{n}}_*&(-1)^{n+1}\imath_0^*\\
(-1)^{n+1}\imath_1^*\bar{t}_{(0,n)}&0&{\imath_{n+1}}_*\end{array}\right)$$
on dimension $n$.

\begin{prop}\label{redu}
If $\bar{t}$ is a natural isomorphism then the cohomology of
$\tilde{G}^*(\bar{t})$ is isomorphic to the translation cohomology
$H^*(sh^3,\bar{t})$.
\end{prop}

\begin{proof}
By Lemma \ref{sur} we can replace the homotopy fiber of the
cochain homomorphism $\bar{t}_*-{sh^3}^*$ by its kernel $\ker
[\bar{t}_*-{sh^3}^*]$. More precisely, the obvious inclusion
$$\ker [\bar{t}_*-{sh^3}^*]\subset\tilde{F}^*(sh^3,\bar{t})$$
induces isomorphisms in all cohomology groups. The reader can
easily check that the cochain homomorphism
$$\kappa\colon\tilde{G}^*(\bar{t})\hookrightarrow\tilde{F}^*(\V,D)$$
defined by $\kappa(c_{-2},c_{-1},c_0)(3m-i)=\bar{t}^m(c_i)$ $(m\in\Z,
0\leq i\leq 2)$ is actually a kernel for $\bar{t}_*-{sh^3}^*$
hence the proposition follows.
\end{proof}

\begin{rem}\label{redu2}
Notice that in fact we have established an explicit
quasi-isomorphism in the proof of Proposition \ref{redu}
$$\upsilon=(\kappa,0)\colon\tilde{G}^*(\bar{t})\To\tilde{F}^*(sh^3,\bar{t}).$$
Moreover, the cochain complexes $\tilde{G}^*$ define a functor on
the category of natural transformations $\bar{t}\colon D\rr
D(sh^3,sh^3)$ between $*$-normalized $\V$-bimodules, the same as
$\tilde{F}^*(sh^3,-)$, see Remark \ref{fun}. Recall that morphisms
$\lambda\colon \bar{t}\r\bar{s}$ between natural transformations
$\bar{t}\colon D\rr D(sh^3,sh^3)$ and $\bar{s}\colon E\rr
E(sh^3,sh^3)$ in this category are natural transformations
$\lambda\colon D\rr E$ such that
$$(\lambda({sh^3}^\op\times sh^3))\vc\bar{t}=\bar{s}\vc\lambda.$$
The induced homomorphism
$$\lambda_*\colon\tilde{G}^*(\bar{t})\To\tilde{G}^*(\bar{s})$$
in dimension $n$ is given by the matrix
$$\lambda_*=\left(\begin{array}{ccc}\lambda_{(-2,n-2)}&0&0\\0&\lambda_{(-1,n-1)}&0\\0&0&\lambda_{(0,n)}\end{array}\right).$$
Therefore the quasi-isomorphism $\upsilon$ and hence the isomorphism
in Proposition \ref{redu} are natural in $\bar{t}$.
\end{rem}

Unfortunately $\upsilon$ in Remark \ref{redu2} is not a chain
homotopy equivalence, so we have not a direct way of representing a
cocycle in $\tilde{F}^*(sh^3,\bar{t})$ by a cocycle coming from
$\tilde{G}^*(\bar{t})$. The following lemma solves this problem. 

\begin{lem}\label{peque}
If $\bar{t}\colon D\rr D(sh^3,sh^3)$ is a natural isomorphism
between $*$-normalized $\V$-bimodules and
$(c,b)\in\tilde{F}^n(sh^3,\bar{t})$ is a cocycle
 then
$$(c(-2),c(-1),c(0)-D(\imath_1,1_n)b(-2))\in\tilde{G}^n(\bar{t})$$ is a
cocycle whose image by $\upsilon$ in Remark \ref{redu2} represents
the same cohomology class as $(c,b)$.
\end{lem}

\begin{proof}
With the notation of the proof of Lemma \ref{sur} it is easy to
check that
$$\upsilon(c(-2),c(-1),c(0)-D(\imath_1,1_n)b(-2))=(c,b)-\delta(\omega(b),0),$$
hence this lemma follows.
\end{proof}

We now consider a special example of a natural transformation as
in Proposition \ref{redu}. Suppose that $(\C{A},t)$ is a
translation category given by an additive category $\C{A}$ and an
additive equivalence of categories $t\colon\C{A}\r\C{A}$. We then
have the abelian group $\hom_\C{A}(X,Y)$ of morphisms in $\C{A}$.
Let $\hom^t$ be the
$\C{A}$-bimodule defined by $t$ 
in $\C{A}$ as follows
\begin{equation}\label{homt}
\hom^t=\hom_\C{A}(t,-),
\end{equation}
and let
\begin{equation}\label{tbar}
\bar{t}\colon\hom^t=\hom_\C{A}(t,-)\rr\hom_\C{A}(t^2,t)=\hom^t(t,t)
\end{equation}
be given by $(-1)t$. By (\ref{indt2}) a candidate triangle
$(X,\alpha)$ in $(\C{A},t)$ yields a natural transformation
$(X,\alpha)^*\bar{t}$ which can be used as a coefficient object
for the computation of translation cohomology groups of
$(\V,sh^3)$. The $3$-dimensional abelian group of the cochain
complex $\tilde{G}^*((X,\alpha)^*\bar{t})$ is
$$\tilde{G}^3((X,\alpha)^*\bar{t})=\hom_\C{A}(tX_{-2},X_1)\oplus\hom_\C{A}(tX_{-1},X_2)\oplus\hom_\C{A}(tX_0,X_3).$$
There is a canonical element in this abelian group determined by
the cochain isomorphism $\alpha$, namely
$$(\alpha_{-2},-\alpha_{-1},\alpha_0)\in\tilde{G}^3((X,\alpha)^*\bar{t}).$$

\begin{lem}
The $3$-cochain $(\alpha_{-2},-\alpha_{-1},\alpha_0)$ is indeed a
cocycle.
\end{lem}

The cohomology class represented by this cochain will be denoted
by
\begin{equation}\label{wi} \idd(X,\alpha)\in
H^3(sh^3,(X,\alpha)^*\bar{t})
\end{equation}
and will be called the \emph{class of
identities}\index{identities!class of -}.

\begin{rem}\label{cla1}
If $(X,\alpha)$ is the candidate triangle corresponding to the
small candidate triangle $$A\st{f}\To B\st{i}\To C\st{q}\To tA$$
under the  equivalence  (\ref{sb'}) then the cochain complex
$\tilde{G}^*((X,\alpha)^*\bar{t})$ is given as follows $(n\geq 0)$
\begin{small}
$$\tilde{G}^{3n+i}((X,\alpha)^*\bar{t})=\left\{%
\begin{array}{ll}
    \hom_\C{A}(tA,t^nA)\oplus\hom_\C{A}(tB,t^nB)\oplus\hom_\C{A}(tC,t^nC), & \hbox{if $i=0$;} \\
    &\\
    \hom_\C{A}(tA,t^nB)\oplus\hom_\C{A}(tB,t^nC)\oplus\hom_\C{A}(tC,t^{n+1}A), & \hbox{if $i=1$;} \\
    &\\
    \hom_\C{A}(tA,t^nC)\oplus\hom_\C{A}(tB,t^{n+1}A)\oplus\hom_\C{A}(tC,t^{n+1}B), & \hbox{if $i=2$.} \\
\end{array}%
\right.    $$
\end{small}
Moreover, the class of identities
$\idd(X,\alpha)$ is represented by the $3$-cocycle
$$(1_{tA},-1_{tB},1_{tC})\in\tilde{G}^3((X,\alpha)^*\bar{t}).$$ For
this reason it is called class of identities.
\end{rem}

We will use later another cocycle representing the class of identities given in the following lemma.

\begin{lem}\label{laotra}
The cochain $b=(b_1,0)\in\tilde{F}(sh^3,(X, \alpha)^*\bar{t})$ defined by $b_1(m)=(-1)^m\alpha_m$ is a cocycle representing the class of identities.
\end{lem}

\begin{proof}
We defined $\idd(X,\alpha)$ as the cohomology class represented by
$(\alpha_{-2},-\alpha_{-1},\alpha_0)\in\tilde{G}^3((X,\alpha)\bar{t})$.
The quasi-isomorphism
$$\upsilon\colon\tilde{G}^*((X,\alpha)^*\bar{t})\To\tilde{F}^*(sh^3,(X,\alpha)^*\bar{t})$$
in Remark \ref{redu2} sends $(\alpha_{-2},-\alpha_{-1},\alpha_0)$
to the cochain
\begin{eqnarray*}
\upsilon(\alpha_0,-\alpha_1,\alpha_2)=(c,0)\in
\tilde{F}^3(sh^3,(X,\alpha)^*\bar{t})&=&
\prod_{m\in\Z}\hom_\C{A}(tX_m,X_{m+3}) \\&&\oplus
\prod_{m\in\Z}\hom_\C{A}(tX_{m+3},X_{m+5})
\end{eqnarray*}
defined as $c(3m-i)=((X,\alpha)^*\bar{t})^m((-1)^i\alpha_i)$ for
all $m\in\Z$ and $0\leq i\leq 2$. We check by induction on $m$
that $c(3m-i)=(-1)^{3m-i}\alpha_{3m-i}$. This is obvious for
$m=0$. If it is true for $m$ then
\begin{eqnarray*}
c(3(m+1)-i)&=&((X,\alpha)^*\bar{t})^{m+1}((-1)^i\alpha_i)\\
&=&((X,\alpha)^*\bar{t})_{(3m-i,3m-i+3)}(c(3m-i))\\
&=&((X,\alpha)^*\bar{t})_{(3m-i,3(m+1)-i)}((-1)^{3m-i}\alpha_{3m-i})\\
&=&\alpha_{3(m+1)-i}(-(-1)^{3m-i}t\alpha_{3m-i})(t\alpha_{3m-i}^{-1})\\
&=&(-1)^{3(m+1)-i}\alpha_{3(m+1)-i}.
\end{eqnarray*}

\end{proof}

A small candidate triangle is said to be \emph{$*$-trivial}\index{triangle!$*$-trivial} 
if it has one of the following forms:
$$A\st{1_A}\To A\To *\To tA,$$
$$A\To *\To tA\st{1_{tA}}\To tA,$$
$$*\To A\st{1_A}\To A\To *.$$
Notice that these are exactly those triangles which are required
to be exact in any (pre)triangulated category. Analogously a
candidate triangle $(X,\alpha)$ is $*$-trivial if there exists
$i\in \Z/3$ such that $d_{n+1}\colon X_n\r X_{n+1}$ is an identity
morphism and $X_{n+2}=*$ provided $n=i \mod 3$. The functors
(\ref{sb}) and (\ref{sb'}) preserve $*$-trivial objects.

\begin{lem}\label{esna}
If $(X,\alpha)$ is a $*$-trivial candidate triangle then
$H^n(sh^3,(X,\alpha)^*\bar{t})=0$ for all $n\geq 0$.
\end{lem}

This lemma is a consequence of the property of the class of
identities stated in the following lemma. Recall that a cochain
complex $X$ in $\C{A}$ is
\emph{contractible}\index{triangle!contractible} if there are
morphisms $h_i\colon X_i\r X_{i-1}$ such that
$1_{X_i}=h_{i+1}d_{i+1}+d_ih_i$ for all $i\in\Z$. A candidate
triangle $(X,\alpha)$ is said to be contractible if $X$ is
contractible and the morphisms $h_i$ can be chosen in such a way
that $\alpha_{i-1}(th_i)=h_{i+3}\alpha_i$ for all $i\in\Z$. The
$*$-trivial candidate triangles are the most basic examples of
contractible candidate triangles. The reader can easily check that
under the equivalences of categories in (\ref{sb}) and (\ref{sb'})
contractible candidate triangles correspond to small contractible
candidate triangles in the sense of \cite{snatc} 1.5.

\begin{lem}\label{esna2}
Let $(X,\alpha)$ be a candidate triangle. The following statements are equivalent:
\begin{enumerate}
\item $\idd(X,\alpha)=0$,

\item $(X,\alpha)$ is contractible,

\item $X$ is contractible,

\item $X\, sh^n$ is contractible for all $n\in\Z$,

\item $H^*(\V,D(X,-))=0$ for any functor
$D\colon\C{A}^\op\times\V\r\Ab$ additive in the first variable and
such that $D(-,*)=0$,

\item $H^*(\V,D(-,X))=0$ for any functor
$D\colon\V\times\C{A}\r\Ab$ additive in the second variable and
such that $D(*,-)=0$,

\item $H^*(sh^3,[(X,\alpha),(Y,\beta)]^*\bar{s})=0$ for any natural transformation $\bar{s}\colon D\rr D(t,t)$ with $D$ a biadditive $\C{A}$-bimodule and any
candidate triangle $(Y,\beta)$,

\item $H^*(sh^3,[(Y,\beta),(X,\alpha)]^*\bar{s})=0$ for any natural transformation $\bar{s}\colon D\rr D(t,t)$ with $D$ a biadditive $\C{A}$-bimodule and any
candidate triangle $(Y,\beta)$,

\item $H^*(sh^3,(X,\alpha)^*\bar{s})=0$ for any natural transformation $\bar{s}\colon D\rr D(t,t)$ with $D$ a biadditive $\C{A}$-bimodule.
\end{enumerate}
\end{lem}

\begin{proof}
Obviously $(2)\rr (3)$, $(4)\rr (3)$, $(7)\rr(9)$, $(8)\rr(9)$ and
$(9)\rr (1)$ are trivial. It is elementary to check $(3)\rr (4)$.
If $X$ is contractible and $D$ is as in $(5)$ then the cochain
complex of abelian groups $\tilde{F}^*(\V,D(X,-))$ is also
contractile. The contraction is given by the homomorphisms
$\bar{h}_i\colon\tilde{F}^i(\V,D(X,-))\r\tilde{F}^{i-1}(\V,D(X,-))$
defined by $(\bar{h}_i(c))(m)=(-1)^ih_{m}^*(c(m-1))$. This proves
$(3)\rr(5)$. Similarly if $X$ is contractible and $D$ is as in
$(6)$ the cochain complex of abelian groups
$\tilde{F}^*(\V,D(-,X))$ is also contractile and contraction is
given by the homomorphisms
$\hat{h}_i\colon\tilde{F}^i(\V,D(-,X))\r\tilde{F}^{i-1}(\V,D(-,X))$
defined by $(\hat{h}_i(c))(m)={h_{m+i}}_*(c(m))$. This proves
$(3)\rr(6)$. Since we have already showed $(3)\Leftrightarrow(4)$
and $(3)\rr(5)$ (resp. $(3)\rr(6)$) the exact sequence
(\ref{exacta}) proves that $(5)\rr(7)$ (resp. $(6)\rr(8)$).

Suppose now that $\idd(X,\alpha)=0$. By 
Lemma \ref{laotra} we see that there are morphisms $k_i\colon
tX_i\r X_{i+2}$ such that
\begin{equation*}\tag{a}
d_{i+3}k_i-k_{i+1}(td_{i+1})=(-1)^i\alpha_i,
\end{equation*}
\begin{equation*}\tag{b}
-\alpha_{i+2}(tk_i)(t\alpha_i^ {-1})-k_{i+3}=0
\end{equation*}
for
all $i\in\Z$. Let us check that the morphisms
$$h_i=(-1)^{i}k_{i-3}\alpha_{i-3}^{-1}\colon X_i\r X_{i-1}$$
determine a contraction of $(X,\alpha)$. On one hand
\begin{eqnarray*}
\alpha_{i-1}(th_i) &=& \alpha_{i-1}(-1)^{i}(tk_{i-3})(t\alpha_{i-3}^{-1})\\
&\st{\text{(b)}}{=}& -(-1)^{i}k_i\\
&=& h_{i+3}\alpha_i,
\end{eqnarray*}
and on the other hand
\begin{eqnarray*}
h_{i+1}d_{i+1}+d_ih_i &=& (-1)^{i+1}k_{i-2}\alpha_{i-2}^{-1}d_{i+1}+d_i(-1)^{i}k_{i-3}\alpha_{i-3}^{-1}\\
&=&(-1)^i(-k_{i-2}(td_{i-2})+d_ik_{i-3})\alpha^{-1}_{i-3}\\
&\st{\text{(a)}}{=}& (-1)^i(-1)^i\alpha_{i-3}\alpha^{-1}_{i-3}\\
&=& 1_{X_i}.
\end{eqnarray*}
Here in the second equality we use that $\alpha$ is a cochain
homomorphism. This proves $(1)\rr(2)$. Hence the proof of this
lemma is finished.
\end{proof}

The implication $(2)\rr(9)$ should be compared, after Section
\ref{nablat}, to the fact that contractible small candidate
triangles are exact in all triangulated categories, see
\cite{snatc} 1.6.



\section{$\nabla$-Triangles}\label{nablat}

Here we obtain a purely cohomological characterization of small
$\nabla$-triangles defined in \cite{ccctc1} 2.13 (Proposition
\ref{reao}).

\begin{rem}
Small $\nabla$-triangles determine the ``underlying'' ordinary
(pre)triangulated structure of a cohomologically (pre)triangulated
category in the sense of \cite{ccctc1} 2.16, see \cite{ccctc1} 2.17.
They are defined in \cite{ccctc1} 2.13 by using the main result in
\cite{ptc}, which says that $3$-dimensional translation cohomology
classes are represented by translation track categories. We already
know by \cite{ccctc1} 2.12 that nevertheless the class of small
$\nabla$-triangles is independent of the choice of a translation
track category representing $\nabla$. However the proof of
\cite{ccctc1} 2.12 is still track-theoretical.
\end{rem}

In this section $(\C{A},t)$ will be a translation category given
by an additive category $\C{A}$ and an additive equivalence $t$. 
Let $\hom^t=\hom_\C{A}(t,-)$ be  the
$\C{A}$-bimodule defined in (\ref{homt}) 
and $
\bar{t}\colon\hom^t=\hom_\C{A}(t,-)\rr\hom_\C{A}(t^2,t)=\hom^t(t,t)$
the natural transformation in (\ref{tbar}) given by $(-1)t$. The
translation cohomology of $(t,\bar{t})$ will be denoted in this paper by
\begin{equation}
H^*(\C{A},t)=H^*(t,\bar{t}).
\end{equation}

Given two candidate triangles $(X,\alpha)$ and $(Y,\beta)$ in the
sense of Section \ref{ct} the translation cohomology groups
$(n\in\Z)$
$$H^n(sh^3,[(X,\alpha),(Y,\beta)]^*\bar{t})$$
are
defined by (\ref{indt2}). This actually defines a
$\C{Cand}(\C{A},t)$-bimodule
\begin{equation}\label{unobi}
H^n(sh^3,[-,-]^*\bar{t}),
\end{equation}
more precisely, given a diagram in $\C{Cand}(\C{A},t)$
$$(X,\alpha)\st{f}\To(Y,\beta)\st{g}\To(Z,\gamma)$$
the induced morphisms $g_*$ and $f^*$ are defined by
$$g_*=\hom^t(X,g)_*\colon H^n(sh^3,[(X,\alpha),(Y,\beta)]^*\bar{t})\To H^n(sh^3,[(X,\alpha),(Z,\gamma)]^*\bar{t}),$$
$$f^*=\hom^t(f,Z)_*\colon H^n(sh^3,[(X,\alpha),(Z,\gamma)]^*\bar{t})\To H^n(sh^3,[(Y,\beta),(Z,\gamma)]^*\bar{t}).$$
Notice that by (\ref{diag})
$$H^n(sh^3,[(X,\alpha),(X,\alpha)]^*\bar{t})=H^n(sh^3,(X,\alpha)^*\bar{t}),$$
therefore the classes of identities in (\ref{wi}) define a $0$-cochain of the category $\C{Cand}(\C{A},t)$
with coefficients in $H^3(sh^3,[-,-]^*\bar{t})$.

\begin{lem}\label{wic}
The classes of identities defined in (\ref{wi}) yield a
cohomology class
$$\idd\in H^0(\C{Cand}(\C{A},t),H^3(sh^3,[-,-]^*\bar{t})).$$
\end{lem}

This follows from Lemma \ref{elpul2}.

Let us fix a translation cohomology class
$$\nabla\in H^3(\C{A},t).$$

\begin{lem}\label{ni}
The translation cohomology class $\nabla$ yields a cohomology
class
$$\bar{\nabla}\in H^0(\C{Cand}(\C{A},t),H^3(sh^3,[-,-]^*\bar{t}))$$
defined over any candidate triangle $(X,\alpha)$ by the formula
$$\bar{\nabla}(X,\alpha)=(X,\alpha)^*\nabla\in H^3(sh^3,(X,\alpha)^*\bar{t}).$$
\end{lem}

This follows from Lemma \ref{elpul}.

A \emph{$\nabla$-triangle}\index{triangle!$\nabla$-} is a
candidate triangle $(X,\alpha)$ in $(\C{A},t)$ which is
annihilated by the cohomology class
$$\idd+\bar{\nabla}\in H^0(\C{Cand}(\C{A},t),H^3(sh^3,[-,-]^*\bar{t})),$$
or equivalently, $(X,\alpha)$ satisfies the equation
\begin{equation}\label{tont}
0=\idd(X,\alpha)+(X,\alpha)^*\nabla\in
H^3(sh^3,(X,\alpha)^*\bar{t}).
\end{equation}

\begin{prop}\label{reao}
Small $\nabla$-triangles 
are small candidate triangles 
corresponding to $\nabla$-triangles by the equivalence of
categories (\ref{sb'}).
\end{prop}

The proof of this proposition is in Section \ref{pruebas}.

This result allows to prove the following proposition, which
completes the proof of the main theorem in \cite{ccctc1}.

\begin{prop}\label{ade}
Let $(\C{A},t)$ be a translation category where $\C{A}$ is additive and $t$ is an additive equivalence and let $\nabla$ be a
translation cohomology class in $H^3(\C{A},t)$. Then the class $\mathcal{E}_\nabla$ of small $\nabla$-triangles is defined and
$(\C{A},t,\mathcal{E}_\nabla)$ satisfies the axioms (Tr0), (Tr2) and
(Tr3) in the definition of a pretriangulated category, see
\cite{ccctc1} 1.1.
\end{prop}


\begin{proof}
Let us prove the first part of (Tr0). For this we consider an
isomorphism of candidate triangles $f\colon(X,\alpha)\r(Y,\beta)$.
Since $\idd+\bar{\nabla}$ is a $0$-dimensional cohomology class we
have that
$$f_*(\idd(X,\alpha)+(X,\alpha)^*\nabla)=f^*(\idd(Y,\beta)+(Y,\beta)^*\nabla)\in
H^3(sh^3,[(X,\alpha),(Y,\beta)]^*\bar{t}),$$ but $f_*$ and $f^*$ are
isomorphisms since $f$ is already an isomorphism, therefore
$\idd(X,\alpha)+(X,\alpha)^*\nabla$ vanishes if and only if
$\idd(Y,\beta)+(Y,\beta)^*\nabla$ vanishes, i. e. $(X,\alpha)$ is a
$\nabla$-triangle if and only if $(Y,\beta)$ is a $\nabla$-triangle.
The second part follows from Lemma \ref{esna}. Axiom (Tr2) follows
from Proposition \ref{tr2na} in Section \ref{tr2}. And finally (Tr3)
is a consequence of Proposition \cite{ccctc1} 4.18, as we already
pointed out in the proof of \cite{ccctc1} 2.17.
\end{proof}

\section{Algebraic homotopy pairs and cohomologically triangulated categories}\label{ahp}

Suppose that 
$\C{A}$ is an additive category, $t\colon\C{A}\r\C{A}$ is an
additive equivalence and $\nabla\in H^3(\C{A},t)$. The homomorphism
$j$ in the exact sequence (\ref{exacta}) determines a cohomology
class
$$j\nabla\in H^3(\C{A},\hom^t).$$

Recall that the objects of the additive category $\C{Pair}(\C{A})$ of pairs in $\C{A}$ are the morphisms in $\C{A}$, and morphisms $(h,k)\colon f\r g$ are commutative squares
$$\xymatrix{A\ar[r]^h\ar[d]_f&C\ar[d]^g\\B\ar[r]^k&D}$$
The $\C{Pair}(\C{A})$-bimodule $(\hom^t)^\#$
is defined by
$$(\hom^t)^\#(f,g)=\hom_\C{A}(tA,D)/(g_*\hom_\C{A}(tA,C)+(tf)^*\hom_\C{A}(tB,D)).$$

There are natural homomorphisms
$$\lambda\colon H^{n+1}(\C{A},\hom^t)\To H^n(\C{Pair}(\C{A}),(\hom^t)^\#),\;\;\;n\geq 0,$$
defined in \cite{tbhp}. These homomorphisms can also be obtained from the K\"unneth spectral sequence for cohomology of categories in
Corollary \ref{ibia}, see Remark \ref{ibib}. In particular we can consider the cohomology class
$$\lambda j\nabla\in H^2(\C{Pair}(\C{A}),(\hom^t)^\#).$$
This class will be called the \emph{characteristic class of
algebraic homotopy pairs}. It corresponds to a well-defined
isomorphism class of linear extensions of the category
$\C{Pair}(\C{A})$ by $(\hom^t)^\#$, see \cite{ah} IV.6. Any
$2$-cocycle $\aleph$ representing $\lambda j\nabla$ determines such
a linear extension
\begin{equation}\label{leahp}
(\hom^t)^\#\st{+}\To\C{AHp}(\aleph)\st{\bar{p}}\To\C{Pair}(\C{A}).
\end{equation}
Here $\C{AHp}(\aleph)$ is the \emph{category of algebraic homotopy
pairs}\index{homotopy pair!algebraic} associated to the cocycle
$\aleph$. It has the same objects as $\C{Pair}(\C{A})$ and morphisms
$(h,k,a)\colon f\r g$ are given by a morphism $(h,k)\colon f\r g$ in
$\C{Pair}(\C{A})$ together with an element $a\in (\hom^t)^\#(f,g)$.
Composition is defined by the following formula
$$(\bar{h},\bar{k},\bar{a})(h,k,a)=(h\bar{h},\bar{k}k,(\bar{h},\bar{k})_*a+(h,k)^*\bar{a}-\aleph((\bar{h},\bar{k}),(h,k))).$$
The functor $\bar{p}$ is the identity on objects and $\bar{p}(h,k,a)=(h,k)$. 
The category $\C{AHp}(\aleph)$ is well defined by $\lambda j\nabla$ up to isomorphism over $\C{Pair}(\C{A})$ equivariant with respect to the action
of $(\hom^t)^\#$ in the third coordinate of morphisms, see \cite{ah} IV.3.4.

By \cite{qcsr} 6.2 we have that $\C{AHp}(\aleph)$ has a unique additive category structure such that $\bar{p}$ is an additive functor and
\begin{equation}\label{sesag}
(\hom^t)^\#(f,g)\hookrightarrow\hom_{\C{AHp}(\aleph)}(f,g)\st{\bar{p}}\twoheadrightarrow\hom_{\C{Pair}(\C{A})}(f,g)
\end{equation}
is a short exact sequence of abelian groups where the first homomorphism is defined by $b\mapsto (0,0,b)$.

There is a full additive inclusion
$$\iota\colon\C{A}\To\C{AHp}(\aleph)$$ given by $\iota(A)\colon *\r
A$. In particular for any morphism $f\colon A\r B$ in $\C{A}$ we can
consider the \emph{algebraic cone functor}
\begin{equation}\label{conalgo}
\acone_f=\hom_{\C{AHp}(\aleph)}(f,\iota)\colon\C{A}\To\C{Ab}.
\end{equation}
Recall that an $\C{A}$-module is just an additive functor
$\C{A}\r\C{Ab}$, therefore algebraic cone functors are
$\C{A}$-modules. Morphisms between $\C{A}$-modules are natural
transformations and the category of $\C{A}$-modules is designated by
$\mo(\C{A})$. There is a Yoneda full inclusion
$\C{A}^\op\subset\mo(\C{A})$ sending an object $A$ to the
$\C{A}$-module $\hom_\C{A}(A,-)$. Algebraic cone functors determine
an additive functor
\begin{equation}\label{conalg}
\acone\colon\C{AHp}(\aleph)^\op\To\mo(\C{A})
\end{equation}
in the obvious way.

The exactness of (\ref{sesag}) implies that for any $f\colon A\r B$ in $\C{A}$ there is an exact sequence of $\C{A}$-modules
\begin{equation}\label{exa5}
\hom_\C{A}(tB,-)\st{(tf)^*}\To\hom_\C{A}(tA,-)\st{\tilde{q}_f}\To\acone_f\st{\tilde{i}_f}\To\hom_\C{A}(B,-)\st{f^*}\To\hom_\C{A}(A,-).
\end{equation}
This exact sequence is natural in the category of algebraic homotopy pairs, i. e. for any morphism $(h,k,a)\colon f\r \bar{f}$ in $\C{AHp}(\aleph)$
the following diagram commutes
\begin{equation}\label{conmuta1}
\xymatrix{\hom_\C{A}(tB,-)\ar[r]^{(tf)^*}&\hom_\C{A}(tA,-)\ar[r]^<(.3){\tilde{q}_f}&
\acone_f\ar[r]^<(.2){\tilde{i}_f}&\hom_{\C{A}}(B,-)\ar[r]^{f^*}&
\hom_{\C{A}}(A,-)\\
\hom_\C{A}(t\bar{B},-)\ar[u]_{(tk)^*}\ar[r]_{(t\bar{f})^*}&\hom_\C{A}(t\bar{A},-)\ar[u]_{(th)^*}\ar[r]_<(.3){\tilde{q}_{\bar{f}}}&
\acone_{\bar{f}}\ar[r]_<(.2){\tilde{i}_{\bar{f}}}\ar[u]_{(h,k,a)^*}&\hom_{\C{A}}(\bar{B},-)\ar[u]_{k^*}\ar[r]_{\bar{f}^*}&
\hom_{\C{A}}(\bar{A},-)\ar[u]_{h^*}}
\end{equation}

Let $\C{AHp}^{rep}(\aleph)\subset\C{AHp}(\aleph)$ be the full
subcategory whose objects are morphisms $f$ in $\C{A}$ such that the
associated algebraic cone functor $\acone_f$ is representable. For
each object $f$ in $\C{AHp}^{rep}(\aleph)$ we choose an
$\C{A}$-module isomorphism
\begin{equation}\label{lanabla}
\chi_f\colon\hom_\C{A}(C_f,-)\cong\acone_f.
\end{equation}
A particular choice determines a functor
$$\bar{C}\colon\C{AHp}^{rep}(\aleph)\To\C{A},\;\; f\mapsto C_f,$$
such that the restriction to $\C{AHp}^{rep}(\aleph)$ of $\acone$ in
(\ref{conalg}) factors up the natural equivalence $\chi$ as
$$\C{AHp}^{rep}(\aleph)^{op}\st{\bar{C}^{op}}\To\C{A}^{op}\subset\mo(\C{A}).$$
By using $\bar{C}$ we can define a functor
\begin{equation}\label{papa}
\bar{\zeta}\colon\C{AHp}^{rep}(\aleph)\To\C{cand}(\C{A},t)
\end{equation}
as follows: given an object $f\colon A\r B$ in
$\C{AHp}^{rep}(\aleph)$ the small candidate triangle
$\bar{\zeta}(f)$ is given by
$$A\st{f}\To B\st{i_f}\To C_f\st{q_f}\To tA$$
where $i_f$ and $q_f$ are determined by the $\C{A}$-module morphisms
$\tilde{i}_f$ and $\tilde{q}_f$ in (\ref{exa5}), respectively, the
isomorphism $\chi_f$ in (\ref{lanabla}) and Yoneda's lemma;
and given a morphism $(h,k,a)\colon f\r \bar{f}$ in
$\C{AHp}^{rep}(\aleph)$ then $\bar{\zeta}(h,k,a)$ is the following
commutative diagram
$$\xymatrix{A\ar[r]^{f}\ar[d]^{h}&B\ar[r]^{i_f}\ar[d]^{k}&C_f\ar[r]^{q_f}\ar[d]^{\bar{C}(h,k,a)}&tA\ar[d]^{th}\\
\bar{A}\ar[r]_{\bar{f}}&\bar{B}\ar[r]_{i_{\bar{f}}}&C_{\bar{f}}\ar[r]_{q_{\bar{f}}}&t\bar{A}}$$
This diagram commutes by (\ref{conmuta1}).
The functor $\bar{\zeta}$ does not depend on the different choices made for its definition, up to natural equivalence.%

In \cite{ccctc1} 4.22 we defined a linear action of the
$\C{Pair}(\C{A})$-bimodule
$$\hom_{\mo(\C{A}^\op)}(\ker t,\coker)$$
on the obvious forgetful functor
$$\varsigma\colon\C{cand}(\C{A},t)\To\C{Pair}(\C{A})$$
defined in \cite{ccctc1} 4.15. Moreover, in \cite{ccctc1} 4.23 we
define a natural transformation
$$\Xi\colon(\hom^t)^\#\rr\hom_{\mo(\C{A}^\op)}(\ker t,\coker).$$
This natural transformation is compatible with the functor
$\bar{\zeta}$ in the following sense.

\begin{prop}\label{pipi}
Given a morphism $(h,k,a)\colon f\r\bar{f}$ in
$\C{AHp}^{rep}(\aleph)$ and
$\alpha\in(\hom^t)^\#(\bar{p}(f),\bar{p}(\bar{f}))$ then
$$\bar{\zeta}((h,k,a)+\alpha)=\bar{\zeta}(h,k,a)+\Xi_{{(\bar{p}(f),\bar{p}(\bar{f}))}}(\alpha).$$
\end{prop}

This proposition follows by the same kind of arguments as in the
proof of \cite{ccctc1} 4.24.

We also defined in \cite{ccctc1} 4.20 a natural transformation
$$\varpi\colon H^1(sh^3,[-,-]^*\bar{t})\rr(\hom^t)^\#(\varsigma,\varsigma).$$
We can use this natural transformation to define the following
natural equivalence relation in $\C{AHp}(\aleph)$: two morphisms
$(h,k,a),(\bar{h},\bar{k},\bar{a})\colon f\r\bar{f}$ in
$\C{AHp}(\aleph)$ are equivalent if
$(h,k,a)+\varpi(b)=(\bar{h},\bar{k},\bar{a})$ for some $b\in
H^1(sh^3,[\bar{\zeta}(f),\bar{\zeta}(\bar{f})]^*\bar{t})$. The
quotient category will be denoted by
$$\C{AHp}(\aleph)/(\varpi,\bar{\zeta}).$$
Moreover, we denote
$$\C{AHp}^{rep}(\aleph)/(\varpi,\bar{\zeta})$$
to the full subcategory of the quotient whose objects are the
morphisms $f$ in $\C{A}$ with representable associated algebraic
cone functor.

By Proposition \ref{pipi} and the first part of \cite{ccctc1} 4.25
the functor $\bar{\zeta}$ factors through this quotient category
\begin{equation}\label{goroo}
\bar{\zeta}\colon\C{AHp}^{rep}(\aleph)\twoheadrightarrow\C{AHp}^{rep}(\aleph)/(\varpi,\bar{\zeta})\st{\hat{\zeta}}\To\C{cand}(\C{A},t).
\end{equation}
Here the first functor is the projection onto the quotient category
and $\hat{\zeta}$ is uniquely determined by this factorization.

The following theorem is one of the main results of this paper.

\begin{thm}\label{ñi}
Let $\C{A}$ be an additive category, $t\colon\C{A}\r\C{A}$ an
additive equivalence and $\nabla\in H^3(\C{A},t)$. Then
$(\C{A},t,\nabla)$ is cohomologically pretriangulated if and only if
all algebraic cone functors defined in (\ref{conalgo}) are
representable and, using the equivalence of categories in
(\ref{sb'}) as an identification, the pull-back of the cohomology
class $\idd+\bar{\nabla}$ in Section \ref{nablat} by the functor
$\hat{\zeta}$ in (\ref{goroo}) vanishes
$$0=\hat{\zeta}^*(\idd+\bar{\nabla})\in H^0(\C{AHp}(\aleph)/(\varpi,\bar{\zeta}),
H^3(sh^3,[\bar{\zeta},\bar{\zeta}]^*\bar{t})).$$
\end{thm}

The proof of this result will be given in Section \ref{pruebas}.

\begin{rem}\label{caca}
At a first glance the vanishing condition in Theorem \ref{ñi} seems
to depend on the choice of the cocycle $\aleph$ representing
$\lambda j \nabla$ and the representations $\chi_f$ in
(\ref{lanabla}) used in the definition of the functor
$\bar{\zeta}\colon\C{AHp}(\aleph)\r\C{cand}(\C{A},t)$. However it is
easy to check that if $\bar{\aleph}$ is another cocycle representing
$\lambda j\nabla$ and
$\bar{\xi}\colon\C{AHp}(\bar{\aleph})\r\C{cand}(\C{A},t)$ is defined
as $\bar{\zeta}$ then for any isomorphism of linear extensions
$$\xymatrix{(\hom^t)^\#\ar[r]^+\ar@{=}[d]&\C{AHp}(\aleph)\ar[r]\ar[d]_\cong&\C{Pair}(\C{A})\ar@{=}[d]\\
(\hom^t)^\#\ar[r]^+&\C{AHp}(\bar{\aleph})\ar[r]&\C{Pair}(\C{A})}$$
the isomorphism of categories in the middle induces isomorphisms
between the algebraic cone functors given by $\aleph$ and
$\bar{\aleph}$, therefore the first ones are representable if and
only if the second ones are. Moreover, the diagram
$$\xymatrix{\C{AHp}^{rep}(\aleph)\ar[dd]_\cong\ar[rd]^{\bar{\zeta}}&\\&\C{cand}(\C{A},t)\\\C{AHp}^{rep}(\bar{\aleph})\ar[ru]_{\bar{\xi}}&}$$
commutes up to natural equivalence. This second diagram induces
another one
$$\xymatrix{\C{AHp}^{rep}(\aleph)/(\varpi,\bar{\zeta}) \ar[dd]_\cong\ar[rd]^{\hat{\zeta}}&\\
&\C{cand}(\nabla)\\\C{AHp}^{rep}(\bar{\aleph})/(\varpi,\bar{\xi})\ar[ru]_{\hat{\xi}}&}$$
which is also commutative up to equivalence. If algebraic cone
functors are representable the behavior of cohomology of categories
with respect to natural transformations, see \cite{fcc} 4.3, implies
that the cohomology groups where the classes
$\hat{\zeta}^*(\idd+\bar{\nabla})$ and
$\hat{\xi}^*(\idd+\bar{\nabla})$ live are isomorphic under an
isomorphism carrying one of these classes to the other, therefore
$$\hat{\zeta}^*(\idd+\bar{\nabla})=0\Leftrightarrow\hat{\xi}^*(\idd+\bar{\nabla})=0,$$
so the vanishing condition in Theorem \ref{ñi} is independent of all
choices.
\end{rem}

Suppose now that $(\C{A},t,\nabla)$ is cohomologically
pretriangulated. By Theorem \ref{ñi} and Proposition \ref{reao} the
functor $\hat{\zeta}$ in (\ref{goroo}) takes values in small
$\nabla$-triangles
$$\hat{\zeta}\colon\C{AHp}(\aleph)/(\varpi,\bar{\zeta})\To\C{cand}(\nabla).$$



By \cite{ah} IV.4.10 (e) the exact sequence for a functor
constructed in \cite{ccctc1} 3.1 gives rise to a well-defined
cohomology class
\begin{equation}\label{conabla}
\theta_\nabla\in
H^1(\C{cand}(\nabla),H^2(sh^3,[i_\nabla,i_\nabla]^*\bar{t})).
\end{equation}%
Here $i_\nabla\colon\C{cand}(\nabla)\subset\C{cand}(\C{A},t)$ is the inclusion of the full subcategory of small $\nabla$-triangles.

The following theorem 
is also a main result in this paper.

\begin{thm}\label{vanish}
Let $(\C{A},t,\nabla)$ be a cohomologically pretriangulated category. Then it is cohomologically triangulated if and only if
$$0=\hat{\zeta}^*\theta_\nabla\in H^1(\C{AHp}(\aleph)/(\varpi, \bar{\zeta}),
H^2(sh^3,[i_\nabla\hat{\zeta},i_\nabla\hat{\zeta}]^*\bar{t})).$$
\end{thm}

The proof of this theorem will be given in Section \ref{pruebas}.

As in the case of small $\nabla$-triangles, the original definition
of $\theta_\nabla$ is based on the fact that $3$-dimensional
translation cohomology classes are represented by translation track
categories. Nevertheless, as we did with small $\nabla$-triangles in
Section \ref{nablat}, in the next section we give a purely
cohomological procedure to obtain $\theta_\nabla$ from $\nabla$ by
using a spectral sequence.

Remark \ref{caca} also applies to Theorem \ref{vanish} with the appropriate modifications.


\section{The main spectral sequence}\label{tmss}

Cohomologically triangulated categories $(\C{A},t,\nabla)$ were
defined in \cite{ccctc1} by using the representability of
$3$-dimensional translation cohomology classes by translation track
categories. By means of the classes of identities in Section
\ref{cotri} we were able to show in Proposition \ref{reao} that
small $\nabla$-triangles (the exact triangles of cohomologically
triangulated categories) can be obtained from $\nabla$ by purely
cohomological methods. Theorem \ref{vanish} characterizes
cohomologically triangulated categories in strictly cohomological
terms by using a $1$-dimensional cohomology class $\theta_\nabla$.
This class is obtained from an exact sequence for functors in
\cite{ccctc1} 3.1 encoding the obstruction theory for the category
of track triangles associated to a good translation track category
representing $\nabla$.
However the reader may be missing a more intrinsical approach to
cohomologically triangulated categories, let us say a construction
of $\theta_\nabla$ from $\nabla$ by using only techniques from
cohomology of categories. We address this issue in this section,
which is mainly expository. The most of the proofs and the technical
constructions appear in the next one. As a consequence of these
results we give a new characterization of cohomologically
(pre)triangulated categories in terms of the filtration associated
to the spectral sequences constructed in the next section. This will
be the main result of this paper.

In the next section we construct some spectral sequences for the computation of the (translation) cohomology of products. There we also
define for any translation category $(\C{A},t)$ with $\C{A}$ additive and $t$ an additive equivalence
a canonical translation cohomology class
\begin{equation*}
\id\in H^3(1_{\C{Cand}(\C{A},t)}\times sh^3,[-,-]^*\bar{t}),
\end{equation*}
which we call the \emph{global class of
identities}\index{identities@global class of -}, see
(\ref{clasi}).

Fixed a translation cohomology class $$\nabla\in H^3(\C{A},t)$$ in
this section we show how the spectral sequences in Section \ref{ss}
together with $\nabla$ and the global class of identities determine
the cohomology classes $\idd+\bar{\nabla}$ and $\theta_\nabla$ which
define $\nabla$-triangles and are used to detect when
$(\C{A},t,\nabla)$ is a cohomologically (pre)triangulated category,
see Theorems \ref{ñi} and \ref{vanish}.

In Remark \ref{prinex} we define an ``evaluation'' translation
category morphism
$$(ev,\epsilon)\colon\C{Cand}(\C{A},t)\times (\V,sh^3)\To(\C{A},t).$$
As we explain in that remark the pull-back of $\bar{t}$ in
(\ref{tbar}) along $(ev,\epsilon)$ is denoted by
$$(ev,\epsilon)^*\bar{t}=[-,-]^*\bar{t},$$ since it is given by the natural transformations between $\V$-bimodules defined
by (\ref{indt2}). This ``evaluation'' morphism can be used to pull back the translation cohomology class $\nabla$ to the same
cohomology group where the global class of identities lives
\begin{equation*}
(ev,\epsilon)^*\nabla \in H^3(1_{\C{Cand}(\C{A},t)}\times
sh^3,[-,-]^*\bar{t}).
\end{equation*}

The following theorem is the main result of this paper.

\begin{thm}\label{mamona}
Let $\C{A}$ be an additive category, $t\colon\C{A}\r\C{A}$ an
additive equivalence and $\nabla\in H^3(\C{A},t)$. Then
$(\C{A},t,\nabla)$ is cohomologically triangulated (resp.
pretriangulated) if and only if all algebraic cone functors defined
in (\ref{conalgo}) are representable and the cohomology class
$$(\hat{\zeta}\times 1_\V)^*(\id+(ev,\epsilon)^*\nabla) \in H^3(1_{\C{AHp}(\aleph)/(\varpi,\bar{\zeta})}\times
sh^3,[\hat{\zeta},\hat{\zeta}]^*\bar{t}),$$ where $\hat{\zeta}$ is defined in
(\ref{goroo}), lies in the third (resp. second) term $D^{2,1}$
(resp. $D^{1,2}$) of the filtration
$$0=D^{4,-1}\subset D^{3,0}\subset D^{2,1}\subset D^{1,2}\subset D^{0,3}=
H^3(1_{\C{AHp}(\aleph)/(\varpi,\bar{\zeta})}\times
sh^3,[\hat{\zeta},\hat{\zeta}]^*\bar{t})$$ given by the spectral sequence in Proposition
\ref{lass}.
\end{thm}

This theorem is a consequence of Theorems \ref{ñi} and \ref{vanish}
together with Lemma \ref{mamo}, (\ref{igualdad}), Proposition
\ref{lader} and the naturality of the spectral sequence in
Proposition \ref{lass} below.

Remark \ref{prinex}  and Proposition \ref{lass} yield a spectral sequence
\begin{equation*}
E^{p,q}_2=H^p(\C{Cand}(\C{A},t),H^q(sh^3,[-,-]^*\bar{t}))\Longrightarrow
H^{p+q}(1_{\C{Cand}(\C{A},t)}\times sh^3,[-,-]^*\bar{t}).
\end{equation*}
Let  $$0=D^{n+1,-1}\subset D^{n,0}\subset\dots\subset
D^{p,q}\subset\cdots\subset D^{0,n}=H^n(1_{\C{Cand}(\C{A},t)}\times sh^3,[-,-]^*\bar{t})$$ be the filtration such that
$E^{p,q}_\infty=D^{p,q}/D^{p+1,q-1}$.
There are homomorphisms $(n\geq 0)$
\begin{equation*}
\xi_0\colon H^n(1_{\C{Cand}(\C{A},t)}\times
sh^3,[-,-]^*\bar{t})\To
H^0(\C{Cand}(\C{A},t),H^n(sh^3,[-,-]^*\bar{t}))
\end{equation*}
given by
$$D^{0,n}\twoheadrightarrow D^{0,n}/D^{1,n-1}=E^{0,n}_\infty\subset E^{0,n}_2.$$
The homomorphisms $\xi_0$ satisfy the property stated in the
following lemma.

\begin{lem}\label{mamo}
In the conditions above, given $b\in H^n(\C{A},t)$ and a candidate
triangle $(X,\alpha)$ in $(\C{A},t)$ we have that
$$\xi_0((ev,\epsilon)^*b)(X,\alpha)=(X,\alpha)^*b\in H^n(sh^3,(X,\alpha)^*\bar{t}).$$
\end{lem}

This is Lemma \ref{elpul}, proved in the next section. Lemma
\ref{ni} follows from this one. Moreover, the global class of
identities determines the classes of identities through $\xi_0$.

\begin{lem}For any candidate triangle $(X,\alpha)$ we have
$(\xi_0(\id))(X,\alpha)=\idd(X,\alpha)$.
\end{lem}

This is Lemma \ref{elpul2}, also proved in the next section. This
implies that the $0$-cochain $\idd$ determined by the classes of
identities in (\ref{wi}) is actually a cohomology class as stated in
Lemma \ref{wic}, and moreover
\begin{equation}\label{igualdad}
\xi_0(\id)=\idd.
\end{equation}

In particular by using these lemmas we notice that
$\nabla$-triangles are the candidate triangles annihilated by the
cohomology class
$$\xi_0(\id+(ev,\epsilon)^*\nabla)\in H^0(\C{Cand}(\C{A},t),H^3(sh^3,[-,-]^*\bar{t})).$$


Let $i_\nabla\colon\C{Cand}(\nabla)\r\C{Cand}(\C{A},t)$ be the
inclusion of the full subcategory of $\nabla$-triangles. By
Proposition \ref{lass} there is also a spectral sequence
\begin{equation*}
E^{p,q}_\infty=H^p(\C{Cand}(\nabla),H^q(sh^3,[i_\nabla,i_\nabla]^*\bar{t}))\Longrightarrow
H^{p+q}(1_{\C{Cand}(\nabla)}\times
sh^3,[i_\nabla,i_\nabla]^*\bar{t}).
\end{equation*}
We also denote by  $$0=D^{n+1,-1}\subset D^{n,0}\subset\dots\subset
D^{p,q}\subset\cdots\subset D^{0,n}=H^n(1_{\C{Cand}(\nabla)}\times sh^3,[i_\nabla,i_\nabla]^*\bar{t})$$ to the filtration such that
$E^{p,q}_\infty=D^{p,q}/D^{p+1,q-1}$
and by $\xi_0$ to the homomorphisms $(n\geq 0)$
\begin{equation*}
\xi_0\colon H^n(1_{\C{Cand}(\nabla)}\times
sh^3,[i_\nabla,i_\nabla]^*\bar{t})\To
H^0(\C{Cand}(\nabla),H^n(sh^3,[i_\nabla,i_\nabla]^*\bar{t}))
\end{equation*}
given by
$$D^{0,n}\twoheadrightarrow D^{0,n}/D^{1,n-1}=E^{0,n}_\infty\subset E^{0,n}_2.$$

The category $\C{Cand}(\nabla)$ is the biggest full subcategory
such that its inclusion functor $i_\nabla$ satisfies
$$0=i_\nabla^*\xi_0(\id+(ev,\epsilon)^*\nabla)=\xi_0(i_\nabla\times 1_\V,0^\vc_{i_\nabla\times sh^3})^*(\id+(ev,\epsilon)^*\nabla).$$
Here for the second equality we use the naturality of the spectral
sequence in Proposition \ref{lass}. This shows that
$$(i_\nabla\times 1_\V,0^\vc_{i_\nabla\times sh^3})^*(\id+(ev,\epsilon)^*\nabla)\in D^{1,2}\subset
H^3(1_{\C{Cand}(\nabla)}\times
sh^3,[i_\nabla,i_\nabla]^*\bar{t}).$$ 
Let $\xi_1$ be the homomorphism given by
$$\xi_1\colon
D^{1,2}\twoheadrightarrow D^{1,2}/D^{2,1}=E^{1,2}_\infty\subset
E^{1,2}_2=H^1(\C{Cand}(\nabla),H^2(sh^3,[i_\nabla,i_\nabla]^*\bar{t})).$$

\begin{prop}\label{lader}
Let $${\theta}_\nabla\in
H^1(\C{Cand}(\nabla),H^2(sh^3,[i_\nabla,i_\nabla]^*\bar{t}))$$ be
the cohomology class determined by $\theta_\nabla$ in
(\ref{conabla}), the equivalence of categories in (\ref{sb}) and
Proposition \ref{reao}. Then the following equality holds
$$\theta_\nabla=\xi_1(i_\nabla\times 1_\V,0^\vc_{i_\nabla\times sh^3})^*(\id+(ev,\epsilon)^*\nabla).$$
\end{prop}

This theorem is proved in Section \ref{pruebas}

In particular we obtain that the cohomology class $\theta_\nabla$,
which was defined by using a translation track category
representing $\nabla$, is indeed well defined by $\nabla$. 

\section{The construction of the spectral sequences and the global class of identities}\label{ss}

The normalized cochain complex $\bar{F}^*(\C{C},D)$ defined in
Section \ref{cc} for the computation of cohomology of categories is
the normalization of a cosimplicial abelian group
$F^\bullet(\C{C},D)$, see \cite{ccglg} Appendix B,
$$\bar{F}^*(\C{C},D)=N F^\bullet(\C{C},D).$$
We refer the reader to the Appendix of this paper for basic facts
and notation concerning (bi)cosimplicial abelian groups,
normalized (bi)complexes, diagonals, and total complexes. Here
$F^\bullet(\C{C},D)$ is given by the abelian groups in
(\ref{producto}). The coboundaries $d_i$ and codegeneracies $s_i$
are given by
\begin{eqnarray*}
(d_ic)(\sigma_1,\dots,\sigma_n)&=&\left\{\begin{array}{ll}
{\sigma_1}_*c(\sigma_2,\dots,\sigma_n),&i=0;
\\c(\sigma_1,\dots,\sigma_i\sigma_{i+1},\dots,\sigma_n),&0<i<n;\\
\sigma_n^*c(\sigma_1,\dots,\sigma_{n-1}),& i=n;
\end{array}\right.\\
(s_ic)(\sigma_1,\dots,\sigma_{n-1})&=&c(\sigma_1,\dots,\sigma_i,1_{X_i},\sigma_{i+1},\dots,\sigma_{n-1}).
\end{eqnarray*}
where $c\in F^n(\C{C},D)$ and $X_i$ is the source of $\sigma_i$ and/or the target of $\sigma_{i+1}$. The ordinary functorial properties of this bicomplex
are the same as for the cochain complexes $\bar{F}^*(\C{C},D)$ and $F^*(\C{C},D)$, see \cite{csc} or \cite{fcc}.

Let $\C{C}$ and $\C{B}$ be categories and $D$ a $(\C{C}\times\C{B})$-bimodule. The functorial properties of $F^\bullet$ yield a cosimplicial $\C{C}$-bimodule
$F^\bullet(\C{B},D)$. Moreover, by applying $F^\bullet(\C{C},-)$ to this simplicial $\C{C}$-bimodule we obtain a bicosimplicial abelian group
$$F^\bullet(\C{C},F^\bullet(\C{B},D)).$$

\begin{lem}\label{ibi}
There is a natural isomorphism of cosimplicial abelian groups
$$\diag F^\bullet(\C{C},F^\bullet(\C{B},D))\cong F^\bullet(\C{C}\times\C{B},D).$$
\end{lem}

\begin{proof}
A cochain $c\in F^n(\C{C},F^n(\C{B},D))$ corresponds to the cochain $\bar{c}\in F^n(\C{C}\times\C{B},D)$ defined by
$$\bar{c}((\sigma_1,\tau_1),\dots,(\sigma_n,\tau_n))=(c(\sigma_1,\dots,\sigma_n))(\tau_1,\dots,\tau_n).$$
The reader can easily check that this identification is compatible with the cosimplicial operators.
\end{proof}

The isomorphism in this lemma will be used as an identification. Moreover, this lemma and the cosimplicial Eilenberg-Zilber-Cartier theorem
described in the Appendix yield the following result.

\begin{cor}[K\"unneth spectral sequence]\label{ibia}
Given categories $\C{C}$ and $\C{D}$ and a $(\C{C}\times\C{B})$-bimodule $D$, there is a natural spectral sequence
$$E^{p,q}_2=H^p(\C{C},H^q(\C{B},D))\Longrightarrow H^{p+q}(\C{C}\times\C{B},D).$$
\end{cor}

\begin{proof}
This is the spectral sequence of the bicomplex
$\bar{F}^n(\C{C},\bar{F}^n(\C{B},D))$. The abutment is indeed
isomorphic to the cohomology of $\bar{F}^*(\C{C}\times\C{B},D)$ by
Lemma \ref{ibi} and Theorem \ref{ezc}.
\end{proof}

This spectral sequence can be regarded as a particular case of the
spectral sequence constructed in \cite{cgc}.

\begin{rem}\label{ibib}
For any category $\C{C}$ the category of pairs $\C{Pair}(\C{C})$ introduced in Section \ref{ahp} is nothing but
the category of functors $\mathbb{I}\r\C{C}$ where $\mathbb{I}$ is the unit interval category, with only two objects $0, 1$ and a unique
non-trivial morphism $0\r 1$. In particular there is an evaluation functor
$$\overline{ev}\colon\C{Pair}(\C{C})\times\mathbb{I}\To\C{C}.$$
For any $\C{C}$-bimodule $D$ there is an associated
$\C{Pair}(\C{C})$-bimodule $D^\#$ constructed in \cite{tbhp}
together with natural homomorphisms $$\lambda\colon
H^{n+1}(\C{C},D)\To H^n(\C{Pair}(\C{C}),D^\#),\;\;\; n\geq 0.$$ For
$n=2$ the homomorphism $\lambda$ plays a role in Section \ref{ahp}.
The reader can check by using the normalized cochain complex
$\bar{F}^*(\mathbb{I},-)$ that there is a natural identification
$$D^\#=H^1(\mathbb{I},\overline{ev}^*D).$$
Moreover, $H^n(\mathbb{I},-)$ vanishes for all $n>1$
and the homomorphisms
$$H^{n+1}(\C{C},D)\twoheadrightarrow E^{n,1}_\infty\hookrightarrow E^{n,1}_2=H^n(\C{Pair}(\C{C}),H^1(\mathbb{I},\overline{ev}^*D))=H^n(\C{Pair}(\C{C}),D^\#)$$
induced by the K\"unneth spectral sequence in Corollary \ref{ibia} coincide with $\lambda$ for all $n\geq 0$.
\end{rem}

Suppose now that $\C{C}$ is a category with zero object,
$(\C{B},s)$ is a translation category, $D$ is a
$(\C{C}\times\C{B})$-bimodule, and $\bar{s}\colon D\rr D(1\times
s,1\times s)$ is a natural transformation. Then there is a
morphism of cosimplicial $\C{C}$-bimodules
$$\bar{s}_*-s^*\colon F^\bullet(\C{B},D)\To F^\bullet(\C{B},D(1\times s,1\times s)).$$
By applying $F^\bullet(\C{C},-)$ we obtain now a bicosimplicial
abelian group homomorphism
$$F^\bullet(\C{C},\bar{s}_*-s^*)\colon F^\bullet(\C{C},F^\bullet(\C{B},D))\To F^\bullet(\C{C},F^\bullet(\C{B},D(1\times s,1\times s))).$$

On the other hand the functorial properties of the translation
cochain complex $\bar{F}(s,-)$ described in Remark \ref{fun} imply
that $\bar{F}^*(s,\bar{s})$ is a cochain complex of
$\C{C}$-bimodules, therefore if we apply $\bar{F}^*(\C{C},-)$ we obtain
a bicomplex
\begin{equation}\label{bc1}
\bar{F}^*(\C{C},\bar{F}^*(s,\bar{s})).
\end{equation}

The translation cochain complex
\begin{equation}\label{bc11}
\bar{F}^*(1_\C{C}\times s,\bar{s})
\end{equation}
is also defined.

In the next two lemmas we identify (\ref{bc11}) and the total
complex of (\ref{bc1}) as the homotopy
fibers of the two possible ways of converting
$F^\bullet(\C{C},\bar{s}_*-s^*)$ into a cochain homomorphism. For
this we recall that given a cochain morphism
$$\xi\colon A^*\To B^*$$
its \emph{homotopy fiber} $\fib(\xi)$ is the cochain complex
defined by
$$\fib^n(\xi)=A^n\oplus B^{n-1}$$
and differential
$$\left(\begin{array}{cc}d&0\\\xi&-d\end{array}\right).$$

\begin{lem}\label{ibi2}
There is a natural identification
$$\fib\tot N F^\bullet(\C{C},\bar{s}_*-s^*)=\tot\bar{F}^*(\C{C},\bar{F}^*(s,\bar{s})).$$
\end{lem}

\begin{proof}
We have the following equalities 
\begin{footnotesize}
\begin{eqnarray*}
\fib^n\tot N F^\bullet(\C{C},\bar{s}_*-s^*)&=& \tot^n N F^\bullet(\C{C},F^\bullet(\C{B},D))\\&&\oplus
\tot^{n-1} N F^\bullet(\C{C},F^\bullet(\C{B},D(1\times s,1\times s)))\\
&=&\bigoplus_{p+q=n}N^{p,q}F^\bullet(\C{C},F^\bullet(\C{B},D))\\&&\oplus\bigoplus_{p+q=n-1}N^{p,q} F^\bullet(\C{C},F^\bullet(\C{B},D(1\times s,1\times s)))\\
&=&\bigoplus_{p+q=n}\bar{F}^p(\C{C},\bar{F}^q(\C{B},D)) \\&&\oplus \bigoplus_{p+q=n-1}\bar{F}^p(\C{C},\bar{F}^q(\C{B},D(1\times s,1\times s)))\\
&=&\bigoplus_{p+q=n}\left(\bar{F}^p(\C{C},\bar{F}^q(\C{B},D)) \oplus \bar{F}^p(\C{C},\bar{F}^{q-1}(\C{B},D(1\times s,1\times s)))\right)\\
&=&\bigoplus_{p+q=n}\bar{F}^p(\C{C},\bar{F}^q(s,\bar{s}))\\&=&\tot^n\bar{F}^*(\C{C},\bar{F}^*(s,\bar{s})).
\end{eqnarray*}
\end{footnotesize}
We leave to the reader to check the compatibility with the differentials.
\end{proof}

\begin{lem}\label{ibi3}
There is a natural identification
$$\fib N\diag F^\bullet(\C{C},\bar{s}_*-s^*)=\bar{F}^*(1_\C{C}\times s,\bar{s}).$$
\end{lem}

\begin{proof}
We have the following identifications
\begin{eqnarray*}
\fib^n N\diag F^\bullet(\C{C},\bar{s}_*-s^*)&=& N^n\diag F^\bullet(\C{C},F^\bullet(\C{B},D))\\&&\oplus
N^{n-1}\diag F^\bullet(\C{C},F^\bullet(\C{B},D(1\times s,1\times s)))\\&\st{\text{(a)}}{=}&
N^n F^\bullet(\C{C}\times\C{B},D)\\&&\oplus N^{n-1} F^\bullet(\C{C}\times\C{B},D(1\times s, 1\times s))\\&=&
\bar{F}^n(\C{C}\times\C{B},D)\oplus  \bar{F}^{n-1}(\C{C}\times\C{B},D(1\times s, 1\times s))\\&=& \bar{F}^n(1_\C{C}\times s,\bar{s}).
\end{eqnarray*}
Here we use Lemma \ref{ibi} in (a). The reader can check that this is compatible with the differentials.
\end{proof}

By using these two lemmas and the Eilenberg-Zilber-Cartier theorem in the Appendix we obtain the following result.

\begin{prop}\label{lass}
Consider a category $\C{C}$ with zero object, a translation
category $(\C{B},s)$, a $(\C{C}\times\C{B})$-bimodule $D$, and a
natural transformation $\bar{s}\colon D\rr D(1\times s,1\times
s)$. Then there is a natural spectral sequence
$$E^{p,q}_2=H^p(\C{C},H^q(s,\bar{s}))\Longrightarrow H^{p+q}(1_\C{C}\times s,\bar{s}).$$
\end{prop}

\begin{proof}
This is the spectral sequence of the bicomplex (\ref{bc1}). Lemmas
\ref{ibi2} and \ref{ibi3} together with the naturality of the
cochain homotopy equivalences in Theorem \ref{ezc}, the long exact
sequence for the cohomology of a fibration of cochain complexes,
and the five lemma imply that the cochain homomorphism given by
the matrices
$$\left(%
\begin{array}{cc}
  \mathfrak{F}_n & 0 \\
  0 & \mathfrak{F}_{n-1} \\
\end{array}%
\right),\;\;\; n\in\Z,$$ induces an isomorphism between the
cohomology of (\ref{bc11}) and the cohomology of the total complex
of (\ref{bc1}), therefore the abutment of the spectral sequence is
as in the statement of this proposition.
\end{proof}

We will write $$0=D^{n+1,-1}\subset D^{n,0}\subset\dots\subset
D^{p,q}\subset\cdots\subset D^{0,n}=H^n(1_\C{C}\times s,\bar{s})$$
for the filtration such that $E^{p,q}_\infty=D^{p,q}/D^{p+1,q-1}$.

\begin{rem}\label{prinex}
We are specially interested in the case where $(\C{B},s)=(\V,sh^3)$ and $\C{C}=\C{Cand}(\C{A},t)$ for some
translation category $(\C{A},t)$
as in Section \ref{nablat}. There is a translation category morphism
$$(ev,\epsilon)\colon\C{Cand}(\C{A},t)\times (\V,sh^3)\To(\C{A},t)$$
defined by the zero-object-preserving ``evaluation'' functor
$$ev\colon\C{Cand}(\C{A},t)\times\V\To\C{A}$$
given by $ev((X,\alpha),n)=X_n$ on objects and given a morphism
$f\colon(X,\alpha)\r(Y,\beta)$ of candidate triangles
$ev(f,1_n)=f_n$, $ev(1_{(X,\alpha)},\imath_n)=d_n$, and the
natural transformation
$$\epsilon\colon t\; ev\rr ev(1_{\C{Cand}(\C{A},t)}\times sh^3)$$
with $\epsilon_{((X,\alpha),n)}=\alpha_n\colon tX_n\r X_{n+3}$.
We therefore can pull back the natural transformation $\bar{t}\colon\hom^t\rr\hom^t(t,t)$, see (\ref{homt}) and (\ref{tbar}), along $(ev,\epsilon)$
$$(ev,\epsilon)^*\bar{t}\colon\hom_\C{A}(t\; ev,ev)\rr\hom_\C{A}(t\; ev(1_{\C{Cand}(\C{A},t)}\times sh^3),ev(1_{\C{Cand}(\C{A},t)}\times sh^3)).$$
If we evaluate this natural transformation of $(\C{Cand}(\C{A},t)\times\V)$-bimodules on two candidate triangles $(X,\alpha)$ and $(Y,\beta)$, leaving
free the second variable, we obtain the following natural transformation of $\V$-bimodules
$$((ev,\epsilon)^*\bar{t})_{(((X,\alpha),-),((Y,\beta),-))}=[(X,\alpha),(Y,\beta)]^*\bar{t},$$
see (\ref{indt2}), hence we can write
$$(ev,\epsilon)^*\bar{t}=[-,-]^*\bar{t}.$$
\end{rem}

This remark and Proposition \ref{lass} yield a spectral sequence
\begin{equation}\label{sscand}
E^{p,q}_\infty=H^p(\C{Cand}(\C{A},t),H^q(sh^3,[-,-]^*\bar{t}))\Longrightarrow
H^{p+q}(1_{\C{Cand}(\C{A},t)}\times sh^3,[-,-]^*\bar{t}).
\end{equation}
In particular there are homomorphisms $(n\geq 0)$
\begin{equation}\label{xi0}
\xi_0\colon H^n(1_{\C{Cand}(\C{A},t)}\times
sh^3,[-,-]^*\bar{t})\To
H^0(\C{Cand}(\C{A},t),H^n(sh^3,[-,-]^*\bar{t}))
\end{equation}
given by
$$D^{0,n}\twoheadrightarrow D^{0,n}/D^{1,n-1}=E^{0,n}_\infty\subset E^{0,n}_2.$$
The homomorphisms $\xi_0$ satisfy the property stated in the
following lemma, see (\ref{diag}).

\begin{lem}\label{elpul}
In the conditions above, given $b\in H^n(\C{A},t)$ and a candidate
triangle $(X,\alpha)$ in $(\C{A},t)$ we have that
$$\xi_0((ev,\epsilon)^*b)(X,\alpha)=(X,\alpha)^*b\in H^n(sh^3,(X,\alpha)^*\bar{t}).$$
\end{lem}

\begin{proof}
The homomorphisms $\xi_0$ are given as follows. Let
\begin{eqnarray*}
\bar{c}=(\bar{c}_1,\bar{c}_2)&\in&\bar{F}^n(1_{\C{Cand}(\C{A},t)}\times sh^3,[-,-]^*\bar{t})\\&=&
\bar{F}^n(\C{Cand}(\C{A},t)\times\V,\hom_\C{A}(t\; ev,ev))\\&&\oplus
\bar{F}^{n-1}(\C{Cand}(\C{A},t)\times\V,\hom_\C{A}(t\; ev(1\times sh^3),ev(1\times sh^3)))
\end{eqnarray*}
be a cocycle representing $c\in H^n(1_{\C{Cand}(\C{A},t)}\times sh^3,[-,-]^*\bar{t})$. We are going to use the identification
in Lemma \ref{ibi} 
as well as the homotopy equivalence $\mathfrak{F}$ in Theorem
\ref{ezc}. The cochains
\begin{eqnarray*}
\mathfrak{F}(\bar{c}_1)&\in& \tot^n\bar{F}^*(\C{Cand}(\C{A},t),\bar{F}^*(\V,\hom_\C{A}(t\; ev,ev)))\\
&=&\bigoplus_{p+q=n}\bar{F}^p(\C{Cand}(\C{A},t),\bar{F}^q(\V,\hom_\C{A}(t\; ev,ev))),\\
\mathfrak{F}(\bar{c}_2)&\in& \tot^{n-1}\bar{F}^*(\C{Cand}(\C{A},t),\bar{F}^*(\V,\hom_\C{A}(t\; ev(1\times sh^3),ev(1\times sh^3))))\\
&=&\bigoplus_{p+q=n-1}\bar{F}^p(\C{Cand}(\C{A},t),\bar{F}^q(\V,\hom_\C{A}(t\; ev(1\times sh^3),ev(1\times sh^3)))),
\end{eqnarray*}
have components
\begin{eqnarray*}
\mathfrak{F}(\bar{c}_1)^{0,n}&\in& \bar{F}^0(\C{Cand}(\C{A},t),\bar{F}^n(\V,\hom_\C{A}(t\; ev,ev))),\\
\mathfrak{F}(\bar{c}_2)^{0,n-1}&\in&
\bar{F}^0(\C{Cand}(\C{A},t),\bar{F}^{n-1}(\V,\hom_\C{A}(t\;
ev(1\times sh^3),ev(1\times sh^3)))).
\end{eqnarray*}
The element $\xi_0(c)(X,\alpha)$ is represented by the cocycle
$$(\mathfrak{F}(\bar{c}_1)^{0,n}(X,\alpha),\mathfrak{F}(\bar{c}_2)^{0,n-1}(X,\alpha))\in \bar{F}^n(sh^3,[-,-]^*\bar{t}).$$
By using the definition of $\mathfrak{F}$ in the statement of
Theorem \ref{ezc} we see that this cochain is given as follows
\begin{eqnarray*}
(\mathfrak{F}(\bar{c}_1)^{0,n}(X,\alpha))(\sigma_1,\dots,\sigma_n)&=&\bar{c}_1((1_{(X,\alpha)},\sigma_1),\dots,(1_{(X,\alpha)},\sigma_n)),\\
(\mathfrak{F}(\bar{c}_2)^{0,n-1}(X,\alpha))(\sigma_1,\dots,\sigma_{n-1})&=&\bar{c}_2((1_{(X,\alpha)},\sigma_1),\dots,(1_{(X,\alpha)},\sigma_{n-1})).
\end{eqnarray*}
Suppose now that $$\bar{b}=(\bar{b}_1,\bar{b}_2)\in
\bar{F}^n(t,\bar{t})=\bar{F}^n(\C{A},\hom_\C{A}(t,-))\oplus
\bar{F}^{n-1}(\C{A},\hom_\C{A}(t^2,t))$$ is a cocycle representing
$b\in H^n(\C{A},t)$. Then $c=(ev,\epsilon)^*b$ is represented by
$$\bar{c}=(\bar{c}_1,\bar{c}_2)=(ev,\epsilon)^*\bar{b}=(ev^*\bar{b}_1,\hom_\C{A}(t\epsilon^{-1},1)_*\epsilon^\#(\bar{b}_1)
+\hom_\C{A}(t\epsilon^{-1},\epsilon)_*ev^*\bar{b}_2),$$ and using
the definition of $(ev,\epsilon)$ in Remark \ref{prinex} we obtain
\begin{small}
\begin{eqnarray*}
(ev^*\bar{b}_1)((1_{(X,\alpha)},\sigma_1),\dots,(1_{(X,\alpha)},\sigma_n))&=&\bar{b}_1(X(\sigma_1),\dots,X(\sigma_n))\\
&=&(X^*\bar{b}_1)(\sigma_1,\dots,\sigma_n),\\
(\hom_\C{A}(t\epsilon^{-1}_{((X,\alpha),k-n+1)},1)_*\epsilon^\#(\bar{b}_1))&&\\((1_{(X,\alpha)},\sigma_1),\dots,
(1_{(X,\alpha)},\sigma_{n-1}))&=&
\sum_{i=0}^{n-1}(t\alpha_{k-n+1}^{-1})^*(-1)^i\bar{b}_1(X
sh^3(\sigma_1),\dots,\\&&X sh^3(\sigma_i),\alpha_{k-i},
tX(\sigma_{i+1}),\dots,tX(\sigma_{n-1}))\\
&=&\hom_\C{A}(t\alpha^{-1}_{k-n+1},1)_*\alpha^\#(\sigma_1,\dots,\sigma_{n-1}),\\
(\hom_\C{A}(t\epsilon^{-1}_{((X,\alpha),k-n+1)},\epsilon_{((X,\alpha),k)})_*
ev^*\bar{b}_2)&&\\((1_{(X,\alpha)},\sigma_1),\dots,(1_{(X,\alpha)},\sigma_{n-1}))&=&
(t\alpha^{-1}_{k-n+1})^*{\alpha_k}_*\bar{b}_2(X(\sigma_1),\dots,X(\sigma_{n-1}))\\
&=&(\hom_\C{A}(t\alpha_{k-n+1},\alpha_k)_*X^*\bar{b}_2)(\sigma_1,\dots,\sigma_{n-1}),\\
\end{eqnarray*}
\end{small}
where $k-i$ is the source of $\sigma_i$ and/or the target of
$\sigma_{i+1}$. Hence we observe that
$$(\mathfrak{F}(\bar{c}_1)^{0,n}(X,\alpha),\mathfrak{F}(\bar{c}_2)^{0,n-1}(X,\alpha))=(X,\alpha)^*\bar{b},$$
i. e. the cocycle representative of
$\xi_0((ev,\epsilon)^*b)(X,\alpha)$ also represents $(X,\alpha)^*b$,
hence the proof is finished.
\end{proof}

We now consider the bicomplex
\begin{equation}\label{bc1*}
\tilde{F}^*(\C{Cand}(\C{A},t),\tilde{F}^*(sh^3,[-,-]^*\bar{t})).
\end{equation}
which is the $*$-normalized version of (\ref{bc1}) for the special
case considered in Remark \ref{prinex}. In bidegree $(p,q)$ for
$q>0$ this bicomplex is
$$\begin{array}{l}
\tilde{F}^p(\C{Cand}(\C{A},t),\prod_{m\in\Z}\hom_\C{A}(t\;
ev(-,m),ev(-,m+q)))
\\
\oplus\tilde{F}^p(\C{Cand}(\C{A},t),\prod_{m\in\Z}\hom_\C{A}(t\;
ev(-,m+3),ev(-,m+q+2))).
\end{array}$$
We now define the $(0,3)$-cochain
$$\widetilde{\id}=(\overline{\id},0)\in\tilde{F}^0(\C{Cand}(\C{A},t),\tilde{F}^3(sh^3,[-,-]^*\bar{t}))$$
by
$$\overline{\id}(X,\alpha)(m)=(-1)^m\alpha_m\in\hom_\C{A}(tX_m,X_{m+3})$$
for any candidate triangle $(X,\alpha)$ and $m\in\Z$.

\begin{prop}
The cochain $\widetilde{\id}$ is a cocycle in the total complex of
(\ref{bc1*}).
\end{prop}

\begin{proof}
The horizontal differential $d_h$ of (\ref{bc1*}) is given by the
cohomology of $\C{Cand}(\C{A},t)$ and the vertical differential
$d_v$ is determined by the translation cohomology of $(\V,sh^3)$.
Therefore given a morphism of candidate triangles
$f\colon(X,\alpha)\r(Y,\beta)$ we have
$$(d_h\widetilde{\id})(f)=((\delta\overline{\id})(f),0)$$
and for any $m\in \Z$
\begin{eqnarray*}
((\delta\overline{\id})(f))(m)&=&{f_{m+3}}_*(\overline{\id}(X,\alpha))(m)-f_{m}^*(\overline{\id}(Y,\beta))(m)\\
&=&f_{m+3}\alpha_m-\beta_m(tf_m)\\&=&0.
\end{eqnarray*}
This vanishes because all morphisms of candidate triangles satisfy
$\beta(t f)=f\alpha$, see Section \ref{ct}.

On the other hand for any candidate triangle $(X,\alpha)$ we have
$$(d_v(\widetilde{\id}))(X,\alpha)=(\delta(\overline{\id}(X,\alpha)),
((X,\alpha)^*\bar{t})_*(\overline{\id}(X,\alpha))-{sh^3}^*(\overline{\id}(X,\alpha))).$$
The first coordinate is
\begin{eqnarray*}
(\delta(\overline{\id}(X,\alpha)))(m)&=&{d_{m+4}}_*\overline{\id}(X,\alpha)(m)+d_{m+1}^*\overline{\id}(X,\alpha)(m+1)\\
&=&(-1)^{m}d_{m+4}\alpha_m+(-1)^{m+1}\alpha_{m+1}(td_{m+1})\\&=&0
\end{eqnarray*}
This vanishes since $\alpha$ is a degree $3$ cochain homomorphism.
For the second coordinate by (\ref{indt2}) we have
\begin{eqnarray*}
(((X,\alpha)^*\bar{t})_*(\overline{\id}(X,\alpha)))(m)&&\\-({sh^3}^*(\overline{\id}(X,\alpha)))(m)
&=&(-1)^{m}\alpha_{m+3}\bar{t}_{(X_m,X_{m+3})}(\alpha_m)(t\alpha^{-1}_m)\\&&-(-1)^{m+3}\alpha_{m+3}\\
&=&-(-1)^{m}\alpha_{m+3}(t\alpha_m)(t\alpha^{-1}_m)-(-1)^{m+3}\alpha_{m+3}\\
&=&0.
\end{eqnarray*}
Now the proof is finished.
\end{proof}

The inclusion of bicomplexes
$$\tilde{F}^*(\C{Cand}(\C{A},t),\tilde{F}^*(sh^3,[-,-]^*\bar{t}))\subset \bar{F}^*(\C{Cand}(\C{A},t),\bar{F}^*(sh^3,[-,-]^*\bar{t}))$$
induces a quasi-isomorphism in total complexes since we know that
the corresponding map of spectral sequences induces isomorphisms
in the $E_2$-term. Hence, by using Lemmas \ref{ibi2} and
\ref{ibi3} and the cosimplicial Eilenberg-Zilber-Cartier theorem
in the Appendix, we see that the cochain $\widetilde{\id}$
represents a cohomology class
\begin{equation}\label{clasi}
\id\in H^3(1_{\C{Cand}(\C{A},t)}\times sh^3,[-,-]^*\bar{t}).
\end{equation}
We will refer to it as the \emph{global class of
identities}\index{identities!global class of -}.

The following lemma proves that the global class of identities determines
all  classes of identities defined in (\ref{wi}) through the
homomorphism $\xi_0$ in (\ref{xi0}).

\begin{lem}\label{elpul2}
For any candidate triangle $(X,\alpha)$ we have
$(\xi_0(\id))(X,\alpha)=\idd(X,\alpha)$.
\end{lem}

\begin{proof}
The cohomology class $(\xi_0(\id))(X,\alpha)$ is represented by the
cochain
$$(\overline{\id}(X,\alpha),0)\in\tilde{F}^3(sh^3,(X,\alpha)^*\bar{t}).$$
This cochain is the same as the cochain in Lemma \ref{laotra},
 hence the proof is
finished.
\end{proof}

Let
$i_\nabla\colon\C{Cand}(\nabla)\hookrightarrow\C{Cand}(\C{A},t)$
be the inclusion of the full subcategory of $\nabla$-triangles. By
definition of $\nabla$-triangle and Lemma \ref{elpul} if
$b=(b_1,b_2)\in\tilde{F}^3(t,\bar{t})$ is a representative of
$\nabla\in H^3(\C{A},t)$ then there exists a cochain
$$e=(e_1,e_2)\in\tilde{F}^0(\C{Cand}(\nabla),\tilde{F}^2(sh^3,[i_\nabla,i_\nabla]^*\bar{t}))$$
such that for any $\nabla$-triangle $(X,\alpha)$ the equation
\begin{equation}\label{aysh}
\delta(e_1(X,\alpha),e_2(X,\alpha))=\tilde{\id}(X,\alpha)+(X,\alpha)^*b\in\tilde{F}^3(sh^3,(X,\alpha)^*\bar{t})
\end{equation}
holds.

\begin{lem}\label{aysh2}
Let $e$ be a cochain as above. The cohomology class
$\xi_1(i_\nabla\times 1_\V,0^\vc_{i_\nabla\times
sh^3})^*(\id+(ev,\epsilon)^*\nabla)$ is represented by the
$1$-cochain
$$\theta\in\tilde{F}^1(\C{Cand}(\nabla),H^2(sh^3,[i_\nabla,i_\nabla]^*\bar{t}))$$
such that for any morphism of candidate triangles
$l\colon(X,\alpha)\r(Y,\beta)$ the cohomology class $\theta(l)\in
H^2(sh^3,[(X,\alpha),(Y,\beta)]^*\bar{t})$ is represented by the
$2$-cocycle
$\bar{\theta}(l)\in\tilde{F}^2(sh^3,[(X,\alpha),(Y,\beta)]^*\bar{t})$
given by
\begin{eqnarray*}
(\bar{\theta}(l))(m)&=&(b_1(d_{m+2},d_{m+1},l_m)-b_1(d_{m+2},l_{m+1},d_{m+1})\\&&+b_1(l_{m+2},d_{m+2},d_{m+1})\\
&&-l_{m+2}(e_1(X,\alpha)(m))+(e_1(Y,\beta)(m))(tl_m),\\&&
b_1(d_{m+4},\beta_m,tl_m)(t\alpha_m^{-1})-b_1(\beta_{m+1},td_{m+1},tl_m)(t\alpha_m^{-1})\\&&+b_1(\beta_{m+1},tl_{m+1},td_{m+1})(t\alpha_m^{-1})\\
&&-b_1(d_{m+4},l_{m+3},\alpha_m)(t\alpha_m^{-1})+b_1(l_{m+4},d_{m+4},\alpha_m)(t\alpha_m^{-1})\\&&-b_1(l_{m+4},\alpha_{m+1},td_{m+1})(t\alpha_m^{-1})\\
&&-\beta_{m+1}b_2(d_{m+1},l_m)(t\alpha_m^{-1})+\beta_{m+1}b_2(l_{m+1},d_{m+1})(t\alpha_m^{-1})\\
&&-l_{m+1}(e_2(X,\alpha)(m))+(e_2(Y,\beta)(m))(tl_m)).
\end{eqnarray*}
\end{lem}

The proof of this lemma is technically just a continuation of the
proof of Lemma \ref{elpul}. One only needs to use again two
elementary but tedious formulas, namely the homotopy equivalence
$\mathfrak{F}$ in Theorem \ref{ezc} and the construction of the
spectral sequence of a bicomplex. We leave the details to the
reader.

\section{The axiom (Tr2)}\label{tr2}

Let $(\C{A},t)$ be a translation category with $\C{A}$ additive
and  $t$ an additive equivalence. There is an isomorphism
$$\ominus\colon\C{Cand}(\C{A},t)\To\C{Cand}(\C{A},t)$$
defined on objects as
$$\ominus(X,\alpha)=(1_\C{A},-0^\vc_{t})(X,\alpha)(sh,0^\vc_{sh^4})=(X
sh,-\alpha \,sh)$$ and on morphisms $f$ as $\ominus(f)=f\,sh$.
This endofunctor is related to the axiom (Tr2) in the definition of pretriangulated categories, see \cite{ccctc1} 1. More precisely, if
$$A\st{f}\To B\st{i}\To C\st{q}\To tA$$
is the small candidate triangle corresponding to $(X,\alpha)$ by
the equivalence of categories in (\ref{sb}), then the small
candidate triangle corresponding to $\ominus(X,\alpha)$ is
isomorphic to
$$ B\st{i}\To C\st{q}\To tA\st{-tf}\To tB.$$
The isomorphism is given by the following commutative diagram in
$\C{A}$
$$\xymatrix{B\ar[r]^i\ar@{=}[d]&C\ar[r]^q\ar@{=}[d]&tA\ar[r]^{-tf}\ar@{=}[d]&tB\ar@{=}[d]\\
X_{-1}\ar[r]^{d_0}\ar@{=}[d]&X_0\ar[r]^{\alpha^{-1}_{-2}d_1}\ar@{=}[d]&tX_{-2}\ar[r]^{-td_{-1}}\ar[d]^{\alpha_{-2}}_\cong&tX_{-1}\ar@{=}[d]\\
X_{-1}\ar[r]^{d_0}&X_{0}\ar[r]^{d_1}&X_1\ar[r]^{-\alpha^{-1}_{-1}d_2}&tX_{-1}}$$

\begin{prop}\label{tr2na}
Let $(X,\alpha)$ be a candidate triangle, then $(X,\alpha)$ is a
$\nabla$-triangle if and only if $\ominus(X,\alpha)$ is a
$\nabla$-triangle.
\end{prop}

For the proof of this proposition we need the following two
lemmas.

\begin{lem}\label{manolo}
The homomorphism induced by $(1_\C{A},-0^\vc_{t})$ on translation
cohomology is the identity
$$1=(1_\C{A},-0^\vc_{t})^*\colon H^*(\C{A},t)\To H^*(\C{A},t).$$
\end{lem}

\begin{proof}
We have to prove that the homomorphism given by the identity matrix is homotopic to
\begin{eqnarray*}
(1_\C{A},-0^\vc_{t})^*&=&\left(\begin{array}{cc}
1&0\\
  (\hom^t(-0^\vc_{t},0^\vc_{t}))_*(-0^\vc_{t})^\# & (\hom^t(-0^\vc_{t},-0^\vc_{t}))_*
\end{array}%
\right)\\&=&\left(\begin{array}{cc}
1&0\\
  -(-0^\vc_{t})^\# & 1
\end{array}%
\right).
\end{eqnarray*}

First of all we restrict ourselves to cochains which are
normalized not only with respect to identities and zero morphisms
but also with respect to products in $\C{A}$ in the sense of
\cite{sncc}. By using \cite{sncc} A.11 we see there is no lose of
generality in this restriction.

We claim that a homotopy is given by
$$\left(\begin{array}{cc}
0&0\\
((0^\vc_t,0^\vc_t),\binom{0^\vc_t}{-0^\vc_t})^\%-((0^\vc_t,0^\vc_t),\binom{0^\vc_t}{0})^\%-((0^\vc_t,0^\vc_t),\binom{0}{-0^\vc_t})^\%& 0
\end{array}\right).$$
Here $(-,-)^\%$ are the cochain homotopies in \cite{ptc} 5 and
$0\colon t\rr t$ is the natural transformation given by zero
morphisms.

We have the following equalities
\begin{scriptsize}
\begin{eqnarray*}
\delta((0^\vc_t,0^\vc_t),\binom{0^\vc_t}{-0^\vc_t})^\%-((0^\vc_t,0^\vc_t),\binom{0^\vc_t}{-0^\vc_t})^\%\delta&&\\
-\delta((0^\vc_t,0^\vc_t),\binom{0^\vc_t}{0})^\%+((0^\vc_t,0^\vc_t),\binom{0^\vc_t}{0})^\%\delta&&\\
-\delta((0^\vc_t,0^\vc_t),\binom{0}{-0^\vc_t})^\%+((0^\vc_t,0^\vc_t),\binom{0}{-0^\vc_t})^\%\delta &\st{\text{(a)}}{=}&
-\hom^t(\binom{0^\vc_t}{-0^\vc_t},0^\vc_t)_*(0^\vc_t,0^\vc_t)^\#\\&&-\hom^t(0^\vc_t,(0^\vc_t,0^\vc_t))_*\binom{0^\vc_t}{-0^\vc_t}^\#+0^\#\\&&
+\hom^t(\binom{0^\vc_t}{0},0^\vc_t)_*(0^\vc_t,0^\vc_t)^\#\\&&+\hom^t(0^\vc_t,(0^\vc_t,0^\vc_t))_*\binom{0^\vc_t}{0}^\#-(0^\vc_t)^\#\\&&
+\hom^t(\binom{0}{-0^\vc_t},0^\vc_t)_*(0^\vc_t,0^\vc_t)^\#\\&&+\hom^t(0^\vc_t,(0^\vc_t,0^\vc_t))_*\binom{0}{-0^\vc_t}^\#-(0^\vc_t)\#\\&\st{\text{(b)}}=&
-\hom^t(0^\vc_t,(0^\vc_t,0^\vc_t))_*\binom{0^\vc_t}{-0^\vc_t}^\#+0^\#\\&&
+\hom^t(0^\vc_t,(0^\vc_t,0^\vc_t))_*\binom{0^\vc_t}{0}^\#-(0^\vc_t)^\#\\&&
+\hom^t(0^\vc_t,(0^\vc_t,0^\vc_t))_*\binom{0}{-0^\vc_t}^\#-(0^\vc_t)\#\\
&\st{\text{(c)}}=&
 -(-0^\vc_t)^\#.
\end{eqnarray*}
\end{scriptsize}
Here we use \cite{ptc} (5.11) for (a) and the biadditivity of
$\hom^t$ in (b). Normalization with respect to identities shows
that $(0^\vc_t)^\#=0$, see \cite{ptc} (5.7), normalization with
respect to zero morphisms implies that $0^\#=0$, and these and
normalization with respect to products show that $$\left(%
\begin{array}{c}
  0^\vc_t \\
  0
\end{array}%
\right)^\#=0, \left(\begin{array}{c}
  0\\
  -0^\vc_t
\end{array}%
\right)^\#= \left(\begin{array}{c}
  0\\
  (-0^\vc_t)^\#
\end{array}%
\right)\text{ and } \left(\begin{array}{c}
  0^\vc_t \\
  -0^\vc_t
\end{array}%
\right)^\#= \left(\begin{array}{c}
  0 \\
  (-0^\vc_t)^\#
\end{array}%
\right),$$ therefore (c) follows easily.

The previous chain of equalities proves the claim, and hence the
lemma.
\end{proof}

\begin{lem}\label{rosario}
For any candidate triangle $(X,\alpha)$ the equality
$$(sh,0^\vc_{sh^4})^*\idd(X,\alpha)=\idd(\ominus(X,\alpha))$$ holds in $ H^3(sh^3,(sh,0^\vc_{sh^4})^*(X,\alpha)^*\bar{t})=H^3(sh^3,(\ominus(X,\alpha))^*\bar{t})$.
\end{lem}

\begin{proof}
By Lemma \ref{laotra} $\idd(X,\alpha)$ and
$\idd(\ominus(X,\alpha))$ are represented by the translation
cocycles $(b,0)$ and $(c,0)$, respectively, defined as
$$b(m)=(-1)^m\alpha_m\text{ and } c(m)=(-1)^m(-\alpha_{m+1}),$$
hence the equality already holds at the level of representative
cocycles $(sh,0^\vc_{sh^4})^*(b,0)=(c,0)$.
\end{proof}

\begin{proof}[Proof of Proposition \ref{tr2na}]
By Lemmas \ref{manolo} and \ref{rosario} we have that
$$\idd(\ominus(X,\alpha))+(\ominus(X,\alpha))^*\nabla=(sh^3,0^\vc_{sh^4})^*(\idd(X,\alpha)+(X,\alpha)^*\nabla)$$
and $(sh^3,0^\vc_{sh^4})^*$ is an isomorphism because
$(sh^3,0^\vc_{sh^4})$ is already an isomorphism, therefore the
proposition follows.
\end{proof}

\section{Proofs of Propositions \ref{reao} and \ref{lader} and Theorems \ref{ñi} and \ref{vanish}}\label{pruebas}

In this section we concentrate all proofs where the use of the
track-theoretical techniques developed in \cite{ptc} and
\cite{ccctc1} is strictly necessary. We will freely use the
definitions and notation in \cite{ccctc1} 2 and 3.

In this section $(\C{B},s)$ will always be a good translation track category representing a translation cohomology class $\nabla\in H^3(\C{A},t)$.

In the proof of the following results we will use a fixed global
section $(\ell,\mu,\nu,\eta)$ in the sense of \cite{ptc} 14. We are
actually interested in $\ell$, $\mu$ and part of $\nu$ only. This is
given by a choice of a morphism $\ell f\colon A\r B$ in $\C{B}$ for
any $f\colon A\r B$ in $\C{A}$ with $p(\ell f)=f$ and
$\ell(1_X)=1_X$, tracks $\mu_{f,g}\colon(\ell f)(\ell g)\rr\ell(fg)$
with $\mu_{f,1}=0^\vc_{\ell f}=\mu_{1,f}$, and tracks
$\nu_{f,s}\colon s(\ell f)\rr\ell(tf)$ with $\nu_{1_X,s}=s_X$. Since
$\C{B}$ has a strict zero object $*$ which is preserved by $s$ and
$s$ is normalized at zero maps we can also impose here $\ell(0)=0$,
$\mu_{f,0}=0^\vc_{0}=\mu_{0,f}$ and $\nu_{0,s}=0^\vc_0$.  We recall
from \cite{ptc} that $\nabla$ is represented by the cocycle
$(c_\C{B},b_s)\in\tilde{F}^3(t,\bar{t})$ with
\begin{eqnarray*}
c_\C{A}({f},{g},{h})&=&\sigma_{\ell (fgh)}^{-1}(\mu_{{fg},{h}}\vc
(\mu_{{f},{g}}(\ell {h}))
\vc((\ell {f})\mu_{{g},{h}}^\vi)\vc\mu_{{f},{gh}}^\vi),
\\
b_s ({f},{g})&=&\sigma^{-1}_{\ell { (t(fg))}}(\mu_{{t (f)},{t
(g)}}\vc (\nu_{{f},s }\nu_{{g},s })
\vc s^\vi_{\ell {f},\ell {g}}\vc s (\mu_{{f},{g}}^\vi)
\vc\nu^\vi_{{fg},s }).
\end{eqnarray*}
The cocycle $(c_\C{B},b_s)$ is in principle only normalized with
respect to identities, but the extra conditions required here to the
global section clearly imply that it is also normalized with respect
to zero morphisms. We will also need a collection $\gamma$ of tracks
$\gamma_f\colon\ell p(f)\rr f$ indexed by all maps $f$ in $\C{B}$
such that $\gamma_{\ell g}=0^\vc_{\ell g}$ for any morphism $g$ in
$\C{A}$.

\begin{proof}[Proof of Proposition \ref{reao}]
Let
\begin{equation*}\tag{a}
\xymatrix{A\ar@/^30pt/[rr]^0_{\;}="a"\ar[r]_{f}&B\ar@/_30pt/[rr]_0^{\;}="c"
\ar[r]^<(-0.22){\;}="b"|{i}&
C\ar[r]_<(-0.17){\;}="d"|{q}\ar@/^30pt/[rr]^0_{\;}="f"&sA\ar[r]_{s(f)}^<(-0.23){\;}="e"&
sB\ar@{=>}"b";"a"_{H_0}\ar@{=>}"d";"c"^{H_{1}}\ar@{=>}"e";"f"_{H_{2}}}
\end{equation*}
be a track triangle and let $(X,\alpha)$ be the candidate triangle
associated to the small candidate triangle
\begin{equation*}\tag{a}
A\st{p(f)}\To
B\st{p(i)}\To C\st{p(q)}\To tA
\end{equation*}
under the equivalence (\ref{sb'}). The reader can check that the
coboundary of the $2$-cochain
\begin{eqnarray*}
(\sigma_0^{-1}(H_0\vc(\gamma_i\gamma_f)\vc\mu^\vi_{p(i),p(f)}),\sigma_0^{-1}(H_1\vc(\gamma_q\gamma_i)\vc\mu^\vi_{p(q),p(i)}),\\
\sigma_0^{-1}(H_2\vc(s(\gamma_f)\gamma_q)\vc(\nu^\vi_{p(f),s}(\ell p(q)))\vc\mu^\vi_{tp(f),p(q)}))&\in&\tilde{G}^2((X,\alpha)^*\bar{t})
\end{eqnarray*}
is
\begin{small}
$$\begin{array}{l}(1_{tA},-1_{tB},1_{tC})\\+(c_\C{B}(p(q),p(i),p(f)),c_\C{B}(tp(f),p(q),p(i)),c_\C{B}(tp(i),tp(f),p(q))-b_s(p(i),p(f))(tp(q))),\end{array}$$
\end{small}
and by Lemma \ref{peque} and Remark \ref{cla1} this represents the cohomology class $\idd(X,\alpha)+(X,\alpha)^*\nabla$, therefore $(X,\alpha)$ is a  $\nabla$-triangle.

On the other hand, if the ``long'' candidate triangle associated to the small one
\begin{equation*}
A\st{f}\To
B\st{i}\To C\st{q}\To tA
\end{equation*}
is a small $\nabla$-triangle we can choose a $2$-cochain $(e_{-2},e_{-1},e_0)\in\tilde{G}^2((X,\alpha)^*\bar{t})$ such that
$$\delta(e_{-2},e_{-1},e_0)=
(1_{tA},-1_{tB},1_{tC})+(c_\C{B}(q,i,f),c_\C{B}(tf,q,i),c_\C{B}(ti,tf,q)-b_s(i,f)(tq)).$$
Now it is not difficult to check that
$$\xymatrix{A\ar@/^30pt/[rr]^0_{\;}="a"\ar[r]_{\ell {f}}&B\ar@/_30pt/[rr]_0^{\;}="c"
\ar[r]^<(-0.22){\;}="b"|{\ell {i}}&
C\ar[r]_<(-0.17){\;}="d"|{\ell {q}}\ar@/^30pt/[rr]^0_{\;}="f"&sA\ar[r]_{s(\ell f)}^<(-0.23){\;}="e"&
sB\ar@{=>}"b";"a"_{H_0}\ar@{=>}"d";"c"^{H_1}\ar@{=>}"e";"f"_{H_2}}$$
with
\begin{eqnarray*}
H_0&=&\sigma_0(e_{-2})\vc\mu_{i,f},\\
H_1&=&\sigma_0(e_{-1})\vc\mu_{q,i},\\
H_2&=&\sigma_0(e_0)\vc\mu_{tf,q}\vc(\nu_{f,s} (\ell q)),
\end{eqnarray*}
is a track triangle.
\end{proof}


\begin{proof}[Proof of Proposition \ref{lader}]

As we remark in the paragraph preceding the statement of Lemma
\ref{aysh2} there is a cochain
$$e=(e_1,e_2)\in\tilde{F}^0(\C{Cand}(\nabla),\tilde{F}^2(sh^3,[-,-]^*\bar{t}))$$
such that for any $\nabla$-triangle $(X,\alpha)$ the equation
\begin{equation*}
\delta(e_1(X,\alpha),e_2(X,\alpha))=\widetilde{\id}(X,\alpha)+(X,\alpha)^*(c_\C{B},b_s)
\end{equation*}
is satisfied. Suppose that $(X,\alpha)$ is the image of the small
$\nabla$-triangle $T$ given by
$$A\st{f}\To B\st{i}\To C\st{q}\To tA$$
by the equivalence of categories in (\ref{sb'}). The previous
equation implies that in the cochain complex
$\tilde{G}^*((X,\alpha)^*\bar{t})$ the coboundary of
$$(e_1(X,\alpha)(-2),e_1(X,\alpha)(-1),e_1(X,\alpha)(0)-e_2(X,\alpha)(-2)(td_1))$$
is
$$(1_{tA},-1_{tB},1_{tC})+(c_\C{B}(q,i,f),c_\C{B}(tf,q,i),c_\C{B}(ti,tf,q)-b_s(i,f)(tq)),$$
therefore the diagram
$$\xymatrix{A\ar@/^30pt/[rr]^0_{\;}="a"\ar[r]_{\ell {f}}&B\ar@/_30pt/[rr]_0^{\;}="c"
\ar[r]^<(-0.22){\;}="b"|{\ell {i}}& C\ar[r]_<(-0.17){\;}="d"|{\ell
{q}}\ar@/^30pt/[rr]^0_{\;}="f"&sA\ar[r]_{s(\ell
f)}^<(-0.23){\;}="e"&
sB\ar@{=>}"b";"a"_{H_0}\ar@{=>}"d";"c"^{H_1}\ar@{=>}"e";"f"_{H_2}}$$
with
\begin{eqnarray*}
H_0&=&\sigma_0(e_1(X,\alpha)(-2))\vc\mu_{i,f},\\
H_1&=&\sigma_0(e_1(X,\alpha)(-1))\vc\mu_{q,i},\\
H_2&=&\sigma_0(e_1(X,\alpha)(0)-e_2(X,\alpha)(-2)(tq))\vc\mu_{tf,q}\vc(\nu_{f,s}
(\ell q)),
\end{eqnarray*}
is a track triangle that we call $\wp T$, see the proof of
Proposition \ref{reao}.

With this choice of a track triangle $\wp T$ for any small
$\nabla$-triangle $T$ we can obtain a $1$-cocycle $\theta$
representing the cohomology class $\theta_\nabla$ in
(\ref{conabla}), see the proof of that corollary. More precisely,
given a morphism $k\colon T\r \bar{T}$ between small
$\nabla$-triangles as in the following diagram
\begin{equation*}
\xymatrix{A\ar[r]^{f}\ar[d]^{k_0}&B\ar[r]^{i}\ar[d]^{k_1}&C\ar[r]^{q}\ar[d]^{k_2}&tA\ar[d]^{tk_0}\\
\bar{A}\ar[r]^{\bar{f}}&\bar{B}\ar[r]^{\bar{i}}&\bar{C}\ar[r]^{\bar{q}}&t\bar{A}}
\end{equation*}
if we write $(Y,\beta)$ for the $\nabla$-triangle corresponding to
$\bar{T}$ under the equivalence of categories in (\ref{sb'}) then
the cohomology class $\theta(k)=\vartheta_{\wp T,\wp
\bar{T}}(k)\in H^2(sh^3,[(X,\alpha),(Y,\beta)]^*\bar{t})$ is
represented by the cocycle
$(b_{-2},b_{-1},b_0)\in\tilde{G}^2([(X,\alpha),(Y,\beta)]^*\bar{t})$
with
\begin{eqnarray*}
b_{-2}&=&\sigma_0^{-1}(((\sigma_0(e_1(Y,\beta)(-2))\vc\mu_{\bar{i},\bar{f}})(\ell
k_0))\vc
((\ell\bar{i})(\mu_{\bar{f},k_0}^\vi\vc\mu_{k_1,f}))\\&&\vc((\mu^\vi_{\bar{i},k_1}\vc\mu_{k_2,i})(\ell
f))\vc((\ell k_2)(\mu_{i,f}^\vi\vc\sigma_0(-e_1(X,\alpha)(-2))))),\\
b_{-1}&=&\sigma_0^{-1}(((\sigma_0(e_1(Y,\beta)(-1))\vc\mu_{\bar{q},\bar{i}})(\ell
k_1))\vc
((\ell\bar{q})(\mu_{\bar{i},k_1}^\vi\vc\mu_{k_2,i}))\\&&\vc((\mu^\vi_{\bar{q},k_2}\vc\mu_{tk_0,q})(\ell
i))\vc((\ell t k_0)(\mu_{q,i}^\vi\vc\sigma_0(-e_1(X,\alpha)(-1))))),\\
b_0&=&\sigma_0^{-1}(((\sigma_0(e_1(Y,\beta)(0)-e_2(Y,\beta)(-2)(t\bar{q}))\vc\mu_{t\bar{f},\bar{q}}\vc(\nu_{\bar{f},s}
(\ell \bar{q})))(\ell
k_2))\\&&\vc(s(\ell\bar{f})(\mu_{\bar{q},k_2}^\vi\vc\mu_{t
k_0,q}\vc(\nu_{k_0,s}(\ell q))))\vc((s^\vi_{\ell \bar{f},\ell
k_0}\vc s(\mu_{\bar{f},k_0}^\vi\vc\mu_{k_1,f})\vc s_{\ell k_1,\ell
f})(\ell q))\\&&\vc(s(\ell k_1)((\nu^\vi_{f,s}(\ell
q))\vc\mu^\vi_{tf,q}\vc\sigma_0(-e_1(X,\alpha)(0)+e_2(X,\alpha)(-2)(tq))))).
\end{eqnarray*}
The reader can check that
\begin{eqnarray*}
b_{-2}&=&c_\C{B}(\bar{i},\bar{f},k_0)-c_\C{B}(\bar{i},k_1,f)+c_\C{B}(k_2,i,f)\\&&-k_2(e_1(X,\alpha)(-2))+(e_1(Y,\beta)(-2))(t
k_0),\\
b_{-1}&=&c_\C{B}(\bar{q},\bar{i},k_1)-c_\C{B}(\bar{q},k_2,i)+c_\C{B}(tk_0,q,i)\\&&-(tk_0)(e_1(X,\alpha)(-1))+(e_1(Y,\beta)(-1))(t
k_2),\\
b_0&=&c_\C{B}(t\bar{f},\bar{q},k_2)-c_\C{B}(t\bar{f},tk_0,q)+c_\C{B}(tk_1,tf,q)\\&&-(tk_1)(e_1(X,\alpha)(0))+(e_1(Y,\beta)(0))(t
k_2)\\&&+b_s(\bar{f},k_0)(t q)-b_s(k_1,f)(t q)\\&&
+k_1(e_2(X,\alpha)(-2))(tq)-(e_2(Y,\beta)(-2))(t (\bar{q}k_2)),
\end{eqnarray*}
therefore this proposition follows from Lemmas \ref{aysh2} and
\ref{peque} and the fact that $(c_\C{B},b_s)$ is normalized.
\end{proof}

\begin{proof}[Proof of Theorem \ref{ñi}]
By \cite{tbhp} 2.4 and 2.5 there is a map from the weak linear
extension of homotopy pairs in \cite{ccctc1} 4.3 to the linear
extension of algebraic homotopy pairs in (\ref{leahp})
\begin{equation*}\tag{a}
\xymatrix{(\hom^t)^\#\ar[r]^+\ar@{=}[d]&\C{Hopair}(\C{B})\ar[d]_\sim^m\ar[r]^{\tilde{p}}&\C{Pair}(\C{A})\ar@{=}[d]\\
(\hom^t)^\#\ar[r]^+&\C{AHp}(\aleph)\ar[r]^{\bar{p}}&\C{Pair}(\C{A})}
\end{equation*}
Here $m$ is an equivalence of categories. We point out that the
category of homotopy pairs considered in \cite{tbhp} does not
coincide strictly with our $\C{Hopair}(\C{B})$ defined in
\cite{ccctc1} 4. However it is 
equivalent to $\C{Hopair}(\C{B})$ in an obvious way which is
compatible with the (weak) linear extension structures above.

The equivalence of categories $m$ in (a) induces isomorphisms
between the track-theoretical cone functors in \cite{ccctc1} 4.4 and
the algebraic ones in (\ref{conalg}), i. e. for any map $f$ in
$\C{B}$ we have
\begin{equation*}\tag{c}
m\colon\cone_f\cong\acone_{p(f)}.
\end{equation*}
In particular $\cone_f$ is representable if and only if
$\acone_{p(f)}$ is. These isomorphisms show that the following
diagram commutes up to equivalence
\begin{equation*}\tag{d}
\xymatrix{\C{Hopair}(\C{B})^\op\ar[dd]_\sim^{m^\op}\ar[rd]^{\cone}&\\&\mo(\C{A})\\\C{AHp}(\aleph)^\op\ar[ru]_{\acone}}
\end{equation*}
As we saw in \cite{ccctc1} 4.8 and Section \ref{ahp} diagram (d)
restricted to (algebraic) homotopy pairs with representable
associated (algebraic) cone functor can be extended in a commutative
way up to equivalence
\begin{equation*}
\xymatrix{\C{Hopair}^{rep}(\C{B})^\op\ar[dd]_\sim^{m^\op}\ar[rd]^{C^\op}\ar@/^20pt/[rrd]^{\cone}&&\\&\C{A}^\op\ar[r]&\mo(\C{A})\\
\C{AHp}^{rep}(\aleph)^\op\ar[ru]_{\bar{C}^\op}\ar@/_20pt/[rru]_{\acone}&&}
\end{equation*}
The functors $C$ and $\bar{C}$ were used to define functors $\zeta$ and $\bar{\zeta}$ that therefore fit in a diagram commuting up to equivalence
\begin{equation*}\tag{e}
\xymatrix{\C{Hopair}^{rep}(\C{B})\ar[dd]_\sim^{m}\ar[rd]^{\zeta}&\\&\C{cand}(\C{A},t)\\\C{AHp}^{rep}(\aleph)\ar[ru]_{\bar{\zeta}}}
\end{equation*}
Now by Proposition \ref{reao} the vanishing condition in the
statement of Theorem \ref{ñi} is equivalent to the fact that
$\hat{\zeta}$ takes values in small $\nabla$-triangles. This is
obviously the same as saying that $\bar{\zeta}$ takes values in
small $\nabla$-triangles and by the commutativity properties of (e)
this happens if and only if $\zeta$ takes values in small
$\nabla$-triangles, therefore the theorem follows from \cite{ccctc1}
4.26.
\end{proof}

\begin{proof}[Proof of Theorem \ref{vanish}]
This proof is technically a continuation of the proof of Theorem
\ref{ñi} above since we need here all arguments already considered
there. In this proof we will pay attention to the functors in the
following extended diagram
\begin{equation*}\tag{f}
\xymatrix{&\C{candt}(\C{B},t)\ar[ld]_{\bar{\varrho}}\ar[dd]^\varrho\ar@{}[lddd]|<(.3){\text{(g)}}\\
\C{Hopair}(\C{B})\ar[dd]_\sim^{m}\ar[rd]^{\zeta}&\\&\C{cand}(\nabla)\\\C{AHp}(\aleph)\ar@{->>}[rd]\ar[ru]_{\bar{\zeta}}\\
&\C{AHp}(\aleph)/(\varpi,\bar{\zeta})\ar@{^{(}->}[uu]_{\hat{\zeta}}}
\end{equation*}
The lower triangle is strictly commutative, see Section \ref{ahp},
and (g) commutes up to equivalence by \cite{ccctc1} 4.16. Moreover,
$\hat{\zeta}$ is faithful here by the second part of \cite{ccctc1}
4.25.

By \cite{ccctc1} 4.27 what we have to prove is the following:
$\hat{\zeta}^*\theta_\nabla=0$ if and only if there exists a full
subcategory $\C{dist}(\C{B},s)\subset\C{candt}(\C{B},s)$ such that
\begin{equation*}\tag{h}
\C{dist}(\C{B},s)\subset\C{candt}(\C{B},s)\st{\bar{\varrho}}\To\C{Hopair}(\C{B})\st{\tilde{p}}\To\C{Pair}(\C{A})
\end{equation*}
is surjective on objects and \
\begin{equation*}\tag{i}
\C{dist}(\C{B},s)\subset\C{candt}(\C{B},s)\st{\bar{\varrho}}\To\C{Hopair}(\C{B})
\end{equation*}
is full.

For this we give the following description of a normalized
$1$-cocycle $\theta$ representing $\theta_\nabla$ in terms of the
obstruction operator $\vartheta$ in the exact sequence for the
functor $\varrho$ in \cite{ccctc1} 3.1. For each small
$\nabla$-triangle $T$ we choose a track triangle $\wp T$ such that
$\varrho\wp T=T$. Then given a morphism of small $\nabla$-triangles
$k\colon T\r\bar{T}$ we have
$$\theta(f)=\vartheta_{\wp T,\wp\bar{T}}(k).$$
Actually all cocycles representing $\theta_\nabla$ can be obtained in this way, compare \cite{ah} IV.4.10 (d).

Suppose first that there exists a full subcategory
$\C{dist}(\C{B},s)$ such that the functors (h) and (i) satisfy the
required properties. For any morphism $f$ in $\C{A}$ we can choose
the representation $$\chi_f\colon\hom_\C{A}(C_f,-)\cong\acone_f$$
used in the definition of $\bar{\zeta}$ to be induced by the
isomorphism in (c) and a track triangle $\wp\bar{\zeta}(f)$ in $\C{dist}(\C{B},s)$ whose
image by $\tilde{p}\bar{\varrho}$ is $f$, see \cite{ccctc1} 4.12.  Now it is
immediate to notice by using \cite{ccctc1} 3.2 and the commutativity
properties of (f) that the pull-back of the $1$-cocycle $\theta$ by
the functor $\hat{\zeta}$ is already trivial
$\hat{\zeta}^*\theta=0$, in particular
$\hat{\zeta}^*\theta_\nabla=0$.

Suppose now that $\hat{\zeta}^*\theta_\nabla=0$. We choose a
$0$-cocycle $c$ such that $\delta(c)=\hat{\zeta}^*\theta$ and a natural isomorphism $\mathfrak{q}\colon\varrho\cong\zeta\bar{\varrho}$. We take
$\C{dist}(\C{B},s)$ to be the full subcategory given by the track
triangles
$$\tilde{T}_f=\wp\hat{\zeta}(f)+(-c(f))$$
in the sense of \cite{ah} IV.4.10 (d), where $f$ runs along all
morphisms in $\C{A}$. With this definition (h) is obviously
surjective on objects, since
$\tilde{p}\bar{\varrho}\tilde{T}_f=\tilde{p}\bar{\varrho}\wp\hat{\zeta}(f)=f$.
Let us check that (i) is full. We have to show that any morphism of
homotopy pairs
$[k_0,k_1,K]\colon\bar{\varrho}\tilde{T}_f\r\bar{\varrho}\tilde{T}_g$
is the image by $\bar{\varrho}$ of a track triangle morphism
$\tilde{k}\colon\tilde{T}_f\r\tilde{T}_g$. Consider the next chain
of equalities
\begin{eqnarray*}
\vartheta_{\tilde{T}_f,\tilde{T}_g}(\mathfrak{q}_{\tilde{T}_g}^{-1}\zeta[k_0,k_1,K]\mathfrak{q}_{\tilde{T}_f})&=&
\vartheta_{\wp\hat{\zeta}(f),\wp\hat{\zeta}(g)}(\mathfrak{q}_{\tilde{T}_g}^{-1}\zeta[k_0,k_1,K]\mathfrak{q}_{\tilde{T}_f})
\\&&+(-(\mathfrak{q}_{\tilde{T}_g}^{-1}\zeta[k_0,k_1,K]\mathfrak{q}_{\tilde{T}_f})_*c(f)
\\&&+(\mathfrak{q}_{\tilde{T}_g}^{-1}\zeta[k_0,k_1,K]\mathfrak{q}_{\tilde{T}_f})^*c(g))
\\&=&\theta(\mathfrak{q}_{\tilde{T}_g}^{-1}\zeta[k_0,k_1,K]\mathfrak{q}_{\tilde{T}_f})
\\&&-(\delta(c))(\mathfrak{q}_{\tilde{T}_g}^{-1}\zeta[k_0,k_1,K]\mathfrak{q}_{\tilde{T}_f})\\&=&0.
\end{eqnarray*}
Here we use \cite{ccctc1} 3.3 for the first equality. By
\cite{ccctc1} 3.2, and since the triangle (g) in diagram (f)
commutes up to the natural equivalence $\mathfrak{q}$, we see that
$\zeta[k_0,k_1,K]\colon\zeta\bar{\varrho}\tilde{T}_f\r\zeta\bar{\varrho}\tilde{T}_g$
in $\C{cand}(\nabla)$ is the image by $\zeta\bar{\varrho}$ of a
morphism $\check{k}\colon\tilde{T}_f\r\tilde{T}_g$ in
$\C{candt}(\C{B},s)$, however $\bar{\varrho}\check{k}$ need not be
the same as $[k_0,k_1,K]$. Nevertheless by the commutativity
properties of (f), since $\hat{\zeta}$ is faithful and $m$ is an
equivalence fitting into the map of linear extensions (a) then there
exists $b\in H^1(sh^3,[\hat{\zeta}(f),\hat{\zeta}(g)]^*\bar{t})$
such that
$$[k_0,k_1,K]=\bar{\varrho}(\check{k})+\varpi_{(\hat{\zeta}(f),\hat{\zeta}(g))}(b),$$
therefore by \cite{ccctc1} 4.21 the morphism
$$\tilde{k}=\check{k}+b\colon\tilde{T}_f\To\tilde{T}_g$$
satisfies $\bar{\varrho}(\tilde{k})=[k_0,k_1,K]$ as desired. Now the proof is finished.
\end{proof}

\appendix

\section{The cosimplicial Eilenberg-Zilber-Cartier theorem}

Let $\Delta$ be the usual simplicial category. A
\emph{bicosimplicial abelian group}\index{bicosimplicial abelian
group} is a functor $A^{\bullet,\bullet}\colon \Delta\times\Delta
\r\C{Ab}$. This functor is completely determined by the abelian
groups $A^{p,q}$ $(p,q\geq 0)$, the horizontal cosimplicial
coboundaries  $d^h_i$ and codegeneracies $s^h_i$ defining the
horizontal cosimplicial abelian groups $A^{\bullet,q}$ $(q\geq
0)$, and the vertical coboundaries  $d^v_i$ and codegeneracies
$s^v_i$ which give rise to the vertical cosimplicial abelian
groups $A^{p,\bullet}$ $(p\geq 0)$. Horizontal and vertical
operators commute.

The \emph{diagonal}\index{diagonal cochain complex}
$\diag(A^{\bullet,\bullet})$ of a bicosimplicial abelian group is
the cosimplicial abelian group given by
$\diag^n(A^{\bullet,\bullet})=A^{n,n}$ with coboundaries
$d_i=d^h_id^v_i=d^v_id^h_i$ and codegeneracies
$s_i=s^h_is^v_i=s^v_is^h_i$.

The \emph{normalized cochain complex}\index{normalized cochain
complex} $N(B^\bullet)$ of a cosimplicial abelian group
$B^\bullet$ is given by
$$N^n(B^\bullet)=\bigcap_{i=0}^{n-1}\ker[s_i\colon B^n\r B^{n-1}]$$
with differential $d=\sum_{i=0}^{n+1}(-1)^id_i\colon
N^n(B^\bullet)\r N^{n+1}(B^\bullet)$.

One can also define the \emph{normalized
bicomplex}\index{normalized bicomplex} $N(A^{\bullet,\bullet})$ of
a bicosimplicial abelian group as
$$N^{p,q}(A^{\bullet,\bullet})=\left(\bigcap_{i=0}^{p-1}\ker[s_i^h\colon
A^{p,q}\r A^{p-1,q}]\right)\cap\left(\bigcap_{i=0}^{q-1}
\ker[s_i^v\colon A^{p,q}\r A^{p,q-1}]\right).$$ The horizontal
differential $d_h\colon N^{p,q}(A^{\bullet,\bullet})\r
N^{p+1,q}(A^{\bullet,\bullet})$ is given by
$d_h=\sum_{i=0}^{p+1}(-1)^id^h_i$ and the vertical differential
$d_v\colon N^{p,q}(A^{\bullet,\bullet})\r
N^{p,q+1}(A^{\bullet,\bullet})$ is
$d_v=\sum_{i=0}^{q+1}(-1)^id^v_i$. Notice that horizontal and
vertical differentials commute.

Recall that for any bicomplex $C^{*,*}$ with commuting horizontal
$d_h$ and vertical $d_v$ differentials $d_hd_v=d_vd_h$ the
\emph{total complex}\index{total complex} $\tot C^{*,*}$ is the
cochain complex defined by
$$\tot^n C^{*,*}=\bigoplus_{p+q=n}C^{p,q}$$
and such that the differential of an element $a\in C^{p,q}$ is
$d(a)=d_h(a)+(-1)^pd_v(a)$.

The following well known theorem is the dual of the
Eilenberg-Zilber-Cartier theorem, see \cite{hnaf} 2.9. For the
statement we recall that a \emph{$(p,q)$-shuffle}\index{shuffle}
$\sigma$ is a permutation of the set $\set{0,\dots,p+q-1}$ such
that
$$\sigma(0)<\cdots<\sigma(p-1)\text{ and }\sigma(p)<\cdots<\sigma(p+q-1).$$
The parity of the $(p,q)$-shuffle $\sigma$ as a permutation is denoted by $\epsilon(\sigma)\in\set{0,1}$. It is known that
$\epsilon(\sigma)=\sum_{i=0}^{p-1}(\sigma(i)-i) \mod 2$.

\begin{thm}\label{ezc}
Let $A^{\bullet,\bullet}$ be a bicosimplicial abelian group. There are mutually inverse natural cochain homotopy equivalences
$$\tot N(A^{\bullet,\bullet})\mathop{\leftrightarrows}_{\mathfrak{G}}^{\mathfrak{F}} N\diag(A^{\bullet,\bullet})$$
defined by
$$\mathfrak{F}(a)=\sum_{p+q=n}\sum_{\sigma\in\set{\text{$(p,q)$-shuffles}}}
(-1)^{\epsilon(\sigma)}s^v_{\sigma(0)}\cdots s^v_{\sigma(p-1)}s^h_{\sigma(p)}\cdots s^h_{\sigma(p+q-1)}(a) ,\;\; a\in A^{n,n};$$
$$\mathfrak{G}(b)=  d^h_{p+q}\cdots d^h_{p+1} d^v_{p-1}\cdots d^v_0 (b),\;\; b\in A^{p,q}.$$
\end{thm}

\bibliographystyle{alpha}
\bibliography{Fernando}